\documentclass[a4paper,12pt]{article}

\usepackage[top=3.5cm, bottom=3.5cm, left=2.5cm, right=2.5cm]{geometry} 
\usepackage{amsmath}
\usepackage{amsthm}
\usepackage{amssymb}
  \usepackage{paralist}
  \usepackage{graphics} 
  \usepackage{epsfig} 
\usepackage{graphicx}  \usepackage{epstopdf}

 \usepackage[colorlinks=false]{hyperref}
\usepackage[nottoc, notlof, notlot]{tocbibind}

\usepackage{float}
\usepackage{tikz}
\usetikzlibrary{arrows}
\usetikzlibrary[patterns]

\usepackage[boxed]{algorithm2e}
\usepackage{setspace}
\usepackage{etoolbox}
\AtBeginEnvironment{algorithm}{\setstretch{1.2}}

\usepackage{array}
\usepackage{tabularx}
\usepackage{ltablex} 
\usepackage{caption}

\usepackage{authblk}
\usepackage{mathrsfs}
\usepackage{thmtools}
\usepackage{bbm}

\usepackage{stmaryrd}

\makeatletter
\newsavebox{\@brx}
\newcommand{\llangle}[1][]{\savebox{\@brx}{\(\m@th{#1\langle}\)}%
  \mathopen{\copy\@brx\kern-0.5\wd\@brx\usebox{\@brx}}}
\newcommand{\rrangle}[1][]{\savebox{\@brx}{\(\m@th{#1\rangle}\)}%
  \mathclose{\copy\@brx\kern-0.5\wd\@brx\usebox{\@brx}}}
\makeatother

\usepackage{cleveref}
\crefname{II}{}{}
\creflabelformat{II}{#2II-#1#3}
\crefname{IIeq}{}{}
\creflabelformat{IIeq}{#2(II-#1)#3}
\crefname{I}{}{}
\creflabelformat{I}{#2#1#3}
\crefname{Ieq}{}{}
\creflabelformat{Ieq}{#2(#1)#3}

\usepackage{comment}

\usepackage{nameref,zref-xr}
\zxrsetup{toltxlabel}
\zexternaldocument*{part2}

\newcommand{\R}{\mathbb{R}}
\newcommand{\N}{\mathbb{N}}
\newcommand{\1}{\mathbbm{1}}
\newcommand{\E}{\mathbb{E}}

\newcommand{\e}{\mathrm{e}}
\newcommand{\dd}{\mathrm{d}}
\newcommand{\pb}{\mathcal{P}}
\DeclareMathOperator{\Tr}{Tr}
\DeclareMathOperator{\Dom}{Dom}

\DeclareMathOperator{\diag}{diag}

\newtheorem{theorem}{Theorem}[section]
\newtheorem{corollary}{Corollary}

\newtheorem{lemma}[theorem]{Lemma}
\newtheorem{proposition}{Proposition}

\newtheorem{assumption}{Assumption}
\theoremstyle{definition}
\newtheorem{definition}[theorem]{Definition}
\newtheorem{remark}{Remark}
\newtheorem{example}{Example}

\title{\bf \LARGE Propagation of chaos: a review of models, methods and applications. \\ I. Models and methods}

\author[1]{Louis-Pierre \textsc{Chaintron}}
\author[2]{Antoine \textsc{Diez}}
\affil[1]{\small
DMA, \'Ecole Normale Sup\'erieure

45 rue d'Ulm

75005 Paris, France

\url{lchaintron@clipper.ens.fr}
\bigskip
}
\affil[2]{\small
Department of Mathematics, 

Imperial College London, South Kensington Campus,

London, SW7 2AZ, UK

\url{antoine.diez18@imperial.ac.uk}
\bigskip
}

\begin{document}
\maketitle

\begin{abstract}
The notion of propagation of chaos for large systems of interacting particles originates in statistical physics and has recently become a central notion in many areas of applied mathematics. The present review describes old and new methods as well as several important results in the field. The models considered include the McKean-Vlasov diffusion, the mean-field jump models and the Boltzmann models. The first part of this review is an introduction to modelling aspects of stochastic particle systems and to the notion of propagation of chaos. The second part presents concrete applications and a more detailed study of some of the important models in the field.
\end{abstract}



\newpage

\tableofcontents

\section{Introduction}

When Boltzmann published his most famous article \cite{boltzmann_weitere_1872} one century and a half ago, the study of large systems of interacting particles was entirely motivated by the microscopic modelling of thermodynamic systems. Although it was far from being an accepted idea at that time, Boltzmann postulated that since a macroscopic volume of gas contains a myriad of elementary particles, it is both hopeless and needless to keep track of each particle and one should rather seek a statistical description. He thus derived the equation that now bears his name and which gives the time evolution of the continuum probability distribution (in the phase space) of a typical particle. With the H-theorem, he also extended and justified the pioneering works of Maxwell and Clausius for equilibrium thermodynamic systems, paving the way alongside Gibbs for a consistent kinetic theory of gases. The Boltzmann equation is derived from first principles under a crucial assumption, called \emph{molecular chaos}. This assumption was already known from Maxwell and is often called the \emph{Stosszahlansatz} since Ehrenfest. Informally, it translates the idea that, despite the multitude of interactions, two particles taken at random should be statistically independent when the total number of particles grows to infinity. It is not so clear how the appearance of probability theory should be interpreted in this context. In the following years, the Stosszahlansatz and its consequences (the H-theorem) were the object of a fierce debate among physicists as they seem to break the microscopic reversibility. Beyond the scientific debate, it has raised metaphysical and philosophical questions about the profound nature of time and randomness. 

The rigorous justification of the work of Boltzmann and the status of molecular chaos became true mathematical questions when Hilbert addressed them in his Sixth Problem at the Paris International Congress of Mathematicians in $1900$. 
Quoting Hilbert, the problem which motivates the present work is to \textit{``[develop] mathematically the limiting processes 
[\ldots] which lead from the atomistic view to the laws of motion of continua''}. Our starting point will be the seminal article of Kac \cite{kac_foundations_1956}. More than half a century after Hilbert, Kac gave the first rigorous mathematical definition of chaos and introduced the idea that for time-evolving systems, chaos should be propagated in time, a property therefore called the \emph{propagation of chaos}. Kac was still motivated by the mathematical justification of the classical collisional kinetic theory of Boltzmann for which he developed a simplified probabilistic model. Soon after Kac, McKean \cite{mckean_propagation_1969} introduced a class of diffusion models which were not originally part of Boltzmann theory but which satisfy Kac's propagation of chaos property. In the classical kinetic theory of Boltzmann, the problem is the derivation of continuum models starting from deterministic, Newtonian, systems of particles. In comparison, the fundamental contribution of Kac and McKean is to have shown that the classical equations of kinetic theory also have a natural stochastic interpretation. This philosophical shift is addressed in the enlightening introduction of Kac \cite{cohen_probabilistic_1973} written for the centenary of the Boltzmann equation. 

Kac and McKean introduced a new mathematical formalism, gave many insights on the stochastic modelling in kinetic theory and proved the two building block theorems (Theorem \cref{thm:mckean} and Theorem \cref{thm:kac}). Their works have stimulated the development of a rich and still active \emph{mathematical kinetic theory}. Keeping strong connections with the original theory of Boltzmann, some fundamental questions raised several decades ago have been answered only recently (see for instance \cite{bodineau_brownian_2016,gallagher_newton_2014,mischler_kacs_2013}). On the other hand, systems of interacting particles are ubiquitous in many applications now and over the last two decades, the tools and concepts developed in kinetic theory have somehow escaped the realm of pure statistical physics. This review paper is motivated by the growing number of models in applied mathematics where the notion of chaos plays a central role. Some recent new domains of applications include the following ones. In mathematical biology and social sciences, self-organization models describe systems of indistinguishable particles (birds, insects, bacteria, crowds\ldots) with a behaviour which can hardly be predicted at the microscopic scale but which are (sometimes) well explained by continuum models derived within the framework of mathematical kinetic theory \cite{bellomo_active_2017, bellomo_active_2019, naldi_mathematical_2010, muntean_collective_2014, degond_mathematical_2018, albi_vehicular_2019, vicsek_collective_2012}. In another context, the recent theory of mean-field games studies the asymptotic properties of games with a large numbers of players \cite{cardaliaguet_master_2019, cardaliaguet_notes_2010, carmona_probabilistic_2018, carmona_probabilistic_2018-1}. Even more recently, systems of particles have been used to model complex phenomena in data sciences, with applications in Markov Chain Monte Carlo theory \cite{del_moral_measure-valued_1998, del_moral_feynman-kac_2004, del_moral_mean_2013}, in optimization \cite{pinnau_consensus-based_2017, totzeck_trends_2021, grassi_particle_2021, carrillo_consensus-based_2021}, or for the training of neural networks \cite{mei_mean_2018, rotskoff_trainability_2019, sirignano_mean_2020, chizat_global_2018, de_bortoli_quantitative_2020}. Compared to the models in statistical physics, many aspects should be reconsidered. To cite a few examples: the basic conservation laws (of momentum, energy\ldots) do not always hold for biological systems and may be replaced by other types of constraints (optimization constraints, geometrical constraints\ldots); the intrinsic randomness (or uncertainty) of the models in applied sciences is often a crucial modelling assumption; the complexity of the interaction mechanisms entails new analytical tools etc. These differences have motivated many new techniques, new insights on the question of propagation of chaos and in the end, new results. 

This review article on propagation of chaos is not the first one on the subject. The course of Sznitman at Saint-Flour \cite{sznitman_topics_1991} studies many of the most important historical probabilistic models. The probabilistic methods are explained in details in the book \cite{graham_probabilistic_1996} (in particular the courses of M\'el\'eard \cite{meleard_asymptotic_1996} and Pulvirenti \cite{pulvirenti_kinetic_1996}). More recently, the review of Jabin and Wang \cite{jabin_mean_2017} focuses on McKean mean-field systems and PDE applications. By its nature, the notion of chaos lies in the interplay between probability theory and Partial Differential Equations. The present review discusses both analytic and probabilistic methods and includes many (very) recent results. We also refer to the article by Hauray and Mischler \cite{hauray_kacs_2014} which is to our knowledge, the most complete reference on Kac's chaos (without propagation of chaos). For deterministic systems which will not be considered here, we refer to the very thorough reviews \cite{jabin_review_2014,golse_dynamics_2016}.

\subsubsection*{Outline.} The article is organised as follows. 

Section \cref{sec:models} introduces the setting and the conventions that will be used throughout the article. A gallery of the models which will be studied is presented; we will distinguish McKean's mean-field jump and diffusion models (Section \cref{sec:meanfieldmodels}) and Boltzmann-Kac models (Section \cref{sec:boltzmann}). 

Section \cref{sec:def} is devoted to the description of the fundamental tools and concepts needed in the study of exchangeable particle systems. The central notions of chaos (Section \cref{sec:kacchaos}) and propagation of chaos (Section \cref{sec:poc}) are defined in this section. 

Section \cref{sec:proving} is a review of the methods used to prove propagation of chaos. Several probabilistic and analytical techniques are described as well as abstract theorems which will be applied to specific models in the second part of this review.  

Finally, for the reader's convenience, we collect in Appendix \cref{appendix:probabilityreminders} useful notions and results in Probability theory regarding stochastic processes in the Skorokhod space, Markov processes, martingale methods, large deviations, the Girsanov transform and the theory of Poisson random measures. Some corollaries of the quantitative Hewitt-Savage theorem are gathered in Appendix \cref{appendix:hewittsavage}. 

The second part of this review will be devoted on the one hand to the classical models introduced by Kac and McKean and their recent developments and on the other hand to a gallery of recent applications in applied mathematics and beyond. Throughout this first part, references to the second part are indicated by ``II-'' (for instance Section \cref{sec:summary} refers to the second section of the second part).

\subsection*{Notations and conventions}

\subsubsection*{Sets}

\noindent\begin{tabularx}{\linewidth}{@{}lX}
    $C( I , E)$ & The set of continuous functions from a time interval $I=[0,T]$ to a set $E$, endowed with the uniform topology.   \\
    $C_b(E)$, $C_b^k(E)$ & Respectively the set of real-valued bounded continuous functions and the set of functions with $k\geq1$ bounded continuous derivatives on a set $E$. \\
    $C_c(E)$ & The set of real-valued continuous functions with compact support on a locally compact space $E$. \\
    $C_0(E)$ & The set of real-valued continuous functions vanishing at infinity on a locally compact space $E$, i.e. $\varphi\in C_0(E)$ when for all $\varepsilon>0$, there exists a compact set $K_\varepsilon\subset E$ such that $|\varphi(x)|<\varepsilon$ for all $x\in E$ outside $K_\varepsilon$.\\ 
    $D( I , E)$ & The space of functions which are right continuous and have left limit everywhere from a time interval $I=[0,T]$ to a set $E$, endowed with the Skorokhod $J1$ topology. This is the space of the \emph{c\`adl\`ag} functions. This space is also called the \emph{Skorokhod space} or the \emph{path space}. \\
    $L^p(E)$ or $L^p_\mu(E)$ & The set of measurable functions $\varphi$ defined almost everywhere on a measured space $(E,\mu)$ such that the $|\varphi|^p$ is integrable for $p\geq1$. When $p=+\infty$, this is the set of functions with a bounded essential supremum. We do not specify the dependency in $\mu$ when no confusion is possible. \\
    $\mathcal{M}_d(\R)$ & The set of $d$-dimensional square real matrices. \\
    $\mathcal{M}(E)$ & The set of signed measures on a measurable space $E$. \\
    $\mathcal{M}^+(E)$ & The set of positive measures on a measurable space $E$. \\
    $\pb(E)$ & The set of probability measures on a space $E$.\\
    $\pb_p(E)$ & The set of probability measures with bounded moment of order $p\geq1$ on a space $E$.\\
    $\widehat{\pb}_N(E)$ & The set of empirical measures of size $N$ over  a set $E$, that is measures of the form $\mu = \frac{1}{N}\sum_{i=1}^N \delta_{x^i}$, where $x^i\in E$. \\
    $\R_+$ & The set $[0,+\infty)$.\\
    $\mathfrak{S}_N$ & The permutation group of the set $\{1,\ldots,N\}$.\\
    $\mathbb{S}^{d-1}$ & The sphere of dimension $d-1$. 
\end{tabularx}

\subsubsection*{Generic elements and operations}

\noindent\begin{tabularx}{\linewidth}{@{}lX}

$C$ & A generic nonnegative constant, the value of which may change from line to line. \\
$C(a_1,\ldots a_n)$ & A generic nonnegative constant which depends on some fixed parameters denoted by $a_1,\ldots,a_n$. Its value may change from line to line.\\
$\diag(x)$ & The $d$-dimensional diagonal matrix whose diagonal coefficients $x_1,\ldots,x_d$ are the components of the $d$-dimensional vector $x$. \\ 
$\nabla\cdot V$ & The divergence of a vector field $V:\R^d\to\R^d$ or of a matrix field $V:\R^d\to\mathcal{M}_d(\R)$, respectively defined by $\nabla\cdot V = \sum_{i=1}^d \partial_{x_i} V_i$ or componentwise by $(\nabla\cdot V)_i = \sum_{j=1}^d \partial_{x_j}V_{ij}$. \\
$A:B$ and $\|A\|$ & The Frobenius inner product of two matrices $A,B\in\mathcal{M}_d(\R)$ defined by $A:B:=\sum_{i=1}^d\sum_{j=1}^d A_{ij}B_{ij}$ and the associated norm $\|A\|:=\sqrt{A:A}$. \\
$\nabla^2 V$ & The Hessian matrix of a scalar field $V:\R^d\to \R$ defined componentwise by $(\nabla^2 V)_{ij} = \partial^2_{x_i, x_j} V$. \\
$I_d$ & The $d$-dimensional identity matrix.\\
$\mathrm{Id}$ & The identity operator on a vector space.\\
$\langle x, y\rangle$ or $x\cdot y$ & The Euclidean inner product of two vectors $x,y\in\R^d$ defined by $\langle x,y\rangle \equiv x\cdot y := \sum_{i=1}^d x^i y^i$. One notation or the other may be preferred for typographical reasons in certain cases. \\ 
$M_{ij}$ & The $(i,j)$ (respectively row and column indexes) component of a matrix $M$. \\ 
$\mathsf{P}(u)$ & The projection matrix $\mathsf{P}(u):=I_d-\frac{u\otimes u}{|u|^2}$ on the plane orthogonal to a vector $u\in\R^d$.\\
$\varphi\in C_b(E)$ & A generic bounded continuous test function on $E$.\\
$\varphi_N\in C_b(E^N)$ & A generic bounded continuous test function on the product space $E^N$. \\
$\Phi\in C_b(\pb(E))$ & A generic bounded continuous test function on the set of probability measures on $E$. \\
$u\otimes v$, $\mu\otimes\nu$ or $\varphi\otimes\psi$ & Respectively, the matrix tensor product of two vectors $u,v\in\R^d$ defined componentwise by $(u\otimes v)_{ij} = u_iv_j$; the product measure on $E\times F$ of two measures $\mu,\nu$ respectively on $E$ and $F$; the product function on $E\times F$ defined by $(\varphi\otimes\psi)(x,y)=\varphi(x)\psi(y)$ for two real-valued function $\varphi,\psi$ respectively on $E$ and $F$. \\
$\Tr M$ & The trace of the matrix $M$. \\
$M^\mathrm{T}$ & The transpose of the matrix $M$.\\
$\mathbf{x}^N=(x^1,\ldots,x^N)$ & A generic element of a product space $E^N$. The components are indexed with a superscript. \\
$\mathbf{x}^{M,N}=(x^1,\ldots,x^M)$ & The $M$-dimensional vector in $E^M$ constructed by taking the $M$ first components of $\mathbf{x}^N$. \\
$x=(x_1,\ldots,x_d)^\mathrm{T}$ and $|x|$ & A generic element of a $d$-dimensional space and its norm. The coordinates are indexed with a subscript. The norm of $x$ denoted by $|x|$ is the Euclidean norm. 
\end{tabularx}

\subsubsection*{Probability and measures}

\noindent\begin{tabularx}{\linewidth}{@{}lX}
$K\star\mu$ & The convolution of a function $K:E\times F\to G$ with a measure $\mu$ on $F$ defined as the function $K\star\mu:x\in E \mapsto \int_F K(x,y)\mu(\dd y)\in G$. When $E=F=G=\R^d$ and $K:\R^d\to\R^d$, we write $K\star\mu(x)=\int_{\R^d}K(x-y)\mu(\dd y)$. \\
$\delta_x$ & The Dirac measure at the point $x$. \\
$\mu_{\mathbf{x}^N}$ & The empirical measure defined by $\mu_{\mathbf{x}^N}=\frac{1}{N}\sum_{i=1}^N\delta_{x^i}$ where $\mathbf{x}^N=(x^1,\ldots,x^N)$. \\
$\E_\mu[\varphi]$ & Alternative expression for $\langle \mu,\varphi\rangle$ when $\mu$ is a probability measure. When $\mu=\mathbb{P}$ on $(\Omega,\mathscr{F},(\mathscr{F}_t)_t,\mathbb{P})$, the expectation is simply denoted by $\mathbb{E}$. \\
$H(\nu|\mu)$ & The relative entropy (or Kullback-Leibler divergence) between two measures $\mu,\nu$, see Defintion \cref{def:entropyfisher}.\\ 
$\langle \mu,\varphi\rangle$ & The integral of a measurable function $\varphi$ with respect to a measure $\mu$. \\
$\mathrm{Law}(X)$ & The law of a random variable $X$ as an element of $\pb(E)$ where $X$ takes its value in the space $E$.\\
$(\Omega,\mathscr{F},(\mathscr{F}_t)_t,\mathbb{P})$ & A filtered probability space. Unless otherwise stated, all the random variables are defined on this set. The expectation is denoted by $\E$. \\
$\sigma(X^1,X^2,\ldots)$ & The $\sigma$-algebra generated by the random variables $X^1, X^2,\ldots$. \\
$T_\#\mu$ & The pushforward of the measure $\mu$ on a set $E$ by the measurable map $T:E\to F$. This is a measure on the set $F$ defined by $T_\#\mu(\mathscr{A}) = \mu(T^{-1}(\mathscr{A}))$ for any measurable set $\mathscr{A}$ of $F$. \\
$\|\cdot\|_{\mathrm{TV}}$ & The Total Variation (TV) norm for measures. \\
$W_p$ & The Wasserstein-$p$ distance between probability measures (see Definition \cref{def:wasserstein}).\\
$X\sim\mu$ & It means that the law of the random variable $X$ is $\mu$. \\ 
$(\mathsf{X}_t)_t$ or $(\mathsf{Z}_t)_t$ & The canonical process on the path space $D(I,E)$ defined by $\mathsf{X}_t(\omega)=\omega(t)$. \\
$(\mathbf{X}^N_t)^{}_t$ or $(\mathbf{Z}^N_t)^{}_t$ & The canonical process on the product space $D(I,E)^N$ with components $\mathbf{X}^N_t=(\mathsf{X}^1_t,\ldots,\mathsf{X}^N_t)$. 
\end{tabularx}

\subsubsection*{Systems of particles and operators}

\noindent\begin{tabularx}{\linewidth}{@{}lX}
$E$ & The state space of the particles, assumed to be at least a Polish space. \\
$f^N_t$ & The $N$-particle distribution in $\pb(E^N)$ at time $t\geq0$.\\
$f^{k,N}_t$ & The $k$-th marginal of $f^N_t$.\\
$f^N_I$ & The $N$-particle distribution on the path space in $\pb(D(I,E^N))$ or $\pb(C(I,E^N))$ for a time interval $I=[0,T]$. We identify $D(I,E^N)\simeq D(I,E)^N$. \\
$f_t$ & The limit law in $\pb(E)$ at time $t\geq0$.\\
$f_I$ & The limit law on the path space in $\pb(D(I,E))$ or $\pb(C(I,E))$. \\
$F^N_t$ & The law of the empirical process in $\pb(\pb(E))$ at time $t\geq0$.\\
$F^{\mu,N}_I$ & The weak pathwise law of the empirical process in $\pb(D(I,\pb(E)))$ on the time interval $I=[0,T]$. \\
$F^N_I$ & The strong pathwise law of the empirical process in $\pb(\pb(D(I,E)))$ on the time interval $I=[0,T]$. \\
$\mathcal{L}_N$ & The $N$-particle generator acting on (a subset of) $C_b(E^N)$. \\
$\mathcal{L}^N$ & The $N$-particle generator acting on $\pb(E^N)$ defined as the formal adjoint of $\mathcal{L}_N$. \\
$L\diamond_i\varphi_N$ & The action of an operator $L$ on (a subset of) $C_b ( E )$ against the $i$-th variable of a function $\varphi_N$ in $C_b ( E^N )$, defined as the function in (a subset of) $C_b ( E^N )$
$L \diamond_i \varphi_N : (x^1,\ldots,x^N) \mapsto L [ x \mapsto \varphi_N ( x^1 , \ldots , x^{i-1} , x , x^{i+1} , \ldots , x^N ) ] ( x^i )$. The definition readily extends to the case of an operator $L^{(2)}$ acting on $C_b(E^2)$ and two indexes $i<j$ in which case we write $L^{(2)}\diamond_{ij}\varphi_N$. \\
$(\mathcal{X}^N_t)^{}_t$ & The $N$-particle process, with components
$\mathcal{X}^N_t=(X^{1,N}_t,\ldots,X^{N,N}_t)\in E^N$. Often we write $X^{i,N}_t\equiv X^i_t$ and $(\mathcal{X}^N_t)^{}_t\equiv\mathcal{X}^N_{[0,T]}$. \\
$(\mathcal{Z}^N_t)^{}_t$ & An alternative notation for the $N$-particle process with $\mathcal{Z}^N_t=(Z^{1,N}_t,\ldots,Z^{N,N}_t)$. Often used for Boltzmann particle systems or kinetic systems. 
\end{tabularx}

\section{Models and properties}\label[I]{sec:models}

\subsection{Particle systems, setting and conventions}\label[I]{sec:particlesystemsconventions}

The starting point of this review is a system of $N$ particles 
\[\mathcal{X}^N_I \equiv (\mathcal{X}^N_t)^{}_{t\in I} \equiv ( X^{1,N}_t, \ldots , X^{N,N}_t )^{}_{t\in I},\]
where each particle $(X^{i,N}_t)^{}_{t\in I}$ is a stochastic process with values in the state space $E$ which is at least Polish (i.e. separable and completely metrizable) and defined on a time interval $I=[0,T]$ with $T\in(0,+\infty]$. When no confusion is possible, we drop the $N$ superscript and only write $X^i_t\equiv X^{i,N}_t$ for the $i$-th particle. The $N$-particles are not independent; they are said to \emph{interact}. 

Throughout this review, $(\mathcal{X}^N_t)^{}_t$ is a nice stochastic process in $E^N$ which satisfies the strong Markov property and which has c\`adl\`ag sample paths (the related topology is the $J1$ Skorokhod topology, see definition \cref{def:SkorokhodTopoJ1}). Several examples will be given in the next sections but to fix ideas, the particle system will either be a Feller diffusion process in $E=\R^d$ (or in a Borel subset or in a manifold) or a jump process in a more general state space which satisfies the $C_b$-Feller property (i.e. the transition operator is strongly continuous and maps $C_b(E^N)$ to $C_b(E^N)$). One may also consider mixed jump-diffusion processes. The $N$-particle process will often be given as the solution of a Stochastic Differential Equation (SDE) but in full generality, it will be described by its generator denoted by $\mathcal{L}_N$ acting on a domain  $\Dom(\mathcal{L}_N)\subset C_b(E^N)$ or $\Dom(\mathcal{L}_N)\subset C_0(E^N)$ (a detailed description of the analytical setting for Markov processes and their generators is provided in Appendix \cref{appendix:timeinhomoegeneousmarkov}, see in particular Definitions \cref{def:feller}, \cref{def:infgene} and Theorem \cref{thm:charactFullGene}). The generator determines what is called the \emph{interaction mechanism}. The only but crucial assumption that is made on this interaction mechanism is the symmetry or \emph{exchangeability}.

\begin{definition}[Exchangeability]
A family $(X^i)_{i \in I}$ of random variables is said to be \emph{exchangeable} when the law of $(X^i)_{i \in I}$ is invariant under every permutation of a finite number of indexes $i \in I$.
\end{definition}

In a dynamical seting, the \emph{pathwise exchangeability} is assumed in the sense that exchangeability holds for the family of processes $(X^{i,N}_I )^{}_{1 \leq i \leq N}$, at the level of trajectories. Taking the time-coordinate (i.e. the push-forward of the family law by the map $\omega \mapsto \omega(t)$), this implies the \emph{pointwise exchangeability} i.e. exchangeability for the position vector $(\mathcal{X}^{1,N}_t,\ldots,\mathcal{X}^{N,N}_t)$ at any time $t\geq0$. Formally, $\mathcal{X}^N_t$ can be seen as an element of $E^N/\mathfrak{S}_N$, where $\mathfrak{S}_N$ denotes the group of all permutations of $\{1,\ldots,N\}$, although for simplicity we will keep using $E^N$ as the state space. Such a particle system will be called \emph{exchangeable}.

\subsubsection*{The statistical description.} In statistical physics, the previous description in terms of stochastic process is sometimes called the \emph{microscopic scale} because the trajectory of each individual particle is recorded. When $N$ is \emph{large}, the microscopic scale contains too much information and a \emph{statistical description} is sought. There are at least three statistical points of view on the particle system, detailed below. 
\begin{enumerate}
\item The easiest one, is simply given by the $N$-particle distribution $f^N_t\in \pb(E^N)$ at time $t\in I$. From the general theory of Markov processes (see Appendix~\cref{appendix:timeinhomoegeneousmarkov}), $f^N_t$ satisfies the forward Kolmogorov equation written in weak form:
\begin{equation}\label[Ieq]{eq:weakliouville} \forall \varphi_N \in \Dom ( \mathcal{L}_N ), \quad \frac{\dd}{\dd  t} \langle f^N_t , \varphi_N \rangle = \langle f^N_t ,  \mathcal{L}_N \varphi_N  \rangle. \end{equation}
Here and throughout this review, the bracket notation $\langle\cdot,\cdot\rangle$ denotes the integral of a test function (here $\varphi_N$) against a probability measure (here $f^N_t$).  
\begin{remark} Note that 
\[\langle f^N_t, \varphi_N\rangle = \langle f^N_0, u_N(t,\cdot) \rangle,\]
where for $\mathbf{x}^N\in E^N$, $u_N \equiv u_N(t,\mathbf{x}^N) := \E\big[\varphi_N(\mathcal{X}^N_t)|\mathcal{X}^N_0 = \mathbf{x}^N\big]$ solves the backward Kolmogorov equation 
\[\partial_t u_N = \mathcal{L}_Nu_N.\]
The equation \cref{eq:weakliouville} thus describes the dynamics of an \emph{observable} of the system.
\end{remark}
Equation \cref{eq:weakliouville} is also called the \emph{master equation} in a probabilistic context and is better known as the \emph{Liouville equation} in classical (deterministic) kinetic theory. In this review, we follow this latter terminology and Equation \cref{eq:weakliouville} will be called the \emph{(weak) Liouville equation}. The forward Kolmogorov equation, or (strong) Liouville equation, reads \[\partial f^N_t = \mathcal{L}^Nf^N_t,\]
where $\mathcal{L}^N\equiv \mathcal{L}_N^\star$ is the dual operator of $\mathcal{L}_N$. In general, no explicit expression for $\mathcal{L}^N$ is available and it is thus easier to focus on the weak point of view. From the weak Liouville equation, it is possible to compute the time evolution of any \emph{observable} of the particle system, that is of any of the averaged quantities $\langle f^N_t ,  \varphi_N \rangle$ for a test function $\varphi_N$. The drawback is that $f^N_t\in\pb(E^N)$ belongs to a high dimensional space (since $N$ is large). However, by the exchangeability assumption, the law $f^N_t$ is a symmetric probability measure and it is thus possible to define for any $k\in\N$ the $k$-th marginal $f^{k,N}_t\in\pb(E^k)$ by:
\[ \forall \varphi_k \in C_b ( E^k ), \quad \langle f^{k,N}_t , \varphi_k \rangle := \langle f^{N}_t , \varphi_k \otimes 1^{\otimes (N - k)} \rangle. \]
The exchangeability ensures that the term on the right-hand side does not depend on the indexes of the $k$ variables; for instance, it would be equivalent to take $1^{\otimes(N-k)}\otimes\varphi_k$ as a test function instead of $\varphi_k \otimes 1^{\otimes (N - k)}$. Each marginal distribution satisfies a Liouville equation obtained from \cref{eq:weakliouville} by taking $\varphi_k$ as a test function. This equation may not be closed in the sense that, depending on $\mathcal{L}_N$, the right-hand side may depend on $f^N_t$ or on the other marginals. 

\begin{remark}
Note that with a slight abuse we take $\varphi_k$ in the space $C_b(E^k)$ although we should say that $\varphi_k$ belongs to a subset of $C_b(E^k)$ and is such that $\varphi_k\otimes 1^{\otimes(N-k)}$ belongs to $\Dom(\mathcal{L}_N)$. We will often keep doing that in the following. 
\end{remark}

\item From the point of view of stochastic analysis, the particle system $\mathcal{X}^N_I$ can be seen as a random element of $D(I,E^N)\simeq D(I,E)^N$ so its law is a probability measure on the path space, denoted by $f^N_I\in\pb(D(I,E^N))$. This \emph{pathwise law} is generally given as the unique solution of the following martingale problem (probability reminders can be found in Appendix \cref{appendix:probabilityreminders}). 

\begin{definition}[Particle martingale problem]\label[I]{def:martingaleproblemparticles}
Let $T\in(0+\infty]$. A pathwise law $f^N_{[ 0 , T ]}\in\pb\big(D\big([0,T],E^N\big)\big)$ is said to be a solution of the martingale problem associated to the particle system issued from $f^N_0\in\pb(E^N)$ whenever for all $\varphi_N \in \Dom(\mathcal{L}_N)$,
\begin{equation}\label[Ieq]{eq:martingaleparticles}M^{\varphi_N}_t := \varphi_N {\left( \mathbf{X}^N_t \right)} - \varphi_N {\left( \mathbf{X}^N_0 \right)} - \int_0^t \mathcal{L}_N \varphi_N {\left( \mathbf{X}^N_s \right)}\dd s, \end{equation}
is a $f^N_{[ 0 , T ]}$-martingale, where $( \mathbf{X}^N_t )_{t \geq 0}$ denotes the canonical process on the Skorokhod space $D \big( [ 0 , T ] , E^N \big)$ defined for $\omega\in D \big( [ 0 , T ] , E^N \big) $ by $\mathbf{X}^N_t(\omega)=\omega(t)$. 
\end{definition}

Note that the time marginal or \emph{pointwise law} is given by $f^N_t = (\mathbf{X}^N_t)_\# f^N_I$. The weak Liouville equation \cref{eq:weakliouville} can be recovered by taking the expectation in \cref{eq:martingaleparticles}. This description is called \emph{pathwise} and the previous one \emph{pointwise}. 

\item Finally an exchangeable particle system can also be described by its empirical measure: 
\begin{equation}\label[Ieq]{eq:empiricalintro} \mu_{\mathcal{X}^N_t} :=\frac{1}{N}\sum_{i=1}^N \delta_{X^{i,N}_t} \in \pb(E). \end{equation}
Contrary to $f^N_t$, the measure $\mu_{\mathcal{X}^N_t}$ is a \emph{random object}: it can be seen as a measure-valued random variable whose law is the push-forward measure of $f^N_t$ by the application $\boldsymbol{\mu}_N : \mathbf{x}^N \in E^N \mapsto \mu_{\mathbf{x}^N} \in\pb(E)$. Thanks to the exchangeability, the $\pb(E)$-valued process $(\mu_{\mathcal{X}^N_t})_{t\geq0}$ is a Markov process. Moreover the expectation of the empirical measure gives the first marginal of $f^N_t$ :
\[\forall \varphi\in C_b(E), \quad \E \langle \mu_{\mathcal{X}^N_t} , \varphi \rangle = \langle f^{1,N}_t , \varphi \rangle. \]
It follows that from the empirical measure, it is possible to reconstruct the law of any individual particle. As we shall see later, in fact, it characterises the full $N$-particle distribution $f^N_t$. Its pathwise version is the empirical measure on the path space: 
\[\mu_{\mathcal{X}^N_I} := \frac{1}{N} \sum_{i=1}^N \delta_{X^i_I}\in \pb(D(I,E)),\]
where each particle $X^i_I$ is seen as a $D(I,E)$-valued process. 
\end{enumerate}

\subsubsection*{The mesoscopic scale.} 

The main concern of this review is the description of the limit dynamics when $N\to+ \infty$. This will be given by a \emph{nonlinear} object which describes the \emph{average} behaviour of the system. Various points of view may be adopted: from the previous discussion, a natural idea is to study the limit $N\to+ \infty$ of the statistical objects $f^{k,N}_t$, for all $k\in\N$, or $\mu_{\mathcal{X}^N_t}$. The central notion of this review is the \emph{Propagation of Chaos} property which states that for all $t\in I$ there exists $f_t\in\pb(E)$ such that, 
\begin{equation}\label[Ieq]{eq:chaosintro}\forall k\in \N,\quad f^{k,N}_t \underset{N\to+\infty}{\longrightarrow} f_t^{\otimes k},\end{equation}
for the weak convergence of probability measures and provided that this property holds true when $t=0$. As we will see in the following, the property \cref{eq:chaosintro} is equivalent to 
\begin{equation}\label[Ieq]{eq:empiricalchaosintro}\mu_{\mathcal{X}^N_t} \underset{N\to+\infty}{\longrightarrow} f_t,\end{equation}
for the convergence in law (note that the limit is a deterministic object). The term \emph{propagation} introduced by Kac \cite{kac_foundations_1956} refers to the idea that the above convergence at $t=0$ is sufficient to prove the convergence at a later time. As we shall see in the following, the status of the time variable is a central question: one may try to quantify the convergence speed with respect to $t$, analyse its behaviour when $I=[0,+\infty)$ is infinite or study the more general question of the pathwise convergence of the trajectories. All these aspects will define as many notions of propagation of chaos.  

In \cref{eq:chaosintro}, the tensor product indicates that any $k$ particles taken among the $N$ become statistically independent when $N\to+\infty$. Any subsystem of the $N$-particle system thus behaves as a system of i.i.d processes with common law $f_t$ (note that the particles are always identically distributed by the exchangeability assumption). This translates the physical idea that for large systems, the correlations between two (or more) given particles which are due to the interactions become negligible. By looking at the whole system, only an averaged behaviour can be observed instead of the detailed correlated trajectories of each particle. This level of description is called the \emph{mesoscopic scale} in statistical physics. 

The central question is therefore the description of the limit law $f_t$. In turns out that in most cases, it is relatively easy to see that $f_t$ is formally the solution of one of the following nonlinear problems. 
\begin{enumerate}
    \item The solution of a (nonlinear) Partial Differential Equation (PDE) obtained from the Liouville equation by some closure assumption. In some cases, $f_t$ can also be seen as the law of a nonlinear Markov process in the sense of McKean, typically the law of a nonlinear SDE (these notions will be properly defined later). Just as in the classical case, the two approaches are linked by It\={o}'s formula.
    \item The solution of a (nonlinear) martingale problem. This description is more detailed than the previous one, since it gives the probability distribution on the path space $f_I\in \pb(D(I,E))$ with time marginals $(f_t)_t$.
\end{enumerate}
For each of the particle systems considered in this review, the program is thus the following: 
\begin{enumerate}[(1)]
    \item Prove that some limit point \cref{eq:chaosintro} or \cref{eq:empiricalchaosintro} exists in a suitable topology.  
    \item Identify the limit as the solution of a nonlinear problem. 
    \item Prove that the nonlinear problem is wellposed so that the limit is uniquely defined. 
\end{enumerate}
Note that the three steps can be carried out in any order. In this review the main concern will be the first step, which is the core of the propagation of chaos property. This also provides an existence result for the nonlinear problem of the second step. The third step is often proved beforehand. Actually, in many cases, the nonlinear problem has its own dedicated literature; many of its properties are known and may be useful to carry out the first step. 

We conclude this introductory section by a brief overview of the models studied in this work and which will be detailed in the following subsections. 

Section \cref{sec:mckeanintro} is devoted to the description of various diffusion processes, starting from the prototypical example introduced in the seminal work of McKean \cite{mckean_propagation_1969,mckean_class_1966}. At the microscopic scale, the particle system is defined as a system of interacting It\=o processes, where the interaction depends only on (observables of) the empirical measure \cref{eq:empiricalintro}. Physically, it means that each particle interacts with a single averaged quantity in which the other particles contribute with a weight of the order $1/N$. This type of interaction is called \emph{mean-field interaction} and the propagation of chaos is a particular instance of a \emph{mean-field limit}. In section \cref{sec:meanfieldjump} the diffusion interaction is replaced by a jump mechanism. 

Section \cref{sec:boltzmann} presents another class of models which extends the work of Kac \cite{kac_foundations_1956} on the Boltzmann equation. The $N$-particle process is driven by time-discrete pairwise interactions which update the state of only two particles at each time. In classical kinetic theory, the particles are said to \emph{collide}.

\begin{remark}[Notational convention]
We will adopt the following notational convention: the $N$-particle mean-field processes are denoted by the letter $\mathcal{X}$ and the Boltzmann $N$-particle processes by the letter $\mathcal{Z}$. Historically, many of the Boltzmann models that we are going to study are spatially homogeneous versions of a Boltzmann kinetic system. We thus also use the letter $\mathcal{Z}$ for kinetic systems, that is systems where each particle $Z^i_t$ is defined by its position and its velocity, respectively denoted by the letters $X^i_t$ and $V^i_t$. 
\end{remark}

\subsection{Mean-field models}\label[I]{sec:meanfieldmodels}

\subsubsection{Abstract mean-field generators and mean-field limits}

A \emph{mean-field particle system} is a system of $N$ particles characterised by a generator of the form
\begin{equation}\label[Ieq]{eq:Nparticlemeanfieldgenerator}\mathcal{L}_N\varphi_N(\mathbf{x}^N) = \sum_{i=1}^N L_{\mu_{\mathbf{x}^N}}\diamond_i\varphi_N(\mathbf{x}^N),\end{equation}
where given a probability measure $\mu\in\pb(E)$, $L_\mu$ is the generator of a Markov process on $E$. Throughout this review, the notation $L\diamond_i\varphi_N$ denotes the action of an operator $L$ defined on (a subset of) $C_b ( E )$ against the $i$-th variable of a function $\varphi_N\in C_b ( E^N )$; in other words, $L\diamond_i\varphi_N$ is defined as the function:
\[L \diamond_i \varphi_N : (x^1,\ldots,x^N) \in E^N \mapsto L [ x \mapsto \varphi_N ( x^1 , \ldots , x^{i-1} , x , x^{i+1} , \ldots , x^N ) ] ( x^i ) \in\R.\]
There are two main classes of mean-field models, depending on the form of the generator $L_\mu$.
\begin{enumerate}
    \item In Section \cref{sec:mckeanintro}, $L_\mu$ is the generator of a diffusion process, and the associated $N$-particle system is called a \emph{McKean-Vlasov diffusion}. 
    \item In Section \cref{sec:meanfieldjump}, $L_\mu$ is the generator of a jump process and the $N$-particle system is called a \emph{mean-field jump process}. When $L_\mu$ is the sum of a pure jump generator and the generator of a deterministic flow, the process is called a mean-field Piecewise Deterministic Markov Process (PDMP for short). 
\end{enumerate}
It is also possible to consider \emph{mixed processes} when $L_\mu$ is the sum of a diffusion generator and a jump generator. 

It is classically assumed that the domain of the generator $L_\mu$ does not depend on $\mu$. This domain will be denoted by $\mathcal{F}\subset C_b(E)$. 

In that case, it is easy to guess the form of the associated nonlinear system obtained when $N\to+\infty$. Taking a test function of the form 
\[\varphi_N(\mathbf{x}^N) = \varphi(x^1),\]
where $\varphi\in \mathcal{F}$, one obtains the one-particle Kolmogorov equation: 
\[\frac{\dd}{\dd t} \langle f^{1,N}_t,\varphi\rangle = \int_{E^N} L_{\mu_{\mathbf{x}^N}}\varphi(x^1)f^N_t(\dd \mathbf{x}^N).\]
Note that the right-hand side depends on the $N$-particle distribution. As already discussed in the introduction, if the limiting system exists then, its law $f_t$ at time $t\geq0$ is typically obtained as the limit of the empirical measure process: when $\mathcal{X}^N_t\sim f^N_t$,
\[\mu_{\mathcal{X}^N_t} \underset{N\to+\infty}{\longrightarrow} f_t.\]
This also implies $f^{1,N}_t\to f_t$. Reporting formally in the previous equation, it follows that $f_t$ should satisfy 
\begin{equation}\label[Ieq]{eq:weakpdemeanfield}\forall\varphi\in \mathcal{F},\quad\frac{\dd}{\dd t} \langle f_t,\varphi\rangle = \langle f_t, L_{f_t}\varphi\rangle.\end{equation}
This is the weak form of an equation that is called \emph{the (nonlinear) evolution equation}. Note that the evolution equation is nonlinear due to the dependency of $L$ on the measure argument $f_t$. With a slight abuse, we will often simply write $C_b ( E )$ instead of $\mathcal{F}$ in the following. 

This is a very analytical derivation. Its probabilistic counterpart is the following nonlinear martingale problem. 
\begin{definition}[Nonlinear mean-field martingale problem]\label[I]{def:nonlinearmeanfieldmartingaleproblem} Let $T\in(0,+\infty]$ and let us write $I=[0,T]$. A pathwise law $f_I\in\pb(D([0,T],E))$ is said to be a solution of the nonlinear mean-field martingale problem issued from $f_0\in\pb(E)$ whenever for all all $\varphi\in\mathcal{F}$,
\begin{equation}\label[Ieq]{eq:mckeanmartnigaleproblem} M^\varphi_t := \varphi(\mathsf{X}_t)-\varphi(\mathsf{X}_0) - \int_0^t L_{f_s}\varphi(\mathsf{X}_s)\dd s,\end{equation}
is a $f_I$ martingale, where $(\mathsf{X}_t)_t$ is the canonical process and for $t\geq0$, $f_t := (\mathsf{X}_t)_\# f_I$. The natural filtration of the canonical process is denoted by $\mathscr{F}$. 
\end{definition}

Note that $f_I$ contains much more information than the evolution equation \cref{eq:weakpdemeanfield} and as the notation implies, $f_t=(\mathsf{X}_t)_\#f_I\in\pb(E)$ solves the evolution equation. If the nonlinear martingale problem is wellposed then the canonical process $(\mathsf{X}_t)_t$ is a time inhomogeneous Markov process on the probability space $(D([0,T],E),\mathscr{F},f_I)$. This may seem a little bit abstract for now but in what follows, we will see that most often, given the usual abstract probability space $(\Omega,\mathscr{F},\mathbb{P})$, one can define a process on $\Omega$ such that its (pathwise) law is a solution of the nonlinear martingale problem. Such process is called \emph{nonlinear in the sense of McKean} or simply \emph{nonlinear} for short.

\subsubsection{McKean-Vlasov diffusion} \label[I]{sec:mckeanintro}

Let be given two functions 
\begin{equation}\label[Ieq]{eq:mckeandriftsiffusionfunctions}
b : \R^d\times\pb(\R^d)\to \R^d,\quad \sigma : \R^d\times\pb(\R^d)\to\mathcal{M}_d(\R)\end{equation}
respectively called the \emph{drift} vector and the \emph{diffusion} matrix. For a fixed $\mu\in\pb(\R^d)$, the following generator is the generator of a diffusion process in $\R^d$:
\begin{equation}\label[Ieq]{eq:FPgenerator}\forall\varphi\in C^2_b(\R^d),\quad L_\mu\varphi(x) := b(x,\mu)\cdot\nabla\varphi + \frac{1}{2}\sum_{i,j=1}^d a_{ij}(x,\mu)\partial_{x_i}\partial_{x_j}\varphi,\end{equation}
where $a(x,\mu):=\sigma(x,\mu)\sigma(x,\mu)^\mathrm{T}$. The $N$-particle generator \cref{eq:Nparticlemeanfieldgenerator} associated to this class of diffusion generators defines a process called a \emph{McKean-Vlasov diffusion process}. The associated $N$-particle process is governed by the following system of SDEs:
\begin{equation}\label[Ieq]{eq:mckeanvlasov}
\forall i\in\{1,\ldots,N\}, \quad \dd X^{i,N}_t = b \big( X^{i,N}_t,\mu_{\mathcal{X}^N_t} \big) \dd t + \sigma \big( X^{i,N}_t,\mu_{\mathcal{X}_t} \big) \dd B^{i}_t
\end{equation}
where $B^{1}_t,\ldots,B^{N}_t$ are $N$ independent Brownian motions.

\begin{remark}[On the assumptions on $b$ and $\sigma$]
Note that for simplicity, here, we have implicitly assumed that Equation \cref{eq:mckeanvlasov} has a strong solution: it may require precise assumptions on the coefficients $b$ and $\sigma$ (typically Lipschitz assumptions as in Proposition \cref{prop:wellposednessmckeanlipschitz} below). This will also be the case for Equation \cref{eq:mckeanvlasov-pde} below. However we stress out that these assumptions are not satisfied in many important cases. In this introductory section, such strong formulations are only used for better readability but in the following and in particular in the second part of this review, the assumptions on $b$ and $\sigma$ will play a crucial role. Example \cref{example:mckeanvlasov} below present some important models where $b$ and $\sigma$ are singular and for which Equations \cref{eq:mckeanvlasov} and \cref{eq:mckeanvlasov-pde} are not defined in strong form. These models will be thoroughly reviewed in the second part.  
\end{remark}

\begin{remark}
Note that in the system \cref{eq:mckeanvlasov} there are actually $dN$ independent one-dimensional Brownian motions. This remark may be helpful in cases where the Brownian motions in the different directions are different. In particular, for \emph{kinetic} particles defined by their positions and velocities, the noise is often added on the velocity variable only (this case is nevertheless covered by \cref{eq:mckeanvlasov} with a block-diagonal matrix $\sigma$ with a vanishing block on the position variable). 
\end{remark}

In this case, the evolution equation \cref{eq:weakpdemeanfield} can be written in a strong form (at least formally) and reads:
\begin{equation}\label[Ieq]{eq:mckeanvlasov-pde}
\partial_t f_t( x ) = -\nabla_x\cdot\{b( x ,f_t)f_t\}+\frac{1}{2}\sum_{i,j=1}^d \partial_{x_i}\partial_{x_j}\{a_{ij}(x,f_t)f_t \} 
\end{equation}
This is a nonlinear Fokker-Planck equation which is used in many important modelling problems (see Example \cref{example:mckeanvlasov}). This equation was obtained (formally) previously using only the generators when $N\to+\infty$. Here, there is an alternative way to derive the limiting system: looking at the SDE system \cref{eq:mckeanvlasov}, the empirical measure can be formally replaced by its expected limit $f_t$. Since all the particles are exchangeable, this can be done in any of the $N$ equations. The result is a process $(\overline{X}_t)_t$ which solves the SDE:
\begin{equation}\label[Ieq]{eq:mckeanvlasov-limit}
\dd \overline{X}_t = b {\left( \overline{X}_t,f_t \right)}\dd t + \sigma {\left( \overline{X}_t,f_t \right)} \dd B_t.
\end{equation}
where $B_t$ is a Brownian motion and $\overline{X}_0\sim f_0$. Moreover, since for all $i$, $X^i_t$ has law $f^{1,N}_t$ and since it is expected that $f^{1,N}_t\to f_t$, the process $\overline{X}_t$ and the distribution $f_t$ should be linked by the relation 
\[f_t=\mathrm{Law}(\overline{X}_t).\]
The dependency on its law of the solution of a SDE is a special case of what is called a nonlinear process in the sense of McKean. It would now be desirable to prove that the process \cref{eq:mckeanvlasov-limit} is well defined or (equivalently) that the PDE \cref{eq:mckeanvlasov-pde} or the martingale problem \cref{eq:mckeanmartnigaleproblem} are wellposed. The following result gives the reference framework in which all these objects are well defined.  

\begin{proposition}\label[I]{prop:wellposednessmckeanlipschitz}
Let us assume that the functions $b$ and $\sigma$ are globally Lipschitz: there exists $C>0$ such that for all $x,y\in\R^d$ and for all $\mu,\nu\in\pb_2(\R^d)$ it holds that:
\[|b(x,\mu)-b(y,\nu)| + |\sigma(x,\mu)-\sigma(y,\nu)|\leq C\big(|x-y|+W_2(\mu,\nu)\big),\]
where $W_2$ denotes the Wasserstein-2 distance (see Definition \cref{def:wasserstein}). Assume that $f_0\in\pb_2(\R^d)$. Then for any $T>0$ the SDE \cref{eq:mckeanvlasov-limit} has a unique strong solution on $[0,T]$ and consequently, its law is the unique weak solution to the Fokker-Planck equation \cref{eq:mckeanvlasov-pde} and the unique solution to the nonlinear martingale problem \cref{eq:mckeanmartnigaleproblem}.
\end{proposition}

The proof of this proposition is fairly classical. In some special linear cases (see below), it can be found in \cite[Section 3]{mckean_propagation_1969}, \cite[Theorem 1.1]{sznitman_topics_1991} or \cite[Theorem~2.2]{meleard_asymptotic_1996}. For the most general case which includes the above proposition, we refer to \cite[Theorem 1.7]{carmona_lectures_2016}. The proof is based on a fixed point argument that is sketched below. 

\begin{proof} Let us define the map: 
\[\Psi : \pb_2\big(C([0,T],\R^d)\big)\to\pb_2\big(C([0,T],\R^d)\big), \quad m\mapsto \Psi(m),\]
where for $m\in\pb_2\big(C([0,T],\R^d)\big)$, $\Psi(m)$ is the law (on the path space) of the solution $(X^m_t)_{0\leq t\leq T}$ of the following SDE: 
\[\dd X^m_t = b(X^m_t,m_t)\dd t + \sigma(X^m_t,m_t)\dd B_t.\]
Note that the map $t\in[0,T]\mapsto m_t\in\pb_2(\R^d)$ is continuous for the $W_2$-distance where~$m_t$ is the time marginal of $m$. The goal is to prove that $\Psi$ admits a unique fixed point. Let $m,m'\in \pb_2\big(C([0,T],\R^d)\big)$ and let $t\in[0,T]$. Then by the Burkholder-Davis-Gundy inequality (see Proposition \cref{prop:bdg}) and the Lipschitz assumptions on $b$ and $\sigma$, one can prove that there exists a constant $C>0$ such that: 
\begin{multline*}
\E{\left[\sup_{0\leq s\leq t} \big|X^m_s-X^{m'}_s\big|^2\right]}\\\leq CT{\left(\int_0^t \E{\left[\sup_{0\leq r\leq s} \big|X^m_r-X^{m'}_r\big|^2\right]}\dd s + \int_0^t W_2^2(m_s,m_s')\dd s \right)}.
\end{multline*}
A similar computation will be detailed in the proof of Theorem \cref{thm:mckean}. By Gronwall lemma, we obtain that for a constant $C(T)$ it holds that:
\[\E{\left[\sup_{0\leq s\leq t} \big|X^m_s-X^{m'}_s\big|^2\right]}\leq C(T)\int_0^t W_2^2(m_s,m'_s)\dd s.\]
Let us denote by ${W}_{d_{\mathcal{C}_t},2}$ the Wasserstein-2 distance on the space of probability measures on a path space of the form $\mathcal{C}_t = C([0,t],\R^d)$ endowed with the uniform topology (see Definition \cref{def:spaceswasserstein}) for a given $t\in[0,T]$. For any $t\in[0,T]$, we also write $m_{[0,t]}\in\pb_2\big(C([0,t],\R^d)\big)$ for the restriction of $m$ on $[0,t]$. In particular, note that if $s\leq t$ then $W_{d_{\mathcal{C}_s},2}(m_{[0,s]},m'_{[0,s]})\leq W_{d_{\mathcal{C}_t},2}(m_{[0,t]},m'_{[0,t]})$. Then by definition of $\Psi$ and ${W}_{d_{\mathcal{C}_t},2}$, we conclude that: 
\begin{align*}{W}_{d_{\mathcal{C}_t},2}^2\big(\Psi\big(m_{[0,t]}\big),\Psi\big(m_{[0,t]}'\big)\big) &\leq C(T) \int_0^t W_2^2(m_s,m'_s)\dd s \\ &\leq C(T) \int_0^t {W}_{d_{\mathcal{C}_s},2}^2\big(m_{[0,s]},m'_{[0,s]}\big) \dd s.\end{align*}
By iterating this inequality, the $k$-th iterate $\Psi^k$ of $\Psi$ satisfies: 
\begin{align*}
{W}_{d_{\mathcal{C}_T},2}^2\big(\Psi^k(m),\Psi^k(m')\big) &\leq c(T)^k \int_0^T \frac{s^k}{(k-1)!} {W}_{d_{\mathcal{C}_s},2}^2\big(m_{[0,s]},m'_{[0,s]}\big) \dd s \\
&\leq \frac{(c(T)T)^k}{k!} {W}_{d_{\mathcal{C}_T},2}^2(m,m'),
\end{align*}
from which it can be seen that $\Psi^k$ is a contraction and thus admits a unique fixed point for $k$ large enough. 
\end{proof}

\begin{example}\label[I]{example:mckeanvlasov}

Depending on the form of the drift and diffusion coefficients, the McKean-Vlasov diffusion can be used in a wide range of modelling problems. Some examples are gathered below and many other will be given in Section \cref{sec:applications}.

\begin{enumerate}
    \item The first case is obtained when $b$ and $\sigma$ depend linearly on the measure argument. Namely, for $n,m\in\N$, let us consider two functions 
    \[K_1:\R^d\times\R^d\to\R^m,\quad K_2:\R^d\times\R^d\to\R^n,\]
    and let us take
    \[b(x,\mu) = \tilde{b}(x,K_1\star\mu(x)),\quad \sigma(x,\mu)=\tilde{\sigma}(x,K_2\star \mu(x)),\]
    where $\tilde{b}:\R^d\times\R^m\to\R^d$ and $\tilde{\sigma}:\R^d\times\R^m\to\mathcal{M}_d(\R)$. When $K_1,K_2$ and $\tilde{b},\tilde{\sigma}$ are Lipschitz and bounded, the propagation of chaos result is the given by McKean's theorem (Theorem \cref{thm:mckean}). 
    
    In many applications, $\sigma$ is a constant diffusion matrix, $K_1(x,y)\equiv K(y-x)$ for a fixed symmetric radial kernel $K:\R^d\to\R^d$ and $b(x,\mu)=K\star\mu(x)$. The case where $K$ has a singularity is much more delicate (see Section \cref{sec:mckeantowardssingular}, Section \cref{sec:jabin}, Section \cref{sec:holding}) but contains many important cases. For instance, in fluid dynamics, when $K$ is the Biot-Savart kernel $K(x) = x^\perp/|x|^2 $ in dimension $d=2$ (defining $(x_1,x_2)^{\perp} = (-x_2,x_1)$) and $\sigma(x,\mu)\equiv \sqrt{2\sigma}I_2$ for a fixed $\sigma>0$, the limit Fokker-Planck equation  reads:
    \begin{equation}\label[Ieq]{eq:mckeanlinearintro}\partial_t f_t + \nabla\cdot (f_t K\star f_t ) = \sigma\Delta f_t,\end{equation}
    By invariance by translation, the quantity $\omega = f_t - 1$ is the solution of the famous \emph{vorticity equation} which can be shown to be equivalent to the 2D incompressible Navier-Stokes system (see \cite{jabin_quantitative_2018}). 
    
    In biology, still in dimension $d=2$ but with $K(x)=x/|x|^2$, the equation \cref{eq:mckeanlinearintro} is an example of the Patlak-Keller-Segel model for chemotaxis (see \cite{bresch_mean-field_2019}). Other examples in mathematical physics and mathematical biology are presented in Section \cref{sec:classicalpdekinetic}. 
    
    Another important model of this form is the Kuramoto model obtained when $E=\R$ and $K(\theta)=K_0\sin(\theta)$ for a given $K_0>0$. In this case, the particles model the frequencies of a system of \emph{oscillators} which tends to synchronize, see for instance \cite{acebron_kuramoto_2005, goncalves_large_2015, bertini_dynamical_2009, bertini_synchronization_2014} and Section \cref{sec:kuramoto}. 
    
    \item The case of \emph{gradient systems} is a sub-case of the previous one when $\sigma(x,\mu)\equiv\sigma I_d$ for a constant $\sigma>0$ and
    \begin{equation}\label[Ieq]{eq:gradientmckean} 
    b(x,\mu) = -\nabla V(x)-\int_{\R^d}\nabla W(x-y)\mu(\dd y),
    \end{equation}
    where $V,W$ are two symmetric potentials on $\R^d$ respectively called the confinement potential and the interaction potential. The limit Fokker-Planck equation 
    \[\partial f_t = \frac{\sigma^2}{2}\Delta f_t + \nabla\cdot(f_t\nabla(V+W\star f_t)),\]
    is called the granular-media equation and will be studied in Section \cref{sec:gradientsystems}. An important issue is the long-time behaviour of gradient systems which is often studied under convexity assumptions on the potentials (see Section \cref{sec:gradientsystems}, Section \cref{sec:reflectioncoupling}, Section \cref{sec:optimalcoupling}). 
    \item A kinetic particle $Z^{i,N}_t = (X^{i,N}_t,V^{i,N}_t)\in\R^d\times\R^d$ is a particle defined by two arguments, its position $X^{i,N}_t$ and its velocity $V^{i,N}_t$ defined as the time derivative of the position. The evolution of a system of kinetic particles is usually governed by Newton's laws of motion. In a random setting, the typical system of SDEs is thus the following: for $i\in\{1,\ldots,N\}$,
    \[
    \left\{\begin{array}{rcl}
    \dd X^{i,N}_t &=& V^{i,N}_t\dd t \\ 
    \dd V^{i,N}_t &=& F\big(X^{i,N}_t,V^{i,N}_t,\mu_{\mathcal{X}^N_t}\big)\dd t + \sigma\big(X^i_t,V^{i,N}_t,\mu_{\mathcal{X}^N_t}\big)\dd B^{i}_t,
    \end{array}
    \right.
    \]
    where $F:\R^d\times\R^d\times\pb(\R^d)\to\R^d$ and $\sigma:\R^d\times\R^d\times\pb(\R^d)\to\mathcal{M}_d(\R)$. Note that it is often assumed that the force field induced by the interactions between the particles depends only on their positions, which is why we have written 
    \[\mu_{\mathcal{X}^N_t} := \frac{1}{N}\sum_{i=1}^N \delta_{X^{i,N}_t}\in\pb(\R^d)\]
    instead of $\mu_{\mathcal{Z}^N_t}\in\pb(\R^d\times\R^d)$. This special case of the McKean-Vlasov diffusion in $E=\R^{d}\times\R^d$ is also often called a \emph{second order system} by opposition to the \emph{first order systems} when $E=\R^d$. Note that when $\sigma=0$, the limit Equation \cref{eq:mckeanvlasov-pde} is the renowned Vlasov equation which is historically one of the first and most important models in plasma physics and celestial mechanics. In the following, we will nevertheless most often consider stochastic models although some of the results still apply in this deterministic case (in particular the important Theorem \cref{thm:mckean}). For a detailed account of the Vlasov equation in this context, we refer to the review article \cite{jabin_review_2014}. Several examples of stochastic kinetic particle systems are given in Section \cref{sec:flocking} which deals with swarming models. For instance, the (stochastic) Cucker-Smale model \cite{cucker_mathematics_2007, ha_emergence_2009, naldi_particle_2010, cattiaux_stochastic_2018} describes a system of bird-like particles which interact by aligning their velocities to the ones of their neighbours: 
    \[\dd V^i_t = \frac{1}{N}\sum_{j\ne i} K(|X^j_t-X^i_t|)(V^j_t-V^i_t) \dd t +\sigma \dd B^i_t,\]
    where $\sigma\geq0$ is a noise parameter and $K:\R_+\to\R_+$ is a smooth nonnegative function vanishing at infinity which models the vision of the particles. Other classical models of this form include the attraction-repulsion models \cite{dorsogna_self-propelled_2006, carrillo_double_2009} or the stochastic Vicsek models \cite{degond_continuum_2008, degond_phase_2015, giacomin_alignment_2019}.
    
    \item The general case \cref{eq:mckeandriftsiffusionfunctions}, where $b$ and $\sigma$ have a possibly nonlinear dependence on $\mu$ can be extended to even more general cases. A simple extension is the case of time-dependent functions $b$ and $\sigma$. They may also be random themselves and the $x$ and $\mu$ arguments may be replaced respectively by a full trajectory on $C([0,T],\R^d)$ and a pathwise probability distribution on $\pb(C([0,T],\R^d))$. The most general setting is thus:
    \[b:[0,T]\times \Omega\times C([0,T],\R^d)\times \pb(C([0,T],\R^d))\to \R^d.\]
    Such cases may be very difficult to handle but have recently been used in the theory of mean-field games \cite{cardaliaguet_master_2019, lacker_strong_2018, carmona_probabilistic_2018,carmona_probabilistic_2018-1, cardaliaguet_notes_2010}. Under strong Lipschitz assumptions (in the appropriate topology), some very general results can be obtained by a relatively simple adaptation of the proofs valid in the linear case (see Section \cref{sec:mckeangeneralinteractionsextendingmckean}). For more general systems, we will only briefly mention some existing results in Section \cref{sec:chaosviagirsanov} and Section \cref{sec:mckeangeneralinteractionsothers}. 
\end{enumerate}

\end{example}

\begin{remark}[Martingale measures]\label[I]{rem:martingalemeasure} Starting from an arbitrary nonlinear Fokker-Planck with operator \cref{eq:FPgenerator}, one may wonder if it can always be written (at least formally) as the limit of a particle system. The McKean-Vlasov diffusion positively answers when the diffusion matrix in the Fokker-Planck equation is of the form $a(x,\mu) = \sigma(x,\mu)\sigma(x,\mu)^\mathrm{T}$. For more general matrices $a$, the situation is more complicated. For instance the Landau equation would correspond to a matrix of the form $a(x,\mu) = \int_{\R^d} \sigma(x,y)\sigma(x,y)^\mathrm{T}\mu(\dd y)$. In this case, the problem has been studied with a stochastic point of view in \cite{funaki_certain_1984} and later in \cite{meleard_systemes_1988} and \cite{fontbona_measurability_2009} where an explicit approximating particle system is given (see also Section \cref{sec:landau}). The $N$-particle system is characterised as the solution of a system of $N$ SDEs similar to \cref{eq:mckeanvlasov} but where the $N$ Brownian motions are replaced by $N$ \emph{martingale measures} with intensity $\mu_{\mathcal{X}^N_t}(\dd y)\otimes\dd t$. The notion of martingale measure which originates in the Stochastic PDE literature is studied for instance in \cite{el_karoui_martingale_1990}. Except for the cases investigated in the aforementioned works and although it seems to generalise many of the models presented in this review, there is, to the best of our knowledge, no general theory of propagation of chaos for particle systems driven by martingale measures. 
\end{remark}

\subsubsection{Mean-field jump processes and PDMPs} \label[I]{sec:meanfieldjump}

In this section $E$ is any Polish space. Let us be given a family of probability measures called the \emph{jump measures}:
\[P:E\times\pb(E)\to\pb(E),\,\, (x,\mu)\mapsto P_\mu(x,\dd y),\]
and a positive function, called the \emph{jump rate}:
\[\lambda:E\times\pb(E)\to\R_+,\,\, (x,\mu)\mapsto \lambda(x,\mu),\]

For a given $\mu\in\pb(E)$ the following generator is the generator of a \emph{pure jump process}:
\[L_\mu\varphi(x) = \lambda(x,\mu)\int_{E} \{\varphi(y)-\varphi(x)\}P_\mu(x,\dd y).\]
We will also consider the case of a PDMP when 
\begin{equation}\label[Ieq]{eq:mainjumpgenerator}L_\mu\varphi(x) = a\cdot \nabla\varphi(x)+\lambda(x,\mu)\int_{E} \{\varphi(y)-\varphi(x)\}P_\mu(x,\dd y),\end{equation}
where, with a slight abuse of notation, $a\cdot\nabla$ denotes the transport flow associated to a function $a:E\to E$. Using the family of generators \cref{eq:mainjumpgenerator}, the $N$-particle system with mean-field generator \cref{eq:Nparticlemeanfieldgenerator} can be constructed as follows (see \cite[Theorem 13.2.5]{del_moral_stochastic_2017}). 
\begin{itemize}
\item To each particle $i\in\{1,\ldots,N\}$ is attached a Poisson clock with jump rate $\lambda(\mu_{\mathcal{X}^N_t},X^i_t)$. The jump times of particle $i$ are denoted by $(T^i_n)^{}_n$. The jump times can be constructed recursively by taking $T^i_0 =0$ and 
\[T^i_{n+1} = \inf\Big\{t\geq T^i_n, \int_{T^i_n}^t \lambda(\mu_{\mathcal{X}^N_s},X^i_s) \dd s \geq E^i_n\Big\},\]
where the random variables $(E^i_n)_{n\geq0}$ are i.i.d. and distributed according to an exponential law with parameter 1. 
\item Between two jump times, the motion of a particle is purely deterministic: 
\begin{equation}\label[Ieq]{eq:pdmpflow}
\forall n\in\N,\quad\forall t\in[T^i_{n}, T^i_{n+1}),\quad \dd X^i_t = a(X^i_t)\dd t, 
\end{equation}
\item At each time $T^i_n$, a new state is sampled from the jump measure on $E$: 
\begin{equation}\label[Ieq]{eq:newstatepdmp}
X^i_{T^i_n} \sim P_{\mu_{\mathcal{X}^N_{T^{i-}_n}}}{\left(X^i_{T^{i-}_n},\dd y\right)}\in\mathcal{P}(E).
\end{equation}
\end{itemize}

One expects that in the limit $N\to+\infty$, the law $f_t$ of a particle will satisfy the evolution equation \cref{eq:weakpdemeanfield} which, in this case, reads:
\begin{equation}\label[Ieq]{eq:pdepdmp}\frac{\dd}{\dd t}\langle f_t,\varphi\rangle = \langle f_t, a\cdot\nabla_x\varphi\rangle+ \iint_{E\times E} \lambda(x,f_t)\{\varphi(y)-\varphi(x)\}P_{f_t}(x,\dd y)f_t(\dd x),\end{equation}
for all $\varphi\in C_b(E)$. Two important cases are given in the following examples. 

\begin{example}[Nanbu particle system]\label[I]{example:nanbumeanfield}
Let us take $a=0$ and $\lambda=1$ for simplicity. When the jump measure is linear in $\mu$, i.e. is of the form:
\[P_\mu(x,\dd y) = \int_{z\in E}\Gamma^{(1)}(x,z,\dd y)\mu(\dd z),\]
where $\Gamma^{(1)}:E\times E\to \pb(E)$, then the mean-field generator \cref{eq:Nparticlemeanfieldgenerator} describes a $N$-particle system where at each jump, a particle with state $x$ chooses uniformly another particle, say which has a state $z$, and sample a new state according to the law $\Gamma^{(1)}(x,z,\dd y)$. In \cite{graham_stochastic_1997}, this particle system is called a \emph{Nanbu particle system} in honour of Nanbu who introduced a similar system in \cite{nanbu_direct_1980} and used it as an approximation scheme for the Boltzmann equation of rarefied gas dynamics \cref{eq:Boltzmannphysics}. This equation will be described more thoroughly in Section \cref{sec:Boltzmannclassicalmodels}. When $\Gamma^{(1)}$ is an abstract law, the associated mean-field jump particle system generalizes the one introduced by Nanbu and the limit equation is the following \emph{general Boltzmann equation} (written in weak form): 
\[\frac{\dd}{\dd t }\langle f_t,\varphi\rangle = \int_{E\times E\times E}\varphi(y)\Gamma^{(1)}(x,z,\dd y)f_t(\dd x)f_t(\dd z) - \langle f_t,\varphi\rangle,\]
for all $\varphi\in C_b(E)$. A more classical derivation of the general Boltzmann equation will be given in Section \cref{sec:boltzmann} and the subsequent Example \cref{example:nanbuboltzmann} provides an alternative point of view on the Nanbu particle system. 
\end{example}

\begin{example}[BGK type model]\label[I]{example:bgk}
    In kinetic theory the state space is $E=\R^d\times\R^d$ and the $N$ particles are given by $Z^i_t=(X^i_t,V^i_t)$ with $X^i_t$ the position and $V^i_t$ the velocity of particle $i$ at time $t$. Without external force, it is natural to expect that the particles evolve deterministically and continuously between two jumps as
    \[\dd X^i_t = V^i_t \dd t,\quad \dd V^i_t = 0.\]
    Moreover, the post-jump distribution and the jump rate often do not depend specifically on the pre-jump velocity of the jumping particle but only on its position and on the distribution of particles. Thus we take: 
    \[P_\mu((x,v),\dd x',\dd v') = \delta_x(\dd x')\otimes \mathscr{M}_{\mu,x}(v')\dd v',\]
    where given $\mu\in\pb(E)$ and $x\in\R^d$, $\mathscr{M}_{\mu,x}$ is a probability density function. In this case, Equation \cref{eq:weakpdemeanfield} becomes: 
    \[\frac{\dd}{\dd t}\langle f_t,\varphi\rangle = \langle f_t, v\cdot\nabla_x\varphi\rangle+ \iint_{\R^d\times \R^d}\lambda(f_t,x)\{\varphi(x,v')-\varphi(x,v)\}\mathscr{M}_{f_t,x}(v')\dd v'f_t(\dd x,\dd v),\]
    and its strong form reads: 
    \[\partial_t f_t(x,v) +v\cdot\nabla_x f_t(x,v)= \lambda(x,f_t)\Big(\rho_{f_t}(x)\mathscr{M}_{f_t,x}(v)-f_t(x,v)\Big),\]
    where the spatial density of the particles at time $t$ is defined by:
    \[\rho_{f_t}(x) := \int_{\R^d} f_t(x,v)\dd v.\]
    When $\mathscr{M}_{f_t,x}$ is the Maxwellian distribution
    \[\mathscr{M}_{f_t,x}(v) = \frac{\rho_f}{(2\pi T)^{d/2}}\exp{\left(\frac{|v-u|^2}{2T}\right)},\]
    with $(\rho_f u, \rho_f |u|^2 + \rho_f T) = \int_{\R^d}(v,|v|^2)f_t(x,v)\dd v$, then this equation is called the Bhatnagar-Gross-Krook (BGK) equation \cite{bhatnagar_model_1954}. It is used in mathematical physics as a simplified model of rarefied gas dynamics (for a detailed account of the subject, we refer the interested reader to the reviews \cite{bellomo_macroscopic_2004} and \cite{villani_review_2002} or to the book \cite{cercignani_mathematical_1994}). 
\end{example}

In this review, we found useful to distinguish a class of mean-field jump models that we call \emph{parametric models} which are defined by a jump measure of the form 
\[P_\mu(x,\dd y) = \Big(\psi(x,\mu,\cdot)_\#\nu\Big)(\dd y),\]
where $\nu\in\pb(\Theta)$ is a probability measure on a fixed parameter space $\Theta$ and 
\[\psi : E\times\pb(E)\times\Theta\to E.\]
In this case, for all test function $\varphi\in C_b(E)$, 
\[\int_E \varphi(y)P_\mu(x,\dd y) = \int_\Theta \varphi\big(\psi(x,\mu,\theta)\big)\nu(\dd \Theta).\]
The $N$-particle process associated to a parametric model admits a SDE representation using the formalism of Poisson randon measure which is briefly recalled in Appendix \cref{appendix:poisson}. 

\begin{example}[SDE representation for parametric models]\label[I]{example:parametricjump} 
Let us assume for all $\theta\in \Theta$, the function $\psi(\cdot,\cdot,\theta):E\times\pb(E)\to E$ is Lipschitz for the distance on $E$ and the Wasserstein-1 distance on $\pb(E)$, with a Lipschitz constant $L(\theta)>0$ and a function $L\in L^1_\nu(\Theta)$. This (classical) hypothesis will ensure the wellposedness of the SDE representations of both the particle system and its nonlinear limit, see \cite[Section 3.1]{andreis_mckeanvlasov_2018} and \cite[Theorem 1.2 and Theorem 2.1]{graham_mckean-vlasov_1992}. 

To each particle $i\in\{1,\ldots,N\}$ is attached an independent Poisson random measures $\mathcal{N}^i(\dd s,\dd u,\dd\theta)$ on $[0,+\infty)\times[0,+\infty)\times \Theta$ with intensity measure $\dd s\otimes\dd u\otimes\nu(\dd\theta)$ where $\dd t$ and $\dd u$ denote the Lebesgue measure. The $N$ independent random measures $\mathcal{N}^i$ play a comparable role to the $N$ independent Brownian motions which define a McKean-Vlasov diffusion in \cref{eq:mckeanvlasov}. In the present case, the mean-field jump $N$-particle process is the solution of the following system of SDEs driven by the measures $\mathcal{N}^i$ 
\begin{align}\label[Ieq]{eq:meanfieldjumpsde}
    &X^i_t = X^i_0 + \int_0^t a(X^i_s)\dd s \nonumber\\
    &+ \int_{0}^t\int_{0}^{+\infty}\int_\Theta {\left\{\psi{\left(X^i_{s^-},\mu_{\mathcal{X}^N_{s^-}},\theta\right)}-X^{i}_{s^-}\right\}}\1_{\big(0,\lambda{\big(X^i_{s^-},\mu_{\mathcal{X}^N_{s^-}}\big)}\big]}(u)\,\,\mathcal{N}^i(\dd s,\dd u,\dd\theta).
\end{align}

\begin{remark}
Parametric models are reminiscent but different from the so-called \emph{disordered models} where random i.i.d. parameters $\theta^i$ (called \emph{disorder}) are attached to each particle. For these models, various notions of propagation of chaos can be defined, for instance one can prove a classical propagation of chaos result for each realization of the disorder (\emph{quenched} behavior) or the convergence of the joint empirical measure of the processes and the disorder (\emph{quenched averaged} behavior). These problems are of particular interest from the point of view of Large Deviations. The main case of application is the renowned Kuramoto model \cite{acebron_kuramoto_2005, goncalves_large_2015, dai_pra_mckean-vlasov_1996} which will be briefly discussed in Section \cref{sec:kuramoto}. 
\end{remark}

In neurosciences, the variable $X^i_t$ represents the membrane potential of a neuron indexed by $i$ at time $t$ and the Poisson random measures model the interactions between the neurons due to the chemical synapses. A random jump is called a \emph{spike}. In the model introduced by \cite{fournier_toy_2016}, the effect of the spikes is to reset the potential of the membrane to a fixed value, fixed to 0. In Equation \cref{eq:meanfieldjumpsde} this corresponds to the simple case where $X^i_t\in\R_+$ and $\psi\equiv0$. Note that in this particular case there is no need to consider a parameter space $\Theta$ and $\mathcal{N}^i$ is a Poisson random measures on $[0,+\infty)\times[0,+\infty)$ only. Note that in \cite{fournier_toy_2016}, the deterministic drift $a \equiv a(x,\mu)$ also depends on the empirical measure of the system: it models the effect of electrical synapses which tends to relax the membrane potential of the neurons towards the average potential of the system. An additional interaction mechanism is described in the subsequent example. 
\end{example}

\begin{example}[Simultaneous jumps]\label[I]{example:simultaneousjumps}  The neurons models \cite{de_masi_hydrodynamic_2015, fournier_toy_2016, andreis_mckeanvlasov_2018} extend the (parametric) mean-field jump model \cref{eq:meanfieldjumpsde} to allow \emph{simultaneous jumps} at each jump time $T^i_n$. It models the effect that at each spiking event of a neuron $i$, the membrane potential of all the other neurons $j\ne i$ is also increased by a small amplitude.

In the parametric setting with $E=\R^d$ (or more generally when $E$ has a vector space structure), the mean-field jump model with simultaneous jumps is a defined by the following objects. 
\begin{itemize}
    \item The jump rate function:
    \[\lambda : (x,\mu)\in E\times\pb(E) \mapsto \lambda(x,\mu)\in \R_+.\]
    \item A symmetric probability measure $\nu_N$ on the $N$-fold product of the parameter space $\Theta^N$. We also assume that there exists a symmetric probability measure on $\Theta^\mathbb{N}$ such that $\nu_N$ coincides with the projection of $\nu$ on the first $N$ coordinates. This assumption is natural to be able to take the limit $N\to+\infty$. In~\cite{andreis_mckeanvlasov_2018}, the parameter space is $\Theta=[0,1]$.   
    \item The main jump measure 
    \[P_\mu(x,\dd y) = \Big(\psi(x,\mu,\cdot)_\#\nu\Big)(\dd y),\]
    where
    \[\psi : E\times\pb(E)\times\Theta\to E,\quad (x,\mu,\theta)\mapsto x + \alpha(x,\mu,\theta),\]
    and $\alpha:E\times\pb(E)\times\Theta\to E$ is the jump amplitude. 
    \item The collateral jump measures
    \[\widetilde{P}^N_\mu(x,z,\dd y) = \Big(\widetilde{\psi}^N(x,z,\mu,\cdot)_\#\nu_2\Big)(\dd y),\]
    where
    \[\widetilde{\psi}^N : E\times E\times\pb(E)\times \Theta^2\to E,\quad (x,z,\mu,\theta_1,\theta_2)\mapsto x +\frac{ \widetilde{\alpha}(x,z,\mu,\theta_1,\theta_2)}{N},\]
    and $\widetilde{\alpha}:E\times E\times\pb(E)\times\Theta^2\to E$ is the collateral jump amplitude. It satisfies  $\widetilde{\alpha}(x,x,\mu,\theta_1,\theta_2)=0$ for all $x\in E$, $\mu\in\pb(E)$ and $\theta_1,\theta_2\in\Theta$. In \cite{fournier_toy_2016}, the amplitude is fixed $\widetilde{\alpha}(x,z,\mu,\theta_1,\theta_2)\equiv1$ for $x\ne z$.  
\end{itemize}
The $N$-particle process can be defined as before by an algorithmic description. At each time $T^i_n$, a parameter $\theta\sim\nu_N$ is drawn and then the state of particle $i$ is updated by adding the jump amplitude
\[\alpha{\left(X^i_{T^{i-}_n},\mu_{\mathcal{X}^N_{T^{i-}_n}},\theta_i\right)}.\]
But in this case, at time $T^i_n$, all the other particles $j\ne i$ also jumps by the amplitude
\[\frac{\widetilde{\alpha}{\left(X^j_{T^{i-}_n},X^i_{T^{i-}_n},\mu_{\mathcal{X}^N_{T^{i-}_n}},\theta_j,\theta_i\right)}}{N}.\]
When the parameters $\alpha,\widetilde{\alpha}$ satisfy the Lipschitz integrability conditions of \cite[Section 3.1]{andreis_mckeanvlasov_2018}, a SDE representation of the particle system can also be given. As before, let $\mathcal{N}^i(\dd s,\dd u, \dd\theta)$ be a set of $N$ independent Poisson random measures on $[0,+\infty)\times[0,+\infty)\times\Theta^\mathbb{N}$ with intensity $\dd s\otimes\dd u\otimes \nu$, where $\dd s$ and $\dd u$ denote the Lebesgue measure. The SDE respresentation of the $N$-particle system is given by the following system of SDEs: 
\begin{align}\label[Ieq]{eq:simultaneousjumps}
&X^i_t = X^i_0 + \int_0^t a(X^i_s)\dd s \nonumber\\
&+\int_{0}^t\int_{0}^{+\infty}\int_{\Theta^\mathbb{N}} {\left\{\psi{\left(X^i_{s^-},\mu_{\mathcal{X}^N_{s^-}},\theta_i\right)}-X^{i}_{s^-}\right\}}\1_{\big(0,\lambda{\big(X^i_{s^-},\mu_{\mathcal{X}^N_{s^-}}\big)}\big]}(u)\,\,\mathcal{N}^i(\dd s,\dd u,\dd\theta)\nonumber\\
&+ \sum_{j\ne i} \int_{0}^t\int_{0}^{+\infty}\int_{\Theta^\mathbb{N}} {\left\{\widetilde{\psi}^N{\left(X^i_{s^-},X^j_{s^-},\mu_{\mathcal{X}^N_{s^-}},\theta_i,\theta_j\right)}-X^{i}_{s^-}\right\}}\times\nonumber\\
&\phantom{+ \sum_{j\ne i} \int_{0}^t\int_{0}^{+\infty}\int_{\Theta^\mathbb{N}} {\Big\{\widetilde{\psi}^N{\Big(X^j_{s^-},X^i_{s^-},\mu_{\mathcal{X}^N_{s^-}}}}}\times\1_{\big(0,\lambda{\big(X^j_{s^-},\mu_{\mathcal{X}^N_{s^-}}\big)}\big]}(u)\,\,\mathcal{N}^j(\dd s,\dd u,\dd\theta).
\end{align}
Compared to the previous framework, in addition to the main jump operator \cref{eq:mainjumpgenerator}, each particle is also subject to the \emph{collateral jump generator} defined for all $\varphi\in~C_b(E)$ and $x\in E$ by:
\begin{equation}\label[Ieq]{eq:simultaneousjumpsoperatorN}
\widetilde{L}^N_\mu\varphi(x) := N\iint_{E\times E}\lambda(z,\mu)\{\varphi(y)-\varphi(x)\}\widetilde{P}_\mu^N(x,z,\dd y)\mu(\dd z).
\end{equation}
Note that this generator depends on $N$ but it satisfies the weak limit: for all $\varphi\in~C^1_b(E)$, $x\in E$ and $\mu\in\pb(E)$, 
\[\widetilde{L}^N_\mu\varphi(x)\underset{N\to+\infty}{\longrightarrow}\widetilde{L}_\mu\varphi(x) := \iint_{E\times \Theta^2}\lambda(z,\mu)\widetilde{\alpha}(x,z,\mu,\theta_1,\theta_2)\cdot\nabla\varphi(x)\mu(\dd z)\nu_2(\dd\theta_1,\dd\theta_2).\]
The $N$-particle system is thus defined by the mean-field generator \cref{eq:Nparticlemeanfieldgenerator} which takes the form:
\[\forall\varphi_N\in C_b(E^N),\quad\mathcal{L}_N\varphi_N := \sum_{i=1}^N \Big\{L_\mu\diamond_i \varphi_N+\widetilde{L}^N_\mu\diamond_i \varphi_N\Big\}.\]
In the limit $N\to+\infty$, the nonlinear evolution equation \cref{eq:weakpdemeanfield} is expected to become:
\[\forall\varphi\in C^1_b(E),\quad\frac{\dd}{\dd t}\langle f_t,\varphi\rangle = \langle f_t,L_{f_t}\varphi\rangle + \langle f_t,\widetilde{L}_{f_t}\varphi\rangle.\] 
Note that this is the equation satisfied by the law of the solution of the following nonlinear SDE:
\begin{align}\label[Ieq]{eq:nonlinearsdesimultaneousjumps}
\overline{X}_t = \overline{X}_0 &+ \int_0^t a(\overline{X}_s)\dd s\nonumber\\
&+\int_{0}^t\int_{0}^{+\infty}\int_{\Theta^\mathbb{N}} {\alpha(\overline{X}_{s^-},f_s,\theta_1)}\1_{\big(0,\lambda{\big(\overline{X}_{s^-},f_s\big)}\big]}(u)\,\,\mathcal{N}(\dd s,\dd u,\dd\theta)\nonumber\\
&+\int_0^t \iint_{E\times\Theta^2} \lambda(z,f_s)\widetilde{\alpha}(\overline{X}_s,z,f_s,\theta_1,\theta_2)f_s(\dd z)\nu_2(\dd\theta_1,\dd\theta_2)\dd s,
\end{align}
with $\mathrm{Law}(\overline{X}_s)=f_s$. In the last equation, $\mathcal{N}(\dd s,\dd u,\dd\theta)$ is a Poisson random measure on $[0,+\infty)\times[0,+\infty)\times\Theta^\mathbb{N}$ with intensity $\dd s\otimes \dd u \otimes\nu(\dd\theta)$ where $\dd s$ and $\dd u$ denote the Lebesgue measure. 
\end{example}

Mean-field jump processes and PDMPs are not so common in the literature compared to the McKean-Vlasov diffusion models or the Boltzmann models (Section~\cref{sec:boltzmann}). The Nanbu particle system serves as a simplified Boltzmann model (see Example \cref{example:nanbuboltzmann}). Mean-field jump processes can also be used as an approximation of a McKean-Vlasov diffusion. Since the dynamics only relies on a sampling mechanism on the state space $E$, compared to diffusion processes, it allows more flexibility and avoids some technicalities for instance when $E$ has a more complex geometrical structure, typically when $E$ is a manifold. In applications, mean-field jump processes model a motion called \emph{run and tumble} which is classical in the study of the dynamics of populations of bacteria. As already mentioned, Example \cref{example:simultaneousjumps} corresponds to a toy example of neuron model. More realistic examples often consider a combination of (simultaneous) jumps and a diffusive behaviour, see \cite{andreis_mckeanvlasov_2018} and the references therein. The nonlinear martingale problem associated to mixed jump-diffusion models is studied in \cite{graham_nonlinear_1992} where the wellposedness is proved under classical Lipschitz and boundedness assumptions on the parameters (see also Section \cref{sec:mckeancompactnessreview}). Other neuron models based on a mean-field jump process will be described in Section \cref{sec:neurons}.

\subsection{Boltzmann models}\label[I]{sec:boltzmann}

\subsubsection{General form}\label[I]{sec:boltzmanngeneral}

Given a Polish space $E$, a Boltzmann model is a $N$-particle system with an infinitesimal generator acting on $\varphi_N\in C_b(E^N)$ of the form:
\begin{equation}\label[Ieq]{eq:boltzmanngenerator}
\mathcal{L}_N\varphi_N = \sum_{i=1}^N L^{(1)}\diamond_i\varphi_N+\frac{1}{N}\sum_{i< j}L^{(2)}\diamond_{ij}\varphi_N,
\end{equation}
where $\varphi_N\equiv\varphi_N(z^1,\ldots,z^N)$ is a test function on the product space $E^N$. The operator $L^{(2)}$ acts on two-variable test functions and stands for binary interactions between particles. The operator $L^{(1)}$ acts on one-variable test functions and describes the individual flow of each particle (and possibly the boundary conditions). More explicitly, let us recall the notations, for $(z^1,\ldots,z^N)\in E^N$ and $i<j$, 
\begin{equation*}
    L^{(1)}\diamond_i \varphi_N(z^1,\ldots,z^n)  =  L^{(1)} \big[u \mapsto \varphi_N(z^1,\ldots,z^{i-1},u,z^{i+1},\ldots, z^N)\big](z^i)
\end{equation*}
and 
\begin{multline*}
    L^{(2)}\diamond_{ij} \varphi_N(z^1,\ldots,z^n) \\ =L^{(2)}\big[(u,v)\mapsto \varphi_N(z^1,\ldots,z^{i-1},u,z^{i+1},\ldots,z^{j-1},v,z^{j+1},\ldots, z^N)\big](z^i,z^j).
\end{multline*}
These models are called \emph{Boltzmann models} in reference to the famous Boltzmann equation of rarefied gas dynamics which is a fundamental equation for mathematicians, physicists and philosophers. It will be explained at the end of this section (see Equation \cref{eq:Boltzmannphysics}) how it can be obtained as the limit of a general particle system of the form \cref{eq:boltzmanngenerator}. The specificity of Boltzmann models is that the particles interact only at random times by pair and not individually with an average of all the other particles as in mean-field models. In full generality, the state space $E$ is an abstract space. In classical kinetic theory, $E=\R^d\times\R^d$ is the phase space of positions and velocities and two particles interact when they are close enough: they are said to \emph{collide} and by analogy, we will keep this terminology to refer to an interaction between two particles even in an abstract space. In addition to these pairwise interactions, each particle is also subject to an individual flow prescribed by the operator $L^{(1)}$. Typical examples in kinetic theory include
\begin{itemize}
    \item (Free transport) $L^{(1)}\varphi(x,v) = v\cdot\nabla_x\varphi$,
    \item (Space diffusion) $L^{(1)}\varphi(x,v) = \Delta_x\varphi$.
    \item (Velocity diffusion) $L^{(1)}\varphi(x,v) = \Delta_v\varphi$.
\end{itemize} 
When two particles collide, the effect of the collision is prescribed by the operator $L^{(2)}$. In kinetic theory, this operator acts on the velocity variable only but in full generality, in an abstract space $E$, it will be assumed to satisfy the following assumptions. 



\begin{assumption}\label[I]{assum:L2} The operator $L^{(2)}$ satisfies the following properties. 
\begin{enumerate}[(1)]
    \item The domain of the operator $L^{(2)}$ is $C_b(E^2)$.
    \item There exist a continuous map called the \emph{post-collisional distribution}
    \[\Gamma^{(2)}:(z_1,z_2)\in E\times E \mapsto \Gamma^{(2)}(z_1,z_2,\dd z_1',\dd z_2')\in\pb(E\times E),\]
    and a symmetric function called the \emph{collision rate}
    \[\lambda : (z_1,z_2)\in E\times E\mapsto \lambda(z_1,z_2)\in \R_+,\]
    such that for all $\varphi_2\in C_b(E^2)$ and all $z_1,z_2\in E$,
    \begin{equation}\label[Ieq]{eq:generalL2}L^{(2)}\varphi_2(z_1,z_2) = \lambda(z_1,z_2)\iint_{E\times E} \{\varphi_2(z_1',z_2')-\varphi_2(z_1,z_2)\}\Gamma^{(2)}(z_1,z_2,\dd z_1',\dd z_2').\end{equation}
    \item For all $z_1,z_2\in E$, the post-collisional distribution is symmetric in the sense that
    \begin{equation}\label[Ieq]{eq:symgamma2}\Gamma^{(2)}(z_1,z_2,\dd z_1',\dd z_2') = \Gamma^{(2)}(z_2,z_1,\dd z_2',\dd z_1').\end{equation}
    \item The function $\lambda$ is measurable on $\{(z_1,z_2)\in E^2,\,\,z_1\ne z_2\}$ and for all $z\in E$, $\lambda(z,z)=0$.
\end{enumerate}
\end{assumption}

\begin{remark}
The assumption that $\lambda$ is a (measurable) function prevents from considering the true classical Boltzmann inhomogeneous case in kinetic theory $\lambda( z_1 , z_2) = \delta_{x_1 = x_2}$ (that is, two particles collide when they are exactly at the same position), which is beyond the scope of this review (see however Example \cref{example:hardsphere}). The collision rate is often assumed to be uniformly bounded 
\begin{equation*}\sup_{z_1,z_2\in E}\lambda(z_1,z_2)\leq \Lambda <\infty.\end{equation*}
This \emph{cutoff} assumption is unfortunately not physically relevant for many models where an infinite number of collisions may happen in finite time.
\end{remark}

Note that when $E$ is a locally compact Polish space, the Riesz-Markov-Kakutani theorem states that any linear operator on the space of two-variable test functions in $C_c(E\times E)$ can be written in the form \cref{eq:generalL2}. The third assumption ensures that the law $f^N_t$ defined by the backward Kolmogorov equation remains symmetric for all time provided that $f^N_0$ is symmetric. This follows from the observation that under \cref{eq:symgamma2}, the action of any transposition $\tau\in \mathfrak{S}_N$ on $C_b(E^N)$ commutes with $\mathcal{L}_N$: 
\[\tau^{-1}\mathcal{L}_N\tau = \mathcal{L}_N.\]

\begin{example}[Jump amplitude]
    When $E=\R^d$ (or more generally when $E$ has a vector space structure), the interaction law $\Gamma^{(2)}$ is often given in terms of jump amplitudes. Given the law $\widehat{\Gamma}^{(2)}$ of the jump amplitudes of the form: 
    \[\widehat{\Gamma}^{(2)}:(z_1,z_2)\in \R^d\times\R^d\mapsto\widehat{\Gamma}^{(2)}(z_1,z_2,\dd h, \dd k)\in \pb(\R^d\times\R^d),\]
    the post-collisional law $\Gamma^{(2)}$ is the image measure of $\widehat{\Gamma}^{(2)}$ by the translation 
    \[(h,k)\in \R^d\times\R^d\mapsto (z_1+h,z_2+k)\in\R^d\times\R^d,\]
    so that 
    \begin{multline*}\iint_{\R^2\times \R^d} \varphi_2(z_1',z_2')\Gamma^{(2)}(z_1,z_2,\dd z_1',\dd z_2') \\= \iint_{\R^d\times \R^d} \varphi_2(z_1+h,z_2+k)\widehat{\Gamma}^{(2)}(z_1,z_2,\dd h,\dd k).\end{multline*}
    This is the case investigated in \cite{meleard_asymptotic_1996,graham_stochastic_1997}.
\end{example}

The $N$-particle generator can be interpreted as a generator of a Piecewise Deterministic Markov Process in the space $E^N$. Consequently, and similarly to the case of a mean-field PDMP process (Section \cref{sec:meanfieldjump}), following \cite[Chapters 12 and 13, Theorem 13.2.5]{del_moral_stochastic_2017}, a simple $N$-particle process with generator $\mathcal{L}_N$ of the form \cref{eq:boltzmanngenerator} is given in the following proposition. 

\begin{proposition}\label[I]{prop:prototypeBoltzmannprocess} Let $\mathcal{Z}^N_t = (Z^1_t,\ldots,Z^N_t)$ be the $N$-particle process defined by the three following rules.
\begin{enumerate}[(i)]  
\item For each (unordered) pair of particles $(i,j)$, consider an independent non homogeneous Poisson process with rate $\lambda(Z^i_t,Z^j_t)/N$. The jump times $(T^{ij}_n)^{}_n$ can be constructed recursively by taking $T^{ij}_0=0$ and 
\[T^{ij}_{n+1} = \inf \Big\{t \geq T^{ij}_{n},\,\,\frac{1}{N}\int_{T^{ij}_n}^t {\lambda(Z^i_s,Z^j_s)} \dd s \geq E^{ij}_n\Big\},\]
where $(E^{ij}_n)_n$ are i.i.d. random variables which follow an exponential law with parameter 1. 
\item Between two jump times, the particles evolve independently according to $L^{(1)}$.
\item At each jump time $T^{ij}_n$ of a pair $(i,j)$, update the states of the particles by:
\begin{equation}\label[Ieq]{eq:newstateboltzmann}\Big(Z^i_{T^{ij}_n},Z^j_{T^{ij}_n}\Big)\sim \Gamma^{(2)}\Big(Z^i_{(T^{ij}_n)^-},Z^j_{(T^{ij}_n)^-},\dd z_1',\dd z_2'\Big).\end{equation}
\end{enumerate}Then the generator of $(\mathcal{Z}^N_t)^{}_t$ is $\mathcal{L}_N$ given by \cref{eq:boltzmanngenerator} under Assumption \cref{assum:L2}.\end{proposition}

Note that if $\lambda(z_1,z_2)$ remains of order 1, the factor $1/N$ on the right-hand side of \cref{eq:boltzmanngenerator} ensures that any particle undergoes on average $\mathcal{O}(1)$ collisions per unit of time, which is crucial to take the limit $N\to+\infty$. Let us now describe what this limit may look like. Ultimately, the goal is to describe the limiting behaviour of the one-particle distribution function $f^{1,N}_t$. More generally, taking a test function of the form 
\[\varphi_N = \varphi_s\otimes 1^{N-s},\]
for $s<N$ and $\varphi_s\in C_b(E^s)$, the weak Liouville equation \cref{eq:weakliouville} becomes:
\begin{multline*}\frac{\dd}{\dd t} \langle f^{s,N}_t,\varphi_s\rangle = \sum_{i=1}^s \langle f^{s,N}_t, L^{(1)}\diamond_i\varphi_s\rangle +\frac{1}{N}\sum_{1\leq i<j\leq s}\langle f^{s,N}_t, L^{(2)}\diamond_{ij}\varphi_s\rangle\\ + \frac{N-s}{N}\sum_{i=1}^s\langle f^{s+1,N}_t, L^{(2)}\diamond_{i,s+1}(\varphi_s\otimes 1)\rangle.\end{multline*}
This equation is not closed and involves the $(s+1)$-marginal. This hierarchy of equations is called the BBGKY hierarchy (see Section \cref{sec:finitesystems}). The nonlinear model associated to the Boltzmann particle system is obtained by taking the closure:
\begin{equation}\label[Ieq]{eq:chaosassumptionintro}\forall t\geq0,\,\exists f_t\in\pb(E),\,\forall s\in\N,\quad f_t^{s,N}\underset{N\to+\infty}{\longrightarrow}f_t^{\otimes s},\end{equation}
which is called the \emph{chaos assumption}. The fundamental question in this review is to justify when this property holds. If the chaos assumption holds, then taking $s=1$ in the Liouville equation shows formally that the one-particle distribution converges towards the weak measure solution $f$ of: 
\[\frac{\dd}{\dd t}\langle f_t,\varphi\rangle = \langle f_t, L^{(1)}\varphi\rangle+\langle f_t^{\otimes 2},L^{(2)}(\varphi\otimes 1)\rangle,\]
which is called the \emph{general Boltzmann equation}. Using Assumption \cref{assum:L2}, this equation can be rewritten 
\begin{multline}\label[Ieq]{eq:Boltzmannequationgeneral}
\frac{\dd}{\dd t}\langle f_t,\varphi\rangle = \langle f_t, L^{(1)}\varphi\rangle\\
+\int_{E^3}\lambda(z_1,z_2)\big\{\varphi(z_1')-\varphi(z_1)\big\}\Gamma^{(2)}(z_1,z_2,\dd z_1',E)f_t(\dd z_1)f_t(\dd z_2),
\end{multline}
or in a more symmetric form, using \cref{eq:symgamma2}:
\begin{multline}\label[Ieq]{eq:symmetricBoltzmannequationgeneral}
\frac{\dd}{\dd t}\langle f_t,\varphi\rangle = \langle f_t, L^{(1)}\varphi\rangle\\
+\frac{1}{2}\int_{E^4}\lambda(z_1,z_2)\big\{\varphi(z_1')+\varphi(z_2')-\varphi(z_1)-\varphi(z_2)\big\}\Gamma^{(2)}(z_1,z_2,\dd z_1',\dd z_2')f_t(\dd z_1)f_t(\dd z_2)
\end{multline}

All the Boltzmann type equations in this review are special instances of this general equation for a specific choice of $\lambda$ and $\Gamma^{(2)}$. Note that the general Boltzmann equation \cref{eq:Boltzmannequationgeneral} is written in weak form. Examples of $\Gamma^{(2)}$ which lead to more classical Boltzmann type equations used in the modelling of rarefied gas dynamics and written in strong form are given in Section \cref{sec:Boltzmannclassicalmodels}. Here, we can only formally write the dual version of \cref{eq:Boltzmannequationgeneral}:
\[\partial_t f_t = L^{(1)\star} f_t + Q(f_t,f_t),\]
where $Q$ is called the \emph{collision operator} which is a (symmetric) quadratic operator on $\pb(E)\times\pb(E)\to\mathcal{M}(E)$, defined weakly, for $\varphi\in C_b(E)$, by: 
\begin{multline*}\langle Q(\mu_1,\mu_2),\varphi\rangle =\\ \frac{1}{2}\int_{E^4}\lambda(z_1,z_2)\big\{\varphi(z_1')+\varphi(z_2')-\varphi(z_1)-\varphi(z_2)\big\}\Gamma^{(2)}(z_1,z_2,\dd z_1',\dd z_2')\mu_1(\dd z_1)\mu_2(\dd z_2).\end{multline*}

\begin{example}[Nanbu particle system, continuation of Example \cref{example:nanbumeanfield}]\label[I]{example:nanbuboltzmann}
The general Boltzmann equation \cref{eq:Boltzmannequationgeneral} only depends on the marginals of $\Gamma^{(2)}$. In other words, the detail of the interaction mechanism at the particle level is lost in the limit. As a consequence, one can construct different mechanisms which lead to the same Boltzmann equation. For instance, let the marginal of $\Gamma^{(2)}$ be denoted by: 
\[\forall(z_1,z_2)\in E^2,\quad\Gamma^{(1)}(z_1,z_2,\dd z_1') := \Gamma^{(2)}(z_1,z_2,\dd z_1',E).\]
Let us consider the new post-collisional law: 
\begin{multline*}\tilde{\Gamma}^{(2)}(z_1,z_2,\dd z_1',\dd z_2') \\= \frac{1}{2}\left(\Gamma^{(1)}(z_1,z_2,\dd z_1')\otimes\delta_{z_2}(\dd z_2')+ \Gamma^{(1)}(z_2,z_1,\dd z_2')\otimes\delta_{z_1}(\dd z_1')\right),\end{multline*}
and let us denote by $\tilde{\mathcal{L}}_N$ the new corresponding $N$-particle generator (with $L^{(1)}$ unchanged). This is the generator associated to a particle system such that when a collision occurs, only one particle among the two updates its state (according to the law $\Gamma^{(1)}$) while the state of the other particle remains unchanged. Such mechanism is called a Nanbu interaction mechansim following the terminology of \cite{graham_stochastic_1997, nanbu_direct_1980}. Nevertheless, one can check that the Boltzmann equation associated to this process is exactly \cref{eq:Boltzmannequationgeneral} with an interaction rate $\lambda$ replaced by $\lambda/2$. In the limit $N\to+\infty$, one cannot distinguish this system from the system where both the particles update their states after a collision. 

Note that as explained in Example \cref{example:nanbumeanfield} the Nanbu particle system is also a special case of a mean-field jump process (see Section \cref{sec:meanfieldjump}) with: 
\[\lambda(z,\mu) = \int_{E}\lambda(z,z')\mu(\dd z') \equiv (\lambda\star\mu)(z),\]
and 
\[P_\mu(z,\dd z') = \frac{\int_{z''\in E} \lambda(z,z'')\Gamma^{(1)}(z,z'',\dd z')\mu(\dd z'')}{\int_{z''\in E}\lambda(z,z'')\mu(\dd z'')}.\]
\end{example}

The two following examples are two variants of the Boltzmann model. 

\begin{example}[External clock] The authors of \cite{carlen_kinetic_2013} consider a model where the time between two collisions is given by a Poisson process with fixed rate $\Lambda N$, $\Lambda>0$, independently of the particles. When a collision occurs, the pair $(i,j)$ of particles which interact is chosen among all the pairs of particles with probability $p_{i,j}(\mathcal{Z}^N_t)$, normalised so that for all $\mathcal{Z}^N\in E^N$
\[\sum_{i<j}p_{i,j}(\mathcal{Z}^N) = 1.\]
In this case, 
\[\mathcal{L}_N\varphi_N(\mathcal{Z}^N) = \Lambda N\sum_{i<j} p_{i,j}(\mathcal{Z}^N) {L}^{(2)}\diamond_{ij}\varphi_N(\mathcal{Z}^N).\]
The situation differs from the previous case where a collision rate is attached to each pair of particles and no normalisation constraint is imposed. Dropping the normalisation constraint and taking $p_{i,j}(\mathcal{Z}^N) = \lambda(Z^i,Z^j)/\Lambda$ would give exactly~\cref{eq:boltzmanngenerator}. In the case 
\[\forall i<j,\quad p_{i,j}(\mathcal{Z}^N) = \frac{2}{N(N-1)},\]
that is when all the pairs are chosen with the same probability, then the situation is equivalent to the previous case with all the collision rates equal to the constant $\Lambda$. In \cite{carlen_kinetic_2013}, propagation of chaos is proved in this case only. It is believable that propagation of chaos does not hold when interactions are driven by a clock independent of the particles. 
\end{example}

\begin{example}[Non cutoff models]\label[I]{example:noncutoffmodels}
In this review, we made the choice to distinguish the jump rate $\lambda$ and the post-collisional distribution $\Gamma^{(2)}$. An alternative choice in the literature (see for instance \cite{meleard_asymptotic_1996}) is to consider $\lambda\equiv1$ and a \emph{collision kernel} 
\[\Gamma^{(2)}:(z_1,z_2)\in E\times E\mapsto \Gamma^{(2)}(z_1,z_2,\dd z_1',\dd z_2')\in \mathcal{M}^+(E\times E),\]
which is a positive measure but not necessarily a probability distribution. The collision rate $\lambda$ is thus directly encoded in the total mass of the collision kernel. Two cases may happen, for given $z_1,z_2\in E$, either
\[\Gamma^{(2)}(z_1,z_2,E,E)<+\infty\]
or
\[\Gamma^{(2)}(z_1,z_2,E,E)=+\infty.\]
The first case is called the cutoff case. This case is equivalent to the previous case~\cref{eq:boltzmanngenerator} with Assumption \cref{assum:L2} and the following post-collisional law and collision rate: 
\[\tilde{\Gamma}^{(2)}(z_1,z_2,\dd z_1',\dd z_2')= \frac{\Gamma^{(2)}(z_1,z_2,\dd z_1',\dd z_2')}{\Gamma^{(2)}(z_1,z_2,E,E)},\quad \tilde{\lambda}(z_1,z_2) = \Gamma^{(2)}(z_1,z_2,E,E).\]
In the second case, called the non-cutoff case, the lack of integrability means that there are an infinite number of collisions in finite time. Such system therefore cannot be simulated by a particle system as in Proposition \cref{prop:prototypeBoltzmannprocess}. Nevertheless it still makes sense to consider the abstract Markov process defined by the generator $\mathcal{L}_N$. Non-cutoff models are historically important as explained in Section \cref{sec:Boltzmannclassicalmodels}. Non-cutoff models are often handled by approximating them by cutoff models. In this review we implicitly consider cutoff models but we will occasionally specify when a technique can be extended to non-cutoff cases.  
\end{example}

The nonlinear limit can also be defined as the solution of a more general martingale problem.

\begin{definition}[Nonlinear Boltzmann martingale problem]\label[I]{def:nonlinearboltzmannmartingaleproblem} Let $T>0$ and $f_0\in\pb(E)$. We write $I=[0,T]$. We say that $f_I\in\pb(D([0,T],E))$ is a solution to the nonlinear Boltzmann martingale problem with initial law $f_0$ when for any test function $\varphi\in \mathrm{Dom}(L^{(1)})$, 
    \[M^\varphi_t = \varphi(\mathsf{Z}_t)-\varphi(\mathsf{Z}_0) - \int_0^t \{L^{(1)}\varphi(\mathsf{Z}_s)+K\varphi(\mathsf{Z}_s,f_s)\}\dd s,\]
    is a $f_I$-martingale, where $(\mathsf{Z}_t)_t$ is the canonical process, $f_s=({Z_s})_{\#} f_I$ and for $\mu\in \pb(E)$ and $z_1\in E$,
    \[K\varphi(z_1,\mu) := \iint_{E\times E}\lambda(z_1,z_2) \{\varphi(z_1')-\varphi(z_1)\}\Gamma^{(1)}(z_1,z_2,\dd z_1')\mu(\dd z_2).\]
\end{definition}

Existence and uniqueness for the nonlinear Boltzmann martingale problem holds under classical Lipschitz and boundedness assumptions on the parameters. It is a special case of the model studied in \cite{graham_nonlinear_1992}. Note that this martingale problem is also a special case of the nonlinear mean-field martingale problem (Definition \cref{def:nonlinearmeanfieldmartingaleproblem}) with 
\[L_\mu\varphi(z) = L^{(1)}\varphi(z) + K\varphi(z,\mu).\]
This translates the fact that the Boltzmann equation is obtained as the limit of both the Boltzmann model and the Nanbu particle system which is a special case of mean-field jump process.

We end this section with a classical useful proposition which states that when the collision rate $\lambda$ is uniformly bounded, then the situation is essentially the same as when it is constant.

\begin{proposition}[Uniform clock trick]\label[I]{prop:acceptreject} Assume that 
\begin{equation}\label[Ieq]{eq:uniformboundlambda}\sup_{z_1,z_2\in E}\lambda(z_1,z_2)\leq \Lambda <\infty,\end{equation}
and let $(\tilde{\mathcal{Z}}^N_t)^{}_t$ be the process defined by the three following rules.
\begin{enumerate}[(i)]
\item To each pair of particles is attached an independent Poisson process with rate $\Lambda/N$.
\item Between two jump times, the particles evolve independently according to $L^{(1)}$.
\item When the clock of the pair $(i,j)$ rings at time $T_{ij}$, then the states of the particles is updated with probability $\lambda(\tilde{Z}^i_t,\tilde{Z}^j_t)/\Lambda$ by: 
\[\Big(\tilde{Z}^i_{T_{ij}^+},\tilde{Z}^j_{T_{ij}^+}\Big)\sim \Gamma^{(2)}\Big(\tilde{Z}^i_{T_{ij}^-},\tilde{Z}^j_{T_{ij}^-},\dd z_1',\dd z_2'\Big),\]
and with probability $(1-\lambda(\tilde{Z}^i_t,\tilde{Z}^j_t)/\Lambda)$, nothing happens (this case is called a \emph{fictitious collision}). 
\end{enumerate}
Then the law of $(\tilde{\mathcal{Z}}^N_t)_t^{}$ is equal to the law of the process constructed in Proposition~\cref{prop:prototypeBoltzmannprocess}.
\end{proposition}

\begin{proof}
Let us compute the generator $\tilde{\mathcal{L}}_N$ of the process $(\tilde{\mathcal{Z}}^N_t)^{}_t$. It holds that
\[\tilde{\mathcal{L}}_N\varphi_N = \sum_{i=1}^N L^{(1)}\diamond_i\varphi_N+\frac{1}{N}\sum_{i<j} \tilde{L}^{(2)}\diamond_{ij}\varphi_N,\]
with, given $\varphi_2\in C_b(E^2)$, 
\begin{align*}\tilde{L}^{(2)}\varphi_2(z_1,z_2) &= \Lambda
\int_0^1 \Big\{ \1_{\eta\leq\frac{\lambda(z_1,z_2)}{\Lambda}}\left(\iint_{E\times E}\varphi_2(z_1',z_2')\Gamma^{(2)}(z_1,z_2,\dd z_1',\dd z_2')\right)\\
&\qquad\qquad\qquad+\1_{\eta\geq\frac{\lambda(z_1,z_2)}{\Lambda}}\varphi(z_1,z_2)\Big\}\dd \eta -\Lambda\varphi_2(z_1,z_2)\\
&= \Lambda\times\frac{\lambda(z_1,z_2)}{\Lambda}\iint_{E\times E}\varphi_2(z_1',z_2')\Gamma^{(2)}(z_1,z_2,\dd z_1',\dd z_2')\\
&\qquad\qquad\qquad+\Lambda\left(1-\frac{\lambda(z_1,z_2)}{\Lambda}\right)\varphi_2(z_1,z_2)-\Lambda\varphi_2(z_1,z_2)\\
&= L^{(2)}\varphi_2(z_1,z_2),
\end{align*}
and thus $\mathcal{L}_N=\tilde{\mathcal{L}}_N$ and the two processes are equal in law.\end{proof}

\subsubsection{Parametric Boltzmann models}\label[I]{sec:boltzmannparametricmodels}

In many applications, the post-collisional distribution is explicitly given as the image measure of a known parameter space $(\Theta,\nu)$ endowed with a probability measure $\nu$ (or a positive measure with infinite mass in the non cutoff case). Analogously to the case of mean-field jump models (see Example \cref{example:parametricjump}), in this review, we distinguish this particular class of models and we call them \emph{parametric Boltzmann models}.  

\begin{definition}[Parametric Boltzmann model]\label[I]{def:Boltzmannparammodel} Let be given two measurable functions
\[\psi_1:E\times E\times\Theta\to E,\quad\psi_2:E\times E\times\Theta\to E,\]
which satisfy the symmetry assumption 
\begin{equation}\label[Ieq]{eq:sympsi1psi2}\forall (z_1,z_2)\in E^2,\quad (\psi_1,\psi_2)(z_1,z_2,\cdot)_{\#}\nu = (\psi_2,\psi_1)(z_2,z_1,\cdot)_{\#}\nu.\end{equation}
Let the function $\psi$ be defined by
\[\psi:E\times E\times\Theta\to E^2,\,(z_1,z_2,\theta)\mapsto \big(\psi_1(z_1,z_2,\theta),\psi_2(z_1,z_2,\theta)\big).\]
A parametric Boltzmann model with parameters $(\Theta,\psi)$ is a Boltzmann model of the form \cref{eq:boltzmanngenerator} with Assumption \cref{assum:L2} and a post-collisional distribution of the form: 
\[\forall (z_1,z_2)\in E^2,\quad \Gamma^{(2)}(z_1,z_2,\dd z_1',\dd z_2') = \psi(z_1,z_2,\cdot)_{\#} \nu.\]
\end{definition}
The symmetry assumption \cref{eq:sympsi1psi2} is the equivalent of \cref{eq:symgamma2} in this special case. In particular, for any two-variable test function $\varphi_2\in C_b(E^2)$ and any $(z_1,z_2)\in E\times E,$
\begin{align*}
\iint_{E\times E} \varphi_2(z_1',z_2')\Gamma^{(2)}(z_1,z_2,\dd z_1',\dd z_2') &= \int_{\Theta} \varphi_2(\psi_1(z_1,z_2,\theta),\psi_2(z_1,z_2,\theta))\nu(\dd\theta)\\
&= \int_{\Theta} \varphi_2(\psi_2(z_2,z_1,\theta),\psi_1(z_2,z_1,\theta))\nu(\dd\theta),
\end{align*}
where the last equality follows from \cref{eq:sympsi1psi2}. A sufficient condition for \cref{eq:sympsi1psi2} to hold is the case investigated in \cite{sznitman_equations_1984} with: 
\[\psi_2(z_1,z_2,\theta) = \psi_1(z_2,z_1,\theta).\]
In terms of particle systems, following Proposition \cref{prop:prototypeBoltzmannprocess}, in a parametric model, when a collision occurs, a parameter $\theta\sim\nu$ is sampled first and then the states of the particle is updated by: 
\[\Big(Z^i_{T_{ij}^+},Z^j_{T_{ij}^+}\Big) = \psi{\left(Z^i_{T_{ij}^-},Z^j_{T_{ij}^-},\theta\right)}.\]

\begin{example}[Symmetrization]\label[I]{example:symmetrization} Wagner \cite{wagner_functional_1996} treats the case of particle systems with a generator of the form: for $\mathbf{z}^N=(z^1,\ldots,z^N)\in E^N$ and $\varphi_N\in C_b(E^N)$, 
\begin{multline}\label[Ieq]{eq:Boltzmannnonsym}\mathcal{L}_N\varphi_N(\mathbf{z}^N) = \sum_{i=1}^N L^{(1)}\diamond_i \varphi_N(\mathbf{z}^N) \\+ \frac{1}{2N}\sum_{i\ne j}\tilde{\lambda}(z^i,z^j)\int_{\tilde{\Theta}}\big\{\varphi_N\big(\mathbf{z}^N\big(i,j,\tilde{\theta}\big)\big)-\varphi_N\big(\mathbf{z}^N\big)\big\}\tilde{\nu}(\dd\tilde{\theta}),\end{multline}
where $\tilde{\lambda}:E\times E\to\R_+$, $\tilde{\Theta}$ is a parameter set endowed with a probability measure $\tilde{\nu}$ and $\mathbf{z}^N(i,j,\tilde{\theta})$ is the $N$ dimensional vector whose $k$ component is equal to 
\[{z}^k(i,j,\tilde{\theta}) = \left\{\begin{array}{rcl} 
z^k & \text{if} & k\ne i,j\\
\tilde{\psi}_1(z^i,z^j,\tilde{\theta}) & \text{if} & k=i\\
\tilde{\psi}_2(z^i,z^j,\tilde{\theta}) & \text{if} & k=j
\end{array}
\right.,\]
for two given functions $\tilde{\psi}_1, \tilde{\psi}_2 :E\times E\times \tilde{\Theta}\to E$. The main difference with the generator \cref{eq:boltzmanngenerator} is that Wagner distinguishes the pairs $(i,j)$ and $(j,i)$ while in \cref{eq:boltzmanngenerator} we consider unordered pairs of particles but add the symmetry assumption \cref{eq:symgamma2}. Consequently, the double sum in \cref{eq:Boltzmannnonsym} runs over all indices $i,j=1,\ldots,N$ while in the sum \cref{eq:boltzmanngenerator}, it runs over the indices $i<j$. Nevertheless, using a simple symmetrization procedure, the model \cref{eq:Boltzmannnonsym} fits into the previous general framework with 
\[\Theta = \tilde{\Theta}\times[0,1],\quad \nu(\dd\theta) = \tilde{\nu}(\dd \tilde{\theta})\otimes \dd\sigma,\]
where $\theta=(\tilde{\theta},\sigma)\in\Theta$, $\dd\sigma$ is the uniform probability measure on $[0,1]$ and for $z_1,z_2\in E$ we define
\begin{align*}
\lambda(z_1,z_2) &= \frac{\tilde{\lambda}(z_1,z_2)+\tilde{\lambda}(z_2,z_1)}{2},\\
\psi_1(z_1,z_2,\theta) &= \1_{\sigma\leq \frac{\tilde{\lambda}(z_1,z_2)}{2\lambda(z_1,z_2)}}\tilde{\psi}_1(z_1,z_2,\tilde{\theta})+\1_{\sigma> \frac{\tilde{\lambda}(z_1,z_2)}{2\lambda(z_1,z_2)}}\tilde{\psi}_2(z_2,z_1,\tilde{\theta}),\\
\psi_2(z_1,z_2,\theta) &= \1_{\sigma\leq \frac{\tilde{\lambda}(z_1,z_2)}{2\lambda(z_1,z_2)}}\tilde{\psi}_2(z_1,z_2,\tilde{\theta})+\1_{\sigma> \frac{\tilde{\lambda}(z_1,z_2)}{2\lambda(z_1,z_2)}}\tilde{\psi}_1(z_2,z_1,\tilde{\theta}).
\end{align*}
One can check that the functions $\psi_1$ and $\psi_2$ satisfy \cref{eq:sympsi1psi2} and that the generator \cref{eq:boltzmanngenerator} of the associated parametric model (Definition \cref{def:Boltzmannparammodel}) is equal to \cref{eq:Boltzmannnonsym}. In this case, the Boltzmann equation \cref{eq:Boltzmannequationgeneral} reads 
\begin{multline*}
\frac{\dd}{\dd t}\langle f_t,\varphi\rangle = \langle f_t, L^{(1)}\varphi\rangle+\frac{1}{2}\int_{\tilde{\Theta}\times E^2}\Big\{\tilde{\lambda}(z_1,z_2)\big[\varphi\big(\tilde{\psi}_1(z_1,z_2,\tilde{\theta})\big)-\varphi(z_1)\big]\\+\tilde{\lambda}(z_2,z_1)\big[\varphi\big(\tilde{\psi}_2(z_2,z_1,\tilde{\theta})\big)-\varphi(z_1)\big]\Big\}\tilde{\nu}(\dd\tilde{\theta})f_t(\dd z_1)f_t(\dd z_2),    
\end{multline*}
or equivalently after the change of variables $(z_1,z_2)\mapsto(z_2,z_1)$, 
\begin{multline}\label[Ieq]{eq:boltzmannpsitilde}
\frac{\dd}{\dd t}\langle f_t,\varphi\rangle = \langle f_t, L^{(1)}\varphi\rangle+\frac{1}{2}\int_{\tilde{\Theta}\times E^2}\tilde{\lambda}(z_1,z_2)\Big\{\varphi\big(\tilde{\psi}_1(z_1,z_2,\tilde{\theta})\big)+\varphi\big(\tilde{\psi}_2(z_1,z_2,\tilde{\theta})\big)\\-\varphi(z_1)-\varphi(z_2)\Big\}\tilde{\nu}(\dd\tilde{\theta})f_t(\dd z_1)f_t(\dd z_2).    
\end{multline}
\end{example}

The introductory section of \cite{wagner_functional_1996} contains many examples of such models, in particular models (in Russian) due to Leontovich in the 30's and Skorokhod in the 80's that we did not manage to find. A more recent example inspired by economic models of wealth distribution \cite{matthes_steady_2008} is given in \cite{cortez_quantitative_2016}. The authors assume $E=\R$ with $L^{(1)}=0$, $\lambda=1$, $\tilde{\Theta}=\R^4$ and 
\[\tilde{\psi}_1\big(z_1,z_2,(L,R,\tilde{L},\tilde{R})\big) = Lz_1+Rz_2,\]
and
\[\tilde{\psi}_2\big(z_1,z_2,(L,R,\tilde{L},\tilde{R})\big) = \tilde{L}z_2+\tilde{R}z_1.\]
In this model, the state of a particle represents the wealth of an individual and the parameters $(L,R,\tilde{L},\tilde{R})$ specify how a trade between two individuals affect their wealth. This model generalises a famous model due to Kac \cite{kac_foundations_1956} which will be discussed in the next section. From the modelling point of view, it is more natural to use generators of the form \cref{eq:Boltzmannnonsym}; several additional examples will be given in particular in Section \cref{sec:socioeconomicmodels} for socio-economic models. On the other hand, the generator form \cref{eq:boltzmanngenerator} will simplify some computations in Section \cref{sec:boltzmannreview}. 

In the parametric framework, the particle system can be advantageously written as the solution of a system of SDEs driven by Poisson measures. In a famous article, Tanaka \cite{tanaka_probabilistic_1978} proposed a SDE approach to study the nonlinear Boltzmann system of rarefied gas dynamics (which will be presented in the next section). He introduced a class of nonlinear SDEs driven by Poisson random measures which will be described in Section \cref{sec:couplingBoltzmann}. As we shall see, although it is relatively easy to write a system of coupled SDEs which describes the particle system, its relationship with Tanaka's SDE is not completely straightfoward. Around the same time, Murata \cite{murata_propagation_1977} tackled the question and proved the propagation of chaos (for a specific model) using a coupling argument between the two systems of SDEs. The idea is of course reminiscent of the well-known McKean's theorem and all the works reviewed in Section \cref{sec:synchronouscouplingreview} for McKean-Vlasov systems. Note however that Murata's work is among the first ones which use the very fruitful idea of coupling to prove propagation of chaos. His argument is based on a clever but not so easy optimal coupling argument which seems to have been largely forgotten in the subsequent literature. A recent series of articles \cite{fournier_rate_2016,cortez_quantitative_2016,cortez_quantitative_2018} has proposed a more contemporary point of view on the question. The arguments are very similar to Murata's but take advantage of the development of the theory of optimal transport. Let us also mention that these articles seem to be based on \cite{fontbona_measurability_2009} which also introduces an optimal coupling argument reminiscent of Murata's but in a different context, namely the derivation of the Landau equation from a system of interacting diffusion processes. We will continue this discussion in Section \cref{sec:couplingBoltzmann}.

\begin{example}[Semi-parametric model]\label[I]{example:semiparametric} A natural extension of the parametric model would consider a measure on $\Theta$ which depends on the state of the particles, for instance one can consider a post-collisional distribution of the form
\begin{multline}\label[Ieq]{eq:semiparametric}\iint_{E\times E} \varphi_2(z_1',z_2')\Gamma^{(2)}(z_1,z_2,\dd z_1',\dd z_2') \\= \int_{\Theta} \varphi_2\big(\psi_1(z_1,z_2,\theta),\psi_2(z_1,z_2,\theta)\big)q(z_1,z_2,\theta)\nu(\dd\theta)\end{multline}
where for all $z_1,z_2\in E$, $q(z_1,z_2,\cdot)$ is a probability density function with respect to the measure ${\nu}\in\mathcal{M}^+(E)$. Wagner \cite{wagner_functional_1996} considered such model that will be called \emph{semi-parametric} in this review. If there exists $M>0$ and $q_0({\theta})$ a probability density function with respect to ${\nu}$ such that
\begin{equation}\label[Ieq]{eq:semiparambound}\forall z_1,z_2\in E,\,\forall {\theta}\in{\Theta},\quad q(z_1,z_2,{\theta})\leq M q_0({\theta}),\end{equation}
then the situation can be reduced to the parametric case thanks to an accept-reject scheme similar to the one in Proposition \cref{prop:acceptreject}. Namely in the extended parameter space 
\[\tilde{\Theta}=\Theta\times [0,1],\]
endowed with the probability measure $q_0(\theta)\dd\theta\otimes d\eta$,
let us define the function
\[\tilde{\psi}(z_1,z_2,(\theta,\eta)) = \left\{\begin{array}{rcl}
\big(\psi_1(z_1,z_2,\theta),\psi_2(z_1,z_2,\theta)\big) & \text{if} & \eta\leq\frac{q(z_1,z_2,\theta)}{Mq_0(\theta)} \\ 
(z_1,z_2) & \text{if} & \eta>\frac{q(z_1,z_2,\theta)}{Mq_0(\theta)}
\end{array}
\right..\]
Then up to a time rescaling $t\to tM$, the parametric model $(\tilde{\Theta},\tilde{\psi})$ is equivalent in law to the semi-parametric model. Note that \cref{eq:semiparambound} automatically holds when $\Theta$ is compact and $q$ bounded.
\end{example}

\subsubsection{Classical models in collisional kinetic theory}\label[I]{sec:Boltzmannclassicalmodels}

The foundations of kinetic theory lie in the seminal work of Boltzmann and Maxwell who attempted to understand the large scale behaviour of a gas of particles defined in the phase space $E=\R^d\times\R^d$ by their position and velocity. Many interactions mechanisms can be considered, depending on the physical assumptions. The starting point is the Newton equations satisfied by the $N$-particle system $(\mathcal{Z}^N_t)^{}_t$, for $i\in\{1,\ldots,N\}$, 

\begin{equation}\label[Ieq]{eq:Newtonpotential}\left\{\begin{array}{rcl}
    \displaystyle{\frac{\dd X^i_t}{\dd t}} & = & V^i_t\\
    \displaystyle{\frac{\dd V^i_t}{\dd t}} & = & \displaystyle{-\sum_{j=1}^{N}\nabla V(|X^j_t-X^i_t|)} 
\end{array}\right.,\end{equation}
where $V$ is a (smooth) repulsive potential, typically an inverse power law. Another important system is the hard-sphere system which will be described in Example~\cref{example:hardsphere}. Without any other assumption, it is not clear that this set of equations defines a binary collision process. In fact, it is more reminiscent of a mean-field system without the (crucial) $1/N$ scaling in front of the sum. Boltzmann and Maxwell considered the case of dilute gases (also called \emph{rarefied gas}), that is gases where the density of particles is so small that in the sum in \cref{eq:Newtonpotential}, there is typically no more than one non-zero term. The dynamics of each particle is therefore mainly driven by the free transport until the particle comes very close to another particle which induces a deviation of its trajectory (as well as the trajectory of the other particle) depending on the potential $V$. During this process, everything is deterministic and the only source of randomness comes from the initial condition. The probabilistic interpretation presented in this section is due to Kac. Boltzmann derived the equation satisfied by the one-particle distribution when $N\to+\infty$. In its most general form, the Boltzmann equation of rarefied gas dynamics reads (in strong form): 
\begin{multline}\label[Ieq]{eq:Boltzmannphysics}
    \partial_t f_t(x,v) + v\cdot\nabla_x f_t\\ = \int_{\R^d}\int_{\mathbb{S}^{d-1}} B(v-v_*,\sigma)\Big(f_t(x,v_*')f_t(x,v')-f_t(x,v_*)f_t(x,v)\Big)\dd v_*\dd\sigma,
\end{multline}
where
\begin{equation}\label[Ieq]{eq:postcollisionsigma}\left\{\begin{array}{rcl}
v' & = & \displaystyle{\frac{v+v_*}{2}+\frac{|v-v_*|}{2}\sigma}\\
v_*' &= & \displaystyle{\frac{v+v_*}{2}-\frac{|v-v_*|}{2}\sigma}
\end{array}
\right.,\end{equation}
are the post-collisional velocities. The parameter $\sigma\in\mathbb{S}^{d-1}$ is often called the scattering angle. This transformation preserves energy and momentum. The function $B:\R^d\times\mathbb{S}^{d-1}\to\R_+$ is of the form
\begin{equation}\label[Ieq]{eq:collisionkernelB}B(u,\sigma) = \Phi(|u|)\Sigma(\theta),\end{equation}
with $\cos\theta=\frac{u}{|u|}\cdot\sigma$, $\theta\in[0,\pi]$. The function $\Phi$ is called the velocity cross-section and the function $\Sigma$ is called the angular cross-section. The function $B$ is referred as the collision kernel (in the literature, it is also sometimes called the cross-section). It is customary to write $B(u,\sigma)\equiv B(|u|,\cos\theta)$. 
Depending on the choice of the potential $V$, some of the most important collision kernels derived by Maxwell are listed below. 
\begin{itemize}
    \item (Hard-sphere) 
    \begin{equation}\label[Ieq]{eq:hardsphereB}\Phi(|u|)=|u|,\quad \Sigma(\theta)=1.\end{equation}
    \item (Inverse-power law potentials) \[\Phi(|u|)=|u|^\gamma, \quad \gamma=\frac{s-(2d-1)}{s-1},\quad s>2,\]
    and $\Sigma$ has a non-integrable singularity when $\theta\to0$, so that 
    \[\int_0^\pi\Sigma(\theta)\dd \theta = +\infty\]
    \item (Maxwell molecules) 
    \begin{equation}\label[Ieq]{eq:truemaxwellmolecules}\Phi(|u|)=1,\quad \int_0^\pi\Sigma(\theta)\dd \theta = +\infty.\end{equation}
    \item (Maxwell molecules with Grad's cutoff)
    \begin{equation}\label[Ieq]{eq:maxwellmoleculescutoff}\Phi(|u|)=1, \quad \int_0^\pi \Sigma(\theta)\dd \theta <+\infty.\end{equation}
\end{itemize}
    
We will not go further into the description of the Boltzmann equation. The interested reader will find a thorough discussion and analysis of these different models in the reviews \cite{villani_review_2002, bellomo_macroscopic_2004} or in the classical books \cite{cercignani_boltzmann_1988,cercignani_mathematical_1994}. We also mention the book \cite{cercignani_ludwig_2006} which contains a very interesting biography of Ludwig Boltzmann as well as a scientific discussion of the physics of his time and of his legacy. 

This review is focused on the rigorous derivation of the Boltzmann equation from a system of particles. On the right-hand side of \cref{eq:Boltzmannphysics}, the variable $x$ (position) only appears as a parameter: this is the limit where collisions between two particles happen only when the two particles are at the same position. From a mathematical point of view, this purely local interaction mechanism makes the derivation very difficult if not impossible with stochastic tools (see \cite{meleard_asymptotic_1996}). A deterministic example (Lanford's theorem) is nevertheless given in Example \cref{example:hardsphere} and Section \cref{sec:lanford}. Apart from this result we will focus on simplified mechanisms: either spatially homogeneous \emph{Kac models} (Example \cref{example:spatiallyhomogeneousboltzmann} and Example \cref{example:kacmodels}) or kinetic mollified models (Example \cref{example:mollifiedmodels}). Both cases have a natural probabilistic interpretation which fits into the framework of Section \cref{sec:boltzmanngeneral}.  

\begin{example}[Spatially homogeneous Boltzmann equation]\label[I]{example:spatiallyhomogeneousboltzmann} 
In the case of a spatially homogeneous problem, the difficulty due to the local interaction does not appear. In this case the spatially-homogeneous Boltzmann equation of rarefied gas dynamics describes a gas of particles defined by their velocity only: 
\begin{equation}\label[Ieq]{eq:spatiallyhomogeneousboltzmann}
\partial_t f(t,v) = \int_{\R^d}\int_{\mathbb{S}^{d-1}} B(v-v_*,\sigma)\Big(f_t(v_*')f_t(v')-f_t(v_*)f_t(v)\Big)\dd v_*\dd\sigma,
\end{equation}
This last equation can be shown to be the strong form of the Boltzmann equation~\cref{eq:Boltzmannequationgeneral} with a parametric post-collisional distribution given by the cross-section~$B$. The $N$-particle stochastic process associated to this equation is given by Proposition \cref{prop:prototypeBoltzmannprocess}. More precisely, the (spatially-homogeneous) hard-sphere model and the (spatially-homogeneous) model of Maxwell molecules with Grad's cutoff fit into the framework of Section \cref{sec:boltzmannparametricmodels} with:
\[\psi_1(v,v_*,\theta) = v',\quad\psi_2(v,v_*,\theta) = v'_*,\]
and 
\begin{align*}\lambda(v,v_*) &= \Phi(|v-v_*|)\int_0^\pi\Sigma(\theta)\dd\theta,\\ \Gamma(v,v_*,\dd z',\dd v_*,\dd v_*') &= \psi(v,v_*,\cdot)_{\#} \left(\frac{\Sigma}{\int_0^\pi\Sigma(\theta)\dd\theta}\right).\end{align*}
To be more precise, with this particular choice of the parameters, the weak-form of the general Boltzmann equation \cref{eq:symmetricBoltzmannequationgeneral} reads: 
\begin{multline*}\frac{\dd}{\dd t}\langle f_t,\varphi\rangle \\= \frac{1}{2}\int_{\R^d\times\R^d\times\mathbb{S}^{d-1}} \Big\{\varphi(v')+\varphi(v_*')-\varphi(v)-\varphi(v_*)\Big\}f_t(v)f_t(v_*)B(v-v_*,\sigma)\dd v\dd v_*\dd\sigma.\end{multline*}
The strong form \cref{eq:spatiallyhomogeneousboltzmann} is obtained thanks to the following classical involutive unit Jacobian changes of variables which allow to exchange $(v,v_*)$ and $(v',v'_*)$ : 
\[(v,v_*,\sigma)\to (v',v_*',\vec{k}),\quad (v,v_*)\to(v_*,v),\]
with $\vec{k}=(v-v_*)/|v-v_*|$ (see \cite[Chapter 1, Section 4.5]{villani_review_2002}). The collision kernel $B$ is invariant by these changes of variables so we can keep its arguments unchanged. 

The non-cutoff cases are more difficult to handle due to the non-integrability of the angular cross-section (see Example \cref{example:noncutoffmodels}). 
\end{example} 

The following example describes the famous Kac model, which is a one dimensional caricature of a gas of Maxwellian molecules. Beyond the relative simplicity of the model, the seminal article of Kac \cite{kac_foundations_1956} is of particular importance because it introduces the mathematical definition of propagation of chaos. 

\begin{example}[Kac model]\label[I]{example:kacmodels}
     Within the framework of Section \cref{sec:boltzmannparametricmodels}, the Kac model is defined in $E=\R$ by $L^{(1)}=0$ and the post-collisional distribution 
    \begin{multline*}\iint_{\R\times\R}\varphi_2(z_1',z_2') \Gamma^{(2)}(z_1,z_2,\dd z_1',\dd z_2') \\:= \frac{1}{2\pi}\int_{-\pi}^{\pi} \varphi_2(z_1\cos\theta+z_2\sin\theta,-z_1\sin\theta+z_2\cos\theta)\dd\theta\end{multline*}
    and the collision rate $\lambda(z_1,z_2)=\nu=\text{constant}$. Then the weak Boltzmann equation becomes 
    \[\frac{\dd}{\dd t}\langle \varphi,f_t\rangle = \frac{\nu}{2\pi}\iint_{\R\times\R} \int_{-\pi}^{\pi}\{\varphi(z_1\cos\theta+z_2\sin\theta)-\varphi(z_1)\}f_t(\dd z_1)f_t(\dd z_2).\]
    With the change of variable (with $\theta$ fixed)
    \[(z_1',z_2') = (z_1\cos\theta+z_2\sin\theta,-z_1\sin\theta+z_2\cos\theta),\]
    followed by $\theta\mapsto -\theta$ (both changes of variable have unit jacobian), Kac obtained the following equation in strong form:
    \[\partial_t f_t(z_1) = \frac{\nu}{2\pi}\int_\R\int_{-\pi}^{\pi} \{f(z_1')f(z_2')-f_t(z_1)f_t(z_2)\}\dd z_2\dd\theta.\]
    We refer the reader to \cite{carlen_entropy_2008,mischler_kacs_2012} for a thorough analysis and discussion of the Kac model and its generalisations in kinetic theory. Keeping a collision rate $\lambda$ constant the authors of \cite{carlen_kinetic_2013} generalised the arguments of the proof of the propagation of chaos to a larger class of models. This generalised result, that we will call \emph{Kac's theorem}, will be discussed in Section \cref{sec:kactheorem}. As already mentioned in Section \cref{sec:boltzmannparametricmodels}, the Kac model is also a special instance of the model studied in \cite{cortez_quantitative_2016} within a framework which will be described in Section \cref{sec:couplingBoltzmann}.
\end{example}

The work of Kac had a very strong influence on the literature so that Boltzmann models with $L^{(1)}=0$ are sometimes called Kac models, or also \emph{homogeneous} Boltzmann models. 

\begin{example}[Mollified Boltzmann models]\label[I]{example:mollifiedmodels} The interaction mechanism of a kinetic particle system is said to be \emph{purely local} when, within the framework of Section \cref{sec:boltzmanngeneral}, the collision rate is taken equal to 
\[\lambda((x_1,v_1),(x_2,v_2)) = \delta_{x_1,x_2}.\]
This indicates that two particles interact if and only if they are exactly at the same position. As explained in \cite{meleard_asymptotic_1996}, the probabilistic interpretation of purely local models is extremely difficult and one can rather consider the smoothened version: 
\[\lambda((x_1,v_1),(x_2,v_2)) = K(|x_1-x_2|),\]
where $K$ is a smooth mollifier with fixed radius (i.e. a non negative radial function which tends to zero at infinity and which integrates to one). The post-collisional distribution is unchanged and acts only on the velocity variable:
\[\Gamma^{(2)}(z_1,z_2,\dd z_1',\dd z_2')\equiv \Gamma^{(2)}(v_1,v_2,\dd v_1',\dd v_2')\otimes\delta_{x_1}(\dd x_1')\otimes\delta_{x_2}(\dd x_2').\]
This model is called a mollified Boltzmann model. Its probabilistic treatment is discussed in \cite{graham_stochastic_1997} and \cite{meleard_asymptotic_1996}. Note that with a general state space $E$, all the models in Section \cref{sec:boltzmanngeneral} and in particular the one in Proposition \cref{prop:prototypeBoltzmannprocess} are implicitly mollified Boltzmann models. Most of the models reviewed in Section \cref{sec:boltzmannreview} are mollified models. Purely local Boltzmann models can be recovered by letting the mollifier converge to a Dirac delta $K\to\delta_0$ (formally or with a quantitative control). Another example of purely local Boltzmann model is the hard-sphere system defined in the next Example \cref{example:hardsphere}. 
\end{example}

\begin{example}[Hard-sphere system]\label[I]{example:hardsphere} A \emph{hard-sphere} is a spherical particle defined by its position, its velocity and its diameter $\varepsilon>0$. Moreover, it is assumed that two hard-spheres cannot overlap. A system of $N$ hard-spheres is thus defined on the domain: 
\[\mathcal{D}_N := \big\{\mathbf{z}^N = (x^i,v^i)_{i\in\{1,\ldots,N\}}\in (\R^{d}\times\R^d)^N,\,\,\forall i\ne j,\,|x^i-x^j|\geq \varepsilon\big\}.\]
The dynamics of the hard-sphere system is a special degenerate case of \cref{eq:Newtonpotential} with a vanishing potential but with an additional boundary condition which tells what happens on the boundary of $\mathcal{D}_N$, that is when two particles are at a distance $\varepsilon$ (the term collision is here self-explanatory). The collision of two hard-spheres is an elastic collision which preserves energy and momentum. Starting with a pair of pre-collisional velocities $(v^i,v^j)$, writing down the conservation laws leads to the following formula for the post-collisional velocities: 
\begin{equation}\label[Ieq]{eq:postcollisionnu}
\begin{array}{ll}
v^{i*}= v^i-\nu^{i,j}\cdot(v^i-v^j)\nu^{i,j}\vspace{0.2cm}\\
v^{j*}= v^j+\nu^{i,j}\cdot(v^i-v^j)\nu^{i,j}
\end{array},
\end{equation}
where $\nu^{i,j}:=(x^i-x^j)/|x^i-x^j|\in\mathbb{S}^{d-1}$. This representation is not the same as the representation \cref{eq:postcollisionsigma} but it can be shown that they are actually equivalent \cite[Chapter 1, Section 4.6]{villani_review_2002}. Pre-collisional means that $(v^i,v^j)$ are such that $(v^i-v^j)\cdot\nu^{i,j}<0$. It can also be checked that the post-collisional velocities satisfy $(v^{i*}-v^{j*})\cdot\nu^{i,j}>0$. Note that this transformation is an involution in the sense that if $(v^i-v^j)\cdot\nu^{i,j}>0$ (that is the $v^i$ and $v^j$ are in a post-collisional configuration), then \cref{eq:postcollisionnu} gives the pre-collisional velocities. Note also that this dynamical system is completely deterministic.

The large scale behaviour when $N\to+\infty$ and $\varepsilon\to0$ is given by the Boltzmann equation \cref{eq:Boltzmannphysics} with the hard-sphere cross-section. Under the chaotic assumption \cref{eq:chaosassumptionintro} at time $t=0$, Lanford's theorem \cite{lanford_time_1975} states that in the Boltzmann-Grad limit $N\varepsilon^{d-1}\to1$, \cref{eq:chaosassumptionintro} also holds for later time. This scaling was introduced by Grad in \cite{grad_asymptotic_1963}. The proof of Lanford's theorem is extremely difficult. We will briefly review the main ideas in Section \cref{sec:lanford}. Our presentation will follow closely \cite{gallagher_newton_2014}. 
\end{example}

\pagebreak
\section{Notions about chaos} \label[I]{sec:def}

\subsection{Topology reminders: metrics and convergence for probability measures} \label[I]{sec:topology}

Since propagation of chaos is about the convergence of probability measures, we first need to present the topological tools that will be constantly used in the following. The content of this section is fairly classical, most of the results specific to our topic can be found in \cite{hauray_kacs_2014}, see also \cite[Section 3.4]{jabin_review_2014}, \cite[Section 2.5]{mischler_kacs_2013}, \cite[Section 3]{mischler_new_2015} or \cite{villani_limite_2001}. A more general overview of the topology of the space of probability measures can be found in the classical books \cite{villani_optimal_2009,billingsley_convergence_1999, parthasarathy_probability_1967}.

\subsubsection{Distances on the space of probability measures}\label[I]{sec:distancesproba}

Let $(\mathscr{E},d_\mathscr{E})$ be a Polish space. For $p \geq 1$, a measure $\mu$ in $\pb ( \mathscr{E} )$ admits a \emph{finite $p$-th moment} when there exists $x_0\in \mathscr{E}$ such that
\[  \E_{X \sim \mu} \big[ d_\mathscr{E}(X,x_0)^p \big] := \int_\mathscr{E} d_\mathscr{E}(x,x_0)^p \mu ( \dd x ) < + \infty. \]
This property does not depend on $x_0$. The space of probability measures with finite $p$-th moment is denoted by $\pb_p ( \mathscr{E} )$. The Wasserstein distance on $\pb_p(\mathscr{E})$ will be the most important one in the following.   

\begin{definition}[Wasserstein distances]\label[I]{def:wasserstein} For $p \geq 1$, the Wasserstein-$p$ distance between the probability measures $\mu$ and $\nu$ in $\pb_p \left( \mathscr{E} \right)$ is defined by
\[W_{d_\mathscr{E},p}(\mu,\nu) := \inf_{\pi\in\Pi(\mu,\nu)} \left(\int_{\mathscr{E}\times \mathscr{E}} d_\mathscr{E}(x,y)^p\pi(\dd x,\dd y)\right)^{1/p} = \inf_{\substack{X\sim \mu\\ Y\sim \nu}} \mathbb{E} \left[ d_\mathscr{E}(X,Y)^p \right]^{1/p} \]
where $\Pi(\mu,\nu)$ is the set of all couplings of $\mu$ and $\nu$, that is to say, the set of probability measures on $E\times E$ with first and second marginals respectively equal to $\mu$ and $\nu$. 
\end{definition}

The total variation distance can be understood as a Wasserstein-1 distance with the trivial distance $d_\mathscr{E}(x,y)=\delta_{x,y}$.

\begin{definition}[Total variation norm]
The total variation distance between two probability measures $\mu$ and $\nu$ in $\pb \left( \mathscr{E} \right)$ is defined by
\[\|\mu-\nu\|_{\mathrm{TV}} = 2\inf_{\substack{X\sim \mu\\ Y\sim \nu}}\mathbb{P}(X\ne Y).\]
\end{definition}

Since $\pb ( \mathscr{E} )$ can be seen as a subset of the dual space $C_b ( \mathscr{E} )^{\star}$, natural strong norms on $\pb ( \mathscr{E} )$ are induced by usual norms on functional spaces. The following proposition links these distances to dual norms:

\begin{proposition}[Duality formulae]\label[I]{prop:dualitywasserstein}
The total variation and Wasserstein-1 distances satisfy: 
\[\|\mu-\nu\|_{\mathrm{TV}} = \sup_{\|\varphi\|_\infty\leq 1}\left\{\int_\mathscr{E} \varphi(x)\mu(\dd x)-\int_\mathscr{E} \varphi(x)\nu(\dd x)\right\}\]
and
\[W_{d,1}(\mu,\nu) = \sup_{\|\varphi\|_{\mathrm{Lip},d_\mathscr{E}}\leq 1}\left\{\int_\mathscr{E} \varphi(x)\mu(\dd x)-\int_\mathscr{E} \varphi(x)\nu(\dd x)\right\}\]
where
\[\|\varphi\|_{\mathrm{Lip},d_\mathscr{E}}:=\sup_{x\ne y} \frac{|\varphi(x)-\varphi(y)|}{d_\mathscr{E}(x,y)}.\]
\end{proposition}

\begin{proof}
See \cite[Theorem 1.14]{villani_topics_2003}
\end{proof}

Another important class of dual norms is given by the negative Sobolev norms $W^{-s,p}$. Let us emphasize two special cases.

\begin{definition}[Some negative Sobolev norms]\label[I]{def:sobolevnorm} When $\mathscr{E} = \R^d$ we define the following norms. 
\begin{itemize}
    \item For $\mu,\nu\in\pb(\mathscr{E})$ and $s>\frac{d}{2}$
    \[\|\mu-\nu\|_{H^{-s}}^2 := \int_{\R^d} |\hat{\mu}(\xi)-\hat{\nu}(\xi)|^2\frac{\dd\xi}{(1+|\xi|^2)^{s}},\]
    where $\hat{\mu}$ is the Fourier transform of $\mu$. 
    \item The dual norm of the Euclidean Lipschitz semi-norm 
    \[ \|\mu-\nu\|_{W^{-1,\infty}} := \sup_{ \| \varphi \|_{W^{1,\infty}} \leq 1} \langle \mu - \nu , \varphi \rangle \]
    where the $W^{1,\infty}$ Sobolev norm is $\| \varphi \|_{W^{1,\infty}} = \| \varphi \|_{\infty} +\| \nabla \varphi \|_{\infty}$.
\end{itemize}
\end{definition}

An important property of the negative Sobolev norm $H^{-s}$ is its polynomial structure (see \cite[Lemma 2.9]{hauray_kacs_2014}), in the sense that the definition of this norm only involves the integration of some functions against tensor products of measures.

\begin{lemma}\label[I]{lemma:sobolevpolynomial}
The negative Sobolev norm $H^{-s}$, $s>d/2$ on $\R^d$ satisfies for any $\mu,\nu\in\pb(\R^d)$, 
\begin{align}\label[Ieq]{eq:polynomialsobolev}
\|\mu-\nu\|_{H^{-s}}^2 &= \int_{\R^{2d}}\Phi_s(x-y)(\mu^{\otimes 2}-\mu\otimes\nu)(\dd x,\dd y)\\
&\qquad\qquad+\int_{\R^{2d}}\Phi_s(x-y)(\nu^{\otimes 2}-\nu\otimes\mu)(\dd x,\dd y)\\
&=\int_{\R^{2d}} \Phi_{s}(x-y)(\mu-\nu)^{\otimes 2}(\dd x,\dd y),
\end{align}
where $\Phi_s(z) := \int_{\R^d} \e^{-iz\cdot\xi}(1+|\xi|^2)^{-s}\dd \xi$.
\end{lemma}

\begin{remark}
A similar distance has also been recently considered in particular in the work of Serfaty \cite{serfaty_systems_2019} where it is called a \emph{modulated energy}. 
\end{remark}

\begin{proposition}[Comparison of distances] \label[I]{comparisondistance}
Assume the distance $d_\mathscr{E}$ to be bounded. The following uniform topological equivalences hold.
\begin{itemize}
    \item The $\mathrm{TV}$ distance dominates the Wassertein-1 distance $W_1$.
    \item $\| \cdot \|_{\mathrm{Lip},d_\mathscr{E}}$ is equivalent to $\| \cdot \|_{W^{1,\infty}}$ (see \cite[Equations (2.4) and (2.5)]{hauray_kacs_2014}) and this implies the same for $W_1$ and $\|\cdot \|_{W^{-1,\infty}}$ distances.
    \item The $W_2$-distance dominates the $W_1$-distance, and for $s > \frac{d+1}{2}$ the $W_1$-distance dominates the square of the $H^{-s}$-distance.
    \item For measures in $\pb_p (\mathscr{E})$ with $p > 0$ and $s \geq 1$, the $H^{-s}$ distance dominates the $W_1$ distance up to a positive exponent.  
    \item For measures in $\pb_p (\mathscr{E})$ with $p > 2$, the $W_1$-distance dominates the $W_2$-distance up to a positive exponent. 
\end{itemize}
\end{proposition}

\begin{proof}
See \cite[Lemma 2.1]{hauray_kacs_2014}, which gives a quantitative version of this.
\end{proof}

Finally, let $\mathcal{F}=\{\varphi^k,\,\, k\in\N\}$ be a countable and separating subset of $C_b(\mathscr{E})$ (see Definition \cref{def:separatingconvergencedetermining} below) and such that $\|\varphi^k\|_\infty\leq 1$ for all $k\in\N$. Then the following expression defines a distance on $\pb(E)$ for any $p\geq1$

\begin{equation}\label[Ieq]{eq:D0expression}
D_p(\mu,\nu) := {\left(\sum_{k=1}^{+\infty} \frac{1}{2^k}| \langle \mu-\nu,\varphi^k \rangle|^p\right)}^{1/p}.
\end{equation}

In the literature, the most encountered distances are $D_1$ and $D_2$ which are used as a convenient tool to metricise the notion of weak convergence defined below (see Example \cref{ex:convergencedeterminingsets}). To conclude, we summarise the main cases of interest for the Wasserstein distance and other related distances.

\begin{definition}\label[I]{def:spaceswasserstein} In problems related to propagation of chaos, the Wasserstein distances are often used in the following cases.
\begin{itemize}
    \item When $\mathscr{E}=E$ is the state space of the particles, endowed with a distance $d_E$, we do not specify the dependency in $d_E$:
    \[W_{d_E,p}\equiv W_p \]
    When $E=\R^d$, the bounded moment assumption can be removed by using the bounded distance $\tilde{d}_E(x,y) := \inf(|x-y|,1)$.
    \item When $\mathscr{E}=E^k$, $k\in\N$ is a product space of the state space $(E,d_E)$, unless otherwise specified we follow \cite{hauray_kacs_2014} and use the normalised distance: for $\mathbf{x}^k=(x^1,\ldots,x^k)$ and $\mathbf{y}^k=(y^1,\ldots,y^k)$,
    \[d_{E^k}(\mathbf{x}^k,\mathbf{y}^k) := \frac{1}{k}\sum_{i=1}^k d_E(x^{i},y^{i}), \]
    and we simply write $W_{d_{E^k},p}\equiv W_p$. When we use the non-normalised distance 
    \[\widetilde{d}_{E^k}(\mathbf{x}^k,\mathbf{y}^k) := \sum_{i=1}^k d_E(x^{i},y^{i}),\]
    we write $W_{\widetilde{d}_{E^k},p}\equiv \widetilde{W}_p$. In the special case $E=\R^d$ endowed with the Euclidean norm, we will rather use the normalised distance, 
    \[d^p_{E^k}(\mathbf{x}^k,\mathbf{y}^k) := \frac{1}{k}\sum_{i=1}^k |x^{i}-y^{i}|^p, \]
    and the non-normalised one
    \[\widetilde{d}^p_{E^k}(\mathbf{x}^k,\mathbf{y}^k) := \sum_{i=1}^k |x^{i}-y^{i}|^p, \]
    so that $\widetilde{W}^p_p(\mu,\nu)=kW_p^p(\mu,\nu)$.
    \item The continuous path space $\mathscr{E}=C\big([0,T],\mathscr{F}\big)$ is endowed with the uniform topology 
    \[d_\mathscr{E}\big((X_t)_t,(Y_t)_t\big) := \sup_{t\in [0,T]} d_\mathscr{F}(X_t,Y_t) \]
    Two important cases are $\mathscr{F}=E$ and $\mathscr{F}=E^k$ for $k\in \N$ and is useful to note that
    \[C\big([0,T],E^k\big)\simeq C\big([0,T],E\big)^k \]
    \item The Skorokhod space $\mathscr{E}=D([0,T],\mathscr{F})$ is endowed with the Skorokhod distance (see Section \cref{appendix:skorokhod}). It is often more convenient to use the uniform topology although it does not make the space complete. However, as the uniform topology is stronger than the Skorokhod topology, any estimate in Wasserstein distance for the uniform topology implies the same estimate for the Skorokhod topology, see \cite[Section 3]{andreis_mckeanvlasov_2018}. 
    \item When $\mathscr{E}=\pb(\mathscr{F})$ is a space of probability measures over a space $\mathscr{F}$ which is typically one of the aforementioned spaces, we will mainly encounter three cases: 
    \[\mathcal{W}_{p} := W_{W_p,p} \qquad \mathcal{W}_{D_1} := W_{D_1,1} \qquad \mathcal{W}_{H^{-s}} := W_{H^{-s},2} \]
\end{itemize}
\end{definition}

\subsubsection{Convergence in the space of probability measures}\label[I]{sec:convergenceproba}

Since $\pb ( \mathscr{E} )$ is a subset of $ C_b ( \mathscr{E})^{\star}$, a weak topology is induced by the weak-$\star$ topology on $C_b ( \mathscr{E} )^{\star}$.

\begin{definition}(Weak convergence)\label[I]{def:convergenceproba} The weak convergence of a sequence of probability measures $( \mu_N )_N$ towards $\mu \in \pb ( \mathscr{E} )$ is defined as the related weak-$\star$ convergence in $C_b ( \mathscr{E} )^{\star}$. More precisely, a sequence of probability measures $(\mu_N)_N$ is said to converge weakly towards $\mu$ when 
\[\forall \varphi \in C_b ( \mathscr{E} ), \quad \langle \mu_N , \varphi \rangle \underset{N\to+\infty}{\longrightarrow} \langle \mu , \varphi \rangle.\]
The corresponding topology is the \emph{weakest} topology which makes the evaluation maps $\nu \mapsto \langle \nu , \varphi \rangle$ measurable. In probability theory, the related convergence for $\mu_N$-distributed random variables is also called \emph{convergence in law} or \emph{convergence in distribution} (see for instance \cite[Section 3]{billingsley_convergence_1999}).
\end{definition}

In many examples, the set of continuous bounded test functions is too large and it is necessary to work with a smaller subspace, for instance the domain of a generator. The minimal needed assumptions on a subspace of test functions are given by the following definition (see \cite[p.112]{ethier_markov_1986}). 

\begin{definition}[Separating and convergence determining class]\label[I]{def:separatingconvergencedetermining} A subset $\mathcal{F}\subset C_b(\mathscr{E})$ is called \emph{separating} whenever for all $\mu,\nu\in\pb(\mathscr{E})$, the condition 
\[\forall \varphi\in\mathcal{F},\quad \langle \mu,\varphi\rangle=\langle \nu,\varphi\rangle,\]
implies that $\mu=\nu$. The subset $\mathcal{F}$ is said to be \emph{convergence determining} whenever for any sequence $(\mu_N)_N$ in $\pb(\mathscr{E})$ and $\mu\in\pb(\mathscr{E})$, the condition 
\[\forall \varphi\in\mathcal{F}, \quad \langle \mu_N,\varphi\rangle\underset{N\to+\infty}{\longrightarrow}\langle \mu,\varphi\rangle,\]
implies that $\mu_N\to \mu$ weakly. Note that a convergence determining set is also separating (the converse is false in general). 
\end{definition}

\begin{example}\label[I]{ex:convergencedeterminingsets} The following sets are convergence determining. 
\begin{itemize}
\item By the Portmanteau theorem \cite[Theorem 2.1]{billingsley_convergence_1999}, the set $UC_b(\mathscr{E})$ of bounded uniformly continuous functions on $\mathscr{E}$ is convergence determining (for any equivalent metric on $\mathscr{E}$). 
\item When $\mathscr{E}$ is locally compact, the space $C_c(\mathscr{E})$ of continuous functions with compact support is convergence determining \cite[Chapter 3, Proposition 4.4]{ethier_markov_1986}. The space $C_0(\mathscr{E})$ of continuous functions vanishing at infinity is thus also convergence determining. Note that the space $\pb(\mathscr{E})$ is not a closed subspace of $C_c(\mathscr{E})^\star$ (nor of $C_0(\mathscr{E})^\star$). 
\item When $\mathscr{E}$ is locally compact, the Stone-Weierstrass theorem implies that $C_0(\mathscr{E})$ is separable. Thus, any dense countable subset $\mathcal{F}=\{\varphi^k,\,\,k\in\N\}\subset C_0(\mathscr{E})$ is convergence determining. Without loss of generality we can assume that $\|\varphi^k\|_\infty\leq 1$ for all $k\in\N$. Consequently, the distance \cref{eq:D0expression} metricises the weak-convergence. Indeed, since each term of the series \cref{eq:D0expression} is bounded by $2^{-k}$, the convergence $D_p(\mu_N,\mu)\to 0$ as $N\to+\infty$ is equivalent to $\langle \varphi^k,\mu_N\rangle\to\langle \mu,\varphi^k\rangle$ for all $k\in\N$. It is also possible to take $\varphi^k$ Lipschitz with a Lipschitz constant bounded by 1 for all $k\in\N$ and vanishing at infinity. 
\item In general $C_b(\mathscr{E})$ is not separable so there is no obvious other countable convergence determining set. There exists nevertheless another classical choice when $\mathscr{E}$ is only separable. By a theorem due to Urysohn, any separable metric space can be topologically imbedded in $[0,1]^\N$ and it is therefore possible to construct on $\mathscr{E}$ an equivalent metric $\tilde{d}_\mathscr{E}$ which makes $(\mathscr{E},\tilde{d}_\mathscr{E})$ a totally bounded set. The completion $\tilde{\mathscr{E}}$ of this space is therefore compact and the set $UC_b(\mathscr{E})$ under this metric is isomorphic to the set $C_b(\tilde{\mathscr{E}})$ which is separable since $\tilde{\mathscr{E}}$ is compact (by Stone-Weierstrass theorem). In conclusion, there exists a countable dense subset  $\mathcal{F}=\{\varphi^k,\,\,k\in\N\}$ in $UC_b(\mathscr{E})$. Up to replacing $\varphi^k$ by $\varphi^k/\|\varphi^k\|_\infty$ one can assume that the $\varphi^k$ are bounded by 1 and the distance \cref{eq:D0expression} thus metricises the weak convergence, see \cite[Theorem 6.6]{parthasarathy_probability_1967} and \cite[Theorem 1.1.2]{stroock_multidimensional_1997}. Note that $(\pb(\mathscr{E}),D_1)$ is separable and $D_1$ is equivalent to a complete metric, see the remark which follows \cite[Theorem 1.1.2]{stroock_multidimensional_1997} and \cite[Remark 3.2.2]{dawson_measure-valued_1993}.
\item Since the space $\mathrm{Lip}(\mathscr{E})$ is dense in the space $C_b(\mathscr{E})$, the functions $\varphi^k$ in the above examples can be taken Lipschitz (with a Lipschitz constant bounded by 1). 
\end{itemize}
\end{example}

The weak convergence is thus metricised by a $D_1$ distance. Since this distance is weaker than the Wasserstein-1 distance (it can be seen by Proposition \cref{prop:dualitywasserstein}), this implies that the topology induced by the Wasserstein distance is stronger than the topology induced by the weak convergence. The topology induced by the Wasserstein distance is described by the following theorem.

\begin{theorem}[Wassertein topology] Let $(\mathscr{E},d_\mathscr{E})$ be a Polish space and $p \geq 1$. The Wasserstein distance $W_{d_\mathscr{E},p}$ metricises the weak convergence in $\pb_p(\mathscr{E})$, defined as the convergence against bounded continuous test functions and the convergence of the $p$-th moments. 
\end{theorem}

\begin{proof} \cite[Theorem 6.9]{villani_optimal_2009}\end{proof}

In the following, an important case is the case $\mathscr{E}=\pb(E)$. Weak convergence of measures in $\pb ( \pb ( E ) )$ is thus defined as the convergence against test functions in $C_b (  \pb ( E ) )$. Their representation is not intuitive, except for \emph{linear} test functions of the kind $\mu \mapsto \langle \mu, \varphi \rangle$ where $\varphi$ belongs to $C_b \left(E\right)$. The following results (stated in a more probabilistic framework) show that these functions are sufficient to prove weak convergence. 

\begin{proposition}[Measure-valued convergence in $\mathcal{W}_{D_1}$] \label[I]{prop:randomCV} Let $D_1$ be a distance given by \cref{eq:D0expression} and Example \cref{ex:convergencedeterminingsets} which metricises the weak convergence on $\pb(E)$. Consider a sequence $( \mu_N )_N$ of $\pb ( E )$-valued random variables and another random probability measure $\mu$. The following properties hold.  
\begin{enumerate}[(i)]
    \item If $\mathcal{W}_{D_1} ( \mathrm{Law} ( \mu_N ) , \mathrm{Law} ( \mu ) ) \to 0$ as $N\to+\infty$ then $(\mu_N)_N$ converges in law towards $\mu$.
    \item If $\E| \langle \mu_N - \mu , \varphi \rangle| \to 0$ as $N\to+\infty$ for all $\varphi\in UC_b(E)$, then it holds that $\mathcal{W}_{D_1} ( \mathrm{Law} (\mu_N) , \mathrm{Law} (\mu) ) \to 0$ and $(\mu_N)_N$ converges in law towards $\mu$. 
\end{enumerate}
\end{proposition}

\begin{proof} Let us recall \cite[Theorem 6.1]{ parthasarathy_probability_1967} that the space of bounded uniformly continuous functions is convergence determining. Thus, let $\Phi \in C_b ( \pb ( E ))$ be a function which is uniformly continuous for the metric $D_1$. For any $\varepsilon>0$, there exists $\delta(\varepsilon)>0$ such that for any $\mu,\nu\in\pb(E)$, 
\[{D_1}(\mu,\nu)\leq \delta(\varepsilon)\Rightarrow |\Phi(\mu)-\Phi(\nu)| \leq \varepsilon.\]
The first point then directly stems from the Markov inequality:
\begin{align*}
|\langle \mathrm{Law} ( \mu_N ) - \mathrm{Law} ( \mu ) , \Phi \rangle | &\leq \mathbb{E} | \Phi ( \mu_N) - \Phi( \mu) | \leq \varepsilon + 2 \| \Phi \|_{\infty} \mathbb{P} {\left( | \Phi ( \mu_N) - \Phi( \mu) | \geq \varepsilon \right)} \\
&\leq \varepsilon + \frac{2 \| \Phi \|_{\infty}}{\delta(\varepsilon)} \mathbb{E}D_1 ( \mu_N , \mu ).    
\end{align*} 
Since this is true for any $\mathrm{Law} ( \mu_N ),\mathrm{Law} ( \mu )$-distributed random variables $\mu_N$, $\mu$ this finally gives
\[ | \langle \mathrm{Law} ( \mu_N ) - \mathrm{Law} ( \mu ) , \Phi \rangle | \leq  \varepsilon + \frac{2 \| \Phi \|_{\infty}}{\delta(\varepsilon)} \mathcal{W}_{D_1} {\left( \mathrm{Law} ( \mu_N ) , \mathrm{Law} ( \mu ) \right)}, \]
and the conclusion follows. For the second point, using the expression \cref{eq:D0expression} for $D_1 \left( \mu_N , \mu \right)$, the monotonic convergence theorem gives
\[ \mathcal{W}_{D_1} ( \mathrm{Law} ( \mu_N ) , \mathrm{Law} ( \mu ) ) \leq \sum_{k=1}^{+\infty} \frac{1}{2^k} \mathbb{E}| \langle \mu_N - \mu , \varphi^k \rangle|, \]
and then the dominated convergence theorem concludes the proof. 
\end{proof}

\begin{remark}[Comparison to $\mathcal{W}_{1}$] For $E$ locally compact, it has been proven at the same time that
\[ \mathcal{W}_{D_1} ( \mathrm{Law} ( \mu_N ) , \mathrm{Law} ( \mu ) ) \leq \sup_{\| \varphi \|_{\mathrm{Lip}} \leq 1} \mathbb{E} | \langle \mu_N - \mu, \varphi \rangle |. \]
Since 
\[ \sup_{\| \varphi \|_{\mathrm{Lip}} \leq 1} \mathbb{E} | \langle \mu_N - \mu, \varphi \rangle | \leq \mathbb{E} {\left[ \sup_{\| \varphi \|_{\mathrm{Lip}} \leq 1} | \langle \mu_N - \mu, \varphi \rangle | \right]} = \mathbb{E} W_1 ( \mu_N , \mu ), \]
this pinpoints, taking the infimum on the $\mathrm{Law} ( \mu_N ),\mathrm{Law} ( \mu )$-distributed random variables $\mu_N$, $\mu$, that $\mathcal{W}_{1}$ is stronger than $\mathcal{W}_{D_1}$ and both are stronger than the weak convergence on $\pb ( \pb ( E ) )$.
\end{remark}

\begin{corollary}[Sufficient conditions in a deterministic case]
With the same assumptions as above, if $\mu$ is a deterministic $\pb ( E )$-valued random variable (i.e. $\mathrm{Law}(\mu) \in \pb ( \pb ( E ) )$ is a Dirac mass), then the following assertions are equivalent.
\begin{enumerate}[(i)]
\item $\mathcal{W}_{D_1} ( \mathrm{Law} ( \mu_N ) , \mathrm{Law} ( \mu ) ) \to 0$ as $N\to+\infty$
\item For all bounded uniformly continuous function $\varphi$ on $E$, $\mathbb{E}| \langle \mu_N - \mu,\varphi \rangle | \to 0$ as $N\to+\infty$.
\end{enumerate}
The second assertion is also equivalent to $\mathbb{E}| \langle \mu_N - \mu,\varphi \rangle |^2 \to 0$ as $N\to+\infty$ for all bounded uniformly continuous function $\varphi$ on $E$. 
\end{corollary}

\begin{proof} The direct implication uses the fact $\mathcal{W}_{D_1}$ metricises the convergence in law of measure-valued random variables and that $\nu \mapsto | \langle \nu - \mu , \varphi \rangle |$ is continuous for the weak-$\star$ topology on $\pb ( E )$ when $\mu$ is deterministic. The converse implication is the second point of the previous proposition.
\end{proof}

The previous results can be found in \cite[Section 2]{del_moral_measure-valued_1998} or in \cite[Section 5, Lemma~10]{villani_limite_2001} for an equivalent argument in $E=\R^d$. In the previous lemma, only \emph{linear} test functions are used. This notion can be generalised by considering the algebra of \emph{polynomials} on $\pb(E)$. Its definition and main properties are stated in the following lemma. 

\begin{lemma}\label[I]{lemma:polynomial} Let $E$ be a Polish space. For $k\in\N$ and $\varphi_k\in C_b(E^k)$, the monomial function of order $k$ on $\pb(E)$ is defined by: 
\[R_{\varphi_k} : \pb(E) \to\R, \quad \mu \mapsto \langle \mu^{\otimes k}, \varphi_k\rangle.\]
The linear span of the set of monomial functions is called the algebra of polynomial functions on $\pb(E)$. The following properties hold.  
\begin{enumerate}[(i)]
    \item Every monomial is bounded and continuous on $\pb(E)$ for the weak topology.
    \item The algebra of polynomial functions is a convergence determining subset of $C_b(\pb(E))$. 
    \item If $E$ is compact then the algebra of polynomials is dense in $C_b(\pb(E))$. 
\end{enumerate}
\end{lemma}

\begin{proof} 
\begin{enumerate}[(i)]
\item First it is clear that every monomial and thus every polynomial is bounded. Let $(\mu_N)_N$ be a sequence in $\pb(E)$ and $\mu\in\pb(E)$ such that $\mu_N\to \mu$ as $N\to+\infty$. For any $k\in\N$, and any tensorized test function $\varphi_k\in\varphi^1\otimes\ldots\otimes\varphi^k\in C_b(E^k)$, it holds that 
\[ \langle \mu^{\otimes k}_N, \varphi_k\rangle = \prod_{j=1}^k \langle \mu_N,\varphi^j\rangle \underset{N\to+\infty}{\longrightarrow} \langle \mu^{\otimes k}, \varphi_k\rangle.\]
Then using \cite[Chapter 3, Proposition 4.6]{ethier_markov_1986}, the set $C_b(E)^{\otimes k}\subset C_b(E^k)$ is convergence determining and thus $\mu^{\otimes k}_N\to \mu^{\otimes k}$. It implies that for all $\varphi_k\in C_b(E^k)$, $R_{\varphi_k}(\mu_N)\to R_{\varphi_k}(\mu)$ and $R_{\varphi_k}$ is therefore continuous. 
\item From the first point, the algebra of polynomials is a subset of $C_b(\pb(E))$. Let $D_1$ be a metric of the form \cref{eq:D0expression} such that $(\pb(E),\widetilde{D}_1)$ is a Polish space for a metric $\widetilde{D}_1$ which is equivalent to $D_1$. A set of functions $\mathcal{F}\subset C_b(\pb(E))$ is said to strongly separates points when for every $\mu\in \pb(E)$, and $\delta>0$, there exists a finite set $\{\Phi^1,\ldots,\Phi^k\}\subset \mathcal{F}$ such that 
\[\inf_{\nu: D_1(\mu,\nu)\geq\delta}\max_{1\leq j\leq k} |\Phi^j(\nu)-\Phi^j(\mu)|>0.\]
Since $D_1$ is equivalent to a complete metric, \cite[Chapter 3, Theorem 4.5]{ethier_markov_1986} states that $\mathcal{F}$ is convergence determining if $\mathcal{F}$ strongly separates points. It is thus enough to prove that the algebra of polynomial functions contains a subset which strongly separates points. Let $(\varphi^k)_k$ be the sequence of functions in $C_b(E)$ which defines $D_1$. Then the set $\{R_{\varphi^k},\,\,k\in\N\}\subset C_b(\pb(E))$ strongly separates points. Indeed, let $\mu\in \pb(E)$, let $\delta>0$ and let $m\in\N$ be such that $2^{-m}<\delta/4$. For any $\nu\in\pb(E)$ such that $D_1(\mu,\nu)\geq\delta$, it holds that
\[\sum_{k=1}^m |\langle \mu,\varphi^k\rangle-\langle \nu,\varphi^k\rangle| \geq \frac{\delta}{2},\]
and hence $\max_{1\leq k\leq m} |\langle \mu,\varphi^k\rangle-\langle \nu,\varphi^k\rangle| \geq \delta/(2m)$. The conclusion follows. 
\item This follows from the Stone-Weierstrass theorem, since $\pb(E)$ is compact in this case. 
\end{enumerate}
\end{proof}

\subsubsection{Entropic convergence}

Powerful tools to compare measures are also given by the Kullback-Leibler divergence, which is traditionally called the \emph{relative entropy} in our context, and the related Fisher information.

\begin{definition}[Entropy, Fisher information and entropic convergence]\label[I]{def:entropyfisher}
Let $\mathscr{E}$ be a Polish space. Given two probability measures $\mu,\nu \in \pb ( \mathscr{E})$ (or more generally two measures), the relative entropy is defined by
\[ H ( \nu | \mu ) := \int_{\mathscr{E}} \frac{\dd  \nu}{\dd  \mu} \log {\left( \frac{\dd \nu}{\dd  \mu} \right)} \dd  \mu , \]
where $\dd\nu/\dd\mu$ is the Radon-Nikodym derivative. When the two measures are mutually singular, by convention, the relative entropy is set to~$+\infty$ (the same holds for the Fisher information below). If moreover $E$ is endowed with a smooth manifold structure, the Fisher information can be defined as
\[ I ( \nu | \mu ) := \int_{\mathscr{E}} {\left| \nabla \log {\left( \frac{\dd  \nu}{\dd \mu} \right)} \right|}^2 \dd  \nu, \]
with the same conventions. These quantities are dimensionally super-additive, equality being achieved only for tensorized distributions, in the sense that given $\nu\in\pb(\mathscr{E}\times \mathscr{E})$ with marginals $\nu_1,\nu_2\in\pb(\mathscr{E})$ and $\mu\in\pb(\mathscr{E})$, then 
\[H(\nu|\mu\otimes\mu) \geq H(\nu_1|\mu)+H(\nu_2|\mu),\]
equality being achieved if only if $\nu=\nu_1\otimes\nu_2$. Moreover $H ( \nu | \mu ) \geq 0$  and $H ( \nu | \mu ) = 0$ if and only if $\mu = \nu$. The \emph{entropic convergence} of a sequence $(\mu_N)_N$ in $\pb (  \mathscr{E} )$ towards~$\mu$ is defined by the convergence of the relative entropy:
\[ H ( \mu_N | \mu ) \underset{N \to+ \infty}{\longrightarrow} 0. \]
\end{definition}

The relative entropy between two probability measures is also called the \emph{Kullback-Leibler divergence}.  

\begin{remark}[Towards dimension free quantities]
For $\mu^N,\nu^N \in \pb ( \E^N )$, the normalized entropy $H_N ( \nu^N | \mu^N ) := \frac{1}{N} H ( \nu^N | \mu^N )$ can be handful, since it leads to estimates which do not depend on $N$ when $\mu^N = \mu^{\otimes N}$ is tensorized; the same holds for $W_1 ( \mu^{\otimes N} , \nu^{\otimes N} )$ when $W_1$ is defined using the normalized distance on $E^N$, see \cite[Proposition 2.6]{hauray_kacs_2014}. An extension to random measures $\pi \in \pb ( \pb ( \mathscr{E}) )$ is provided in \cite{hauray_kacs_2014} setting $\mathcal{H} ( \pi ) = \mathbb{E}_{\nu \sim \pi} H ( \nu | \mu )$ for a given $\mu \in \pb ( \mathscr{E} )$.
\end{remark}

The entropic convergence is stronger than the strongest distance.

\begin{proposition}[Pinsker inequality]
The following inequality implies that the entropy convergence is stronger than the convergence in total variation norm:
\begin{equation}\label[Ieq]{eq:pinsker}  \| \mu - \nu \|_{\mathrm{TV}}^2 \leq 2 H ( \nu | \mu ). \end{equation}
\end{proposition}

The link with the Wassertein-2 distance can be recovered through the following proposition.

\begin{proposition}[HWI inequality \cite{otto_generalization_2000}]
If $\dd \nu = \e^{-\Psi(x)} \dd x$ is a probability measure on $\R^d$ with finite second order moment, $C^2$-regularity for $\Psi$ and $D^2 \Psi \geq \lambda I_d$ for some real $\lambda$, then for any $\mu \in \pb(\R^d)$ which is absolutely continuous with respect to $\nu$, it holds that
\[ H ( \nu | \mu ) \leq W_2 ( \mu,\nu ) \sqrt{ I ( \nu | \mu ) } - \frac{\lambda}{2} W_2^2 ( \mu,\nu ). \]
\end{proposition}

Further results which link the relative entropy and the distances on $\pb(\mathscr{E})$ will be given in Section \cref{concentrationI}. The relative entropy will play an important role in Section~\cref{sec:chaosfromentropy}.

\subsection{Representation of symmetric particle systems}

This section introduces the various points of view to describe a system of particles. So far, we have mainly discussed the case of finite systems described by the $N$-particle distribution function $f^N_t$ at time $t$. The first Section \cref{sec:finitesystems} will detail more of its properties. Then, as the goal is to deal with the limit $N\to+\infty$, a framework for infinite particle systems is needed, this will be described in Section \cref{sec:infinitesystems}. 

\subsubsection{Finite particle systems}\label[I]{sec:finitesystems}

Let $N\in\N$ be a fixed finite number of particles. In full generality, there is only one property of the $N$-particle distribution function that is always true: at any time and for any of the models considered, it is a symmetric probability measure on $E^N$ (the particle system is said to be exchangeable). Let us therefore consider in this section a symmetric probability measure $f^N\in\pb_\mathrm{sym}(E^N)$ (in a static framework, it does not depend on the time). There exist two main representations of $f^N$ which are based on this symmetry assumption. 

\subsubsection*{The marginal distributions and the BBGKY hierarchy.} 

The symmetry assumption implies that for any $k\leq N$, we can define the $k$-th marginal distribution on $E^k$ by:
\[\forall\varphi_k\in C_b(E^k),\quad \langle f^{k,N}, \varphi_k\rangle = \langle f^N, \varphi_k\otimes 1^{\otimes (N-k)}\rangle,\]
and $f^{k,N}\in \pb_{\mathrm{sym}}(E^k)$ is itself a symmetric probability measure. The $N$-th marginal is of course the measure $f^N$ itself. However, keeping in mind that the final goal is to take $N\to+\infty$, one can consider for any fixed $k\in\N$ the limit of $f^{k,N}$ in $\pb(E^k)$, which is not possible for $f^N$ directly since it belongs to a space which depends on~$N$. As we shall see in the following, it is often enough to treat the case $k=2$.

In a dynamic framework, when $f^N_t$ solves the Liouville equation \cref{eq:weakliouville}, for each given $k\in\N$, a natural idea is to derive an equation for the $k$-th marginal distribution by considering a test function in \cref{eq:weakliouville} of the form $\varphi_N=\varphi_k\otimes 1^{\otimes (N-k)}$ with $\varphi_k\in C_b(E^k)$. For Boltzmann models, this computation as already been sketched in Section \cref{sec:boltzmanngeneral} and gave:
\begin{multline}\label[Ieq]{eq:bbgkyboltzmannweak}\frac{\dd}{\dd t} \langle f^{k,N}_t,\varphi_k\rangle = \sum_{i=1}^k \langle f^{k,N}_t, L^{(1)}\diamond_i\varphi_k\rangle +\frac{1}{N}\sum_{1\leq i<j\leq k}\langle f^{k,N}_t, L^{(2)}\diamond_{ij}\varphi_k\rangle\\ + \frac{N-k}{N}\sum_{i=1}^k\langle f^{k+1,N}_t, L^{(2)}\diamond_{i,k+1}(\varphi_k\otimes 1)\rangle.
\end{multline}
For mean-field systems, let us look at the special linear case: 
\[\forall\varphi\in \mathcal{F},\quad L_\mu\varphi(x) = \int_E \widetilde{L}_x\varphi(y)\mu(\dd y),\]
where for any $x\in E$, $\widetilde{L}_x$ is a Markov operator on $\mathcal{F}$ (such that for all $\varphi\in\mathcal{F}$ and $y\in E$, the map $x\mapsto \widetilde{L}_x\varphi(y)$ is measurable). Then one can check similarly (using the symmetry of $f^N_t$) that the $k$-th marginal of the Liouville equation satisfies: 
\begin{multline}\label[Ieq]{eq:bbgkymeanfieldweak}
\frac{\dd}{\dd t} \langle f^{k,N}_t, \varphi_k\rangle = \frac{1}{N}\sum_{1\leq i,j\leq k} \int_{E^k} \widetilde{L}_{x^i}\diamond_i \varphi_k\big(\widehat{\mathbf{x}}^{i,j}\big)f^{k,N}_t\big(\dd\mathbf{x}^k\big) \\+
\frac{N-k}{N} \sum_{i=1}^k \int_{E^{k+1}} \widetilde{L}_{x^i}\diamond_i \varphi_k\big(\widehat{\mathbf{x}}^{i,k+1}\big)f^{k+1,N}_t\big(\dd\mathbf{x}^{k+1}\big),
\end{multline}
where we recall the notation $\mathbf{x}^k=(x^1,\ldots,x^k)\in E^k$ and for $i\leq k$, $\widehat{\mathbf{x}}^{i,j}$ denotes the vector in $E^k$ where the $i$-th element is replaced by $x^j$. 

In both equations \cref{eq:bbgkyboltzmannweak} and \cref{eq:bbgkymeanfieldweak}, the important point to notice is that the leading term (in $N$) on the right-hand side depends on the $(k+1)$-th marginal. Since $f^{k,N}_t$ depends on $f^{k+1,N}_t$ for any $k<N$, this gives a hierarchy of $N$ equations, the $N$-th one being the Liouville equation itself. This hierarchy is called the BBGKY hierarchy, from the names of the mathematicians Bogolioubov, Born, Green, Kirkwood and Yvon. It is more customary to write the BBGKY in the forward form. For the classical Boltzmann model of Section \cref{sec:Boltzmannclassicalmodels}, using the notations of \cref{eq:Boltzmannphysics} and \cref{eq:postcollisionsigma}, the first equation of the hierarchy reads,
\begin{multline*}
    \partial f^{1,N}_t(x,v) + v\cdot\nabla_x f^{1,N}_t \\
=\frac{N-1}{N}\int_{\R^d}\int_{\mathbb{S}^{d-1}} B(v-v_*,\sigma)\Big(f^{2,N}_t(x,v_*',x,v')-f^{2,N}_t(x,v_*,x,v)\Big)\dd v_*\dd\sigma. 
\end{multline*}
We refer to the classical reference \cite{cercignani_mathematical_1994} for a more detailed derivation of the BBGKY hierarchy associated to this model and to \cite{carlen_kinetic_2013-1} for another class of Boltzmann models. 

For the mean-field case, let us consider the diffusion operator in $E=\R^d$ given by:
\[\widetilde{L}_x\varphi(y) = K(x,y)\cdot\nabla_y\varphi(y) + \Delta_y\varphi(y),\]
where $K:\R^d\times \R^d\to\R^d$ is a symmetric function with $K(x,x)=0$ for all $x\in\R^d$. Then the first equation of the BBGKY hierarchy in forward form reads: 
\[\partial_t f^{1,N}(x) = -\frac{N-1}{N}\nabla_x\cdot {\left\{\int_{\R^d} K(x,z)f^{2,N}_t(x,z)\dd z\right\}}+\Delta_x f^{1,N}_t(x).\]
Note that in both cases, if 
\begin{equation}\label[Ieq]{eq:closuref2}
f^{2,N}_t=f^{1,N}_t\otimes f^{1,N}_t,
\end{equation}
then, up to the factor $(N-1)/N$, the first marginal $f^{1,N}_t$ solves the nonlinear limit problem, respectively the Boltzmann equation \cref{eq:Boltzmannphysics} and the Fokker-Planck equation \cref{eq:mckeanvlasov-pde} (with $b(x,\mu)=K\star\mu(x)$ and $\sigma=\sqrt{2}I_d$). The relation \cref{eq:closuref2} is called a \emph{closure assumption} because under this assumption, the marginals (here the first one) satisfy a closed equation. The question of Kac's chaos and the propagation of chaos is precisely to justify this closure assumption in the asymptotic limit $N\to+\infty$. Indeed, the relation \cref{eq:closuref2} is never true since it means that any two particles are statistically independent (which is not the case since they interact). 

The BBGKY hierarchy is useful only in the linear cases described above. For general mean-field models with an operator $L_\mu$ which has a more complicated dependence in $\mu$, it is not possible to derive a BBGKY hierarchy: this procedure would only say that $f^{k,N}_t$ depends on the whole distribution $f^N_t$ which is not informative. For the Boltzmann model described above, the proof of the propagation of chaos and the justification of the closure assumption \cref{eq:closuref2} are reviewed in Section~\cref{sec:lanford} (this result is the renown Lanford's theorem). We also mention that, beyond propagation of chaos, different closure assumptions than \cref{eq:closuref2} can be considered as an approximating procedure in a numerical perspective, see for instance \cite{berlyand_continuum_2019} for the mean-field model described above. 

\subsubsection*{The empirical measure.} With a more probabilistic point of view, a symmetric measure $f^N\in\pb(E^N)$ means that any system $\mathcal{X}^N=(X^1,\ldots,X^N)\in E^N$ of $f^N$-distributed random variables is invariant in law under any permutation of the indexes. Such an exchangeable system is equivalently described by its (random) empirical measure 
\begin{equation}\label[Ieq]{eq:empiricalmeasuresymmetric}
\mu_{\mathcal{X}^N} = \frac{1}{N}\sum_{i=1}^N \delta_{X^i}\in \pb(E),\end{equation}
as this measure contains all the statistical information up to the particle numbering (a quantitative version is stated in the Lemma \cref{lemma:grunbaumlemma} below). One can immediately see the advantage of such representation: it is possible to work with only one element which belongs to the fixed space $\pb(E)$, in contrast to $f^N\in\pb(E^N)$ or to the $N$ marginal distributions. To be completely rigorous, one should work in the quotient space $E^N / \mathfrak{S}_N$, whose elements $\bar{\mathbf{x}}^N$ gather all the permutations of the vector $\mathbf{x}^N \in E^N$. There is the one to one mapping:
\begin{equation}\label[Ieq]{eq:empiricalmapEN}
\boldsymbol{\mu}_N : E^N / \mathfrak{S}_N \rightarrow \widehat{\pb}_N ( E ),\quad\bar{\mathbf{x}}^N \mapsto \mu_{\mathbf{x}^N},
\end{equation}
where $\widehat{\pb}_N(E)$ denotes the space of empirical measures of size $N$ on $E$. Since $\mu_{\mathcal{X}^N}\in\widehat{\pb}_N(E)\subset\pb(E)$ is a random element, a somehow unfortunate complication arises for the space of observables $C_b(\pb(E))$: in this framework, test functions are continuous bounded functions on (a subset of) the set of probability measures (endowed with the weak topology). This is clearly more difficult to handle than usual test functions on $E^N$ or $E^k$. 

\begin{remark} Note that the two sets $C_b ( E^N / \mathfrak{S}_N )$ and $C_b ( \widehat{\pb}_N ( E ) )$ are naturally identified by taking the composition with the previous map. Moreover, since all the measures considered are symmetric, integration on $E^N / \mathfrak{S}_N$ is equivalent to integration on $E^N$. This is why, with a slight abuse, the test functions always belong to $C_b \left( E^N \right)$.
\end{remark}

From the point of view of measure theory, the representation \cref{eq:empiricalmeasuresymmetric} means that the law $f^N$ is replaced by its push-forward by the map \cref{eq:empiricalmapEN} (seen as a map $E^N\to\pb(E)$), defined by: 
\[F^N := (\boldsymbol{\mu}_N)_\# f^N \in \pb(\pb(E)).\]
The following lemma shows that $F^N$ is enough to characterise $f^N$. 

\begin{lemma}[Approximation rate of marginals] \label[I]{lemma:grunbaumlemma} For $1 \leq k\leq N$, let the \emph{moment measure} $F^{k,N}\in\pb(E^k)$ be defined by:
\[\forall \varphi_k\in C_b(E^k),\quad \langle F^{k,N}, \varphi_k\rangle = \int_{\pb(E)} \langle \nu^{\otimes k},\varphi_k \rangle F^{N}(\dd\nu).\]
Then it holds that:
\begin{equation}\label[Ieq]{eq:approxratemarginals} \big\| f^{k,N} - F^{k,N} \big\|_{\mathrm{TV}} \leq 2\frac{k(k-1)}{N}. \end{equation}
\end{lemma}

Coming back to the probabilistic point of view, $F^N$ is the law of the random measure $\mu_{\mathcal{X}^N}$ where $\mathcal{X}^N\sim f^N$. The moment measures can thus be written 
\[F^{k,N} = \E{\left[\mu_{\mathcal{X}^N}^{\otimes k}\right]},\]
where this expression is understood in the weak sense, for all $\varphi_k\in C_b(E^k)$, 
\[\langle F^{k,N},\varphi_k \rangle = {\left\langle\E{\left[\mu_{\mathcal{X}^N}^{\otimes k}\right]},\varphi_k\right\rangle}= \E{\left[{\left\langle\mu_{\mathcal{X}^N}^{\otimes k},\varphi_k\right\rangle}\right]}.\]

\begin{proof} Given a test function $\varphi_k\in C_b ( E^k )$, using the symmetry of $f^{k,N}$, it holds that:
\[\big\langle \E_{\mathcal{X}^N} \mu_{\mathcal{X}^N}^{\otimes k} , \varphi_k \big\rangle = \int_{E^N} \big\langle \mu_{\mathbf{x}^N}^{\otimes k},\varphi_k\big\rangle f^N\big(\dd\mathbf{x}^N\big),\]
and
\[\langle f^{k,N} , \varphi_k \rangle = \int_{E^N}\frac{1}{N!}\sum_{\sigma\in\mathfrak{S}_N} \varphi_k{\left( x^{\sigma (1)} , \ldots , x^{\sigma (k)} \right) }f^N\big(\dd\mathbf{x}^N\big).\]
Consequently,
\[ \big| \big\langle f^{k,N} - \E_{\mathcal{X}^N} \mu_{\mathcal{X}^N}^{\otimes k} , \varphi_k \big\rangle \big| \leq \sup_{\mathbf{x}^N \in E^N} {\left| \frac{1}{N!} \sum_{\sigma \in \mathfrak{S}_N} \varphi_k {\left( x^{\sigma (1)} , \ldots , x^{\sigma (k)} \right)} - \big\langle \mu_{\mathbf{x}^N}^{\otimes k} , \varphi_k \big\rangle \right|}.  \]
Moreover,
\[\frac{1}{N!} \sum_{\sigma \in \mathfrak{S}_N} \varphi_k {\left( x^{\sigma (1)} , \ldots , x^{\sigma (k)} \right)} = \frac{1}{A^k_N} \sum_{ \substack{ i_1, \ldots , i_k \\ \text{pairwise distinct} } }\varphi_k {\left( x^{i_1} , \ldots , x^{i_k} \right)},\]
and
\begin{align*}
\big\langle \mu_{\mathbf{x}^N}^{\otimes k} , \varphi_k \big\rangle &= \frac{1}{N^k} \sum_{i_1, \ldots , i_k} \varphi_k {\left( x^{i_1} , \ldots , x^{i_k} \right)}\\
&= \frac{1}{N^k} \sum_{ \substack{ i_1, \ldots , i_k \\ \text{pairwise distinct} } } \varphi_k {\left( x^{i_1} , \ldots , x^{i_k} \right)} + R_{k,N},
\end{align*}
where $A^k_N:= N!/(N-k)!$. The number of pairwise distinct tuples $\left( x^{i_1} , \ldots , x^{i_k} \right)$ of integers between $1$ and $N$ being $A^k_N$, this leads $| R_{N,k} | \leq \|\varphi_k\|_{\infty}(1-A^k_N/N^k)$. Consequently, using the same argument
\[ \left| \frac{1}{N!} \sum_{\sigma \in \mathfrak{S}_N} \varphi_k {\left( x^{\sigma (1)} , \ldots , x^{\sigma (k)} \right)} - \big\langle \mu_{\mathbf{x}^N}^{\otimes k} , \varphi_k \big\rangle \right| \leq 2 \| \varphi_k \|_{\infty} \left( 1 - \frac{A^k_N}{N^k} \right) \, .\]
The conclusion follows by noticing that
\[   1 - \frac{A^k_N}{N^k} \leq 1 - {\left(1-\frac{k-1}{N}\right)}^k\leq \frac{k(k-1)}{N},  \]
and using Proposition \cref{prop:dualitywasserstein}.
\end{proof}

This (elementary) lemma is known at least since \cite{grunbaum_propagation_1971} where it was used to prove propagation of chaos (see Section \cref{sec:abstractmischlermouhot}). This lemma can also be seen as finite system version of the de Finetti theorem, see \cite[Theorem 13]{diaconis_finite_1980}. The case of infinite systems is discussed in the following section. Note that the result of \cite[Theorem 13]{diaconis_finite_1980} is actually an \emph{existence} result for a measure $F^N\in\pb(\pb(E))$ which satisfies \cref{eq:approxratemarginals}. Finally the empirical measure map is an isometry for the Wasserstein distance as shown in the following proposition.  

\begin{proposition}[Proposition 2.14 in \cite{hauray_kacs_2014}] \label[I]{prop:W1measure}
Let $f^N,g^N$ be two symmetric probability measures on $E^N$ and let $F^N=(\boldsymbol{\mu}_N)_\# f^N$ and $G^N=(\boldsymbol{\mu}_N)_\# g^N$ be the associated empirical law in $\pb(\pb(E))$. Then it holds that
\[ W_1 \big( f^N , g^N \big) = \mathcal{W}_1 \big( F^N,G^N \big).\]
\end{proposition}

This result also holds for the Wasserstein-2 distance \cite[Lemma 11]{carrillo_-convexity_2020}.

\subsubsection{Infinite particle systems and random measures}\label[I]{sec:infinitesystems}

In the previous section, finite exchangeable particle systems are described either by the marginal distributions or by the empirical measure. In this section the framework to take the limit $N\to~+\infty$ is presented. An infinite set of exchangeable random variables $(\overline{X}{}^1,\overline{X}{}^2,\ldots)$ is described by one of the two following objects. 
\begin{enumerate}
    \item The infinite hierarchy of marginals distributions $f^k\in\pb_\mathrm{sym}(E^k)$, $k\in\N$ such that 
    \[f^k = \mathrm{Law}{\left(\overline{X}{}^1,\ldots,\overline{X}{}^k\right)}.\]
    They satisfy the compatibility relation: for every $1 \leq j \leq k$, 
    \begin{equation}\label[Ieq]{eq:compatibilitymeasures}
        \forall \varphi_j\in C_b(E^j),\quad \langle f^k,\varphi_j\otimes 1^{\otimes (k-j)}\rangle = \langle f^j,\varphi_j\rangle.
    \end{equation}
    In other words, the $j$-particle marginal of $f^k$ is $f^j$.
    \item The infinite sequence of random empirical measures of size $N$, $N\in\N$,
    \begin{equation}\label[Ieq]{eq:empiricalmeasureinfinite}\mu_{\overline{\mathcal{X}}{}^N} = \frac{1}{N}\sum_{i=1}^N\delta_{\overline{X}{}^i}.\end{equation}
\end{enumerate}
The two representations are linked by the de Finetti and Hewitt-Savage theorems stated below. Let us first state some preliminary useful results.

Given an infinite system of exchangeable particles $( \overline{X}{}^i )_{i \geq 1}$, important measurable events are given by two particular $\sigma$-algebras.

\begin{definition}[Symmetric and asymptotic $\sigma$-algebras] Let $C_{\mathrm{sym}} ( E^N )$ denote the set of symmetric continuous $\R$-valued functions on $E^N$ which are invariant under permutations of their arguments. 
\begin{itemize}
    \item The $\sigma$-algebra of \emph{exchangeable} events (i.e. events which do not depend on any finite permutation of the $\overline{X}{}^i$) is defined by:
\[ \mathcal{S}_\infty := \bigcap_{k \geq 1} \sigma \Big( \sigma \big( \varphi_k ( \overline{X}{}^1 , \ldots , \overline{X}{}^k ) , \varphi_k \in C_{\text{sym}} ( E^k ) \big) , \overline{X}{}^{k+1} , \overline{X}{}^{k+2} , \ldots \Big), \]
where we recall that $\sigma(X^1,X^2,\ldots)$ is the $\sigma$-algebra generated by the random variables $X^1, X^2,\ldots$.
    \item The \emph{asymptotic} $\sigma$-algebra (whose events do not depend on any finite number of the $\overline{X}{}^i$) is defined by:
\[ \mathcal{A}_\infty := \bigcap_{k \geq 1} \sigma {\left( \overline{X}{}^{k+1} , \overline{X}{}^{k+2} , \ldots \right)}.\]
\end{itemize}
\end{definition}

The fundamental result for exchangeable systems is the following proposition. 

\begin{proposition}[\cite{letta_sur_1989}]
For exchangeable systems, the following equality holds
\[ \mathcal{S}_\infty = \mathcal{A}_\infty. \]
\end{proposition}

\begin{corollary}[Hewitt-Savage 0-1 law]\label[I]{thm:hewittsavage01}
In the special case where the $\overline{X}{}^i$ are i.i.d. variables (and then automatically exchangeable), then any event in the $\sigma$-algebra $\mathcal{S}_\infty$ or in the $\sigma$-algebra $\mathcal{A}_\infty$ has measure 0 or 1. This is known as the Kolmogorov~0-1 law for $\mathcal{A}_\infty$ and the Hewitt-Savage 0-1 law for $\mathcal{S}_\infty$.
\end{corollary}

Since the empirical measures \cref{eq:empiricalmeasureinfinite} are random measures, a criteria for the convergence in law in $\pb(\pb(E))$ is often needed; this motivate the following results. A thorough discussion of the theory of random measures can be found in \cite{dawson_measure-valued_1993}. An important notion is the notion of moment measure already introduced earlier and properly defined below. 

\begin{definition}[Moment measures]\label[I]{def:momentmeasure}
For $k\in\N$, the $k$-th moment measure of a measure $\pi\in\pb (\pb ( E ) )$ is defined by:  
\[ \pi^k := \int_{\pb ( E )} \nu^{\otimes k} \pi ( \dd  \nu ) = \mathbb{E}_{\nu \sim \pi} {\left[ \nu^{\otimes k} \right]} \in\pb(E^k). \]
This definition is understood in the weak sense, so that $\langle \pi^k , \varphi_k \rangle = \mathbb{E}_{\nu \sim \pi} \langle \nu^{\otimes k} , \varphi_k \rangle$ for any $\varphi_k$ in $C_b ( E^k )$. 
\end{definition}

Note that the sequence of moment measures $(\pi^k)_k$ satisfies the compatibility property. They also characterise the convergence in $\pb(\pb(E))$.

\begin{lemma}[Convergence of random measures] \label[I]{lemma:CVrandomM}
A sequence $( \pi_N )_N$ of random measures in $\pb ( \pb ( E ) )$ converges weakly towards $\pi \in \pb ( \pb ( E ) )$ if and only if 
\[ \forall k \geq 1, \quad \pi^k_N \underset{N\to+\infty}{\longrightarrow} \pi^k, \]
where the convergence is the weak convergence in $\pb(E^k)$. 
\end{lemma}

\begin{proof} The direct implication stems from the fact the maps $\pi \mapsto \pi^k$ are continuous for the respective weak-$\star$ topologies. For the converse, the weak convergence of $( \pi^k_N )^{}_N$ towards $\pi^k$ implies that for all $\varphi_k\in C_b(E^k)$, $\langle \pi_N,R_{\varphi_k}\rangle\to\langle \pi,R_{\varphi_k}\rangle$ where $R_{\varphi_k}$ is the monomial function:
\[ R_{\varphi_k}:\nu \in \pb ( E ) \mapsto \int_{E^k} \varphi_k ( x^1 , \ldots , x^k ) \nu^{\otimes k} ( \dd x^1, \ldots , \dd x^k ) \in\R.\]
The conclusion follows from Lemma \cref{lemma:polynomial}.
\end{proof} 

The following lemma is a useful tightness criterion in $\pb ( \pb ( E ) )$; it can be found in \cite[Proposition 2.2 (2.5)]{sznitman_topics_1991}, where the first moment measures $\pi^1$ is referred as the \emph{intensity measure} related to $\pi$ (this terminology reminiscent of the intensity of a Poisson random measure).

\begin{lemma}[Tightness for random measures] \label[I]{lemma:empiricaltightness}
The tightness of a sequence $( \pi_N )_N$ in $\pb ( \pb ( E ) )$ is equivalent to the tightness of the sequence $( \pi^1_N )^{}_N$ in $\pb ( E )$.
\end{lemma}

\begin{proof} The direct implication stems from the fact the map $\pi \mapsto \pi^1$ is continuous for the respective weak-$\star$ topologies. For the converse, assume the $( \pi^1_N )_N$ is tight. For every $\varepsilon > 0$, there exists a compact subset $K^c_{\varepsilon}\subset E$  such that $ \pi^1_N ( K^c_{\varepsilon} ) \leq \varepsilon$ for every $N$. By the Markov inequality, for every $k \geq 1$ and every $N\geq1$, it holds that
\[ \pi_N {\left(  {\left\{ \nu \in \pb ( E ),\, \, \, \nu {\left( K^c_{\varepsilon ( k 2^k )^{-1}} \right)} \geq \frac{1}{k} \right\}} \right)} \leq k \pi^1 {\left( K^c_{\varepsilon ( k 2^k )^{-1}} \right)} \leq \frac{\varepsilon}{2^k},  \]
so that
\[ \pi_N {\left( \bigcap_{k \geq 1} {\left\{ \nu \in \pb ( E ),\,\, \, \nu {\left( K^c_{\varepsilon ( k 2^k )^{-1}} \right)} \leq \frac{1}{k} \right\}}  \right)} \geq 1 - \sum_{k \geq 1} \frac{\varepsilon}{2^k} = 1 - \varepsilon. \]
Since the intersection at the last line is a compact subset of $\pb(E)$, the sequence $(\pi_N)_N$ is tight. 
\end{proof}

\begin{example}[The case of empirical measures]
This lemma is particularly interesting for random empirical measures $\mu_{\mathcal{X}^N}$, since it reduces the question of tightness of $( \mu_{\mathcal{X}^N} )_N$ in $\pb(\pb(E))$ to tightness of $( X^{1,N} )_N$ in $\pb(E)$.
\end{example}

The following theorems are the two main results of this section. The first one states that (the law of) a random measure can always be represented by an infinite exchangeable particle system. This theorem is due to de Finetti and can be found in \cite[Theorem 11.2.1]{dawson_measure-valued_1993}. 

\begin{theorem}[De Finetti representation theorem for random measures] \label[I]{thm:deFinetti}
Let $\pi\in\pb ( \pb ( E ) )$. Then there exists a sequence $( \overline{X}{}^i )_{i \geq 1}$ of $E$-valued exchangeable random variables such that the following properties hold. 
\begin{enumerate}[(1)]
\item For any $k \geq 1$, $( \overline{X}{}^1, \ldots , \overline{X}{}^k )$ has joint distribution $\pi^k$.
\item The weak limit 
\[\mu = \lim_{k\to+\infty}\frac{1}{k} \sum_{i = 1}^k \delta_{\overline{X}{}^i}\in\pb(E),\]
exists almost surely and $\mu$ is $\pi$-distributed.
\item The random measure $\mu$ is $\mathcal{S}_\infty$-measurable, and conditionally on $\mathcal{S}_\infty$ the random variables $\overline{X}{}^i$ are independent and $\mu$-distributed.
\end{enumerate}
\end{theorem}

\begin{example} A famous example is given in \cite{dawson_measure-valued_1993}: a de Finetti representation of the Fleming-Viot measure-valued process is given by the Moran particle system.
\end{example}

Note that the last property says that exchangeability implies \emph{conditional} independence and thus exchangeable particles are not so far from i.i.d. variables.

Conversely, an infinite exchangeable particle system is always associated to a unique element in $\pb(\pb(E))$. The following theorem is also due to de Finetti in the case of Bernoulli random variables. It has been generalised to any exchangeable Borel measurable variables in a Polish space by Hewitt and Savage. The following quantitative version of the Hewitt-Savage theorem is due to \cite[Theorem 5.1]{hauray_kacs_2014}. 

\begin{theorem}[De Finetti, Hewitt-Savage] \label[I]{thm:quantitativehewitt} Let $E$ be a locally compact Polish space. Let $( f^N )_N$ be an infinite sequence of symmetric probability measures on $E^N$, $N\in\N$, which satisfy the compatibility relation \cref{eq:compatibilitymeasures}. Then the following properties hold.
\begin{enumerate}[(1)]
    \item There exists a unique $\pi\in\pb(\pb(E))$ such that: 
    \[f^N = \pi^N := \int_{\pb(E)}\nu^{\otimes N}\pi(\dd\nu).\]
    \item When $E\subset \R^d$ is a Borel set, for any $s>d/2$, the sequence $(  \mathrm{Law} ( \mu_{\mathcal{X}^N} ) )_{N \geq 1}$ is a Cauchy sequence in $\pb(\pb(E))$ for the distance  $\mathcal{W}_{H^{-s}}$ (see Definition \cref{def:spaceswasserstein}) : for any $N,M\geq1$,
    \[ \mathcal{W}^2_{H^{-s}} \big( \mathrm{Law} ( \mu_{\mathcal{X}^N} ) , \mathrm{Law} ( \mu_{\mathcal{X}^M} ) \big) \leq 2 \| \Phi_s \|_{\infty} {\left( \frac{1}{N} + \frac{1}{M} \right)}, \]
    where $\Phi_s$ is defined by \cref{eq:polynomialsobolev} and $\mathcal{X}^N\sim f^N$. The limit of this sequence is the measure $\pi$ characterised above. 
 \end{enumerate}
\end{theorem}

\begin{proof}[Proof (some ideas)] The original argument of Hewitt and Savage is based on the Krein-Milman theorem and the fact that tensorised measures are extreme points of the convex set $\pb_{\mathrm{sym}}(E)$. A constructive quantitative approach due to Diaconis and Freedman is based on the approximation Lemma \cref{lemma:grunbaumlemma}, see \cite[Theorem 14]{diaconis_finite_1980} in the compact case. An alternative argument based on the density of polynomial functions in $\pb(E)$ (thanks to the Stone-Weierstrass theorem) is due to Pierre-Louis Lions. We refer the interested reader to \cite[Section 2.1]{rougerie_finetti_2015} and the references therein. For the second point proved in \cite[Theorem 5.1]{hauray_kacs_2014}, the Cauchy-estimates relies on the polynomial structure of the $H^{-s}$-norm \cref{eq:polynomialsobolev} combined with the observation that $f^{N+M}$ is a transference plan between $f^M$ and $f^N$ (by the compatibility property). It turns the problem into controlling $\mathbb{E} \| \mu_{\mathcal{X}^N} - \mu_{\mathcal{X}^M} \|^2_{H^{-s}} $ for $( \mathcal{X}^N , \mathcal{X}^M )\sim f^N \otimes f^M$. Once convergence is shown, the limit is identified by the moment measures and Lemma \cref{lemma:grunbaumlemma}. Convergence can be obtained in stronger metrics than $\mathcal{W}^2_{H^{-s}}$, see Corollary \cref{CauchyStrongHewit} in the appendix.
\end{proof}

\subsection{Kac's chaos}\label[I]{sec:kacchaos}

\subsubsection{Definition and characterisation}

The notion of chaos was introduced in the seminal article of Mark Kac \cite{kac_foundations_1956}.

\begin{definition}[Kac's chaos]
Let $f\in\pb ( E )$. A sequence $( f^N )_{N \geq 1}$ of symmetric probability measures on $E^N$ is said to be \emph{$f$-chaotic} when for any $k\in\N$ and any function $\varphi_k \in C_b ( E^k )$, 
\begin{equation}\label[Ieq]{eq:kacchaos} 
\lim_{N \to+ \infty} \langle f^N , \varphi_k \otimes 1^{\otimes N - k} \rangle = \langle f^{\otimes k} , \varphi_k \rangle.
\end{equation}
It means that for all $k\in\N$, the $k$-th marginal satisfies $f^{k,N}\to f^{\otimes k}$ for the weak topology. Kac's chaos can be equivalently defined by considering only tensorized test functions $\varphi_k = \varphi^1 \otimes \ldots \otimes \varphi^k$, since the algebra of tensorized functions in $C_b ( E )$ is a convergence-determining class according to \cite[Chapter 3, Theorem 4.5 and Proposition 4.6, pp.113-115]{ethier_markov_1986}.
\end{definition}

Interpreting $f^N$ as the law of an exchangeable system of $N$ particles, the property~\cref{eq:kacchaos} means that for any group of $k$ particles, the particles become statistically independent as $N$ tends to $+\infty$, hence the terminology of chaos. The results of the previous sections on finite and infinite exchangeable systems lead to the following useful characterization of Kac's chaos.  

\begin{lemma} \label[I]{lemma:caractchaos}
Each of the following assertions is equivalent to Kac's chaos.
\begin{enumerate}[(i)]
\item There exists $k \geq 2$ such that  $f^{k,N}$ converges weakly towards $f^{\otimes k}$.
\item The random measure $\mu_{\mathcal{X}^N}$ with $\mathcal{X}^N\sim f^N$ converges in law towards the deterministic measure $f$.
\end{enumerate}
\end{lemma}

This classical result can be found in \cite[Proposition 2.2]{sznitman_topics_1991}.

\begin{proof} Clearly, Kac's chaos implies $(i)$. Then using Proposition \cref{prop:randomCV}, it can be proved that $(i)\Rightarrow(ii)$. Let $\mathcal{X}^N\sim f^N$. It is enough to prove that for any $\varphi\in C_b(E)$,
\[\mathbb{E}\big| \big\langle \mu_{\mathcal{X}^N} - f , \varphi \big\rangle \big|^2\underset{N\to+\infty}{\longrightarrow}0.\]
Assume \cref{eq:kacchaos} with $k=2$; it then also holds for $k=1$. Using the symmetry of $f^N$, it holds that
\begin{align*}
\mathbb{E}\big| \big\langle \mu_{\mathcal{X}^N} - f , \varphi \big\rangle \big|^2 &=  \frac{1}{N^2} \sum_{i,j=1}^N \mathbb{E}{\left[ \varphi ( X^i ) \varphi ( X^j ) \right]} - \frac{2}{N} \langle f , \varphi \rangle \sum_{i=1}^N \mathbb{E}{ \left[ \varphi ( X^i ) \right]} + \langle f , \varphi \rangle^2 \\
&= \frac{1}{N} \mathbb{E}{\left[ \varphi ( X^1 )^2 \right]} + \frac{N-1}{N} \mathbb{E} {\left[ \varphi ( X^1 ) \varphi ( X^2 ) \right]}  - 2 \langle f , \varphi \rangle \mathbb{E} {\left[ \varphi ( X^1 ) \right]} \\ &\quad+ \langle f , \varphi \rangle^2,
\end{align*} 
where the symmetry of $f^N$ has been used. Since $\varphi$ is bounded, the first term goes to $0$ as $N \to \infty$. The remaining expression vanishes using \cref{eq:kacchaos} with $k=1,2$. This proves $(ii)$. Then the condition $(ii)$ implies Kac's chaos. If $F^N:=\mathrm{Law} {\left( \mu_{\mathcal{X}^N} \right)}\to\delta_f$, then according to Lemma \cref{lemma:CVrandomM}, the $k$-th moment measure $F^{k,N}$ converges weakly towards $f^{\otimes k}$ for every $k \geq 1$. The approximation Lemma \cref{lemma:grunbaumlemma} implies that $f^{k,N}$ converges weakly towards $f^{\otimes k}$ for every $k \geq 1$.
\end{proof}

\begin{remark}[Chaos as a limit of de Finetti representations] Kac's chaos tells that the marginals $f^{k,N}$ converge towards the marginals of an infinite system $( \overline{X}{}^i )_{i \geq 1}$ of i.i.d. $f$-distributed particles. By the de Finetti and Hewitt-Savage theorems, the sequence of empirical measures of this latter system converges towards a random measure $\mu$ which is $\mathcal{S}_{\infty} = \mathcal{A}_{\infty}$-measurable. By the Hewitt-Savage 0-1 law, this $\sigma$-algebra is trivial so that $\mu$ is a deterministic measure. The last part of de Finetti representation Theorem \cref{thm:deFinetti} tells that conditionally on $\mathcal{S}_{\infty}$, the $\overline{X}{}^i$ are $\mu$-distributed so this allows to conclude $\mathrm{Law} ( \mu ) = \delta_f$. 
\end{remark}

\begin{remark}[Chaos as a law of large numbers] \label[I]{rem:chaosasLGN}
Fix $\varphi$ in $C_b ( E )$. Given a bounded continuous function $\theta : \R \to \R$, the function $\nu \mapsto \theta ( \langle \nu , \varphi \rangle )$ is still bounded and weakly-$\star$ continuous on $\pb ( E )$. The convergence in law of $\mu_{\mathcal{X}^N}$ towards $f$ thus implies that
\[ \frac{\varphi ( X^{1,N} ) + \ldots + \varphi ( X^{N,N} )}{N} - \mathbb{E}{\left[\varphi ( X^{1,N} )\right]} = \big\langle \mu_{\mathcal{X}^N} , \varphi \big\rangle - \big\langle f^{1,N} , \varphi \big\rangle \underset{N \to +\infty}{\longrightarrow} 0, \]
where the convergence is the convergence in law. This relation is reminiscent of the law of large numbers. If the $X^i$ were moreover i.i.d. (in this case no need to write $X^{i,N}$, $X^i$ is enough) the law of large numbers would state 
\[ \big\langle \mu_{\mathcal{X}^N} , \varphi \big\rangle \underset{N \to +\infty}{\longrightarrow} \mathbb{E} {\left[ \varphi ( X^1 ) \right]}\quad\text{a.s.},\]
so that almost surely $\mu_{\mathcal{X}^N}  \to \mathrm{Law} ( X^1 )$ weakly. In the general case where particles $X^{i,N}$ are only exchangeable (no more i.i.d.), Kac's chaos states an analogous but weaker result since the convergence of $\mu_{\mathcal{X}^N}$ towards $f$ is only weak; but it however differs since $f$ is the law of a typical particle in the limit system, and not the law of $X^{1,N}$ as in the i.i.d. case, because $X^{1,N}$ still depends on $N$ (i.e. on the other particles). Fluctuations of $\langle \mu_{\mathcal{X}^N} , \varphi \rangle$ in the law of large numbers are described through the central limit theorem; the same can be done for chaos with concentration inequalities and large deviation principles (see Section \cref{sec:fluctuations}). 
\end{remark}

\begin{remark}[Chaos, limit hierarchy and moment measures]\label[I]{rem:chaosbbgky}
Taking (formally) the limit $N \to \infty$ in the BBGKY hierarchy (Section \cref{sec:finitesystems}) gives an infinite set of coupled equations on $( f^k_t )_{k \geq 1}$ which satisfy the compatibility relation \cref{eq:compatibilitymeasures}. This system is Kac's chaotic when this limit hierarchy has the \emph{factorisation property}, that is to say $f^k_t = f_t^{\otimes k}$ for every $k \geq 1$; this implies that the related $\pb ( \pb ( E ) )$-representation of the system is $\delta_{f_t}$. This infinite hierarchy is called the Boltzmann hierarchy in kinetic theory (see \cite{cercignani_mathematical_1994} and Section \cref{sec:lanford}). By the Hewitt-Savage Theorem \cref{thm:quantitativehewitt} it is uniquely associated to an element $\pi\in\pb(\pb(E))$. We thus point out that the Boltzmann hierarchy coincides with the system of moment measures of $\pi$: this object is also commonly used, in another context, in the study of measure-valued processes \cite{dawson_measure-valued_1993}. 
\end{remark}

The following property will be useful for time-dependent systems, its proof is straightforward, see \cite[Proposition 2.4]{sznitman_topics_1991}.

\begin{proposition}[Chaos transportation] \label[I]{prop:chaostransportation}
Let $( f^N )_N$ be a $f$-chaotic sequence and let $\mathsf{T} : E \rightarrow F$ be a $f$-almost surely continuous map between Polish spaces. Then the sequence $( \mathsf{T}_{\#}f^N )_N$ is $\mathsf{T}_{\#}f$-chaotic. 
\end{proposition}

\subsubsection{Quantitative versions of Kac's chaos}\label[I]{sec:quantitativechaos}

Kac's chaos is a non quantitative property which relies only on the weak convergence. Quantitative (stronger) versions can naturally be defined using the topological framework of Section \cref{sec:topology}. The following definitions of quantitative chaos can be found in \cite{hauray_kacs_2014}. A strating point is the notion of chaos in Wasserstein distance. In practise, $W_2$ is well-adapted to the study of diffusion processes, while $W_1$ is often used for jump processes.

\begin{definition}[Chaos in Wasserstein-$p$ distance]\label[I]{def:wassersteinchaos} Let $p\in\N$, let $(f^N)_N$ be a sequence of symmetric measures on $E^N$ and let $f\in\pb(E)$. The following three notions of chaos in Wasserstein-$p$ distance were introduced in \cite{hauray_kacs_2014}:
\begin{itemize}
    \item \emph{(Wasserstein-$p$ Kac's chaos).} For all $k\in\N$, 
    \begin{equation}\label[Ieq]{eq:Omegak}\Omega_k \left( f^N , f \right) := W_p{\left( f^{k,N} , f^{\otimes k} \right)}\underset{N\to+\infty}{\longrightarrow}0.\end{equation}
    \item \emph{(Infinite dimensional Wasserstein-$p$ chaos).} 
    \begin{equation}\label[Ieq]{eq:OmegaN}\Omega_N \left( f^N , f \right) := W_p{\left( f^{N} , f^{\otimes N} \right)}\underset{N\to+\infty}{\longrightarrow}0.\end{equation}
    \item \emph{(Wasserstein-$p$ empirical chaos).} For $\mathcal{X}^N\sim f^N$, 
    \begin{equation}\label[Ieq]{eq:Omegainfty}\Omega_\infty \left( f^N , f \right) := \mathcal{W}_p {\left( \mathrm{Law}(\mu_{\mathcal{X}^N}) , \delta_f \right)}\underset{N\to+\infty}{\longrightarrow}0,\end{equation}
    where $\mathcal{W}_p$ is a Wasserstein-$p$ distance on $\pb(\pb(E))$ (see Defintion \cref{def:spaceswasserstein} for the conventions used when $p=1,2$).
\end{itemize}
\end{definition}

When moment bounds are available and $p=1$, the three notions of chaos \cref{eq:Omegak}, \cref{eq:OmegaN} and \cref{eq:Omegainfty} are actually equivalent. Such a result would not hold if the Wasserstein-1 distance were replaced by another Wasserstein-$p$ distance.

\begin{theorem}[Equivalence in Wasserstein-1 distance]\label[I]{thm:equivalencechaosW1} Let $E=\R^d$ and let $q \geq 1$ such that the sum of moments of order $q$ of $f$ and $f^{1,N}$ are bounded by a constant $\mathscr{M}_q\in(0,\infty)$. Then for any constant $\gamma<(d+1+d/q)^{-1}$, there exists $C=C(d,q,\gamma)\in(0,\infty)$ such that for any $k,\ell\in\{1,\ldots,N\}\cup\{\infty\}$ with $\ell\ne 1$:
\[\Omega_k(f^N,f)\leq C\mathscr{M}_q^{1/q}\left(\Omega_\ell(f^N,f)+\frac{1}{N}\right)^{\gamma},\]
where $\Omega_k$ and $\Omega_\ell$ are defined in Definition \cref{def:wassersteinchaos} with $p=1$. 
\end{theorem}

\begin{proof}[Proof (some ideas)]
See \cite[Theorem 1.2]{hauray_kacs_2014} and \cite[Theorem 2.4]{hauray_kacs_2014}. In particular, the link between \cref{eq:Omegainfty} and \cref{eq:OmegaN} stems from Proposition \cref{prop:W1measure} which states that 
\[W_1 \big( f^{N} , f^{\otimes N} \big) = \mathcal{W}_1 \big( \mathrm{Law}(\mu_{\mathcal{X}^N}) , \mathrm{Law}(\mu_{\overline{\mathcal{X}}{}^N}) \big),\]
where $\overline{\mathcal{X}}{}^N\sim f^{\otimes N}$-distributed. Given such $\overline{\mathcal{X}}{}^N$, $\mathcal{W}_1 \big( \mathrm{Law}(\mu_{\bar{\mathcal{X}}^N}) , \delta_f \big) \leq \E W_1 \big( \mu_{\bar{\mathcal{X}}} , f \big)$ and the quantitative laws of large numbers from \cite{fournier_rate_2015} can be applied.
\end{proof}

\begin{definition}[Strong entropic and TV chaos]\label[Ieq]{eq:strongentropicchaos}
Stronger notions can also be defined using stronger norms. 
\begin{itemize}
    \item $(f^N)_N$ is \emph{$f$-TV chaotic} when for every $k \geq 1$, $\| f^{k,N} - f^{\otimes k} \|_{\mathrm{TV}} \to 0$ as $N \to~+ \infty$. 
    \item $(f^N)_N$ is \emph{$f$-strong entropic chaotic} when for every $k \geq 1$, $H \big(  f^{k,N} | f^{\otimes k} \big) \to 0$ as $N \to + \infty$. 
\end{itemize}
The second one is stronger than the first one by Pinsker's inequality \cref{eq:pinsker}.
\end{definition}

In Definition \cref{eq:strongentropicchaos} and in \cref{eq:Omegak}, a stronger convergence can be obtained when the fixed $k\in\N$ is replaced by a function $k\equiv k(N)$ which depends on $N$. In that case, the chaos is said to hold \emph{for blocks of size $k(N)$}. The infinite dimensional chaos \cref{eq:Omegainfty} corresponds to the case $k(N)=N$. 

When $E=\R^d$ (or $E\subset\R^d$) is endowed with the Lebesgue measure denoted by~$\sigma$ below, other stronger versions of Kac's chaos can also be defined using the notions of entropy and Fisher information. The following notions can be found in \cite{hauray_kacs_2014}. 

\begin{definition}[Entropy and Fisher chaos]\label[I]{def:entropychaotic}
Let $\sigma^N$ denote the Lebesgue measure on $E^N$. 
\begin{itemize}
\item $(f^N)_N$ is $f$\emph{-entropy chaotic} when $f^{1,N} \to f$ weakly and $\frac{H ( f^N | \sigma^N )}{N} \to H ( f | \sigma ) $.
\item $(f^N)_N$ is $f$\emph{-Fisher chaotic} when $f^{1,N} \to f$ weakly and $\frac{I ( f^N | \sigma^N )}{N} \to I ( f | \sigma ) $.
\end{itemize}
\end{definition}

Using sharp versions of the HWI inequality, these notions are classified in a quantitative way in \cite{hauray_kacs_2014}.

\begin{proposition}
Each of the below assertions implies the following. 
\begin{itemize}
\item $(f^N)_N$ is $f$-Fisher chaotic.
\item $(f^N)_N$ is $f$-Kac chaotic with $\Big( \frac{I ( f^N | \sigma^N )}{N} \Big)_N$ bounded.
\item $(f^N)_N$ is $f$-entropy chaotic.
\item $(f^N)_N$ is $f$-Kac chaotic.
\end{itemize}
\end{proposition}

In classical kinetic theory, another important notion of quantitative chaos arises when the $N$ particles are constrained to evolve on the Kac's sphere:
\[E_N=\big\{\mathbf{v}^N\in \R^N,\,\, |v^1|^2+\ldots |v^N|^2 = N\big\} \subset \R^N,\]
or on the Boltzmann's sphere: 
\[E_N=\big\{\mathbf{v}^N\in (\R^3)^N,\,\, |v^1|^2+\ldots |v^N|^2 = N,\,\,v^1+\ldots+v^N=0\big\}\subset(\R^3)^N.\]
In these cases, the adapted notions of entropy chaos and Fisher chaos using a dedicated sequence of reference measures $(\sigma^N)_N$ are defined in \cite{hauray_kacs_2014} and \cite{carlen_entropy_2008}.  

\subsection{Propagation of chaos}\label[I]{sec:poc}

This section finally presents the central concept of this review, the notion of propagation of chaos, which is a dynamical version of Kac's chaos. Let us fix a final time $T\in[0,+\infty]$ and let us write $I=[0,T]$. Let $\mathcal{X}^N_I = ( \mathcal{X}^N_t )^{}_{t \in I}$ be a time-evolving (stochastic) c\`adl\`ag system of $N$-exchangeable particles in $E$ with a $f_0$-chaotic initial distribution $f^N_0\in\pb(E^N)$ where $f_0\in\pb(E)$. One aims to compare the law of a typical particle with a limit flow of measures $( f_t )_{t \in I}$, where $f_t\in\pb(E)$. The propagation of chaos property is said to hold when the initial chaos is propagated at later times. This property can hold either at the level of the law or at the level of trajectories. 

\begin{definition}[Pointwise and pathwise propgation of chaos]\label[I]{def:poc}
Let $f^N_0 \in\pb(E^N)$ be the initial $f_0$-chaotic distribution of $\mathcal{X}^N_0$ at time $t=0$. 
\begin{itemize}
\item \emph{Pointwise propagation of chaos} holds towards a flow of measures $(f_t)_t\in C(I,\pb(E))$ when the law $f^N_t\in\pb(E^N)$ of $\mathcal{X}^N_t$ is $f_t$-chaotic for every time $t\in I$. Note that the flow of measures is continuous in time as it is the solution of a PDE, but the (random) trajectories of the particles are c\`adl\`ag.
\item \emph{Pathwise propagation of chaos} holds towards a distribution $f_I\in\pb(D(I,E))$ on the path space when the law $f^N_I\in\pb\big(D(I,E)^N\big)$ of the process $\mathcal{X}^N_I$ (seen as a random element in $D(I,E)^N$) is $f_I$-chaotic.
\end{itemize}
\end{definition} 

The pointwise level is the analytical point of view where $(f_t)_t$ is the solution of a PDE. At the pathwise level, the limit distribution $f_{[0,T]}$ is often identified as the solution of a nonlinear martingale problem. 

\subsubsection{Quantitative and uniform in time propagation of chaos.} As in Section \cref{sec:quantitativechaos}, it is possible to define quantitative versions of the propagation of chaos by using any of the quantitative notions of Kac's chaos. Then, one can wonder if the propagation of chaos holds uniformly in time, i.e. independently on $T$ and with which rate of convergence in $N$. For instance, a typical quantitative pointwise propagation of chaos estimate reads: 
\begin{equation}\label[Ieq]{eq:typicalpoc}\delta\big(f^{k,N}_t,f^{\otimes k}_t\big) \leq \varepsilon(N,k,T)\Big(1+\delta\big(f^{k,N}_0,f_0^{\otimes k}\big)\Big),\end{equation}
where $\delta$ is any of the distances on $E^k$ defined in Section \cref{sec:topology}, $k\in\N$ is fixed, $t\in[0,T]$ and $\varepsilon(N,k,T)\to0$ as $N\to+\infty$. 

Propagation of chaos is said to hold uniformly in time when $\varepsilon(N,k,T)$ does not depend on $T$. As we shall see, it is usually possible to prove propagation of chaos uniformly in time only for physical models which enjoy some conservation properties. A closely related question when propagation of chaos holds on $I = [0,+\infty)$ is the ergodicity of the process as $t\to+\infty$. For instance, one may wonder if it possible to take the double limit $N\to+\infty$ and $t\to+\infty$ in~\cref{eq:typicalpoc}. It would be possible for instance if propagation of chaos held uniformly in time and $f_t$ converged towards an equilibrium $f_\infty$ as $t\to+\infty$. This would give a relaxation estimate on $\delta(f^{k,N}_t,f_\infty^{\otimes k})$. This question is of particular importance in the study of the Boltzmann equation \cref{eq:Boltzmannphysics} in view of the famous H-theorem. Relaxation towards equilibrium at the particle level, will be mentioned in Section~\cref{sec:kacprogram} for the spatially homogeneous Boltzmann equation \cref{eq:spatiallyhomogeneousboltzmann} and in Section \cref{sec:gradientsystems} for diffusion processes associated to the granular media and Vlasov-Fokker-Planck equations. On the other hand, it is sometimes possible to prove that propagation of chaos does \emph{not} hold uniformly in time, an example is given in \cite{carlen_kinetic_2013-1,carlen_kinetic_2013}.

Finally, for any quantitative version of propagation of chaos, it has been recently shown in \cite{lacker_hierarchies_2021} that the optimal rate of convergence is 
\[\forall k\in\N,\,\,\forall \varphi_k\in C_b(E^k), \quad |\langle f^{k,N}_t,\varphi_k\rangle - \langle f_t^{\otimes k},\varphi_k\rangle| = \mathcal{O}(k/N).\]
In particular, as we shall see in the following, many proofs and in particular the widely known proof of McKean's Theorem \cref{thm:mckean} via coupling arguments (see Section~\cref{sec:coupling}) only yields the sub-optimal rate of convergence $\sqrt{k/N}$. 

\subsubsection{Pointwise from pathwise.} Pathwise propagation of chaos is more general since it keeps tracks of the whole trajectory of the particles. When pathwise propagation of chaos holds, it implies pointwise propagation of chaos: since the coordinate maps are continuous, this directly stems from Proposition \cref{prop:chaostransportation} (it is also a consequence of the convergence of the finite dimensional distributions \cite[Chpater 3, Theorem~7.8]{ethier_markov_1986}). Note however, that this does not preserve the convergence rates. The converse does not always hold (see the counterexample below). In general, pathwise results are more difficult to obtain and can be proved only on finite time intervals. Pointwise propagation of chaos provides also more flexibility since it allows to work on $C ( I , \mathcal{S} )$, where $\mathcal{S}$ may be a subset of $\pb( E )$ or a larger topological space. For instance, useful spaces to study fluctuations are the class of tempered distributions in \cite{dawson_large_1987} or negative weighted Sobolev spaces in \cite{meleard_asymptotic_1996}. In a more analytical perspective, when $(f_t)_t$ solves a known PDE, $\mathcal{S}$ is more naturally identified to a functional space, for instance a Sobolev space \cite{jourdain_propagation_1998,mischler_new_2015}.

\subsubsection{Propagation of chaos via the empirical process.}
The characterisation of Kac's chaos via the empirical measure given in Lemma \cref{lemma:caractchaos} implies that pointwise propagation of chaos is equivalent to the convergence in law of the random measure $\mu_{\mathcal{X}^N_t}$ towards the deterministic measure $f_t$ for any $t\in[0,T]$. At the pathwise level, there are two notions of pathwise empirical propagation of chaos which are presented below. To begin with, a slightly more general definition of the empirical measure map is needed. Given a set $\mathscr{E}$ the empirical measure map is defined by: 
\[\boldsymbol{\mu}^\mathscr{E}_N : \mathscr{E}^N\to \pb(\mathscr{E}),\quad \mathbf{x}^N\mapsto\mu_{\mathbf{x}^N}. \]
In the following, $\mathscr{E}$ is a Polish space; it will be either the state space $E$ (in which case we may omit the superscript $E$ as in \cref{eq:empiricalmapEN}) or the path space. The following maps link the pathwise and pointwise properties.
\begin{itemize}
    \item \textbf{(The evaluation map).} For any $t\in[0,T]$, 
    \[\mathsf{X}^\mathscr{E}_t : D([0,T],\mathscr{E}) \to \mathscr{E},\quad \omega\mapsto \omega(t).\]
    \item \textbf{(The projection map).}
    \[\Pi^\mathscr{E} : \pb\big(D([0,T],\mathscr{E})\big)\to D\big([0,T],\pb(\mathscr{E})\big),\quad \mu\mapsto \big((\mathsf{X}^\mathscr{E}_t)_\#\mu\big)_{0\leq t\leq T}.\]
\end{itemize}

The pathwise and pointwise $N$-particle distributions are linked by
\[f^N_t = \big(\mathsf{X}^{E^N}_t\big){}_\# f^N_{[0,T]} = \Pi^{E^N}{\left(f^N_{[0,T]}\right)}(t).\]
The \emph{empirical measure process} is the measure-valued process defined by: 
\[\big(\mu_{\mathcal{X}^N_t}\big)_t \equiv {\left(\boldsymbol{\mu}^E_N(\mathcal{X}^N_t)\right)_t}.\]
Since this is the image of a process this readily defines the laws for all $t\in[0,T]$,  
\[F^N_t := (\boldsymbol{\mu}^E_N)_\# f^N_t \in \pb(\pb(E)),\]
and the pathwise version: 
\[F^{\mu,N}_{[0,T]} := (\boldsymbol{\mu}^E_N\circ)_\# f^N_{[0,T]}\in \pb(D([0,T],\pb(E))),\]
where $\boldsymbol{\mu}^E_N\circ$ is the natural extension of $\boldsymbol{\mu}^E_N$ on the path space defined by:
\[\boldsymbol{\mu}^E_N\circ : D([0,T],E^N)\to D([0,T],\pb(E)),\quad \omega\mapsto \boldsymbol{\mu}^E_N\circ\omega.\]
Hence, it holds that
\[F^N_t = \Pi^{\mathcal{P}(E)}{\left(F^{\mu,N}_{[0,T]}\right)}(t)= (\mathsf{X}_t^{\pb(E)})_\#F^{\mu,N}_{[0,T]}.\]
But there is another choice: it is also possible to define the pathwise empirical distribution as the push-forward of the $N$-particle pathwise distribution by the empirical map:
\[{F}^N_{[0,T]} := (\boldsymbol{\mu}^\mathcal{D}_N)_\# f^N_{[0,T]} \in \pb(\pb(D([0,T],E))),\]
where we write $\mathcal{D}=D([0,T],E)$. This probability distribution is linked to $F^N_t$ and $F^{\mu,N}_{[0,T]}$ by: 
\[F^{\mu,N}_{[0,T]} = (\Pi^{E})_\#{F}^N_{[0,T]},\quad F^N_t=(\mathsf{X}^{\pb(E)}_t\circ \Pi^E)_\# {F}^N_{[0,T]}.\]
In summary, there are three levels of description of the empirical process. In the following diagram, each space on the top row is a probability space endowed with a probability measure which is the law of the specific version of the random empirical process in the bottom row. The spaces are linked by the maps $\Pi^E$ and $\mathsf{X}_t^{\pb(E)}$.
\begin{equation*}
\begin{array}{rcccl}
\Big(\pb(D([0,T],E)),{F}^N_{[0,T]}\Big) & \overset{\Pi^E}{\longrightarrow} & \Big(D([0,T],\pb(E)),F^{\mu,N}_{[0,T]}\Big)&\overset{\mathsf{X}_t^{\pb(E)}}{\longrightarrow} & \Big(\pb(E),F^{N}_t\Big)\\
\mu_{\mathcal{X}^N_{[0,T]}} &\longmapsto & \big(\mu_{\mathcal{X}^N_t}\big)_{0\leq t\leq T} &\longmapsto & \mu_{\mathcal{X}^N_t}
\end{array}.
\end{equation*}

This gives three notions of empirical propagation of chaos.
\begin{definition}[Empirical propagation of chaos]\label[I]{def:empiricalchaos}
Let $(f_t)_t\in C([0,T],\pb(E))$ be a flow of measures and let $f_{[0,T]}\in\pb(D([0,T],E))$ be such that for any $t\in[0,T]$,
\[f_t = \Pi^E(f_{[0,T]})(t).\]
There are three notions of empirical propagation of chaos defined below (the convergence is the weak convergence). 
\begin{enumerate}[1.]
    \item \emph{(Pointwise empirical propagation of chaos).} For all $t\in [0,T]$, the law $F^N_t$ satisfies
    \begin{equation}\label[Ieq]{eq:pointwiseempirical}F^N_t \underset{N\to+\infty}{\longrightarrow}\delta_{f_t}\in \pb(\pb(E)).\end{equation}
    \item \emph{(Functional law of large numbers).} The law $F^{\mu,N}_{[0,T]}$ satisfies
    \begin{equation}\label[Ieq]{eq:weakpathwise}F^{\mu,N}_{[0,T]}\underset{N\to+\infty}{\longrightarrow}\delta_{\Pi^E(f_{[0,T]})}\equiv \delta_{(f_t)_t}\in \pb\big(D([0,T],\pb(E))\big).\end{equation}
    \item \emph{((Strong) pathwise empirical propagation of chaos).} The law $F^N_{[0,T]}$ satisfies
    \begin{equation}\label[Ieq]{eq:strongpathwise}{F}^N_{[0,T]}\underset{N\to+\infty}{\longrightarrow}\delta_{f_{[0,T]}}\in\pb\big(\pb(D([0,T],E))\big).\end{equation}
\end{enumerate}
\end{definition}

The (strong) pathwise property \cref{eq:strongpathwise} is stronger than the functional law of large numbers \cref{eq:weakpathwise} which is stronger than the pointwise property \cref{eq:pointwiseempirical}. 

\begin{remark}
The functional law of large numbers \cref{eq:weakpathwise} is also a pathwise property which is weaker than \cref{eq:strongpathwise}. To distinguish it with the pathwise empirical propagation chaos \cref{eq:strongpathwise} we may occasionally call the latter property \emph{strong pathwise empirical propagation of chaos}. 
\end{remark}

The implication \cref{eq:strongpathwise}$\Rightarrow$\cref{eq:weakpathwise} is not straightforward because the map $\Pi^E$ is not continuous everywhere. The result holds because the limit is a Dirac mass, this is proved in \cite[Theorem 4.7]{meleard_asymptotic_1996} using a result of L\'eonard \cite[Lemma 2.8]{leonard_large_1995}. The implication \cref{eq:weakpathwise}$\Rightarrow$\cref{eq:pointwiseempirical} is more classical, this is the convergence of the finite dimensional distributions stated in \cite[Chapter 3, Theorem 7.8]{ethier_markov_1986}. 

The converse implications do not hold in general. In fact, since the map $\Pi^E$ is not injective the strong pathwise property is meaningless when only the flow of measures $(f_t)_t$ is known and $f_{[0,T]}$ is not specified. A counterexample is given in \cite[Chapter 1, Section 3(e)]{sznitman_topics_1991}. Sznitman builds two one-dimensional processes $\mathcal{X}^N_{[0,T]}$ and $\widetilde{\mathcal{X}}^N_{[0,T]}$ such that for both processes, the strong pathwise empirical propagation chaos holds but with two different limits $f_{[0,T]}$ and $\widetilde{f}_{[0,T]}$ and with the equality of the flows of time-marginals $\Pi^E(f_{[0,T]})=\Pi^E(\widetilde{f}_{[0,T]})$. The process $\widetilde{\mathcal{X}}^N_{[0,T]}$ is obtained by a re-ordering procedure from $\mathcal{X}^N_{[0,T]}$ which thus ensures the equality of the empirical measure processes: 
\[\forall t\in[0,T],\quad \mu_{\mathcal{X}^N_t} = \mu_{\widetilde{\mathcal{X}}^N_t}.\]
Finally, by Lemma \cref{lemma:caractchaos}, the pointwise empirical propagation of chaos is equivalent to the pointwise propagation of chaos in the sense of Definition \cref{def:poc}. Similarly, the strong pathwise empirical propagation of chaos is equivalent to the pathwise propagation of chaos. On the other hand the functional law of large numbers is an intermediate notion in between pointwise and pathwise propagation of chaos.

\section{Proving propagation of chaos} \label[I]{sec:proving}

Several methods are available to prove propagation of chaos. The choice of the method depends on several aspects, including the following ones.  
\begin{itemize}
    \item \textbf{How is the particle system defined?} A SDE representation allows to control directly the trajectory of each particle. If the system is an abstract Markov process defined by its generator only, one has to pass to the limit inside a ``statistical object'': for instance the Liouville equation \cref{eq:weakliouville} or a martingale problem (Definition \cref{def:martingaleproblemparticles}). 
    \item \textbf{How well is the limit process known?} As mentioned in Section \cref{sec:particlesystemsconventions}, it is sometimes possible to prove at the same time the propagation of chaos and an existence result for the limit object. Many ``historical'' proofs are based on this idea and exploit at the particle level a property of completeness (as in McKean's original proof, Section \cref{sec:mckeancoupling}), a compactness criterion based on a martingale formulation (Section \cref{sec:mckeancompactnessreview}) or an explicit series expansion for the solution of the Liouville equation in the case of a Boltzmann problem (as in Kac's original proof, Section \cref{sec:kactheorem}). Over the years, the study of the limit problem, in parallel to the question of propagation of chaos, has stimulated the development of new techniques where wellposedness results or regularity properties of the limit problem are used to control the particle system. Ultimately, a trade-off has always to be done between regularity of the $N$-particle system and regularity of the limit process. The key idea is to write the $N$-particle system and the limit process in a common framework which allows to compare them.
    \item \textbf{Which kind of propagation of chaos?} In Sections \cref{sec:kacchaos} and \cref{sec:poc}, several notions of chaos and propagation of chaos are introduced. The first distinction to keep in mind is between pathwise and pointwise properties. Then one may seek quantitative estimates. Pathwise chaos is stronger and it is often simpler to get pointwise quantitative estimates.  
\end{itemize}

Keeping these aspects in mind, the present section is organised as follows. Section \cref{sec:coupling} is devoted to an introduction of coupling methods in several cases. These methods (or most of them) exploit a SDE representation of the particle system (it therefore requires  some regularity at the microscopic level and often a wellposedness result for the limit system) and lead to the quantitative pathwise or pointwise propagation of chaos. Section \cref{sec:provingcompactness} introduces some ideas to prove the tightness (and thus the compactness) of the law of the empirical process. This leads to non-quantitative pathwise propagation of chaos results, but as it is only based on the properties of the generator of the $N$-particles process, it remains valid for a wide class of models. A pointwise study of the empirical process via the asymptotic analysis of its generator is described in Section \cref{sec:abstractmischlermouhot}. This leads to a quantitative abstract theorem with a comparable range of applications as the compactness methods. In Section~\cref{sec:largedeviations}, some ideas related to large deviations are presented. It leads to strong (non quantitative) abstract results which go beyond but include the propagation of chaos. Although these results are often too strong or too abstract to be used in practise, the ideas can be reinterpreted to prove propagation of chaos for particle systems with a very weak regularity or with a complex interaction mechanism which are difficult to handle with other methods. Finally, in the case of Boltzmann models, specific tools can be used as described in Section \cref{sec:boltzmanntools}. 

We recall that several applications of these methods will be presented in the second part of this review, in particular in Sections \cref{sec:mckeanreview} and \cref{sec:boltzmannreview}.

\subsection{Coupling methods}\label[I]{sec:coupling}

\subsubsection{Definition}

\begin{definition}[Chaos by coupling the trajectories]\label[I]{def:chaosbycouplingtrajectories}
Let be given a final time $T\in (0,\infty]$, a distance $d_E$ on $E$ and $p\in\N$. Propagation of chaos holds by coupling the trajectories when for all $N\in\N$ there exist
\begin{itemize}
    \item a system of particles $(\mathcal{X}^N_t)^{}_t$ with law $f^N_t\in\pb(E^N)$ at time $t\leq T$, 
    \item a system of \emph{independent} processes $\big(\overline{\mathcal{X}}{}^N_t\big){}^{}_t$ with law $f^{\otimes N}_t\in\pb(E^N)$ at time $t\leq T$, 
    \item a number $\varepsilon(N,T)>0$ such that $\varepsilon(N,T)\underset{N\to+\infty}{\longrightarrow}0$, 
\end{itemize}
such that (pathwise case)
\begin{equation}\label[Ieq]{eq:chaoscouplingpathwise}\frac{1}{N}\sum_{i=1}^N \E{\left[\sup_{t\leq T}d_E\big(X^i_t,\overline{X}{}^i_t\big)^p\right]}\leq \varepsilon(N,T),\end{equation}
or (pointwise case)
\begin{equation}\label[Ieq]{eq:chaoscouplingpointwise}\frac{1}{N}\sum_{i=1}^N \sup_{t\leq T} \E{\left[d_E\big(X^i_t,\overline{X}{}^i_t\big)^p\right]}\leq \varepsilon(N,T).\end{equation}
\end{definition}

By definition of the Wasserstein-$p$ distance (Definition \cref{def:wasserstein}) and Jensen inequality, the bounds \cref{eq:chaoscouplingpathwise} and \cref{eq:chaoscouplingpointwise} imply the infinite dimensional chaos (Definition \cref{def:wassersteinchaos}), respectively: 
\[W_p{\left(f^N_{[0,T]},f^{\otimes N}_{[0,T]}\right)}\underset{N\to+\infty}{\longrightarrow}0,\quad \sup_{t\leq T} W_p\big(f^N_{t},f^{\otimes N}_{t}\big)\underset{N\to+\infty}{\longrightarrow}0,\]
where we recall that the Wasserstein-$p$ distance is associated to the normalised distance on $E^N$ (Definition \cref{def:spaceswasserstein}). The definition can also be weakened by assuming that only the first $k(N)<N$ particles are coupled instead of the whole system of the $N$ particles. This would imply propagation of chaos of the form \cref{eq:Omegak}. 

\begin{remark}
Note that by exchangeability of the particles all the expectations in the sums \cref{eq:chaoscouplingpathwise} and \cref{eq:chaoscouplingpointwise} are equal and the assertions therefore imply the convergence of the one-particle distribution in Wasserstein-$p$ distance on $\pb(E)$. 
\end{remark}

The coupling between the trajectories also implies the quantitative empirical chaos stated in the following lemma in the pointwise case. This is a simple application of \cite{fournier_rate_2015}.

\begin{lemma} Let $E=\R^d$ and assume that $f_t$ has a bounded moment of order $q>p$. If there exists a (pointwise) coupling as in Definition \cref{def:chaosbycouplingtrajectories} then 
\[\sup_{t\leq T}\mathcal{W}^p_p{\big(F^{N}_t,\delta_{f_t}\big)}\leq \varepsilon(N,T)+\beta_d(N),\]
where $\beta_d(N)$ is given by:
\begin{multline*}\beta_d(N) \\= C(p,q)\left\{\begin{array}{ll}
N^{-1/2}+N^{-(q-p)/q} & \text{if}\,\,p>d/2\,\,\text{and}\,\,q\ne 2p\\ 
N^{-1/2}\log(1+N)+N^{-(q-p)/q} & \text{if}\,\,p=d/2\,\,\text{and}\,\,q\ne 2p\\
N^{-p/d}+N^{-(q-p)/q}&\text{if}\,\,p<d/2\,\,\text{and}\,\,q\ne d/(d-p),
\end{array}
\right.
\end{multline*}
for a constant $C(p,q)>0$ which depends only on $p$ and $q$. 
\end{lemma}

\begin{proof} By the triangle inequality, for all $t\leq T$,
\begin{align*}\mathcal{W}^p_p{\big(f^{\mu,N}_t,\delta_{f_t}\big)}&\leq \E W_p^p{\big(\mu_{\mathcal{X}^N_t},f_t\big)} \\
&\leq \E W_p^p{\Big(\mu_{\mathcal{X}^N_t},\mu_{\overline{\mathcal{X}}{}^N_t}\Big)}+\E W_p^p{\Big(\mu_{\overline{\mathcal{X}}{}^N_t},f_t\Big)}.
\end{align*}
The second term on the right-hand side of the last inequality is bounded by $\beta_d(N)$ using \cite[Theorem 1]{fournier_rate_2015}. Moreover the bound \cref{eq:chaoscouplingpointwise} implies that
\[\E W_p^p{\Big(\mu_{\mathcal{X}^N_t},\mu_{\overline{\mathcal{X}}{}^N_t}\Big)}\leq \varepsilon(N,T).\]
\end{proof}

Note that the convergence rate depends on the dimension. When moments of sufficiently high order are available, $\beta_d(N)$ is of the order $N^{-1/2}$ for $p>d/2$ and $N^{-p/d}$ for $p<d/2$. 

We now summarise the most common methods used to construct a coupling in the sense of Definition \cref{def:chaosbycouplingtrajectories}. 

\subsubsection{Synchronous coupling.}\label[I]{sec:synchronouscouplingproving}

When the particle system can be written as the solution of a system of SDEs as in \cref{eq:mckeanvlasov} or \cref{eq:meanfieldjumpsde}, a simple coupling choice consists in constructing $N$ nonlinear processes by taking respectively the same Brownian motions and Poisson random measures as those defining the particle system. This choice is called \emph{synchronous coupling}. 

For the McKean-Vlasov diffusion, the synchronous coupling is thus based on the~$N$ independent processes $\overline{\mathcal{X}}^N_t = (\overline{X}^1_t,\ldots,\overline{X}^N_t)$ defined as the solutions of the~$N$ SDEs: 
\begin{equation}\label[Ieq]{eq:mckeanvlasov-limitsynchronouscoupling}
\dd \overline{X}{}^i_t = b \big( \overline{X}{}^i_t,f_t \big)\dd t + \sigma \big( \overline{X}{}^i_t,f_t \big) \dd B^i_t,
\end{equation}
for $i\in\{1,\ldots,N\}$ where $(B^i_t)^{}_t$ is the same Brownian motion as in \cref{eq:mckeanvlasov} and where we recall that $f_t=\mathrm{Law}\big(\overline{X}{}^i_t\big)$. Since the Brownian motions $B^i_t$ are independent, this gives $N$ independent copies of \cref{eq:mckeanvlasov-limit}. The Theorem \cref{thm:mckean} (and Theorem \cref{thm:mckeangeneral}) in Section \cref{sec:mckeancoupling} will show that this coupling choice leads to an optimal convergence rate (in $N$) for any $T>0$ but with a constant which depends exponentially in~$T$. This comes from the fact that comparing the trajectories of \cref{eq:mckeanvlasov-limitsynchronouscoupling} and \cref{eq:mckeanvlasov} is similar to a stability analysis of the $N$ processes $X^i_t-\overline{X}^i_t$. When the coefficients are globally Lipschitz, the classical Gronwall-based methods imply the stability on any time interval but with an constant which grows exponentially with the time variable. 

The synchronous coupling is by far the most popular choice of coupling method in the literature since \cite{sznitman_topics_1991}. We point out that this was not the original choice of McKean: in the seminal work \cite{mckean_propagation_1969}, McKean uses a synchronous coupling between the $N$-particle system and the subsystem of the $N$ first particles of a system of $M>N$ particles. The proof thus does not necessitate to prove the well-posedness of the nonlinear SDE \cref{eq:mckeanvlasov-limit} as a preliminary step (since it is never used). It is actually a proof of existence which constructs a solution of \cref{eq:mckeanvlasov-limitsynchronouscoupling} by a completeness argument together with a probabilistic reasoning (based on Hewitt-Savage 0-1 law) to recover the independence. 

Similar ideas can be applied for parametric mean-field jump processes (with or without simultaneous jumps). Given the $N$ SDEs \cref{eq:meanfieldjumpsde}, the synchronous coupling is defined by: 
\begin{align}\label[Ieq]{eq:meanfieldjumpsdesynchronouscoupling}
    \overline{X}{}^i_t = X^i_0 &+ \int_0^t a\big(\overline{X}{}^i_s\big)\dd s \nonumber\\
    &+ \int_{0}^t\int_{0}^{+\infty}\int_\Theta {\left\{\psi{\big(\overline{X}{}^i_{s^-},f_s,\theta\big)}-\overline{X}{}^{i}_{s^-}\right\}}\1_{\big(0,\lambda{\big(\overline{X}{}^i_{s^-},f_s\big)}\big]}(u)\,\,\mathcal{N}^i(\dd s,\dd u,\dd\theta),
\end{align}
for $i\in\{1,\ldots,N\}$, where $f_t = \mathrm{Law}\big(\overline{X}{}^i_t\big)$ and where the Poisson random measures~$\mathcal{N}^i$ are the same as in \cref{eq:meanfieldjumpsde}. Equation \cref{eq:meanfieldjumpsdesynchronouscoupling} can be extended straightforwardly to the case of simultaneous jumps (Example \cref{example:simultaneousjumps}). It is then possible to prove similar results as in the case of the McKean-Vlasov diffusion. A complete analysis can be found in \cite{andreis_mckeanvlasov_2018} (see also Section \cref{sec:mckeangeneralinteractions}).

To end this section, let us also mention the recent coupling method introduced in \cite{holding_propagation_2016} which reverses the role of the empirical particle system and the nonlinear law. The author introduces the particle system defined conditionally on $\mathcal{X}^N_t$ by:
\[\dd \widetilde{X}^i_t = b \big( \widetilde{X}^i_t,\mu_{\mathcal{X}^N_t} \big)\dd t + \sigma \big( \widetilde{X}^i_t,\mu_{\mathcal{X}^N_t}\big) \dd \widetilde{B}^i_t.\]
Note that the processes $\widetilde{X}^i_t$ are not independent. More details on how to close the argument (using a generalised Glivenko-Cantelli theorem for SDEs) will be given in Section \cref{sec:holding}. The main advantage is that it allows more singular interactions, namely only H\"older instead of Lipschitz. 

\begin{remark} Coupling methods should still be possible when no SDE is available, since it is always possible to write an evolution equation for an observable $\varphi {\left( \mathcal{X}^N_t \right)}$ using the general It\={o}'s formula for Markov processes (see Section \cref{appendix:Dsemimartingales}). The analog of an SDE can then be recovered taking for $\varphi$ some coordinate functions.
\end{remark}

\subsubsection{Reflection coupling for the McKean-Vlasov diffusion.}\label[I]{sec:reflectioncouplingproving}

Except in specific cases (see Section \cref{sec:gradientsystems}), it is usually not possible or difficult to get uniform in time estimates using a synchronous coupling. One reason is that the strategy can be seen as a stability analysis for the nonlinear system \cref{eq:mckeanvlasov-limit} which classically leads to Gronwall type estimates with a constant which depends exponentially in $T$. The recently developed \emph{reflection coupling} \cite{eberle_reflection_2016,eberle_quantitative_2019} tries to make a better use of the diffusion part to get (hopefully) uniform in times estimates under mild assumptions. This type of coupling was originally introduced in \cite{lindvall_coupling_1986} to study the ergodic properties of diffusion processes in total variation norm. In this context, the coupling time of two diffusion processes with same transition probabilities but different initial distributions is defined as the first time where the two processes hit the same state. The reflection coupling method provides an explicit construction of the processes which (hopefully) ensures that their coupling time is almost surely finite. By letting the two processes being equal after the coupling time, this fact shows that the two processes are asymptotically equal in law. The idea of reflection coupling is easily understood for two 1D Brownian motions starting from the initial conditions $a>0$ and $-a$ : it is sufficient to take one Brownian motion equal to the opposite of the other one. More generally, the solution $(X_t,Y_t)$ of the following system of SDEs in $\R^d$ is called a coupling by reflection: 
\begin{align*}
    \dd X_t &= b(X_t)\dd t + \sigma \dd B_t\\
    \dd Y_t &= b(Y_t)\dd t + \sigma(I_d-2e_t e_t^\mathrm{T})\dd B_t,
\end{align*}
where $b:\R^d\to\R^d$ is the (locally Lipschitz) drift function, $\sigma>0$ is a constant and 
\[e_t = (X_t-Y_t)/|X_t-Y_t|.\]
Moreover, after the coupling time $T:=\inf\{t\geq0, X_t=Y_t\}$, for $t\geq T$ the processes are set to $X_t=Y_t$. In a series of recent works, \cite{eberle_reflection_2016, eberle_quantitative_2019}, Eberle et al. have proved a quantitative contraction result in Wasserstein distance for the law of the two processes. Let $b$ satisfy for all $x,y\in\R^d$,
\[\langle x-y, b(x)-b(y)\rangle \leq-\frac{\sigma^2}{2}\kappa(|x-y|)|x-y|^2,\]
where $\kappa:\R_+\to\R$ satisfies $\liminf_{r\to+\infty} \kappa(r)>0$. In order to measure the discrepancy between the two processes $r_t:=|X_t-Y_t|$, for any fixed smooth function $f$, It\={o}'s formula gives:
\[\dd f(r_t) = r_t^{-1}\langle X_t-Y_t,b(X_t)-b(Y_t)\rangle f'(r_t)\dd t + 2\sigma^2 f''(r_t)\dd t + \dd M_t,\]
where $M_t$ is a martingale. Compared to what the synchronous coupling would give, thanks to It\=o's correction, the coupling by reflection adds a new term in the drift. Using the assumptions on $b$, the drift term is now bounded by
\[2\sigma^2 \left(f''(r_t)-\frac{1}{4}r_t\kappa(r_t)f'(r_t)\right).\]
If $f$ is such that there exists $c>0$ such that for all $r\geq0$
\begin{equation}\label[Ieq]{eq:reflectioneqfintro}f''(r)-\frac{1}{4}r\kappa(r)f'(r)\leq -\frac{c}{2\sigma^2}f(r),\end{equation}
then it gives:
\begin{equation}\label[Ieq]{eq:reflectioncontractexp}\E[f(r_t)] \leq \e^{-ct} \E[f(r_0)].\end{equation}
The idea of \cite{eberle_reflection_2016} is to introduce a positive concave function $f$ so that $d_f(x,y):=f(|x-y|)$ defines a distance on $\R^d$ and such that the bound \cref{eq:reflectioneqfintro} holds. From~\cref{eq:reflectioncontractexp}, this finally gives the following exponential ``contraction bound'' in Wasserstein distance: 
\begin{equation}\label[Ieq]{eq:exponentialcontractionreflectioncoupling}W_{1,d_f}(\mu_t,\nu_t)\leq \e^{-ct} W_{1,d_f}(\mu_0,\nu_0),\end{equation}
where $\mu_t,\nu_t\in\pb(E)$ are the laws of $X_t,Y_t$ at time $t\geq0$. This strategy is successfully applied in \cite{eberle_reflection_2016,eberle_quantitative_2019} to get quantitative contraction and convergence rates for linear and nonlinear gradient McKean-Vlasov systems, even in non convex settings.

Coming back to particle systems, in \cite{durmus_elementary_2020,liu_long-time_2021}, the authors have shown that this idea can be applied to McKean-Vlasov systems \cref{eq:mckeanvlasov} (which can be seen as a classical diffusion equation in the high-dimensional space $\R^{dN}$). They use a componentwise reflection coupling between a particle system and a system of $N$ independent nonlinear McKean-Vlasov systems. The analog of the exponential contraction rate~\cref{eq:exponentialcontractionreflectioncoupling} thus provides a proof of the uniform in time propagation of chaos for gradient systems with milder assumptions than the ones in Section \cref{sec:gradientsystems} obtained with a synchronous coupling. This will be reviewed in Section \cref{sec:reflectioncoupling}. 

\begin{remark}[Extension to more general diffusions] This idea is more natural but not restricted to the case where the diffusion matrix is constant. It can be extended to more general diffusion matrices by ``twisting'' the metric in $\R^d$ to recover a constant diffusion matrix in a modified metric. See \cite{eberle_reflection_2016} for additional details as well as \cite{malrieu_logarithmic_2001} for a similar reasoning in a different context.  
\end{remark}

\subsubsection{Optimal jumps}

For mean-field jump processes and Boltzmann models, the particles interact only at discrete (random) times. They update their state according to a sampling mechanism with respect to a known measure which depends on the empirical measure of the system (see \cref{eq:newstatepdmp} and \cref{eq:newstateboltzmann}). The strategy adopted in \cite{diez_propagation_2020} for mean-field jump processes (see Section \cref{sec:meanfieldjump}) consists in constructing a trajectorial representation of the particle and the nonlinear systems in which the jumps are coupled optimally. Taking the same sequence of jump times $(T^i_n)^{}_n$ for the particle $X^i_t$ and its coupled nonlinear version $\overline{X}{}^i_t$, a post-jump state is sampled for the nonlinear process first: 
\[\overline{X}{}^i_{T^i_n} \sim P_{f_{T^i_n}}{\left(\overline{X}{}^i_{T^{i-}_n},\dd y\right)},\]
and then the post-jump state for the particle is defined as the image:
\begin{equation}\label[Ieq]{eq:optimalcouplingjumpT}X^i_{T^i_n} = \mathsf{T}{\Big(\overline{X}{}^i_{T^i_n}\Big)},\end{equation}
where $\mathsf{T}$ is an optimal transfer map for the $W_1$ distance between the jump measures:
\[\mathsf{T}_\#P_{f_{T^i_n}}{\left(\overline{X}{}^i_{T^{i-}_n},\dd y\right)}= P_{\mu_{\mathcal{X}^N_{T^{i-}_n}}}{\left({X}^i_{T^{i-}_n},\dd y\right)},\]
and
\[\E{\left[\big|\overline{X}{}^i_{T^i_n}-X^i_{T^i_n}\big| \Big| \mathcal{F}_{T^{i-}_n}\right]}=W_1{\Big(P_{f_{T^i_n}}{\left(\overline{X}{}^i_{T^{i-}_n},\dd y\right)},P_{\mu_{\mathcal{X}^N_{T^{i-}_n}}}{\left({X}^i_{T^{i-}_n},\dd y\right)}\Big)},\]
where $\mathcal{F}_t$ denotes the filtration generated by all the processes $(X^i_t)^{}_t$ and $(\overline{X}^i_t)^{}_t$ up to time $t$ and $\mathcal{F}_{T^{i-}_n} = \sigma(\bigcup_{t<T^i_n}\mathcal{F}_t)$. The existence of the optimal transfer map $\mathsf{T}$ is a classical result in optimal transport theory; it holds under mild assumptions, see \cite{champion_monge_2011,feldman_monges_2001} and the references therein. Under Lipschitz assumptions on the jump measure, it can be deduced that 
\begin{equation}\label[Ieq]{eq:optimaljumplipschitz}
    \E{\left[\big|\overline{X}{}^i_{T^i_n}-X^i_{T^i_n}\big| \Big| \mathcal{F}_{T^{i-}_n}\right]}\leq C{\left(\Big|\overline{X}{}^i_{T^{i-}_n}-{X}^i_{T^{i-}_n}\Big| +W_1{\Big(f_{T^i_n},\mu_{\mathcal{X}^N_{T^{i-}_n}}\Big)}\right)}.
\end{equation}
This crude estimate may be refined when the jump measure has a known expression. Note that in the case of a parametric model (Example \cref{example:parametricjump}), the synchronous coupling of the Poisson random measures \cref{eq:meanfieldjumpsdesynchronouscoupling} gives an alternative coupling and an explicit transfer map:
\begin{equation}\label[Ieq]{eq:optimalcouplingjumpparam}\overline{X}{}^i_{T^i_n} = \psi\big(\overline{X}{}^i_{T^{i,-}_n},f_{T^i_n},\theta\big),\quad X^i_{T^i_n} = \psi\Big(X^{i}_{T^{i,-}_n},\mu_{X^i_{T^{i,-}_n}},\theta\Big),\end{equation}
where $\theta\sim \nu(\dd\theta)$ is the same random variable for the two post-jump states. This coupling is not necessarily optimal but under Lipschitz assumptions on $\psi$ it still implies 
\begin{align*}
\E{\left[\big|\overline{X}{}^i_{T^i_n}-X^i_{T^i_n}\big| \Big| \mathcal{F}_{T^{i-}_n}\right]}&=\int_{\Theta}{\Big|\psi{\left(\overline{X}{}^i_{T^{i,-}_n},f_{T^i_n},\theta\right)}-\psi\Big(X^{i}_{T^{i,-}_n},\mu_{X^i_{T^{i,-}_n}},\theta\Big)\Big|}\nu(\dd\theta)\\
&\leq C{\left(\Big|\overline{X}{}^i_{T^{i-}_n}-{X}^i_{T^{i-}_n}\Big| +W_1{\Big(f_{T^i_n},\mu_{\mathcal{X}^N_{T^{i-}_n}}\Big)}\right)}.
\end{align*}

Both couplings \cref{eq:optimalcouplingjumpT} and \cref{eq:optimalcouplingjumpparam} ensure that the $N$ nonlinear processes $(\overline{X}{}^i_t)^{}_t$ remain independent, which is crucial. At each jumping time $T^i_n$, the error between $X^i_t$ and $\overline{X}{}^i_t$ due to the jump is controlled by \cref{eq:optimaljumplipschitz}. A discrete stability analysis then ensures the propagation of the error exponentially in time. Between the jumps, the trajectories are either deterministic or can be controlled by a standard synchronous coupling. This proves the propagation of chaos on any time interval but with a (very) bad behaviour with respect to time. 

For Boltzmann models, the situation is more difficult because two particles ``jump'' at the same time. As discussed in Section \cref{sec:boltzmannparametricmodels}, it is also not completely straightforward to build a SDE representation of the nonlinear process. Moreover, contrary to the mean-field jump processes where the jump measures are typically assumed to have a smooth density, for Boltzmann models, the jumps are obtained by sampling directly from the empirical measure of the system. Since this measure is singular but the solution of the Boltzmann equation is not, an optimal transfer map may not exist. The strategy adopted in \cite{murata_propagation_1977} and then in \cite{cortez_quantitative_2016,cortez_quantitative_2018,fournier_rate_2016} is based on the analysis of the optimal transfer plan (which exists) between the empirical measure of the particle system and the law of the nonlinear system at each jump time. This strategy also needs a kind of synchronous coupling for the jump times. As in the previous cases, propagation of chaos then results from a Gronwall type estimate with a bad behaviour in time, and in this case only for marginals (or block) of size $k(N)=o(N)$. This will be discussed more thoroughly in Section \cref{sec:couplingBoltzmann}.

\subsubsection{Analysis in Wasserstein spaces: optimal coupling}\label[I]{sec:optimalcouplingintro}

The Liouville equation associated to a McKean-Vlasov process with interaction parameters
\[b(x,\mu) \equiv \int_{\R^d}b(x,y)\mu(\dd y), \quad \sigma(x,\mu)\equiv \sigma I_d,\quad \sigma>0,\]
can be written:
\begin{equation}\label[Ieq]{eq:liouvillemckean}
\partial_t f^N_t = -\nabla\cdot(\mathbf{b}^N f^N_t) + \frac{\sigma^2}{2}\Delta f^N_t,
\end{equation}
where
\[\mathbf{b}^N : \mathbf{x}^N\in \R^{dN}\mapsto \left(\frac{1}{N}\sum_{i=1}^N b(x^1,x^i),\ldots,\frac{1}{N}\sum_{i=1}^N b(x^N,x^i)\right)\in\R^{dN}.\]
This equation can be rewritten as a continuity equation:
\begin{equation}\label[Ieq]{eq:mckeanvlasovcontinuity}\partial_t f^N_t = -\nabla\cdot \left(\left(\mathbf{b}^N-\frac{\sigma^2}{2}\nabla \log f^N_t\right)f^N_t\right).\end{equation}
Similarly, the associated nonlinear McKean-Vlasov equation is a continuity equation in $\R^d$ with velocity vector 
\[b\star f_s - \frac{\sigma^2}{2}\nabla \log f_s.\]
Continuity equations are strongly linked to the theory of gradient flows \cite{ambrosio_gradient_2008,pajot_lecture_2014}, which in turn provides new insights on the study of McKean-Vlasov processes. In addition to the present section, see also Section \cref{sec:gradientflows}. 

The two recent works \cite{salem_gradient_2020,del_moral_uniform_2019} are based on a classical result in gradient flow theory (see for instance \cite[Theorem 23.9]{villani_optimal_2009}) which gives an explicit dissipation rate between the solutions of two continuity equations. In the present case, it gives an explicit control of the time derivative of ${W}_2(f^N_t,f^{\otimes N}_t)$ in terms of the well-known maximizing Kantorovich potential $\psi^N_t$ which links the two laws by: 
\[(\nabla\psi^N_t)_\# f_t^{\otimes N} = f^N_t, \quad W_2^2(f^N_t,f_t^{\otimes N}) = \int_{\R^{dN}}|\nabla\psi^N_t(\mathbf{x}^N)-\mathbf{x}^N|^2 f^{\otimes N}_t(\dd\mathbf{x}^N).\]
The existence of $\psi^N_t$ is ensured by Brenier's theorem \cite[Theorem 9.4]{villani_optimal_2009}. This approach of propagation of chaos follows the work of \cite{bolley_convergence_2012,bolley_uniform_2013} where the authors have derived explicit contraction rates in Wasserstein-2 distance for linear Fokker-Planck equations and for the nonlinear granular media equation (which is the nonlinear mean-field limit associated to the gradient system \cref{eq:gradientmckean}). Starting from the same result \cite[Theorem 23.9]{villani_optimal_2009}, the authors of \cite{bolley_convergence_2012,bolley_uniform_2013} introduced a new transportation inequality, the so-called WJ inequality which is then exploited at the particle level in \cite{salem_gradient_2020,del_moral_uniform_2019}. This analysis provides uniform in time propagation of chaos and convergence to equilibrium results for gradient systems in non globally convex settings. The work of \cite{salem_gradient_2020} also provides a new unifying analytical vision of previous coupling approaches. These techniques will be detailed in Section \cref{sec:optimalcoupling}.

\subsection{Compactness methods}\label[I]{sec:provingcompactness}

In this section the main ideas to prove propagation of chaos via compactness arguments are presented. The main advantage of this approach is its wide range of applicability as it can be adapted to jump, diffusion, Boltzmann or even mixed models. The main drawback is that it does not provide any convergence rate.

Compactness methods are based on the empirical representation of the process described in Section \cref{sec:poc}. It therefore reduces the problem of the convergence of a sequence of probability measures on a space which does not depend on $N$ but which in turn is much more delicate to handle than $E$ as it is itself a probability space. The first approach described below is the stochastic analysis approach which is by now classical; the second approach is a more recent analytical approach based on the theory of gradient-flows. 

\subsubsection{Martingale methods.}

Starting from the martingale characterisation of the particle system (Definition \cref{def:martingaleproblemparticles}), it is possible to prove the functional law of large numbers \cref{eq:weakpathwise} and the strong pathwise empirical propagation of chaos \cref{eq:strongpathwise} using the traditional sequence of arguments in stochastic analysis (see for instance \cite{joffe_weak_1986}). 
\begin{enumerate}[(1)]
    \item First prove the tightness of the sequence $\big(F^{\mu,N}_{[0,T]}\big)_N$ (functional law of large numbers) or $\big(F^N_{[0,T]}\big)_N$ (strong pathwise). This will come from usual tightness criteria (see Section \cref{appendix:skorokhod} and Section \cref{appendix:tightness}) but requires some care regarding the spaces (namely, the path space with values in a set of probability measures). The classical and more advanced tools which are used are reminded in Appendix \cref{appendix:probabilityreminders}. By Prokhorov theorem, it is then possible to extract a converging subsequence towards a limit, respectively $\pi\in\pb(D([0,T],\pb(E)))$ (functional law of large numbers) or $\pi\in\pb(\pb(D([0,T],E)))$ (strong pathwise).  
    \item Then identify the $\pi$-distributed limit points as solutions respectively of the weak limit PDE or the limit martingale problem. Again, this may require some care, in particular for c\`adl\`ag processes due to the topology of the Skorokhod space.  
    \item Finally prove the uniqueness of the solution of the previous problem. Usually well-posedness (that is existence and uniqueness) can be proved beforehand although existence is not required (it is automatically provided by the tightness). In conclusion, $\pi$ is a Dirac delta at the desired limit point. 
\end{enumerate}
Strong pathwise empirical chaos \cref{eq:strongpathwise} is not significantly harder to prove than the weaker functional law of large numbers \cref{eq:weakpathwise} but it requires the uniqueness property of the limit martingale problem which is a stronger assumption than the corresponding one for the weak limit PDE. The strong pathwise case is detailed in M\'el\'eard's course \cite[Section 4]{meleard_asymptotic_1996}. 

This approach is historically linked to the study of the spatially homogeneous Boltzmann equation of rarefied gas dynamics \cref{eq:spatiallyhomogeneousboltzmann}: Tanaka \cite{tanaka_probabilistic_1983} proved weak pathiwse empirical chaos for hard-spheres and inverse power Maxwellian molecules. Strong pathwise empirical chaos is proved in \cite{sznitman_equations_1984} for a class of parametric Boltzmann models which includes the stochastic hard-sphere model. See also \cite{wagner_functional_1996} for a functional law of large numbers applied to a large class of parametric Boltzmann models. The martingale method is exploited to treat the case of the McKean-Vlasov diffusion in \cite{sznitman_nonlinear_1984} (with boundary conditions) and in \cite{oelschlager_martingale_1984, gartner_mckean-vlasov_1988} (functional law of large numbers with general interaction functions), see also \cite{leonard_loi_1986}. Strong pathwise empirical chaos is proved for a mixed jump-diffusion mean-field model in \cite{graham_stochastic_1997,meleard_asymptotic_1996}. See also \cite[Theorem 4.1]{fournier_markov_2001} for a strong pathwise result on a cutoff approximation of a non-cutoff Boltzmann model.

The corresponding results will be detailed in Section \cref{sec:martingalecompactness} in the mean-field case and in Section \cref{sec:martingaleboltzmannreview} for Boltzmann models. 

\subsubsection{Gradient flows.}\label[I]{sec:gradientflowsintro}

This second approach gives a pointwise version of the empirical propagation of chaos and is restricted to the McKean-Vlasov gradient system~\cref{eq:gradientmckean}. It is entirely analytical and exploits recent results of the theory of gradient flows \cite{ambrosio_gradient_2008,pajot_lecture_2014,villani_optimal_2009}. 
We briefly recall below one definition of gradient-flows in a metric space.
\begin{definition}[Gradient flows in $(\mathscr{E},d_\mathscr{E})$] Let $(\mathscr{E},d_\mathscr{E})$ be a geodesic metric space.
\begin{enumerate}
    \item \textbf{(Absolutely continuous curves and metric derivative).} A $\mathscr{E}$-valued continuous curve $\mu:(a,b)\subset\R\to\mathscr{E}$ is said to be absolutely continuous whenever there exists $m\in L^1_\mathrm{loc}(a,b)$ such that 
    \[\forall a<s\leq t< b,\quad d_\mathscr{E}(\mu_s,\mu_t)\leq\int_s^t m(r)\dd r.\]
    In this case, the limit 
    \[|\mu'|(t) = \lim_{s\to t} \frac{d_\mathscr{E}(\mu_s,\mu_t)}{|t-s|},\]
    exists for almost every $t\in(a,b)$ and is called the metric derivative of $\mu$ at the point $t$. 
    \item \textbf{(Gradient flow).} Let $T\in(0,+\infty]$ and let $\mu\in C([0,T),\mathscr{E})$ be an absolutely continuous curve. Let us consider $\lambda\in\R$ and a $\lambda$-convex map $\mathcal{F}:\mathscr{E}\to\R\cup\{+\infty\}$. Then $\mu$ is called a $\lambda$-gradient flow associated to the energy $\mathcal{F}$ whenever $\mathcal{F}(\mu_t)<+\infty$ for all $t\in[0,T)$ and $\mu$ satisfies the following \textit{Evolution Variational Inequality} (EVI) :
    \begin{equation}\label[Ieq]{eq:evi}
    \forall \mathrm{a.e.}\, t\in[0,T),\,\,\forall \nu\in\mathscr{E},\quad \frac{1}{2}\frac{\dd}{\dd t} d_\mathscr{E}^2(\mu_t,\nu) + \frac{\lambda}{2}d_\mathscr{E}^2(\mu_t,\nu) \leq \mathcal{F}(\nu)-\mathcal{F}(\mu_t).
    \end{equation}
    \end{enumerate}
\end{definition}

In the following, gradient flows will be considered in the two cases $(\mathscr{E},d_\mathscr{E})=(\pb_2(\R^d),W_2)$ and $(\mathscr{E},d_\mathscr{E})=(\pb_2(\pb_2(\R^d)),\mathcal{W}_2)$. The fundamental result to keep in mind is that evolutionary PDEs, as defined below, have a unique distributional solution which is a gradient-flow. 

\begin{definition}[Evolutionary PDEs]\label[I]{def:evolutionarypde} An evolutionary PDE is a PDE of the form 
\[\partial_t \rho = \nabla\cdot\left(\rho\nabla\frac{\delta\mathcal{F}(\rho)}{\delta \rho}\right),\]
where $\mathcal{F}:\pb_2(\R^d)\to\R\cup\{+\infty\}$ and where the first variation of $\mathcal{F}$ is defined as the unique (up to an additive constant) measurable function $\frac{\delta\mathcal{F}(\rho)}{\delta \rho} :\R^d\to\R$ such that the equality
\[\frac{\dd}{\dd h}\mathcal{F}(\rho+h\chi)\Big|_{h=0} = \int \frac{\delta\mathcal{F}(\rho)}{\delta \rho} \dd \chi,\]
holds for every measure $\chi$ such that $\rho+h\chi\in\pb(\R^d)$ for small enough $h$. 
\end{definition}

In a recent article \cite{carrillo_-convexity_2020}, the authors prove the pointwise empirical propagation of chaos for gradient systems \cref{eq:gradientmckean} using a gradient flow characterisation of $(f^N_t)^{}_t$ and $(F^N_t)^{}_t$ seen as time continuous curves with values respectively in $\pb_2(\R^{dN})$ and $\pb_2(\pb_2(\R^d))$. The argument is based on a compactness criterion in the space $C([0,T],\pb_2(\pb_2(\R^d)))$ which follows from Ascoli's theorem. Under the initial chaos hypothesis, the limit of $(F^N_t)^{}_t$ as $N\to+\infty$ is also identified as a gradient-flow and is shown to be the curve $(\delta_{f_t})_t\in C([0,T],\pb_2(\pb_2(\R^d)))$. Gradient flows are identified by the EVI \cref{eq:evi} which plays a comparable role as the martingale characterisation in a stochastic context. 

The results of \cite{carrillo_-convexity_2020} will be summarised in Section \cref{sec:gradientflows}.

\subsection{A pointwise study of the empirical process} \label[I]{sec:abstractmischlermouhot}

This section is devoted to an analytical pointwise study of the empirical process within the framework developed in \cite{mischler_kacs_2013,mischler_new_2015} following an idea of \cite{grunbaum_propagation_1971}. We also refer to \cite{mischler_kacs_2012} for a review of these results. 

The goal is to obtain a quantitative control of the evolution of the law $F^N_t$ of the empirical process seen as the law of a $\pb(E)$-valued process. This control is obtained via a careful asymptotic analysis of the infinitesimal generator of the empirical process which is shown to converge (in a sense to define) towards the generator of the flow of the limit PDE starting from a random initial condition and also seen as a measure-valued process. This method is intrinsically very abstract and can be applied to a wide range of models (in theory, at least to all the models studied in the present review). The main idea traces back to Gr\"unbaum and his study of the spatially homogeneous Boltzmann hard-sphere model \cref{eq:spatiallyhomogeneousboltzmann} \cite{grunbaum_propagation_1971}. However, the seminal article of Gr\"unbaum was incomplete and based on an unproven assumption (which happened to be false in some cases). The study of measure-valued Markov processes is in general very delicate. This is mainly due to the fact that $\pb(E)$ is only a metric space and not a vector space which causes several technical problems, the most important one being the precise definition of the notion of infinitesimal generator. A probabilistic point of view can be found in \cite{dawson_measure-valued_1993}. In the framework introduced by \cite{mischler_kacs_2013,mischler_new_2015}, a new notion of differential calculus in $\pb(E)$ is defined in order to give a rigorous definition of the limit generator. The question of propagation of chaos is then stated in a very abstract framework which leads to an abstract theorem (Theorem \cref{thm:abstractMischler}) which can be applied to various models after a careful check of a set of five assumptions (Assumption \cref{assum:mischlermouhot}). Most of these assumptions are related to the regularity of the nonlinear solution operator semigroup of the limit PDE. A notable application of this method is the answer to many of the questions raised by Kac in his seminal article \cite{kac_foundations_1956} related to the spatially homogeneous Boltzmann equation \cref{eq:spatiallyhomogeneousboltzmann}. This will be reviewed in Section \cref{sec:kacprogram}. 

The rest of this section is organised as follows. The generator and transition semigroup of the empirical process are defined next. An introductory toy example using this formalism is presented in Subsection \cref{sec:toyexampleempiricalprocess}. This leads to propagation of chaos for only a very small class of linear models. The next two Subsections \cref{sec:limitsemigroup} and \cref{sec:limitgenerator} present the core of the abstract framework and the main difficulties of the approach, in particular the notion of differential calculus needed to define the limit generator. In the last Subsection \cref{sec:abstracttheorem} the five assumptions and the abstract theorem of \cite{mischler_new_2015} are stated.

\subsubsection{The empirical generator.} 
In Section \cref{sec:poc}, the empirical process has been defined as the image of the $N$-particle process by the empirical measure map $\boldsymbol{\mu}_N$. This process is a $\widehat{\pb}_N(E)$-valued Markov process and it is possible to define its transition semi-group and generator by pushing-forward those of the $N$-particle process. More precisely, the empirical transition semi-group is given by: 
\begin{equation}\label[Ieq]{eq:TNhat}\widehat{T}_{N,t}\Phi(\mu_{\mathbf{x}^N}) = T_{N,t}[\Phi\circ \boldsymbol{\mu}_N](\mathbf{x}^N),\end{equation}
and the empirical generator by:
\begin{equation}\label[Ieq]{eq:LNhat}\widehat{\mathcal{L}}_N\Phi(\mu_{\mathbf{x}^N}) = \mathcal{L}_N[\Phi\circ\boldsymbol{\mu}_N](\mathbf{x}^N).\end{equation}
This is a consequence of the identity: 
\begin{equation*}\widehat{T}_{N,t} \Phi {\left( \mu_{\mathcal{X}^N_s} \right)}=\E {\left[ \Phi \big( \mu_{\mathcal{X}^N_{t + s}} \big) | \mathcal{X}^N_s \right]} = \E {\left[ \left( \Phi \circ \boldsymbol{\mu}_N \right) \left( \mathcal{X}^N_{t + s} \right) | \mathcal{X}^N_s \right]} = T_{N,t} {\left[ \Phi \circ \boldsymbol{\mu}_N \right]} {\left( \mathcal{X}^N_s \right)}.\end{equation*}
Note that the empirical semi-group and generator are well defined as operators $C_b(\pb(E))\to C_b(\widehat{\pb}_N(E))$ and the initial law $F^N_0=(\boldsymbol{\mu}_N)_\# f^N_0$ is a probability measure on $\widehat{\pb}_N(E)$. Nonetheless, in order to take the limit $N\to+\infty$, it is more convenient to look at the empirical process as a $\pb(E)$-valued process since $\widehat{\pb}_N(E)\subset \pb(E)$.

\begin{example}[The empirical generator for a mean-field jump process] For mean-field generators, this simply reads
\[ \widehat{\mathcal{L}}_N \Phi {\left( \mu_{\mathbf{x}^N} \right)} = \sum_{i=1}^N L_{\mu_{\mathbf{x}^N}} \diamond_i [ \Phi \circ \mu_N ] \big( \mathbf{x}^N \big). \]
In the special case of a mean-field jump process, one has a more explicit formula:
\begin{align*}
&\widehat{\mathcal{L}}_N \Phi {\left( \mu_{\mathbf{x}^N} \right)} = \sum_{i=1}^N \int_E \lambda {\left( x^i , \mu_{\mathbf{x}^N} \right)} {\left[ \Phi{\left(\mu_{\mathbf{x}^N}-\frac{1}{N}\delta_{x^i} +\frac{1}{N}\delta_y \right)}- \Phi {\left( \mu_{\mathbf{x}^N} \right)} \right]} P_{\mu_{\mathbf{x}^N}} \left( x^i,\dd y \right) \\
&= N {\left\langle \mu_{\mathbf{x}^N} , \lambda ( \cdot , \mu_{\mathbf{x}^N} ) \int_E {\left[ \Phi{\left(\mu_{\mathbf{x}^N} -\frac{1}{N}\delta_{\cdot} +\frac{1}{N}\delta_y \right)}-\Phi {\left( \mu_{\mathbf{x}^N} \right)} \right]} P_{\mu_{\mathbf{x}^N}} {\left( \cdot , \dd y \right)} \right\rangle}.
\end{align*}
\end{example}

\subsubsection{A toy-example using the measure-valued formalism}\label[I]{sec:toyexampleempiricalprocess}

The model presented in this section is a simple introduction to the measure-valued formalism. We follow an idea which was originally suggested to us by P.-E. Jabin. As we shall see, simple computations can lead to propagation of chaos in the limited case of linear models. Let us consider $\mathcal{X}^N_t = (X^1_t,\ldots,X^N_t)_t$ be a PDMP (see Section \cref{sec:meanfieldjump}) defined by
\begin{itemize}
    \item a deterministic flow:
    \[\dd X^i_t = a(X^i_t)\dd t,\]
    where $a$ is $C^1$ and globally Lipschitz,
    \item a jump transition kernel $P : E\times \pb(E)\to \pb(E), (x,\mu)\mapsto P_\mu(x,\dd y)$,
    \item a constant jump rate $\lambda\equiv 1$.
\end{itemize}

We take $E=\R^d$ for simplicity. From \cref{eq:LNhat}, the associated empirical process is a measure-valued Markov process with generator: 
\begin{equation*}
    \widehat{\mathcal{L}}_N\Phi(\mu) = (a\cdot\nabla\Phi)(\mu) + N\iint_{E\times E} {\left\{\Phi{\left(\mu-\frac{1}{N}\delta_x+\frac{1}{N}\delta_y \right)}-\Phi(\mu)\right\}}P_\mu(x,\dd y)\mu(\dd x),
\end{equation*}
where $\Phi\in C_b(\pb(E)) \subset C_b( \widehat{\pb}_N (E) )$ is a test function on $\pb(E)$ and $(a\cdot\nabla\Phi)$ is well defined when $\Phi$ is a polynomial. We recall the notation $F_0^N(\dd\mu)\in\pb(\pb(E))$ for the initial law (supported on the set of empirical measures) and  $F^N_t(\dd\mu)\in\pb(\pb(E))$ for the law at time $t$ of the measure-valued Markov process with generator $\widehat{\mathcal{L}}^N$. For all test functions $\Phi\in C_b(\pb(E))$, one can write the evolution equation for the observables of the empirical process:
\begin{multline}\label[Ieq]{eq:weakpdelaw}
\frac{\dd}{\dd t}\int_{\pb(E)}\Phi(\mu)F^N_t(\dd \mu)+\int_{\pb(E)}\Phi(\mu)\nabla\cdot(aF^N_t)(\dd \mu) \\= N\int_{E\times E\times\pb(E)}{\left\{\Phi{\left(\mu-\frac{1}{N}\delta_x+\frac{1}{N}\delta_y \right)}-\Phi(\mu)\right\}}P_\mu(x,\dd y)\mu(\dd x)F^N_t(\dd \mu).
\end{multline}

The right-hand side is the jump operator and on the left-hand side there is a transport operator, again well defined for $\Phi$ polynomial. Let us also recall the associated nonlinear jump operator acting on $C_b(E)$ : 
\begin{equation}\label[Ieq]{eq:nonlinearjumpop}L_\mu\varphi(x) := \int_E \{\varphi(y)-\varphi(x)\}P_\mu(x,\dd y).\end{equation}
Its associated carr\'e du champ operator is denoted by $\Gamma_\mu$. Our goal is to try to prove pointwise empirical propagation of chaos: namely that for any $t>0$,
\[F^N_t \underset{N\to+\infty}{\longrightarrow}\delta_{f_t},\]
where $f_t\in\pb(E)$ is the solution of the following weak PDE:
\[\forall \varphi\in C_b(E),\quad \frac{\dd}{\dd t} \langle f_t,\varphi\rangle +\langle a\cdot\nabla\varphi,f_t\rangle = \langle f_t, L_{f_t}\varphi\rangle.\]
Sufficient conditions for well-posedness are $C^1$ regularity for $a$ and a global Lispchitz bound, together with Lipschitz-continuity for $P_{\mu}(x , \dd y )$ in $x$ and $\mu$ (for the Wassertein topology). To prove propagation of chaos, we will try a ``direct analytical approach'' and compare directly $F^N_t$ to its limit with a weak distance. Note that it is possible to do that because we work in $\pb(\pb(E))$ which is a space which does not depend on $N$ (it is one of the main clear advantage of the approach). 

Let us consider the distance $\mathcal{W}_{D_2}$, which is the Wasserstein-1 distance on the space $\pb(\pb(E))$ associated to the distance $D_2$ on $\pb(E)$ (see Definition \cref{def:spaceswasserstein}). Since the limit is a Dirac mass, it holds that 
\[\mathcal{W}^2_{D_2}\big(F^N_t,\delta_{f_t}\big) = \int_{\pb(E)} D^2_2(\mu,f_t)F^N_t(\dd\mu) = \sum_{n=1}^{+\infty} \frac{1}{2^n} \int_{\pb(E)}\langle \mu-f_t, \varphi_n\rangle^2F^N_t(\dd\mu),\]
and it is then enough to bound the quantity 
\begin{equation}\label[Ieq]{eq:GNsupphi}g_N(t) := \sup_{\|\varphi\|_\mathrm{Lip}\leq 1} \int_{\pb(E)} \langle \mu-f_t, \varphi\rangle^2 F^N_t(\dd\mu).\end{equation}
Let us fix a test function $\varphi\in \mathrm{Lip}_1(E)$. In order to control the deterministic flow, we define the (time-dependent) ``modified test function'' $\tilde{\varphi}\equiv\tilde{\varphi}(s,x)$ as the solution of the backward transport equation:  
\begin{equation}\label[Ieq]{eq:modifiedphi}
    \left\{
    \begin{array}{rcl}
    \partial_s\tilde{\varphi}+a\cdot\nabla_x\tilde{\varphi}&=&0\\
    \tilde{\varphi}(s=t,x)&=&\varphi(x)
    \end{array}
    \right.
\end{equation}
Since $a$ is a globally Lipschitz vector-field, the function $\tilde{\varphi}$ is Lipschitz for all $s\leq t$ with Lipschitz semi-norm:
\begin{equation}\label[Ieq]{eq:lipschitztildephi}\|\tilde{\varphi}\|_\mathrm{Lip} \leq \e^{(t-s)\|a\|_\mathrm{Lip}}.\end{equation}
We define 
\begin{equation}\label[Ieq]{eq:gphi}\tilde{g}_\varphi(s) = \int_{\pb(E)}\langle \mu-f_s,\tilde{\varphi}\rangle^2 F^N_s(\dd \mu),\end{equation}
where we do not specify the dependency in $N$ for notational simplicity. In order to apply a Gronwall-like argument, we fix $t>0$ and for $s<t$, thanks to \cref{eq:weakpdelaw}, we compute the time derivative: 
\begin{align*}
    &\tilde{g}'_\varphi(s)=\frac{\dd}{\dd s}\int_{\pb(E)} \langle\mu-f_s,\tilde{\varphi}\rangle^2 F^N_s(\dd\mu)\\
    &= \frac{\dd}{\dd s}{\left[\int_{\pb(E)} \langle\mu,\tilde{\varphi}\rangle^2 F^N_s(\dd\mu)-2\langle f_s,\tilde{\varphi}\rangle\int_{\pb(E)}\langle\mu,\tilde{\varphi}\rangle F^N_s(\dd\mu)+\langle f_s,\tilde{\varphi}\rangle^2\right]}\\
    & = \int_{\pb(E)} \langle \mu,\tilde{\varphi}\rangle^2\partial_s F^N_s(\dd\mu)-2\langle f_s,\tilde{\varphi}\rangle\int_{\pb(E)}\langle\mu,\tilde{\varphi}\rangle\partial_sF^N_s(\dd\mu)\\
    &\qquad-2\langle\partial_s f_s,\tilde{\varphi}\rangle\int_{\pb(E)}\langle\mu,\tilde{\varphi}\rangle F^N_s(\dd\mu)+2\langle f_s,\tilde{\varphi}\rangle\langle\partial_s f_s,\tilde{\varphi}\rangle\\
    &\qquad + 2\int_{\pb(E)}\langle\mu,\tilde{\varphi}\rangle\langle \mu,\partial_s\tilde{\varphi}\rangle F^N_s(\dd\mu)\\
    &\qquad-2\langle f_s,\partial_s\tilde{\varphi}\rangle\int_{\pb(E)}\langle\mu,\tilde{\varphi}\rangle F^N_s(\dd\mu)\\
    &\qquad-2\langle f_s,\tilde{\varphi}\rangle\int_{\pb(E)}\langle\mu,\partial_s\tilde{\varphi}\rangle F^N_s(\dd\mu) + 2\langle f_s,\tilde{\varphi}\rangle\langle f_s,\partial_s\tilde{\varphi}\rangle.
\end{align*}
Using the fact that for two test functions $\varphi_1,\varphi_2\in C_b(E)$:
\begin{multline*}\int_{\pb(E)}\langle\mu,\varphi_1\rangle\langle\mu,\varphi_2\rangle\partial_s F^N_s(\dd\mu) = \int_{\pb(E)}\Big\{\langle \mu,a\cdot\nabla_x\varphi_1\rangle\langle\mu,\varphi_2\rangle+\langle\mu,\varphi_1\rangle\langle \mu,a\cdot\nabla_x\varphi_2\rangle
\\+\langle\mu,\varphi_1\rangle\langle \mu,L_{\mu}\varphi_2\rangle+\langle \mu,L_\mu\varphi_1\rangle\langle\mu,\varphi_2\rangle+\frac{1}{N}\langle\mu,\Gamma_\mu(\varphi_1,\varphi_2)\rangle\Big\} F^N_s(\dd\mu),
\end{multline*}
a direct computation shows that: 
\begin{multline*}
    \tilde{g}_\varphi'(s) = \int_{\pb(E)} \Big(\langle\mu,\tilde{\varphi}\rangle-\langle f_s,\tilde{\varphi}\rangle\Big)\Big(\langle \mu,L_\mu\tilde{\varphi}\rangle-\langle  f_s,L_{ f_s}\tilde{\varphi}\rangle\Big) F^N_s(\dd\mu)\\
    +\frac{1}{N}\int_{\pb(E)}\langle \mu,\Gamma_\mu(\tilde{\varphi},\tilde{\varphi})\rangle F^N_s(\dd\mu)
    \end{multline*}
By the arithmetic-geometric mean inequality, we obtain: 
\begin{multline*}
\tilde{g}'_\varphi(s) \leq \frac{1}{2}\int_{\pb(E)} \langle\mu-f_s,\tilde{\varphi}\rangle^2 F^N_s(\dd\mu)+\frac{1}{2}\int_{\pb(E)}\big(\langle\mu,L_\mu\tilde{\varphi}\rangle-\langle f_s,L_{f_s}\tilde{\varphi}\rangle\big)^2 F^N_s(\dd\mu)\\+\frac{1}{N}\int_{\pb(E)}\langle \mu,\Gamma_\mu(\tilde{\varphi},\tilde{\varphi})\rangle F^N_s(\dd\mu).
\end{multline*}
The last term on the right-hand side can be controlled using a mild moment assumption: 
\begin{equation*}\exists\gamma>0,\,\,\forall x\in E,\quad \int_{E} |x-y|^2 P_\mu(x,\dd y) \leq \gamma,\end{equation*}
which implies
\[\Gamma_\mu(\tilde{\varphi},\tilde{\varphi})\leq \gamma\|\tilde{\varphi}\|_\mathrm{Lip}^2 \leq \gamma \e^{2(t-s)\|a\|_\mathrm{Lip}}.\]
we deduce that: 
\begin{equation}\label[Ieq]{eq:Fphibeforelastterm}\tilde{g}_\varphi'(s) \leq \frac{1}{2}\tilde{g}_\varphi(s) + \frac{\gamma\e^{2(t-s)\|a\|_\mathrm{Lip}}}{N}+\frac{1}{2}\int_{\pb(E)}\big(\langle\mu,L_\mu\tilde{\varphi}\rangle-\langle f_s,L_{f_s}\tilde{\varphi}\rangle\big)^2 F^N_s(\dd\mu).\end{equation}
Unfortunately it is not possible in general to close the argument. One would like to bound the last term on the right-hand side of \cref{eq:Fphibeforelastterm} by a quantity of the form~\cref{eq:gphi} (possibly with a different $\varphi$). This is hopeless in general because \cref{eq:gphi} depends on (the square of) a linear quantity in $\mu$ but due to the nonlinearity, the collision operator 
\[\langle Q(\mu),\tilde{\varphi}\rangle := \langle \mu, L_\mu\tilde{\varphi}\rangle\]
is at least quadratic in $\mu$ (or controlled by a quadratic quantity as soon as $P_\mu$ is Lipschitz in $\mu$ for a Wasserstein distance). This fact is analogous to the BBGKY hierarchy at the level of the empirical process. It is possible to close the argument only in linear cases, for instance when 
\[P_\mu(x,\dd y) = \int_{z\in E} K(y,z)\mu(\dd z)\dd y,\]
where $K:E\times E\to\R_+$ is a fixed symmetric interaction kernel with $\int_E K(y,z)\dd y=1$ for all $z\in E$ (note that Lipschitz continuity for $K$ is sufficient to ensure the well-posedness). In this case the collision operator is a linear operator: 
\begin{equation}\label[Ieq]{eq:linearPmu}
\langle Q(\mu),\tilde{\varphi}\rangle  = \langle \mu, K\star \tilde{\varphi}\rangle - \langle \mu, \tilde{\varphi}\rangle.  \end{equation}
Reporting into \cref{eq:Fphibeforelastterm}, we get, 
\begin{multline*}\tilde{g}_\varphi'(s) \leq \frac{1}{2}\tilde{g}_\varphi(s) + \frac{\gamma\e^{2(t-s)\|a\|_\mathrm{Lip}}}{N}+\int_{\pb(E)}\langle\mu-f_s,\tilde{\varphi}\rangle^2 F^N_s(\dd\mu)\\+\int_{\pb(E)}\langle \mu-f_s,K\star\tilde{\varphi}\rangle^2 F^N_s(\dd\mu).\end{multline*}
Since the Lipschitz norm of $\tilde{\varphi}$ is controlled by \cref{eq:lipschitztildephi} and that this bound is preserved by the convolution with $K$, using \cref{eq:GNsupphi}, we conclude that: 
\[\tilde{g}'_\varphi(s) \leq \frac{5}{2}g_N(s)\e^{2(t-s)\|a\|_\mathrm{Lip}} + \frac{\gamma\e^{2(t-s)\|a\|_\mathrm{Lip}}}{N}.\]
Integrating between 0 and $t$, and taking the supremum over $\varphi$ on the left-hand side, we can apply Gronwall lemma and conclude that
\begin{equation}\label[Ieq]{eq:chaosempiricallinear}\mathcal{W}_{D_2}^2\big(F^N_t,\delta_{f_t}\big) \leq g_N(t) \leq C(\gamma,a,t){\left(\mathcal{W}_{D_2}^2\big(F^N_0,\delta_{f_0}\big) + \frac{1}{N}\right)}.\end{equation}
The first term on the right-hand side depends on the initial condition and can be controlled using e.g. \cite{fournier_rate_2015} (using their estimates, it will determine the final rate of convergence since their best convergence rate is $\mathcal{O}(N^{-1})$); sharper estimates can be available though.

\begin{remark}
Note that despite its relative simplicity, the model \cref{eq:linearPmu} has its own interest. In \cite{carlen_kinetic_2013,carlen_kinetic_2013-1} it is called the ``choose the leader'' model. In population dynamics, it is also a time-continuous version of the well-known Moran model \cite{dawson_measure-valued_1993}. It describes the ``neutral'' evolution of a population of individuals where the death and reproduction events happen simultaneously. The kernel $K$ plays the role of a mutation kernel. A different scaling which leads to a different limit is presented in Section \cref{sec:flemingviot}.
\end{remark}

\subsubsection{The limit semi-group and the nonlinear measure-valued process}\label[I]{sec:limitsemigroup}

The previous approach is too coarse as it tries to compare directly the empirical law $F^N_t$ to its limit $\delta_{f_t}$. By looking only at the expectation of some fixed (though infinitely many) observables, we do not keep track of the detailed dynamics of the particle system. In this section we give some insights on an approach which is originally due to Gr\"unbaum \cite{grunbaum_propagation_1971} but which is has been made rigorous in \cite{mischler_new_2015,mischler_kacs_2013}. Rather than looking only at the law $F^N_t$, the idea is to compare the empirical process and a ``nonlinear'' measure-valued process through their semi-groups and generators, thus effectively keeping track of the dynamics. Note that we could define an ``obvious'' nonlinear empirical process by taking the image by $\boldsymbol{\mu}_N$ of $N$ i.i.d. $f_t$-distributed processes. This would lead us to compare $\widetilde{F}^N_t = (\boldsymbol{\mu}_N)_\# f^{\otimes N}_t$ to $F^N_t=(\boldsymbol{\mu}_N)_\# f^{N}_t$ but this would not be simpler than comparing directly $f^N_t$ to $f^{\otimes N}_t$. This could be handled by the coupling approach (Section \cref{sec:coupling}). Instead, the approach of \cite{mischler_kacs_2013, mischler_new_2015} considers the nonlinear dynamics in $\pb(E)$ from the PDE point of view. In the most abstract setting, the limiting nonlinear law $f_t$ is the solution of 
\begin{equation}\label[Ieq]{eq:abstractpde} \partial_t f_t = Q ( f_t ), \end{equation}
where $Q$ is a nonlinear operator. Assuming that this PDE is wellposed, this gives rise to a nonlinear time-continuous semi-group $(\overline{S}_t)_{t\geq0}$ acting on $\pb(E)$ such that the solution of \cref{eq:abstractpde} is given as:
\[ f_t = \overline{S}_t ( f_0 ),\]
where $f_0\in\pb(E)$ is the initial condition and
\[\overline{S}_{t+s} = \overline{S}_t \circ \overline{S}_s = \overline{S}_s \circ \overline{S}_t,\quad \partial_t \overline{S}_t = Q \circ \overline{S}_t.\]
From the stochastic point of view, the operator $\overline{S}_t$ is the dual of the transition operator $T^{f_0}_{t}$ of a nonlinear Markov process in the sense of McKean (see Appendix~\cref{appendix:timeinhomoegeneousmarkov}). The main observation is that the deterministic dynamics \cref{eq:abstractpde} can be seen as a stochastic process in $\pb(E)$, for instance, it is possible to choose a random initial condition: although it is expected to be a given $f_0\in \pb(E)$, it is also natural to take as initial condition the same as the one of the empirical process, that is a random empirical measure with $N$ points sampled from $f^N_0$. Remember that $f^N_0$ is assumed to be initially $f_0$-chaotic. With this choice, the goal is to compare the empirical process $\big(\mu_{\mathcal{X}^N_t}\big)_t$ and the nonlinear process $\big(\overline{S}_t(\mu_{\mathcal{X}^N_0})\big)_t$. The laws of these processes in $\pb(\pb(E))$ at time $t>0$ are respectively given by: 
\begin{equation}\label[Ieq]{eq:PPlawempiricalprocess}F^N_t = (\boldsymbol{\mu}_N)_\# S^N_t(f^N_0),\quad \overline{F}{}^N_t := (\overline{S}_t\circ \boldsymbol{\mu}_N)_\# f^N_0,\end{equation}
where $S^N_t$ denotes the $N$-particle semigroup acting on $\pb(E^N)$, so that $f^N_t =S^N_t(f^N_0)$. The semigroups $(S^N_t)_t$ and $(\overline{S}_t)_t$ describe the forward dynamics of the probability distributions. The dynamics of the observables is described by the dual operators acting on the space of test functions on $C_b(\pb(E))$. As explained at the beginning of this section, for the particle dynamics, everything is given in terms of $(T_{N,t})_t$ which is the semigroup acting on $C_b(E^N/\mathfrak{S}_N)$. For the nonlinear system, the following operator is defined in \cite{mischler_new_2015, mischler_kacs_2013}:
\begin{equation}\label[Ieq]{eq:Tinfty}\forall \Phi\in C_b(\pb(E)),\,\,\forall \nu\in\pb(E),\quad T_{\infty,t}\Phi(\nu) := \Phi{\left(\overline{S}_t(\nu)\right)},\end{equation}
so that
\[\int_{E^N} T_{\infty,t}\Phi\big(\mu_{\mathbf{x}^N}\big)f^N_0\big(\dd\mathbf{x}^N\big)= \big\langle \overline{F}{}^N_t,\Phi\big\rangle = \int_{E^N} \Phi{\left(\overline{S}_t\big(\mu_{\mathbf{x}^N}\big)\right)}f^N_0\big(\dd\mathbf{x}^N\big).\]
Note that the dependence on $N$ is only in the initial condition. On the other hand, it holds that: 
\[\int_{E^N} \widehat{T}_{N,t}\Phi\big(\mu_{\mathbf{x}^N}\big)f^N_0\big(\dd\mathbf{x}^N\big) = \big\langle F^N_t,\Phi\big\rangle,\]
so very loosely speaking, the goal is to prove the convergence of the semi-groups:
\[\widehat{T}_{N,t} \underset{N\to+\infty}{\longrightarrow} T_{\infty,t}.\]
In a classical setting, the convergence of a sequence of semi-groups acting on a set of test functions over a Banach space is solved by Trotter \cite{trotter_approximation_1958} by proving the convergence of the generators. The generator $\widehat{\mathcal{L}}_N$ associated to $\widehat{T}_{N,t}$ is defined by~\cref{eq:LNhat}, although its image is restricted to the subdomain $C_b(\widehat{\pb}_N(E))\subset C_b(\pb(E))$. The generator of $(T_{\infty,t})_t$ is much more delicate to define because $\pb(E)$ is only a metric space and not a Banach space. Its precise and rigorous definition is one of the main contributions of \cite{mischler_new_2015,mischler_kacs_2013}. We will give insights on this later, but for now let us assume that it is possible to define a generator $\mathcal{L}_\infty$ on a sufficiently large subset of $C_b(\pb(E))$ and such that for $\Phi$ in this subset,
\begin{equation}\label[Ieq]{eq:LinftyTinfty}\frac{\dd}{\dd t} T_{\infty,t}\Phi=\mathcal{L}_\infty[T_{\infty,t}\Phi] = T_{\infty,t}\mathcal{L}_\infty\Phi.\end{equation}
We now briefly explain how generator estimates will give an estimate on the discrepancy between $F^N_t$ and $\overline{F}{}^N_t$. First, using Lemma \cref{lemma:grunbaumlemma}, it is sufficient to look at the moment measures $F^{k,N}_t,\overline{F}{}^{k,N}_t\in\pb(E^k)$ for all $k\in\N$. Let $\varphi_k\in C_b(E^k)$ be a test function and let
\[R_{\varphi_k}(\nu)\equiv \Phi_k(\nu) := \langle \nu^{\otimes k},\varphi_k\rangle,\]
be the associated polynomial function. We recall that the operators $T_{N,t}$ and $\mathcal{L}_N$ are directly linked to their empirical versions \cref{eq:TNhat}, \cref{eq:LNhat} by $\boldsymbol{\mu}_N$ and the linear transpose map: 
\[ \boldsymbol{\mu}_N^\mathrm{T} :  C_b ( \pb ( E ) ) \rightarrow C_b {\left( E^N \right)},\quad \Phi \mapsto \Phi \circ \boldsymbol{\mu}_N. \]
By definition, it holds that: 
\[\big\langle F^{k,N}_t,\varphi_k\big\rangle = \big\langle F^N_t, \Phi_k\big\rangle = \big\langle f^N_0, T_{N,t}[\Phi_k\circ\boldsymbol{\mu}_N] \big\rangle = \big\langle f^N_0, T_{N,t}[\boldsymbol{\mu}_N^\mathrm{T}\Phi_k]\big\rangle,\]
and
\[\big\langle\overline{F}{}^{k,N}_t,\varphi_k\big\rangle = \big\langle \overline{F}{}^N_t,\Phi_k\big\rangle = \big\langle f^N_0, T_{\infty,t}\Phi\circ \boldsymbol{\mu}_N\big\rangle = \big\langle f^N_0, \boldsymbol{\mu}_N^\mathrm{T}[T_{\infty,t}\Phi_k]\big\rangle.\]
The difference between these two quantities is controlled by using the formula for $0 \leq t \leq T$ :
\begin{align*} 
T_{N,t} \boldsymbol{\mu}_N^\mathrm{T} - \boldsymbol{\mu}_N^\mathrm{T} T_{\infty,t} &= - \int_0^t \frac{\dd}{\dd s} {\left[ T_{N,t-s} \boldsymbol{\mu}_N^\mathrm{T} T_{\infty,s} \right]} \dd s \\
&= \int_0^t T_{N,t-s} {\left[ \mathcal{L}_N \boldsymbol{\mu}_N^\mathrm{T} -  \boldsymbol{\mu}_N^\mathrm{T} \mathcal{L}_{\infty} \right]} T_{\infty,s} \dd s,\end{align*}
which leads to the following bound: 
\begin{align}
{\left| \big\langle F^{k,N}_t\varphi_k\big\rangle -\big\langle \overline{F}{}^{k,N}_t , \varphi_k \big\rangle \right|} &\leq \int_0^t {\left| \big\langle f_0^N , T_{N,t-s} {\left[ \mathcal{L}_N\boldsymbol{\mu}_N^\mathrm{T} -  \boldsymbol{\mu}_N^\mathrm{T} \mathcal{L}_{\infty} \right]} T_{\infty,s} \Phi_k \big\rangle \right|} \dd s\nonumber \\
&= \int_0^t {\left| \big\langle f^N_{t-s} , {\left[ \mathcal{L}_N \boldsymbol{\mu}_N^\mathrm{T} -  \boldsymbol{\mu}_N^\mathrm{T} \mathcal{L}_{\infty} \right]} T_{\infty,s}\Phi_k \big\rangle \right|} \dd s\nonumber \\
&= \int_0^t {\left| \big\langle f^N_{t-s} , \mathcal{L}_N {\left[ T_{\infty,s} \Phi_k \circ \boldsymbol{\mu}_N \right]} - \mathcal{L}_{\infty} {\left[ T_{\infty,s} \Phi_k \right]} \circ \boldsymbol{\mu}_N \big\rangle \right|} \dd s\nonumber \\
&\leq T \sup_{0 \leq t \leq T} {\left\|{\left( \widehat{\mathcal{L}}_N {\left[ T_{\infty,t} \Phi_k \right]} - \mathcal{L}_{\infty} {\left[ T_{\infty,t} \Phi_k \right]}\right)} \circ \boldsymbol{\mu}_N \right\|}_{\infty}. \label[Ieq]{eq:MMgeneratorestimate}
\end{align}
At this point, in order to obtain a convergence bound in $N$, a generator estimate is needed to compare the behaviour of the empirical generator $\widehat{\mathcal{L}}_N$ to the one of $\mathcal{L}_{\infty}$ against $T_{\infty,t} \Phi_k$ (the map $\boldsymbol{\mu}_N$ is just an artefact to write this comparison in $E^N$).
\begin{remark} The estimate \cref{eq:MMgeneratorestimate} can be made more uniform in time as soon as the particle system preserves some quantity $m:E^N\to\R_+$, in which case \cref{eq:MMgeneratorestimate} becomes: 
\begin{multline*}{\left| \big\langle F^{k,N}_t\varphi_k\big\rangle -\big\langle \overline{F}{}^{k,N}_t , \varphi_k \big\rangle \right|} \\ \leq \sup_{0 \leq t \leq T} {\Big\|\frac{1}{m}{\left( \widehat{\mathcal{L}}_N {\left[ T_{\infty,t} \Phi_k \right]} - \mathcal{L}_{\infty} {\left[ T_{\infty,t} \Phi_k \right]}\right)} \circ \boldsymbol{\mu}_N \Big\|}_{\infty}\int_0^T \langle f^N_{t-s},m\rangle \dd s. \end{multline*}
\end{remark}

In Section \cref{sec:abstracttheorem} we will review the abstract theorem of \cite{mischler_new_2015} which shows how to recover propagation of chaos in the usual framework from an estimate on \cref{eq:MMgeneratorestimate}. Before doing that, we give more insights on the definition of the generator $\mathcal{L}_\infty$ within the framework of \cite{mischler_new_2015}. 

\subsubsection{More on the limit generator}\label[I]{sec:limitgenerator}

As an introductory example, let us consider a mean-field generator of the form \cref{eq:Nparticlemeanfieldgenerator} and a tensorized test function of order two:
\[\varphi_2 = \varphi^1\otimes \varphi^2,\]
as well as the associated polynomial function on $\pb(E)$ defined by $\Phi_2(\nu) = \langle \nu^{\otimes 2},\varphi_2\rangle$. We recall that it is not a restriction to consider tensorized test functions \cite[Chapter~3, Theorem 4.5 and Proposition 4.6, pp.113-115]{ethier_markov_1986}. Then, a direct computation which is detailed in Lemma \cref{lemma:quadraticmeanfield} (a similar computation is also used in the proof of Theorem \cref{thm:martingaleweakpathwiseempiricalchaos}) gives: 
\begin{multline*}\mathcal{L}_N[\Phi_2\circ \boldsymbol{\mu}_N](\mathbf{x}^N) = \langle \mu_{\mathbf{x}^N},L_{\mu_{\mathbf{x}^N}}\varphi^1\rangle\langle \mu_{\mathbf{x}^N}, \varphi^2\rangle + \langle \mu_{\mathbf{x}^N}, \varphi^1\rangle\langle \mu_{\mathbf{x}^N}, L_{\mu_{\mathbf{x}^N}}\varphi^2\rangle \\+ \frac{1}{N}\big\langle \mu_{\mathbf{x}^N}, \Gamma_{L_{\mu_{\mathbf{x}^N}}}(\varphi^1,\varphi^2)\big\rangle,\end{multline*}
where $\Gamma$ is the carr\'e du champ operator. In order to have 
\[\|(\widehat{\mathcal{L}}_N[\Phi_2]-\mathcal{L}_\infty[\Phi_2])\circ\boldsymbol{\mu}_N\|_\infty\underset{N\to+\infty}{\longrightarrow}0,\]
the operator $\mathcal{L}_\infty$ should necessarily satisfy: 
\[\mathcal{L}_\infty\Phi_2(\nu)  = \langle \nu,L_\nu\varphi^1\rangle\langle \nu,\varphi^2\rangle+\langle \nu,\varphi^1\rangle\langle \nu, L_\nu\varphi^2\rangle.\]
This computation is generalised (see Lemma \cref{lemma:largedegreemeanfield}) to any $k$-fold tensorized test function
\[\varphi_k = \varphi^1\otimes\ldots\otimes\varphi^k,\]
and the operator $\mathcal{L}_\infty$ is defined against monomial functions by: 
\begin{equation}\label[Ieq]{eq:Linftyagainstmonomial}\mathcal{L}_\infty\Phi_k(\nu) = \sum_{i=1}^k \langle Q(\nu),\varphi^i\rangle\prod_{j\ne i} \langle \nu,\varphi^j\rangle,\end{equation}
where $\Phi_k(\nu)=\langle \nu^{\otimes k},\varphi_k\rangle$ and with a slight abuse of notation, we write
\[\langle Q(\nu),\varphi\rangle\equiv \langle \nu,L_\nu\varphi\rangle,\]
for the integral of a test function $\varphi$ against $Q(\nu)$ (which is not a probability measure). The relation \cref{eq:Linftyagainstmonomial} can be extended to any polynomial by linearity. We would also get the same relation (with a different operator $Q$) for a Boltzmann operator (see Lemma \cref{lemma:quadraticboltzmann}).

A natural idea would be to use the relation \cref{eq:Linftyagainstmonomial} as a \emph{definition} of an operator acting on polynomials and then extend it to the completion of the space of polynomials, which is a large Banach subset of $C_b(\pb(E))$. However, this would not necessarily imply that the right-hand side of \cref{eq:MMgeneratorestimate} goes to zero for any $\Phi$ (since for a given polynomial, the convergence implied by Lemma \cref{lemma:largedegreemeanfield} may not be uniform in the degree or number of monomials) and proving this convergence does not seem to be an easy task. We recall that the final goal is to apply $\mathcal{L}_\infty$ to the test functions 
\begin{equation}\label[Ieq]{eq:Tinftytphik}T_{\infty,t}\Phi_k(\nu) = \big\langle\overline{S}_t(\nu)^{\otimes k},\varphi_k\big\rangle,\end{equation}
which are not polynomial in general. In particular, the relation \cref{eq:Linftyagainstmonomial} needs to be extended to be able to write 
\[\mathcal{L}_\infty T_{\infty,t} \Phi_k(\nu) = \sum_{i=1}^k \big\langle Q\big(\overline{S}_t(\nu)\big),\varphi^i\big\rangle\prod_{j\ne i} \big\langle \overline{S}_t(\nu),\varphi^j\big\rangle,\]
in order to have \cref{eq:LinftyTinfty}.

\begin{remark} \label[I]{rem:markovlimitappendix}
The stochastic point of view gives more insight on the form of the nonlinearity in \cref{eq:Tinftytphik}. Under the assumption that $f_t$ is the law of a nonlinear Markov process in the sense of McKean (see Appendix \cref{appendix:timeinhomoegeneousmarkov}), one can write the dual form when $\varphi_k = \varphi^1\otimes\ldots\otimes\varphi^k$: 
\begin{equation}\label[Ieq]{eq:TinftytphikMcKean}T_{\infty,t}\Phi_k(\nu) = \prod_{i=1}^k \big\langle \nu,T^{\nu}_t\varphi^i\big\rangle,\end{equation}
where $T^{\nu}_t$ is the nonlinear transition operator of the process. Up to the dependency of $T^{\nu}_t$ on the measure argument $\nu$, the test function $T_{\infty,t}\Phi_k$ is thus close to be a polynomial. When $T^{\nu}_t$ does not depend on $\nu$, $f_t=\overline{S}(f_0)$ is the law of a classical time homogeneous Markov process and thus satisfies a linear equation. In this case, everything is much simpler because the operator $T_{\infty,t}$ acts on the space of polynomial and it is not necessary to extend $\mathcal{L}_\infty$ to a larger subspace of $C_b(\pb(E))$. One could actually bypass the definition of the limit generator. Note also that the linear case is the framework of our toy example in Section \cref{sec:toyexampleempiricalprocess}. Unfortunately, in the nonlinear case, the dual backward form \cref{eq:TinftytphikMcKean} does not really seem to be more helpful since $T^{\nu}_t$ has no reason to behave well with respect to its measure argument. We finally point out that all the argument in \cite{mischler_new_2015,mischler_kacs_2013} is more general as it is entirely based on \cref{eq:Tinftytphik} and does not use the fact that $f_t$ is the law of a nonlinear Markov process in the sense of McKean (even though it is the underlying application case). 
\end{remark}

In his seminal article \cite{grunbaum_propagation_1971}, Gr\"unbaum originally relies on a clever completion of the space of polynomial functions. He then identifies $\mathcal{L}_\infty$ and proves the convergence of the generators on a class $C'$ of ``continuously differentiable functions'' on $\pb(E)$. In order to apply Trotter's result on $T_{\infty,t}$, Gr\"unbaum uses an unproven smoothness assumption on the nonlinear operator $\overline{S}_t$ which ensures that \cref{eq:Tinftytphik} belongs to the class $C'$. This assumption has since been proved to be false for some models. 

The generator $\mathcal{L}_\infty$ is rigorously defined in \cite{mischler_new_2015,mischler_kacs_2013} using a new notion of differential calculus on the space $\pb(E)$ that is very briefly sketched below. The fundamental idea is to consider a distance on $\pb(E)$ inherited from a normed vector space~$\mathcal{G}$ to be chosen. Let $m_\mathcal{G}:E\to\R_+$ be given together with the weighted subspace of probability measures: 
\[\pb_\mathcal{G}(E):=\{f\in\pb(E),\,\,\langle f,m_\mathcal{G}\rangle<+\infty\}.\]
The weight function $m_\mathcal{G}$ may typically be a polynomial function (in which case $\pb_\mathcal{G}$ is the space of probability measures with a bounded moment) but may also depend on the normed vector space $\mathcal{G}$ which is assumed to contain the space of increments: 
\begin{equation}\label[Ieq]{eq:spaceincrements}\mathcal{I}\pb_\mathcal{G}(E) := \Big\{f_1-f_2,\,\,f_1,f_2\in\pb_\mathcal{G}(E)\Big\}\subset \mathcal{G}.\end{equation}
This naturally defines a distance on $\pb_\mathcal{G}(E)$ by: 
\[\forall f_1,f_2\in \pb_\mathcal{G}(E),\quad d_\mathcal{G}(f_1,f_2) := \|f_1-f_2\|_\mathcal{G}.\]
\begin{example}
A simple example is obtained when taking for $\mathcal{G}$ the Banach space of finite mass signed measures endowed with the total variation distance. $\pb_\mathcal{G}(E)$ is then a subspace of $\pb(E)$ endowed with the total variation distance, while the space $\mathcal{I}\pb_\mathcal{G}(E)$ of its increments is seen as a subset of the finite mass signed measures. 
\end{example}
Other examples and their relations with the distances defined in Section \cref{sec:topology} are detailed in \cite[Section 3.2]{mischler_new_2015}. With such notion of distance, it is possible to define a differential structure on the space of probability measures which is inherited from the classical differentiable structure of the Banach space $\mathcal{G}$. Namely, a test function $\Phi: \pb_\mathcal{G}(E)\to \R$ is said to be continuously differentiable at $f\in\pb_\mathcal{G}(E)$ when there exists a continuous linear application $\dd\Phi[f] : \mathcal{G}\to\R$ and a constant $C>0$ such that: 
\begin{equation}\label[Ieq]{eq:diffmeasure}\forall g\in\pb_\mathcal{G}(E),\quad\big|\Phi(g)-\Phi(f)-\langle \dd\Phi[f],g-f\rangle_{\mathcal{G}',\mathcal{G}}\big|\leq Cd_\mathcal{G}(f,g).\end{equation}
Note that the bracket in the inequality is the duality bracket between $\mathcal{G}'$ and $\mathcal{G}$. The main difference with the usual notion of differentiability in a Banach space is that the space of increments \cref{eq:spaceincrements} has no vectorial structure. More details on this notion of differential calculus is given in \cite[Section 3.3 and Section 3.4]{mischler_new_2015} with a special focus on polynomial functions. This definition can be extended to higher order differentiability and to functions with values in $\pb_{\widetilde{\mathcal{G}}}(E)$ instead of $\R$ (which is the case of the operator $\overline{S}_t$). 

The definition of $\mathcal{L}_\infty$ then directly comes from the differentiation of the definition of the pullback semigroup \cref{eq:Tinfty}: let $\Phi$ be a continuously differentiable function on $\pb_\mathcal{G}(E)$, then for all $\nu\in\pb_\mathcal{G}(E)$, by the composition rule (see \cite[Lemma 3.12]{mischler_new_2015}), it holds that:
\begin{align}\label[Ieq]{eq:Linftydef}
    \mathcal{L}_\infty\Phi(\nu) = \frac{\dd}{\dd t}T_{\infty,t}\Phi(\nu)\Big|_{t=0} =\frac{\dd}{\dd t}\Phi\big(\overline{S}_t(\nu)\big)\Big|_{t=0}&={\left\langle \dd\Phi[\nu],\frac{\dd}{\dd t} \overline{S}_t(\nu)\Big|_{t=0}\right\rangle}\nonumber\\
    &= \langle \dd\Phi[\nu], Q(\nu)\rangle. 
\end{align}
This computation is almost rigorous up to the assumption that $Q(\nu)\in\mathcal{G}$. The precise assumptions on $(\overline{S}_t)_t$ which make this computation fully rigorous are given by \cite[Assumption (A2)]{mischler_new_2015} and will be summarised in the next section. 

\begin{example}[Generators estimate for jump processes] 
Quite formal computations may also motivate the introduction of a proper notion of differential calculus on $\pb(E)$ and lead to generators estimates. Taking the example of jump processes, we have seen that: 
\[\widehat{\mathcal{L}}_N \Phi {\left( \mu_{\mathbf{x}^N} \right)} = \sum_{i=1}^N \int_E \lambda{\left(x^i,\mu_{\mathbf{x}^N}\right)}{\left[ \Phi{\left(\mu_{\mathbf{x}^N}-\frac{1}{N}\delta_{x^i} +\frac{1}{N}\delta_y \right)}- \Phi {\left( \mu_{\mathbf{x}^N} \right)} \right]} P_{\mu_{\mathbf{x}^N}} {\left( x^i,\dd y \right)}.\]
The term in the integral is precisely of the form \cref{eq:diffmeasure} with an increment of size $1/N$. Assuming that it is possible to differentiate $\Phi$, we get: 
\begin{align*}
&\widehat{\mathcal{L}}_N \Phi {\left( \mu_{\mathbf{x}^N} \right)}\\
&= \sum_{i=1}^N \int_E \lambda{\left(x^i,\mu_{\mathbf{x}^N}\right)}{\left[ {\left\langle \dd \Phi {\left( \mu_{\mathbf{x}^N} \right)} , -\frac{1}{N}\delta_{x^i} + \frac{1}{N}\delta_y \right\rangle} + o {\left( \frac{1}{N} \right)} \right]} P_{\mu_{\mathbf{x}^N}} {\left( x^i,\dd y \right)} \\
&= {\left\langle \dd \Phi {\left( \mu_{\mathbf{x}^N} \right)}  , Q {\left( \mu_{\mathbf{x}^N} \right)} \right\rangle} + o(1) = \mathcal{L}_{\infty} \Phi {\left(  \mu_{\mathbf{x}^N} \right)} + o( 1 ),
\end{align*}
where $Q$ is given in the weak form by the left-hand side of \cref{eq:pdepdmp} (with $a=0$).
\end{example}

\subsubsection{The abstract theorem}\label[I]{sec:abstracttheorem}

The main theorem \cite[Theorem 2.1]{mischler_new_2015} is based on the following set of assumptions. The first one is the only one which concerns the particle system (it is always implicitly assumed). The second and third ones are motivated by the previous sections. The fourth and fifth ones are stated more informally, more details on their role will be given in the sketch of the proof of the main theorem.

\begin{assumption}\label[I]{assum:mischlermouhot} The following assumptions are respectively numbered (A1) to (A5) in \cite{mischler_new_2015}. 
\begin{enumerate}[(1)]
    \item \textbf{(On the particle system).} The $N$-particle semigroup $( T_{N,t} )_{t \geq 0}$ is a strongly continuous semigroup on $C_b(E^N)$ with generator $\mathcal{L}_N$.
    \item \textbf{(Existence of the pull-back semigroup).} There exists a Banach space $\mathcal{G}$ such that the nonlinear semigroup $(\overline{S}_t)_t$ on $\pb_\mathcal{G}(E)$ is Lipschitz for the distance $d_\mathcal{G}$ uniformly in time. The operator $Q:\pb_\mathcal{G}(E)\to\mathcal{G}$ is bounded and $\delta$-H\"older for a $\delta\in(0,1]$. This implies the existence of the limit generator $\mathcal{L}_\infty$ defined by \cref{eq:Linftydef} (see \cite[Lemma 4.1]{mischler_new_2015}).
    \item \textbf{(Generators estimate).} There exists a sequence $\varepsilon(N)$ such that $\varepsilon(N)\to0$ as $N\to+\infty$ and such that for sufficiently regular test functions $\Phi$ on $\pb_\mathcal{G}(E)$, 
    \begin{equation} \label[Ieq]{eq:controlgenerator}
    \big\| (\widehat{\mathcal{L}}_N \Phi  -  \mathcal{L}_{\infty} \Phi ) \circ \boldsymbol{\mu}_N \big\|_{\infty} \leq \varepsilon(N) \|\Phi\|_{C^{k,1}(\pb_\mathcal{G}(E))},
    \end{equation}
    where $\|\cdot\|_{C^{k,1}(\pb_\mathcal{G}(E))}$ is a norm related to the notion of higher order differentiability. 
    \item \textbf{(Differential stability of the nonlinear semigroup).} The nonlinear semi-group $( \overline{S}_{t} )_{t}$ is differentiable (for the generalised version of the notion of differentiability mentioned above), and its derivatives can be bounded uniformly in time using adapted $\mathcal{G}$-dependent norms.
    \item \textbf{(Weak stability of the nonlinear semigroup).} For a Banach space $\widetilde{\mathcal{G}}$ possibly different from $\mathcal{G}$, the nonlinear semi-group is Lipschitz for the distance~$d_{\widetilde{\mathcal{G}}}$. More generally this can be replaced by the existence of a concave modulus of continuity $\Theta_T:\R_+\to\R_+$ such that for all $f_0,g_0\in \pb_{\widetilde{\mathcal{G}}}(E)$, 
    \[ \sup_{0 \leq t \leq T} d_{\widetilde{\mathcal{G}}}{\left( \overline{S}_{t} {\left( f_0 \right)} ,  \overline{S}_{t} {\left( g_0 \right)} \right)} \leq \Theta_T {\left( d_{\widetilde{\mathcal{G}}} {\left( f_0 , g_0 \right)} \right)}. \]
\end{enumerate}
\end{assumption}

The last assumption (5) quantifies the continuity of the operator $\overline{S}_{t}$ with respect to its measure argument: together with the assumption (4), it should be understood as a regularity assumption with respect to the initial condition for the limit nonlinear PDE. The fact that it holds for a possibly different Banach space $\widetilde{G}$ is essentially technical: it allows more flexibility in the choice of a well-suited topology to deal with the nonlinear problem.

The following abstract theorem is stated and proved in \cite[Theorem 2.1]{mischler_new_2015}.

\begin{theorem}[The abstract theorem in \cite{mischler_new_2015}] \label[I]{thm:abstractMischler} Let Assumption \cref{assum:mischlermouhot} hold true. Let $T>0$, $k\in\N$ and $N\geq 2k$. Then there exist a continuously embedded subset $\mathcal{F}\subset C_b(E)$ and some absolute constants $C, C(T)>0$ such that for any tensorized test function
\[\varphi_k = \varphi^1\otimes\ldots\otimes\varphi^k\in\mathcal{F}^{\otimes k},\]
it holds that:
\begin{align}
\sup_{0 \leq t \leq T} \big| \big\langle f^{k,N}_t - f_t^{\otimes k} , \varphi_k \big\rangle \big| \leq C\frac{k^2 \| \varphi_k \|_{\infty}}{N} &+ k^2 C(T) \| \varphi_k \|_{\mathcal{F}_1} \varepsilon ( N ) \nonumber \\
&+ k \| \varphi_k \|_{\mathcal{F}_2} \Theta_T \Big(\mathcal{W}_{d_{\widetilde{\mathcal{G}}}} \big( \overline{F}{}^N_0 , \delta_{f_0} \big) \Big),\label[Ieq]{eq:abstractthm}
\end{align}
where $\varepsilon(N)$ is defined in Assumption \cref{assum:mischlermouhot}(3), $\mathcal{W}_{d_{\widetilde{\mathcal{G}}}}$ is a Wassertein distance on the space $\pb ( \pb ( E ) )$ related to $\widetilde{\mathcal{G}}$ given by Assumption \cref{assum:mischlermouhot}(5) and $\| \cdot \|_{\mathcal{F}_1}$ and $\| \cdot \|_{\mathcal{F}_2}$ are some norms on $C_b ( E^k )$ which are defined in the complete version of Assumption \cref{assum:mischlermouhot} (see \cite[Section 4]{mischler_new_2015}).
\end{theorem}

\begin{proof}[Proof (main ideas)]

The proof in \cite{mischler_new_2015} relies on three main steps:
\begin{itemize}
    \item Approximate $f^{k,N}_t$ by $F^{k,N}_t$ thanks to Lemma \cref{lemma:grunbaumlemma}.
    \item Approximate $f^{\otimes k}_t = \overline{S}{}_t^{\otimes k} ( f_0 )$ by $\E_{\mathcal{X}^N_0} \overline{S}{}_t^{\otimes k} ( \mu_{\mathcal{X}_0}^N )=\overline{F}{}^{k,N}_t$. Since initial chaos $ \mu_{\mathcal{X}_0}^N \to f_0 $ is assumed, this should stem from the regularity on the limit equation.
    \item Compare the time evolution of $F^{k,N}_t$ to the one of $\overline{F}{}^{k,N}_t$, which motivates the content of Section \cref{sec:limitsemigroup} and Section \cref{sec:limitgenerator}. 
\end{itemize}

We recall that $F^{k,N}_t$ and $\overline{F}{}^{k,N}_t$ are the moment measures (Definition \cref{def:momentmeasure}) associated to the laws $F^N_t$ and $\overline{F}{}^N_t$ defined by \cref{eq:PPlawempiricalprocess}. Each of the three terms on the right-hand side of \cref{eq:abstractthm} thus comes from the splitting:
\begin{align}\label[Ieq]{eq:abstracttheoremsplitting}
\big\langle f^{k,N}_t - f_t^{\otimes k} , \varphi_k \big\rangle = \big\langle f^{k,N}_t - F^{k,N}_t , \varphi_k \big\rangle &+ \big\langle F^{k,N}_t - \overline{F}{}^{k,N}_t , \varphi_k \big\rangle + \big\langle \overline{F}{}^{k,N}_t - f_t^{\otimes k} , \varphi_k \big\rangle.
\end{align}
The first term on the right-hand side is handled with rate $\mathcal{O} ( k^2 N^{-1} )$ by the approximation Lemma \cref{lemma:grunbaumlemma}, using purely combinatorial arguments. The second term is technically the most difficult one. Assumption \cref{assum:mischlermouhot}(2) gives a precise meaning of the relation \cref{eq:MMgeneratorestimate} formally derived earlier. The role of Assumption \cref{assum:mischlermouhot}(3) is thus self-explanatory and Assumption \cref{assum:mischlermouhot}(4) ensures that $\Phi=T_{\infty,t}\Phi_k$ has enough regularity to be taken as a test function in \cref{eq:controlgenerator}. The third term contains two approximations: the first one is how well the initial data is approximated by a (random) empirical measure and then how well this error is propagated in time, which is Assumption~\cref{assum:mischlermouhot}(5).   
\end{proof}

Applications of the abstract Theorem \cref{thm:abstractMischler} to classical models can be found in~\cite{mischler_new_2015}. The assumptions are rigorously justified for Maxwell molecules with cut-off (see Section \cref{sec:Boltzmannclassicalmodels}), the classical McKean-Vlasov diffusion (with a non-optimal convergence rate) and a mixed jump-diffusion model. The main advantage of this abstract method is its wide range of applicability, although each model requires a careful and dedicated verification of the five assumptions. The choice of the different spaces $\mathcal{G}$ indeed strongly depends on the structure of the model. In the companion paper~\cite{mischler_kacs_2013}, the abstract method is developed in a more general framework: the five assumptions are modified to include conservation relations in order to treat the case of Boltzmann models with unbounded collision rates, possibly uniformly in time. The results will be summarised in Section \cref{sec:kacprogram}. 

\begin{remark}[BBGKY hierarchy, statistical solution and limit generator]
We previously made the remark that taking the limit $N\to+\infty$ in the BBGKY hierarchy~\cref{eq:bbgkyboltzmannweak} or \cref{eq:bbgkymeanfieldweak} leads to an infinite hierarchy of equations called the Boltzmann hierarchy (Remark \cref{rem:chaosbbgky}). By the Hewitt-Savage theorem, at every time $t>0$, the Boltzmann hierarchy is associated to a unique $\pi_t\in\pb(\pb(E))$ which is sometimes called a \emph{statistical solution} of the limit problem. In \cite[section 8]{mischler_kacs_2013}, the authors show that for cutoff Boltzmann models, given an initial $\pi_0\in\pb(\pb(E))$, the statistical solution $\pi_t$ is the unique solution to the evolution equation 
\[\partial_t\pi = \mathcal{L}^\infty\pi_t,\]
where $\mathcal{L}^\infty$ is the formal adjoint of the limit generator $\mathcal{L}_\infty$ on $C_b(\pb(E))$. It means that $\pi_t$ satisfies the weak equation 
\[\forall\Phi \in C_b(\pb(E)),\quad \frac{\dd}{\dd t}\langle \pi_t,\Phi\rangle = \langle \pi_t, \mathcal{L}_\infty\Phi\rangle.\]
When $\pi_0=\delta_{f_0}$, there exists a unique statistical solution which is \emph{chaotic} in the sense that this solution is equal to $\delta_{f_t}$ where $f_t$ solves the nonlinear PDE. As already explained many times, this fact is equivalent to the propagation of chaos. Moreover it is proven in \cite[section 8]{mischler_kacs_2013} that the operators which generate (in a sense which is made rigorous) the BBGKY-hierarchy converge towards the generators of the processes related to moment measures of $\pi_t$. However it should be noted that in general there exists many other statistical solutions. The notion of statistical solution is an important notion in fluid mechanics, where it originates. 
\end{remark}

\subsection{Large Deviation Related Methods}\label[I]{sec:largedeviations}

Various approaches related to large deviations theory are investigated here. It is possible to motivate them by looking back at Remark \cref{rem:chaosasLGN}, which suggests that chaos can be seen as a kind of weak law of large numbers, as it implies the weak convergence:
\[ \langle \mu_{\mathcal{X}^N} , \varphi \rangle - \mathbb{E}_{\mathcal{X}^N}{\left[ \varphi{\left( X^{1,N} \right)}\right]} \underset{N \to +\infty}{\longrightarrow} 0. \]
When a strong law of large numbers holds, it is natural to look at the fluctuations of $\langle \mu_{\mathcal{X}^N} , \varphi \rangle$ by establishing some weak central limit theorem. Nonetheless one can look at this issue the other way round, trying to deduce some weak law of large numbers from a fluctuation result. Indeed, the usual central limit theorem implies a weak version of the law of large numbers, although the latter is classically proven using quite different tools. Actually, it should be noted that, if a typical large deviation estimate of the form $\mathbb{P}(D(\mu_{\mathcal{X}^N_t},f_t)>\varepsilon)\leq \e^{-N I(\varepsilon)}$ holds for a distance $D$ on $\pb(E)$ and $I$ a rate function which has a unique zero, then propagation of chaos immediately follows from the Borel-Cantelli lemma. Note however that if one is only interested in the law of large numbers, a large deviation result may be quite overworked. As we shall see, the large deviation theory provides new tools and useful insights on propagation of chaos. In Section \cref{sec:chaosthroughldp} we give a mostly historical description of large deviation results which imply as a byproduct a weak form propagation of chaos in some specific cases. These results are related to Laplace's theory of fluctuations, which has been widely used in statistical physics to study out of equilibrium systems. The relative entropy functional (Definition \cref{def:entropyfisher}) plays a crucial role in this analysis: in Section \cref{sec:chaosfromentropy}, we gather classical results which link propagation of chaos and entropy bounds. Section \cref{concentrationI} is devoted to the study of (quantitative) concentration inequalities which will be useful in the following sections to strengthen propagation of chaos results. We will later give a brief overview of ``pure'' large deviation results which go beyond propagation of chaos in Section~\cref{sec:fluctuations}. Some classical material on large deviation theory can be found in Appendix~\cref{appendix:largedeviations}.

\subsubsection{Chaos through Large Deviation Principles}\label[I]{sec:chaosthroughldp}

In the seminal article \cite{ben_arous_methode_1990}, the authors improve results from \cite{kusuoka_gibbs_1984} and \cite{bolthausen_laplace_1986} on Large Deviation Principles (LDP) for Gibbs measure and obtain as a byproduct a pathwise propagation of chaos result for the McKean-Vlasov diffusion. Firstly, \cite[Theorem A]{ben_arous_methode_1990} below states a large deviation principle for Gibbs measures with a polynomial potential. 

\begin{theorem}[Polynomial Potential]\label[I]{thm:ldpgibbspolynomial} Let $\mathscr{E}$ be a Polish measurable space. Let $\mu_0\in\pb(\mathscr{E})$. Let us
consider a random vector $\mathcal{X}^N$ in $\mathscr{E}^N$, distributed according to the Gibbs measure:
\begin{equation}\label[Ieq]{eq:gibbsmeasure} \mu^N \big( \dd \mathbf{x}^N \big)= \frac{1}{Z_N} \exp {\left[ NG ( \mu_{\mathbf{x}^N} ) \right]} \mu_0^{\otimes N} \big( \dd \mathbf{x}^N \big),\end{equation}
where $Z_N$ is a normalization constant and $G$ is a polynomial function on $\mathcal{P} ( \mathscr{E} )$ (called the energy functional) of the form:
\[G(\mu) = \sum_{k=2}^r \langle \mu^{\otimes k}, V_k\rangle,\]
for some symmetric continuous bounded functions $V_k$ on $\mathscr{E}^k$. Then the laws of $\mu_{\mathcal{X}^N}$ satisfy a large deviation principle in  $\mathcal{P} ( \mathcal{P} ( \mathscr{E} ) )$ with speed $N^{-1}$ and rate function $\mu \mapsto H ( \mu | \mu_0 ) - G(\mu) - \inf_{\mathcal{P} ( \mathscr{E} )} ( H ( \cdot | \mu_0 ) - G )$. 
\end{theorem}

Denote by $m_0$ the infimum of $H ( \cdot | \mu_0 ) - G $ in $\mathcal{P} ( \mathscr{E} )$ and $\mathcal{P}^{m_0}$ the set of probability measures which achieve it. The study of $\mathcal{P}^{m_0}$ is related to the study of the quadratic form $\Theta_{\nu}$ on $L^2_0 \left( \mathscr{E} , \dd  \nu \right)$ (the space of centered square $\nu$-integrable functions on $E$) defined for any $\nu$ in $\mathcal{P}^{m_0}$ by:
\[\forall f,g\in L^2_0(\mathscr{E},\dd \nu),\quad \langle \Theta_\nu f,g\rangle := \sum_{k=2}^r k(k-1)\langle \nu^k, f\otimes g\otimes 1^{k-2}V_k\rangle.\]
The following \cite[Theorem B]{ben_arous_methode_1990} quantifies the fluctuations of $\mu_{\mathcal{X}^N}$ in the \emph{non-degenerate case}. Analogous results for the degenerate case are given in \cite[Theorem~C]{ben_arous_methode_1990}.

\begin{theorem}[Chaos and Fluctuations]\label[I]{thm:chaosfromLDP}
Assume that $\pb^{m_0}$ is non degenerate in the sense that for all $\nu\in\mathcal{P}^{m_0}$,
\[\mathrm{Ker} ( \mathrm{Id} - \Theta_{\nu} ) = \{ 0 \}.\] 
In this case, let us consider the quantities:
\[ d ( \nu ) := [ \det ( \mathrm{Id} - \Theta_{\nu} ) ]^{-1/2}, \quad  \bar{d} ( \nu ) := \frac{d ( \nu )}{\sum_{\nu' \in \mathcal{P}^{m_0}} d ( \nu' )}.\]
Then the following properties hold. 
\begin{enumerate}
    \item The set $\mathcal{P}^{m_0}$ is finite.
    \item $\lim_{N \rightarrow \infty} \e^{N m_0} Z_N = \sum_{\nu \in \mathcal{P}^{m_0}} d ( \nu )$.
    \item For any integer $k \geq 1$ and any $\varphi$ in $C_b ( \mathscr{E}^k )$,
    \[ \langle \mu^N , \varphi \otimes 1^{\otimes N-k} \rangle \xrightarrow[N \rightarrow \infty]{} \sum_{\nu \in \mathcal{P}^{m_0}} \bar{d} ( \nu ) \langle \nu^{\otimes k} , \varphi \rangle. \]
    \item The random measures $\mu_{\mathcal{X}^N}\in\pb(\mathscr{E})$ satisfy local and global Central Limit Theorems. 
\end{enumerate}
\end{theorem}

When $\mathcal{P}^{m_0}$ reduces to a single non-degenerate minimizer $f$, then the third assertion exactly tells that the sequence $( \mu^N )_N$ is $f$-chaotic. 

Going back to an interacting particle system, let $f^N_0\in \pb(E^N)$ be an initial law and $f^N_{[0,T]} \in \mathcal{P} ( C ( [ 0 , T ] , E^N ) ) \simeq  \mathcal{P} ( C ( [ 0 , T ] , E )^N )$ be the pathwise law of the particle system with initial law $f^N_0$. In the same way, let $f_{[0,T]} \in \mathcal{P} ( C ( [ 0 , T ] , E ) )$ be the law of the targeted limit nonlinear process with initial law $f_0 \in \mathcal{P} ( E )$. Following~\cite{ben_arous_methode_1990} and~\cite{ben_arous_increasing_1999}, pathwise chaos on $[0,T]$ can be recovered from the above theorem, essentially by taking $\mathscr{E}=C([0,T],E)$. 

\begin{corollary}[Pathwise chaos from LDP] 
Assume that the following properties hold. 
\begin{enumerate}
    \item $f^N_0$ is a Gibbs measure of the form \cref{eq:gibbsmeasure} with respect to $f_0^{\otimes N}$ for a polynomial energy functional $G\in C_b(\pb(E))$.
    \item The functional $\mu \in \mathcal{P} ( E ) \mapsto H ( \mu | f_0 ) - G ( \mu )$ admits a unique minimizer $\mu_{\star}$ which is non-degenerate. 
    \item $f^N_{[0,T]}$ is a Gibbs measure of the form \cref{eq:gibbsmeasure} with respect to $f^{\otimes N}_{[0,T]} $ for a polynomial energy functional $\mathcal{G} \in C_b ( \pb(C ( [ 0 , T ] , E )) )$. 
    \item The functional $\nu \in \mathcal{P} ( C ( [ 0 , T ] , E ) ) \mapsto H ( \nu | f_{[0,T]} ) - \mathcal{G} ( \nu )$ has a unique minimizer $f_{[0,T]}^\star$ which is non-degenerate. Moreover, $f_{[0,T]}^\star$ is the pathwise law of the nonlinear process with initial condition $\mu_\star$. 
\end{enumerate}
Then the sequence $\big( f^N_{[0,T]} \big)_N$ is $f_{[0,T]}^\star$-chaotic. 
\end{corollary}

The two first two assumptions are related to the initial data. The \emph{propagated} property is more the LDP than the chaoticity since the initial measure is assumed to be Gibbsian and no more chaotic as usual. To recover the usual setting, the first assumption has to be replaced by the $f_0$-chaoticity of $f^N_0$, that is to say $G = 0$. The unique minimizer of $H ( \cdot | \mu_0 )$ is thus in this case $\mu_{\star} = f_0$ and $f_{[0,T]}^\star  = f_{[0,T]} $ is the desired law for the limit nonlinear process. The third assumption tells that the Gibbs form of the density is also valid at the pathwise level. For a McKean-Vlasov diffusion with regular coefficients (typically Lipschitz \cite{dai_pra_mckean-vlasov_1996},\cite{malrieu_logarithmic_2001}), the Gibbs density $\frac{\dd  f^N_{[0,T]}}{\dd  f^{\otimes N}_{[0,T]}}$ can typically be computed using Girsanov's formula (see Appendix~\cref{appendix:girsanov}). Thus the remaining difficulty often lies in the fourth point. 

\begin{example}[Application to several models] Checking that the above assumptions hold can be very technical. To give a flavour of the possible applications, we mention here a few examples.
\begin{itemize}
    \item \emph{(McKean-Vlasov system with regular gradient forces and constant diffusion).} The assumptions are exhaustively checked in the original paper \cite{ben_arous_methode_1990}, leading to  the desired pathwise chaos on finite time intervals.
    \item \emph{(McKean-Vlasov system with only continuous drift and H\"older position dependent diffusion).} Pathwise chaos is proved in the seminal work \cite{dawson_large_1987} by establishing a LDP principle and by showing that the limit law is the only minimizer of the related rate function. The method is close to the one which is described above, but it is driven in some abstract dual spaces in order to weaken the regularity conditions: the diffusion can depend on the position with H\"older-regularity (but does not depend on its law), but no such regularity is needed for the non-linear drift.
    \item \emph{(Hamiltonian systems with random medium interactions).} Pathwise chaos is proved in \cite{dai_pra_mckean-vlasov_1996} by extending the above method. The main difficulty in this system comes from the control of the random jumps and from the random medium. 
    \item \emph{(Curie-Weiss and Kuramoto models).} Once again, it is an application of the above method. The Curie-Weiss model is obtained as a corollary in \cite{dawson_large_1987} and \cite{dai_pra_mckean-vlasov_1996}, while the Kuramato model is the last part of \cite{dai_pra_mckean-vlasov_1996}. The method is applied to analogous jump processes with random interactions in \cite{leonard_large_1995-1}. See also Section \cref{sec:kuramoto}. 
\end{itemize}
\end{example}

\subsubsection{Chaos from entropy bounds}\label[I]{sec:chaosfromentropy}

The theorems stated in the previous Section \cref{sec:chaosthroughldp} strongly suggest that the relative entropy (between the $N$-particle distribution and the tensorised limit law) is an important quantity to look at. In fact, Pinsker inequality \cref{eq:pinsker} implies that 
\[{\left\|f^N_t-f^{\otimes N}_t\right\|}^2_\mathrm{TV} \leq 2H {\left( f^N | f^{\otimes N} \right)},\]
so if the right-hand side goes to zero as $N\to+\infty$, it implies propagation of chaos in Total Variation norm. But as it can be expected, it is very demanding to prove that the relative entropy vanishes. The following lemma shows that a simple bound may be sufficient for a slightly weaker result. 

\begin{lemma}[Dimensional bounds on entropy, \cite{csiszar_sanov_1984}]\label[I]{lemma:dimensionalboundsentropy}
Let $\mathscr{E}$ be a measurable space. For every symmetric probability measure $f^N$ on $\mathscr{E}^N$, and every nonnegative integer $k(N) \leq N$, it holds that
\begin{equation}\label[Ieq]{eq:csiszarmarginal} H {\left( f^{k(N),N} \big| f^{\otimes k(N)} \right)} \leq \frac{k(N)}{N}  H {\left( f^N | f^{\otimes N} \right)}. \end{equation}
\end{lemma}

A bound on $H ( f^N | f^{\otimes N} )$ thus implies propagation of chaos in Total Variation norm for blocks of size $k(N)=o(N)$. This technique is by now classical and various applications will be presented in the following sections.
\begin{remark}
Note that a bound on $H ( f^N | f^{\otimes N} )$ implies that the normalised entropy goes to zero as $N\to+\infty$ :
\[\widetilde{H}( f^N | f^{\otimes N} ):=\frac{1}{N}H ( f^N | f^{\otimes N} )\underset{N\to+\infty}{\longrightarrow}0.\]
\end{remark}

A first historical example of entropy bound for Gibbs measures (with a continuous bounded but non necessarily polynomial potential) can be found in the article \cite{ben_arous_increasing_1999} subsequent to \cite{ben_arous_methode_1990}.  

\begin{theorem}[Entropy bound for Gibbs measures, \cite{ben_arous_increasing_1999}] \label[I]{Gibbs_bound} Let $\mu^N$ be a non degenerate Gibbs measure of the form \cref{eq:gibbsmeasure}. Then, with the same notations as in Theorem~\cref{thm:chaosfromLDP}, the measure $\mu^N_{\star} = \sum_{\nu \in \mathcal{P}^{m_0}} \bar{d} ( \nu ) \nu^{\otimes N}$ satisfies the entropy bound
\[ \limsup_{N \rightarrow + \infty} H {\left( \mu^N | \mu^N_{\star} \right)} < + \infty.\]
\end{theorem}

Note that this result strengthens \cite[Theorem B]{ben_arous_increasing_1999} (Theorem \cref{thm:chaosfromLDP}). As before, this readily implies pathwise propagation of chaos.  

\begin{corollary}[Pathwise McKean-Vlasov, $C^2_b$ potentials \cite{ben_arous_increasing_1999}]
If the rate function has a unique zero \cite[Assumption (A1)]{ben_arous_increasing_1999}, the pathwise entropy bound on $\left[ 0,T \right]$ holds for McKean-Vlasov gradient systems with $C^2_b$ coefficients.
\end{corollary}

Theorem \cref{Gibbs_bound} is very strong and quite general. It is mainly intended for true Gibbs measures and in our case, it may look too powerful (besides, the hypothesis may not be easily checked). In the literature, there are more direct approaches to bound the entropy.   

The two following lemmas are valid respectively in the pathwise and the pointwise cases for the McKean-Vlasov diffusion. They do not directly prove a propagation of chaos result but suggest a way to bound the relative entropy  under very weak assumptions on the drift. They are a direct consequence of Girsanov theorem (Appendix \cref{appendix:girsanov}) which provides an explicit expression of the relative entropy as an observable of the particle system. The classical application is a strengthening result: if this observable can be controlled by a weak form of propagation of chaos, then the entropy bound strengthens the weak propagation of chaos result into strong (pathwise) propagation of chaos in TV norm, see for instance Corollary \cref{coro:mckeantvchaos}. Moreover, in some cases where the drift is singular and when no easy propagation of chaos result is available, computing a bound on only one observable is often a quite sensible strategy. This will be reviewed in particular in Section \cref{sec:jabin}.

\begin{lemma}[Pathwise entropy bound]\label[I]{lemma:entropyboundgirsanov}
Let $T>0$ and $I=[0,T]$. Assume that the nonlinear martingale problem associated to the McKean-Vlasov diffusion (Definition \cref{def:nonlinearmeanfieldmartingaleproblem} and \cref{eq:FPgenerator}) with 
\[b:\R^d\times\pb(\R^d)\to\R^d,\quad\sigma = I_d,\]
is wellposed and let $f_I\in \pb(C([0,T],\R^d))$ be its solution. For $N\in\N$, let $f^N_I\in\pb(C([0,T],(\R^{d})^N))$ be the law of the associated particle system $(\mathcal{X}^N_t)^{}_t$. Then, for any $k\leq N$ it holds that
\begin{equation}\label[Ieq]{eq:entropyboundgirsanov}H{\big(f_I^{k,N}|f_I^{\otimes k}\big)}\leq\frac{k}{2}\E\left[\int_0^T \big|b\big({X}^1_t,\mu_{\mathcal{X}^N_t}\big)-b({X}^1_t,f_t)\big|^2\dd t\right].\end{equation}
\end{lemma}

This lemma is a mere application of Girsanov's theorem; for simplicity, the result is stated in the case of a constant diffusion matrix but it is also valid in the case of a diffusion matrix which depends on the positional argument but not on the measure argument (see Remark \cref{rem:generaldiffusionmatrix}). 

\begin{proof}
Since the nonlinear martingale is well-posed, it is well-known (see \cite[Chapter 5, Proposition 4.6]{karatzas_brownian_1998} or \cite[Chapter 5, Proposition 3.1]{ethier_markov_1986}) that we can construct a filtration and $N$ independent adapted $f_I$-Brownian motions $(\overline{B}^i_t)^{}_t$ on the path space such that 
\begin{equation}\label[Ieq]{eq:mckeansdecanonical}\dd\mathsf{X}^i_t = b(\mathsf{X}^i_t,f_t)\dd t + \dd \overline{B}^i_t,\end{equation}
where we recall that $\mathbf{X}^N_t = (\mathsf{X}^1_t,\ldots,\mathsf{X}^N_t)$ is the canonical process on $C([0,T],\R^{dN})\simeq C([0,T],\R^d)^N$. In other word, the canonical process is a weak solution of the nonlinear McKean-Vlasov SDE on the path space $(C([0,T],\R^d),\mathscr{F},f_I)$. For $i\in\{1,\ldots,N\}$ let us define the processes:
\[\Delta^i_t = b\big(\mathsf{X}^i_t,\mu_{\mathbf{X}^N_t}\big)-b(\mathsf{X}^i_t,f_t).\]
and
\[H^N_t := \sum_{i=1}^N {\left[\int_0^t \Delta^i_s\cdot\dd \overline{B}^i_s - \frac{1}{2}\int_0^t |\Delta^i_s|^2\dd s\right]}.\]
It is classical to check that $\exp(H^N)$ defines a $f_I$-martingale (see \cite[Chapter 3, Corollary 5.16]{karatzas_brownian_1998}). Then by Girsanov theorem (see Appendix \cref{appendix:girsanov}) on the product space $(C([0,T],\R^d)^N,\mathscr{F}^{\otimes N},f_I^{\otimes N})$, it is possible to define a probability measure $f_I^{N}$ on $C([0,T],\R^d)^N$ such that for $i\in\{1,\ldots,N\}$, the processes
\begin{equation}\label[Ieq]{eq:browniangirsanov}B^i_t := \overline{B}^i_t - \int_0^t \mathsf{X}^i_s\dd s\end{equation}
are $N$ independent $f^{N}_I$-Brownian motions. Reporting \cref{eq:browniangirsanov} into \cref{eq:mckeansdecanonical} we see that
\[\dd \mathsf{X}^i_t = b\big(\mathsf{X}^i_t,\mu_{\mathbf{X}^N_t}\big)\dd t + \dd B^i_t.\]
In other words, $(\mathbf{X}^N_t)^{}_t$ is a weak solution of the McKean-Vlasov particle system on the path space $(C([0,T],\R^d)^N,\mathscr{F},f^N_I)$ and, as the notation implies, $f^{N}_I$ is the $N$-particle distribution. Moreover, the Girsanov theorem gives a formula for the Radon-Nikodym derivative: 
\[\frac{\dd f^N_I}{\dd f^{\otimes N}_I} = \exp(H^N_T) = \exp{\left(\sum_{i=1}^N \left[\int_0^t \Delta^i_s\cdot\dd \overline{B}^i_s - \frac{1}{2}\int_0^t |\Delta^i_s|^2\dd s\right]\right)}.\]
As an immediate consequence, it is possible to compute the relative entropy as follows: 
\begin{align*}
    H(f_I^{N}|f_I^{\otimes N}) &:= \E_{f_I^{N}}{\left[\log\frac{\dd f^{N}_I}{\dd f^{\otimes N}_I}\right]}\\
    &= \E_{f_I^{N}}{\left[H^N_T\right]}\\
    &= \E_{f_I^{N}}{\left[\sum_{i=1}^N {\left[\int_0^t \Delta^i_t\cdot\dd \overline{B}^i_t - \frac{1}{2}\int_0^T |\Delta^i_t|^2\dd t\right]}\right]}\\
    &= \E_{f_I^{N}}{\left[\sum_{i=1}^N {\left[\int_0^T \Delta^i_t\cdot\dd {B}^i_t + \frac{1}{2}\int_0^T |\Delta^i_t|^2\dd t\right]}\right]}\\
    &= \frac{1}{2}\E_{f_I^{N}}{\left[\sum_{i=1}^N\int_0^T|\Delta^i_t|^2\dd t\right]},
\end{align*}
which, by exchangeability, eventually gives: 
\[H(f_I^{N}|f_I^{\otimes N})=\frac{N}{2}\E_{f_I^{N}}{\left[\int_0^T |b(\mathsf{X}^1_t,\mu_{\mathbf{X}^N_t})-b(\mathsf{X}^1_t,f_t)|^2\dd t\right]}.\]
Coming back to our usual notations on the abstract probability space $(\Omega,\mathscr{F},\mathbb{P})$ on which a particle system $(\mathcal{X}^N_t)^{}_t\sim f^N_I$ is defined, it simply means that 
\[H(f_I^{N}|f_I^{\otimes N})=\frac{N}{2}\E{\left[\int_0^T |b({X}^1_t,\mu_{\mathcal{X}^N_t})-b({X}^1_t,f_t)|^2\dd t\right]}.\]
The conclusion follows from Lemma \cref{lemma:dimensionalboundsentropy}.
\end{proof}

\begin{remark} \label[I]{rem:generaldiffusionmatrix}
The inequality \cref{eq:entropyboundgirsanov} is actually an equality (see the proof of \cite[Theorem 2.6(3)]{lacker_strong_2018}). This relatively direct computation can be seen as a very special case of \cite[Theorem 2.4]{leonard_girsanov_2012}. The result readily extends to the case of time-dependent parameters $b,\sigma$ and to the case of a non constant diffusion matrix $\sigma\equiv \sigma(t,x)$ which does not depend on the measure argument and which is assumed to be invertible everywhere. The only difference in \cref{eq:entropyboundgirsanov} is that $b$ should be replaced by $\sigma^{-1}b$. An even more general setting is the one given in \cite{lacker_strong_2018} where
\[b:[0,T]\times C([0,T],\R^d)\times \pb(C([0,T],\R^d)) \to \R^d, \,\, \sigma : [0,T]\times C([0,T],\R^d) \to \mathcal{M}_d(\R),\]
are assumed to be jointly measurable. This does not affect the final result \cref{eq:entropyboundgirsanov} nor the argument. 
\end{remark}

\begin{remark}
It is worth noticing that this approach does not seem to be restricted to diffusion processes. On the one hand, the full Girsanov theory can be applied to jump processes as well (see \cite{leonard_girsanov_2012} and the references therein) and it is actually a very powerful and general result in the theory of stochastic integration \cite[Section 5.5]{le_gall_brownian_2016}. On the other hand, any model presented in this review can be written as the solution of a very general martingale problem. To the best of our knowledge, an analogous generalised result does not seem to exist in the literature yet. For the Nanbu system, it may be contained in \cite[Theorem 2.11]{leonard_loi_1986}.

\end{remark}

In a pointwise setting, the time derivative of the relative entropy can be directly computed using the generator of the particle system. 

\begin{lemma}[General bound on the time-derivative entropy] \label[I]{lemma:computeH}
Let $f_t\in\pb(\R^d)$ be the solution of \cref{eq:mckeanvlasov-pde} at time $t$ with 
\[b:\R^d\times\pb(\R^d)\to\R^d,\quad \sigma=I_d,\]
and let $f^N_t\in\pb((\R^d)^N)$ be the law of the associated particle system. Then for every $\alpha > 0$ it holds that
\begin{equation}\label[Ieq]{eq:computeH}\frac{\dd}{\dd t} H \big( f^N_t | f_t^{\otimes N} \big) \leq \frac{\alpha - 1}{2} I \big( f^N_t | f_t^{\otimes N} \big) + \frac{N}{2 \alpha} \E \left[ \big| b \big( X^1_t , \mu_{\mathcal{X}^N_t} \big) - b ( X^1_t ,f_t ) \big|^2 \right]. \end{equation} 
\end{lemma}

The following proof is mostly formal as we assume that the limit $f_t$ and $\log f_t$ are regular enough to be taken as test functions. The computations can be fully justified in the cases where the lemma will be applied. 

\begin{proof} Let us recall that the generator of the $N$-particle system is defined by 
\[\mathcal{L}_N \varphi_N(\mathbf{x}^N) = \sum_{i=1}^N L_{\mu_{\mathbf{x}^N}}\diamond_i \varphi_N(\mathbf{x}^N),\]
where, given $\mu\in\pb(E)$,
\[L_\mu\varphi(x) := \langle b(x,\mu),\nabla\varphi\rangle+\frac{1}{2}\Delta\varphi.\]
The Kolmogorov equation for the $N$-particle system reads: 
\[ \frac{\dd}{\dd t} \left\langle f^N_t , \varphi_N \right\rangle = \left\langle f^N_t , \mathcal{L}_N \varphi_N \right\rangle.\]
The Kolmogorov equation associated to a system of $N$ independent $f_t$-distributed particles reads:
\[ \frac{\dd }{\dd t} \left\langle f^{\otimes N}_t , \varphi_N \right\rangle = \left\langle f^{\otimes N}_t , L^{\diamond N}_{f_t} \varphi_N\right\rangle, \]
where we define the generator 
\[L^{\diamond N}_{f_t} \varphi_N := \sum_{i=1}^N L_{f_t}\diamond_i \varphi_N.\]
The relative entropy is defined by:
\[ H \big( f^N_t | f_t^{\otimes N} \big) = \E_{\mathcal{X}^N_t \sim f^N_t} \left[ \log \frac{\dd  f^N_t}{\dd  f_t^{\otimes N}} \left( \mathcal{X}^N_t \right) \right] \equiv \left\langle f^N_t , \log \frac{f^N_t}{f_t^{\otimes N}} \right\rangle. \]
In the last term, $\frac{\dd  f^N_t}{\dd  f_t^{\otimes N}}$ has been replaced by $\frac{f^N_t}{f_t^{\otimes N}}$, which makes sense since $f^N_t$ and $f_t$ are probability density functions. Using the product derivation rule:
\[ \frac{\dd }{\dd t} H \big( f^N_t | f_t^{\otimes N} \big) = \left\langle f^N_t , \mathcal{L}_N \left( \log \frac{f^N_t}{f_t^{\otimes N}} \right) \right\rangle + \left\langle f^N_t , \frac{\dd }{\dd t} \log \frac{f^N_t}{f_t^{\otimes N}} \right\rangle \]
The last term can be written
\[ \left\langle f^N_t , \frac{\dd }{\dd t} \log \frac{f^N_t}{f_t^{\otimes N}} \right\rangle = \left\langle f^N_t , \frac{f^{\otimes N}_t}{f_t^N} \frac{\dd }{\dd t} \frac{f^N_t}{f_t^{\otimes N}} \right\rangle = \left\langle f_t^{\otimes N} , \frac{\dd }{\dd t} \frac{f^N_t}{f_t^{\otimes N}} \right\rangle. \]
The mass conservation for $f_t^N$ gives
\[ \left\langle f_t^{\otimes N} , \frac{f^N_t}{f_t^{\otimes N}} \right\rangle = \left\langle f_t^N , 1 \right\rangle = 1, \]
and therefore
\[ \left\langle f_t^{\otimes N} , \frac{\dd }{\dd t} \frac{f^N_t}{f_t^{\otimes N}} \right\rangle = - \left\langle f_t^{\otimes N} , L^{\diamond N}_{f_t} \frac{f^N_t}{f_t^{\otimes N}} \right\rangle.\]
Using the definition of the generator, we get:
\[L^{\diamond N}_{f_t} \frac{f^N_t}{f_t^{\otimes N}}(\mathbf{x}^N) = \sum_{i=1}^N \left\langle b(x^i,f_t),\frac{f^N_t}{f_t^{\otimes N}}(\mathbf{x}^N)\nabla_{x^i} \log\frac{f^N_t}{f_t^{\otimes N}}(\mathbf{x}^N)\right\rangle+\frac{1}{2}\Delta\left(\frac{f^N_t}{f_t^{\otimes N}}\right).\]
And thus it holds that 
\begin{align}\frac{\dd }{\dd t} H \big( f^N_t | f_t^{\otimes N} \big) &= \E_{f^N_t}\left[\sum_{i=1}^N \left\langle b(X^i_t,\mu_{\mathcal{X}^N_t})-b(X^i_t,f_t), \nabla_{x^i}\log\frac{f^N_t}{f_t^{\otimes N}}(\mathcal{X}^N_t)\right\rangle\right] \nonumber\\
&\qquad\qquad+ \frac{1}{2}\left\langle f^N_t, \Delta\left( \log\frac{f^N_t}{f_t^{\otimes N}}\right)\right\rangle-\frac{1}{2}\left\langle f^{\otimes N}_t,\Delta\left(\frac{f^N_t}{f_t^{\otimes N}}\right)\right\rangle\nonumber\\
&= \E_{f^N_t}\left[\sum_{i=1}^N \left\langle b(X^i_t,\mu_{\mathcal{X}^N_t})-b(X^i_t,f_t), \nabla_{x^i}\log\frac{f^N_t}{f_t^{\otimes N}}(\mathcal{X}^N_t)\right\rangle\right] \nonumber\\
&\qquad\qquad- \frac{1}{2}\left\langle f^{\otimes N}_t,\left|\nabla \frac{f^N_t}{f_t^{\otimes N}}\right|^2\right\rangle.\label[Ieq]{eq:laststepcomputeH}
\end{align}
and the last term involves
\[\left\langle f^{\otimes N}_t,\left|\nabla \frac{f^N_t}{f_t^{\otimes N}}\right|^2\right\rangle = \left\langle f^N_t, \left|\nabla \log\frac{f^N_t}{f_t^{\otimes N}} \right|^2\right\rangle=:I(f^N_t|f^{\otimes N}_t).\]
Therefore Cauchy-Schwarz inequality and Young inequality give for any $\alpha>0$
\[\frac{\dd}{\dd t} H \big( f^N_t | f_t^{\otimes N} \big) \leq \frac{\alpha - 1}{2} I \big( f^N_t | f_t^{\otimes N} \big) + \frac{1}{2 \alpha}\sum_{i=1}^N \E \left[ \big| b \big( X^i_t, \mu_{\mathcal{X}^N_t} \big) - b ( X^i_t ,f_t ) \big|^2 \right].\]
The conclusion follows since particles are exchangeable. 
\end{proof}

\begin{remark} Several points should be noticed :
\begin{itemize}
    \item It is possible to take $\alpha = 1$ in order to get rid of the Fisher information as in~\cite{guillin_uniform_2020}, but a further control on $ W_2 {\left( \mu_{\mathcal{X}^N_t} ,f_t \right)}$ is then needed, see \cite{guillin_uniform_2020}.
    \item Before the splitting which introduces $\alpha$, Cauchy-Schwarz's inequality would have lead to a bound close to the $HWI$ inequality in our special case.
\end{itemize}
\end{remark}

To end this section, we would like to emphasize the fact that these results do not require any particular regularity on the drift (this is a well-known but remarkable property of Girsanov's transform). As a general rule, entropy related methods are well suited to handle cases with singular interactions. An example with exceptionally weak regularity assumptions will be given in Section \cref{sec:jabin}. Another example with a complex abstract interaction mechanism will be presented in Section \cref{sec:chaosviagirsanov}.

\subsubsection{Tools for concentration inequalities} \label[I]{concentrationI}

Large Deviation principles imply propagation of chaos, but they do not always give a way to quantify it since their result is often purely asymptotic (for instance, Sanov theorem is non-quantitative). In this section, we gather some results which quantify the deviation of an empirical measure of $N$ samples around its mean. These results are valid for any fixed (sufficiently large) $N$. We first state a classical concentration inequality, obtained as the consequence of a Log-Sobolev inequality.
This inequality and other related functional inequalities are deep structural properties of the system which will also be used to study ergodic properties and long-time propagation of chaos (Section \cref{sec:gradientsystems}). Then we state quantitative versions of Sanov theorem which strengthen the above concentration inequality. In this section, $E$ is assumed to be $\R^d$ (or at least a smooth finite dimensional manifold).

The following Log-Sobolev inequality is another kind of entropy bound.

\begin{definition}[Log-Sobolev Inequality]
For $\lambda >0$, a probability measure $\mu$ with finite second moment satisfies a Logarithmic-Sobolev Inequality $LSI( \lambda )$ when for all $\nu$ in $\mathcal{P} ( E )$
\[ H ( \nu | \mu ) \leq \frac{1}{2\lambda} I ( \nu | \mu ), \]
where $I$ is the Fisher information (Definition \cref{def:entropyfisher}).
\end{definition}

An important consequence is the following lemma. 

\begin{lemma}[Concentration, see Ledoux \cite{ledoux_concentration_1999}] \label[I]{concledoux}
If a probability measure $\mu$ satisfies a $LSI( \lambda )$ then for any Lipschitz test function $\varphi$ with Lipschitz constant bounded by~1 and for any $\varepsilon > 0$, it holds that
\begin{equation}\label[Ieq]{eq:concentrationledoux} \mathbb{P}_{X \sim \mu} \big( \big| \varphi ( X ) - \E [ \varphi ( X ) ] \big| \geq \varepsilon \big) \leq 2 \e^{-\frac{ \lambda \varepsilon^2}{2}}. \end{equation}
\end{lemma}
This lemma is typically applied in $E^N$ for $\mathcal{X}^N_t\sim f^N_t$ with the function 
\[\overline{\varphi}\big(\mathcal{X}^N_t\big)=\frac{1}{N}\sum_{i=1}^N\varphi(X^i_t),\]
where $\varphi$ is 1-Lipschitz on $E$. The Lipschitz norm of the function $\overline{\varphi}$ is bounded by $1/\sqrt{N}$.

A classical result \cite[Theorem 1]{otto_generalization_2000} shows that under mild assumptions, the Log-Sobolev inequality also implies the following Talagrand inequality. 

\begin{definition}[Talagrand Inequalities]
For any real $p \geq 1$ and $\lambda > 0$, a probability measure $\mu$ with finite $p$-th moment satisfies a Talagrand inequality $T_p( \lambda )$ when for all $\nu\in\mathcal{P} ( E )$,
\begin{equation}\label[Ieq]{eq:talagrand} W_p \left( \nu , \mu \right) \leq \sqrt{\frac{2 H \left( \nu | \mu \right)}{\lambda}}. \end{equation}
\end{definition}

This inequality is all the more strong as $\lambda$ and $p$ are big (it is a consequence of Jensen's inequality). It is known that $T_2$ implies some Poincar\'e inequality and a handful characterization is available for $T_1$ inequalities, see \cite{bolley_weighted_2005}. Talagrand inequalities are also useful to quantify ergodicity with respect to the Wassertein distance. Here is another characterization.
\begin{lemma}[Square-exponential moment]
A probability measure $\mu$ with finite expectation satisfies a $T_1$ inequality if and only if there exist $\alpha > 0$ and $x\in E$ such that $\int_E \e^{\alpha |x-y|^2} \mu( \dd  y) < + \infty $. 
\end{lemma}

Note that Talagrand inequalities allow to recover usual Wassertein convergence (and then convergence in law) from entropic convergence. Concentration inequalities can also stem from Talagrand inequalities, although the stronger Logarithmic Sobolev inequality is more often used in this context. The Logarithmic Sobolev inequality can be established thanks to the following criterion. 

\begin{proposition}[Bakry, Emery \cite{azema_diffusions_1985,bakry_hypercontractivite_1994}]\label[I]{prop:bakryemery}
Consider a diffusion process with semi-group $( P_t)_{t \geq 0}$, generator $\mathcal{L}$ and carr\'e du champ operator $\Gamma : (\varphi,\psi) \mapsto \frac{1}{2} [ \mathcal{L} ( \varphi \psi ) - \varphi \mathcal{L} (\psi) - \psi \mathcal{L} ( \varphi ) ]$. Assume that there exists a real $\lambda$ such that for every regular function $\varphi$
\[\Gamma_2 (\varphi) \geq \lambda \Gamma (\varphi)\]
where $\Gamma_2 (\varphi) = \frac{1}{2} [ \mathcal{L} ( \Gamma ( \varphi ) ) - 2 \Gamma ( \varphi \mathcal{L} (\varphi) ) ]$. Then the following properties hold. 
\begin{itemize}
    \item For all $x\in\R^d$ and all $t>0$, if $P_0=\delta_x$, then  $P_t$ satisfies $LSI {\left( \frac{\lambda}{1 - \e^{-\lambda t}} \right)}$ (where the semi-group $P_t$ is identified with the transition probability measure that it generates). 
    \item If $\lambda > 0$, the semi-group is ergodic and $P_t$ converges towards the invariant measure $\mu$ with rate \[H(P_t|\mu)\leq Ce^{-2 \lambda t}.\]
\end{itemize}
\end{proposition}

In the previous concentration result \cref{eq:concentrationledoux}, the test function is fixed before computing the probability. One could take the supremum over all 1-Lipschitz test functions, this would correspond to a weak chaos for the $D_1$ distance. To get stronger estimates on the Wasserstein distance between $\mu_{\mathcal{X}^N_t}$ and $f_t$, the supremum needs to come inside the probability. This is not an easy task, it requires a quantitative version of Sanov theorem, which is proved in \cite[Theorem 2.1]{bolley_quantitative_2006}.

\begin{theorem}[Pointwise Quantitative Sanov \cite{bolley_quantitative_2006}] \label[I]{thm:PQS}
Consider a probability measure $\mu$ on $\R^d$ which satisfies $T_p(\lambda)$ for $p\in[1,2]$ and $\lambda>0$ and which has a bounded square-exponential moment. Let $\mathcal{X}^N\sim\mu^{\otimes N}$ be a system of $N$ i.i.d  $\mu$-distributed random variables. Then for any $\lambda' <\lambda$ and $\varepsilon >0$, there exists a constant $N_{\varepsilon}$ (which depends also on $d$ and the square-exponential moment of $\mu$) such that for all $N\geq N_\varepsilon$,
\[\mathbb{P} {\left( W_p {\left( \mu_{\mathcal{X}^N} , \mu \right)} > \varepsilon \right)} \leq  \e^{-\gamma_p\frac{\lambda'}{2} N \varepsilon^2},\]
where $\gamma_p>0$ is an explicit constant which depends only on $p$.
\end{theorem}

The following pathwise generalization is proved in \cite[Theorem 1]{bolley_quantitative_2010}.

\begin{theorem}[Pathwise Quantitative Sanov and Pathwise chaos \cite{bolley_quantitative_2010}]
Under the same assumptions, the above theorem holds for a measure $\mu$ on the H\"older space $C^{0,\alpha} ( [0,T] , \R^d )$ with $\alpha\in(0,1]$.
\end{theorem}

Except the last one, the results in this section are stated in a static framework. Examples of time dependent systems and applications to propagation of chaos are detailed in Section \cref{sec:concentrationineqgradient}.

\subsection{Tools for Boltzmann interactions}\label[I]{sec:boltzmanntools}

\subsubsection{Series expansions}

Let us consider first the homogeneous Boltzmann system with $L^{(1)}=0$. When the collision rate $\lambda$ satisfies the uniform bound \cref{eq:uniformboundlambda}, then it can be directly checked that the operator $\mathcal{L}_N$ is continuous for the $L^\infty$ norm: 
\[\forall \varphi_N\in C_b(E^N),\quad \|\mathcal{L}_N\varphi_N\|_\infty\leq \Lambda(N-1)\|\varphi_N\|_\infty.\]
Without loss of generality (see Proposition \cref{prop:acceptreject}), we will assume here that $\lambda\equiv\Lambda$ is constant. As a consequence, the exponential series $\e^{t\mathcal{L}_N}$ is absolutely convergent for $t<1/(\Lambda(N-1))$ and there is a semi explicit formula for an observable $\varphi_N\in C_b(E^N)$ at any time $t\geq0$: 
\[\E{\left[\varphi_N(\mathcal{Z}^N_t)\right]}=\sum_{k=0}^{+\infty}\frac{t^k}{k!}\langle f^N_0, \mathcal{L}_N^k\varphi_N\rangle,\]
where $(\mathcal{Z}^N_t)^{}_t$ is the particle process with initial law $f^N_0\in \pb(E^N)$. Then, considering a test function $\varphi_N\equiv\varphi_s\otimes 1^{\otimes(N-s)}$ which depends only on $s$ variables for a fixed $s\in\N$, the term on the right-hand side depends only on $N$ through known quantities, namely the initial law $f^N_0$ and the operator $\mathcal{L}_N$. The initial law $f^N_0$ is assumed to be $f_0$-chaotic so there is an asymptotic control of all its marginals when $N\to+\infty$.
In its seminal article \cite{kac_foundations_1956}, Kac managed to pass to the limit directly in the series on the right-hand side using a dominated convergence theorem argument. This necessitates in particular to prove the absolute convergence of the series on a time interval independent of $N$ when $s$ is fixed. The argument has been generalised in \cite{carlen_kinetic_2013} and will be thoroughly discussed in Section \cref{sec:kactheorem}. As a byproduct it will show the existence of a solution of the Boltzmann equation in the form of an explicit series expansion. The final formula \cref{eq:convergeEvarphisZN} will be a direct extension of the ``exponential formula'' obtained by McKean in \cite{mckean_exponential_1967} for the solution of a simpler Boltzmann model in $E=\{-1,1\}$ (the famous 2-state Maxwellian gas). In a famous work, Wild~\cite{wild_boltzmanns_1951} showed that the solution of the Boltzmann equation for cut-off Maxwellian molecules has a semi-explicit representation in the form of an infinite sum (see also \cite[Chapter 4, Section 1]{villani_review_2002}, \cite{carlen_central_2000} and the references therein). McKean showed that for the 2-state Maxwellian gas, the Formula \cref{eq:convergeEvarphisZN} obtained by propagation of chaos can be interpreted as the dual version of a Wild sum. 

This approach is only focused on the evolution of observables oof the form $\langle f^N_t, \varphi_N\rangle$ for a fixed $\varphi_N$ and does not study directly the evolution of the $N$-particle law $f^N_t$. The evolution of $f^N_t$ is given by the forward Kolmogorov equation. This dual point of view on Kac's theorem is studied in \cite{pulvirenti_kinetic_1996} and will also be reviewed in Section \cref{sec:kactheorem}. The starting point is the BBGKY hierarchy: 
\[\partial_t f^{s,N}_t = \frac{s}{N}\mathcal{L}^s f^{s,N}_t + \frac{N-s}{N}\mathcal{C}_{s,s+1}f_t^{s+1,N}, \]
where $\mathcal{C}_{s,s+1}:\pb(E^{s+1})\to \pb(E^s)$ is defined for $f^{(s+1)}\in \pb(E^{s+1})$ and $\varphi_s\in C_b(E^s)$ by
\[\big\langle \mathcal{C}_{s,s+1}f^{(s+1)},\varphi_s\big\rangle := \sum_{i=1}^s \int_{E^{s+1}} L^{(2)}\diamond_{i,s+1}[\varphi_s\otimes 1](\mathbf{z}^{s+1})f^{(s+1)}(\dd\mathbf{z}^{s+1}).\]
Let $\mathbf{T}_N^{(s)}(t)=\exp(t\frac{s}{N}\mathcal{L}^s)$ denote the semi-group acting on $\pb(E^s)$ generated by $\frac{s}{N}\mathcal{L}^s$. Then, interpreting the last term on the right-hand side as a perturbation of a linear differential equation, Duhamel's formula reads: 
\[f^{s,N}_t = \mathbf{T}_N^{(s)}(t)f^{s,N}_0 + \frac{N-s}{N}\int_0^t \mathbf{T}_N^{(s)}(t-\tau)\mathcal{C}_{s,s+1}f^{s,N}_\tau\dd \tau.\] 
Iterating this formula gives a semi-explicit series expansion in terms of the initial condition: 
\begin{multline*}f_{t}^{s,N} = \sum_{k=0}^{+\infty} \alpha_N^{(s,k)} \int_0^t\int_0^{t_1}\ldots\int_0^{t_{k-1}}\mathbf{T}_N^{(s)}(t-t_1)\mathcal{C}_{s,s+1}\mathbf{T}_N^{(s+1)}(t_1-t_2)\mathcal{C}_{s+1,s+2}\ldots\\
\mathcal{C}_{s+k-1,s+k}\mathbf{T}_N^{(s+k)}(t_k)f_0^{s+k,N}\dd t_1\ldots\dd t_k,\end{multline*}
where $\alpha^{(s,k)}_N=(N-s)\ldots(N-s-k+1)/N^k$ if $s+k\leq N$ and $\alpha^{(s,k)}_N=0$ otherwise. Taking the limit $N\to+\infty$ in this series is the dual viewpoint of the previous approach. If the limit exists, propagation of chaos holds whenever the result can be identified as the infinite hierarchy of tensorised laws $f^{\otimes s}_t$ where $f_t$ solves the Boltzmann equation. Note that this approach can be easily extended to inhomogeneous Boltzmann systems where $L^{(1)}\ne 0$, in which case the semi-group $\mathbf{T}_N^{(s)}$ should be replaced by the semi-group generated by $\sum_{i=1}^s L^{(1)\star}\diamond_i+\frac{s}{N}\mathcal{L}^s$. A famous example is given by Lanford's theorem (see Section \cref{sec:lanford}). The study of the evolution of observables is however more natural for abstract systems when there is no known explicit formula for the dual operators $\mathcal{L}^s$ and $\mathcal{C}_{s,s+1}$ acting on the particle probability distributions.

\subsubsection{Interaction graph}\label[I]{sec:interactiongraph}

In an abstract Boltzmann model given by the generator~\cref{eq:boltzmanngenerator} in Section~\cref{sec:boltzmann}, the binary interactions can be represented by graph structures. Given a trajectorial realisation of the particle system, the \emph{interaction graph} of a particle (or a group of particles) is built backward in time and retain the genealogical interactions which determine the particle at the current time (i.e. the history of the collisions). Before building graphs from particle realisations, the minimal structure of such a possible graph is detailed in the following definition.

\begin{definition}[Interaction graph]
Consider an index $i\in\{1,\ldots,N\}$ (it will stand later for the index of a particle). An interaction graph for $i$ at time $t>0$ is the data of
\begin{enumerate}
    \item a $k$-tuple $\mathcal{T}_k=(t_1,\ldots, t_k)$ of \emph{interaction times} $t>t_1>t_2>\ldots>t_k>0$,
    \item a $k$-tuple $\mathcal{R}_k=(r_1,\ldots,r_k)$ of pairs of indexes, where for $\ell\in\{1,\ldots,k\}$, the pair denoted by $r_\ell=(i_\ell,j_\ell)$ is such that $j_\ell\in\{i_0,i_1,\ldots,i_{\ell-1}\}$ with the convention $i_0=i$ and $i_\ell\in\{1,\ldots,N\}$.
\end{enumerate}
Such an interaction graph is denoted by $\mathcal{G}_i(\mathcal{T}_k,\mathcal{R}_k)$.
\end{definition}

Given a trajectorial realisation of a Boltzmann particle system, the interaction graph of the particle $i$ retains the minimal information needed to compute the state of particle~$i$ at time $t>0$. It is constructed as follows.
\begin{itemize}
    \item The set $(i_1,\ldots,i_k)$ is the set of indexes of the particles which interacted directly or indirectly with particle $i$ during the time interval $(0,t)$ (an indirect interaction means that the particle has interacted with another particle which interacted directly or indirectly with particle $i$) -- note that the $i_\ell$'s may not be all distinct.
    \item The times $(t_1,\ldots,t_k)$ are the times at which an interaction occurred.
    \item For $\ell\in\{1,\ldots,k\}$, the indexes $(i_\ell,j_\ell)$ are the indexes of the two particles which interacted together at time $t_\ell$. 
\end{itemize}

\begin{figure}
    \centering
    \begin{tikzpicture}[>=latex]
    \draw[->] (0,0)--(0,5.5);
    \draw[->] (0,0)--(4.5,0);
    
    \coordinate (t) at (0,5);
    \coordinate (t1) at (0,4);
    \coordinate (t2) at (0,2.5);
    \coordinate (t3) at (0,1.7);
    \coordinate (t4) at (0,0.6);
    
    \node[anchor=east] at (t) {$t$};
    \node[anchor=east] at (t1) {$t_1$};
    \node[anchor=east] at (t2) {$t_2$};
    \node[anchor=east] at (t3) {$t_3$};
    \node[anchor=east] at (t4) {$t_4$};
    
    \coordinate (i) at (1,5);
    \coordinate (j1) at (1,4);
    \coordinate (i1) at (2,4);
    \coordinate (j2) at (1,2.5);
    \coordinate (i2) at (3,2.5);
    \coordinate (j3) at (2,1.7);
    \coordinate (i3) at (3,1.7);
    \coordinate (j4) at (2,0.6);
    \coordinate (i4) at (4,0.6);
    
    \draw[-,color=red] (j3) -- (i3) ;
    
    \node[anchor=north] at (1,0) {$i$};
    \node[anchor=north] at (2,0) {$i_1$};
    \node[anchor=north] at (3,0) {$i_2$};
    \node[anchor=north] at (4,0) {$i_4$};

    \filldraw[black] (j1) circle (2pt);
    \filldraw[black] (i1) circle (2pt);
    \filldraw[black] (j2) circle (2pt);
    \filldraw[black] (i2) circle (2pt);
    \filldraw[black] (j3) circle (2pt);
    \filldraw[black] (i3) circle (2pt);
    \filldraw[black] (j4) circle (2pt);
    \filldraw[black] (i4) circle (2pt);
    
    \draw[-] (1,0) -- (i) ;
    \draw[-] (2,0) -- (i1) ;
    \draw[-] (3,0) -- (i2) ;
    \draw[-] (4,0) -- (i4) ;
    
    \draw[-] (j1) -- (i1) ; 
    \draw[-] (j2) -- (i2) ; 
    \draw[-] (j4) -- (i4) ;
    
    \draw[-,dotted] (t) -- (i)  ;
    \draw[-,dotted] (t1) -- (j1)  ;
    \draw[-,dotted] (t2) -- (j2)  ;
    \draw[-,dotted] (t3) -- (j3)  ;
    \draw[-,dotted] (t4) -- (j4)  ;

    \end{tikzpicture}
    \caption{An interaction graph. The vertical axis represents time. Each particle is represented by a vertical line parallel to the time axis. The index of a given particle is written on the horizontal axis. The construction is done backward in time starting from time $t$ where only particle $i$ is present. At each time $t_\ell$, if $i_\ell$ does not already belong to the graph, it is added on the right (with a vertical line which starts at $t_\ell$). The couple $r_\ell=(i_\ell,j_\ell)$ of interacting particles at time $t_\ell$ is depicted by an horizontal line joining two big black dots on the vertical line representing the particles $i_\ell$ and $j_\ell$. for instance, on the depicted graph, $r_2=(i_2,i)$. Note that at time $t_3$, $r_3=(i_1,i_2)$ (or indifferently $r_3=(i_2,i_1)$) where $i_1$ and $i_2$ were already in the system. Index $i_3$ is skipped and at time $t_4$, the route is $r_4=(i_4,i_1)$. The \emph{recollision} occurring at time $t_3$ is depicted in red.}
    \label[I]{fig:interactiongraph}
\end{figure}
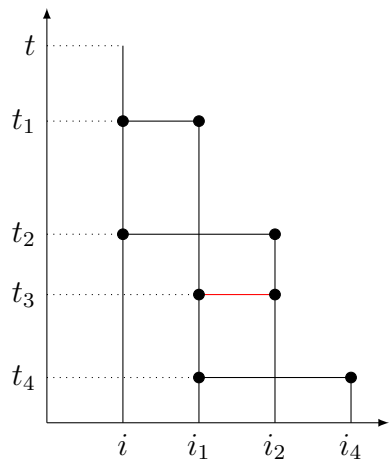

Following the terminology of \cite{graham_stochastic_1997}, a \emph{route} of size $q$ between $i$ and $j$ is the union of $q$ elements $r_{\ell_k}=(i_{\ell_k},j_{\ell_k})$, $k=1,\ldots,q$ such that $i_{\ell_1}=i$, $i_{\ell_{k+1}}=j_{\ell_k}$ and $j_{\ell_q}=j$. A route of size 1 (\emph{i.e} a single element $r_\ell$) is simply called a route. A route which involves two indexes which were already in the graph before the interaction time (backward in time) is called a \emph{recollision}. This construction is more easily understood with the graphical representation of an interaction graph shown on Figure \cref{fig:interactiongraph}. 

The definition of interaction graphs can be extended straightforwardly starting from a group of particles instead of only one particle. This representation does not take into account the physical trajectories of the particles, it only retains the history of the interactions among a group of particles. Note that the graph is not a tree in general since the $i_\ell$'s are not necessarily distinct. It is a tree when no recollision occurs. 

Interaction graphs are a classical tool in the study of Boltzmann particle systems. As we shall see in Section \cref{sec:lanford} they are particularly useful to give a physical interpretation of the series expansions discussed in the previous section. The connections between interaction graphs and series expansions are more thoroughly discussed in~\cite{mckean_exponential_1967} and \cite{carlen_central_2000}. 

In a more probabilistic setting, the following definition extends the construction of an interaction graph to the case of random parameters.

\begin{definition}[Random interaction graph]
Let $\Lambda>0$, $N\in\N$, $i\in\{1,\ldots,N\}$ and $t>0$. Let $(T^{m,\ell})_{1\leq k<\ell\leq N}$ be $N(N-1)/2$ independent Poisson processes with rate $\Lambda/N$. For each Poisson process $T^{m,\ell}$ we denote by $(T^{m,\ell}_n)^{}_n$ its associated increasing sequence of jump times. The sets of times $\mathcal{T}_k=(t_1,\ldots,t_k)$ and routes $\mathcal{R}_k=(r_1,\ldots,r_k)$ are defined recursively as follows. Initially, $t_0=t$ and $i_0=i$ and for $k\geq0$,
\begin{equation}\label[Ieq]{eq:recursivejumpingtimes}t_{k+1} = \max_{\ell,p,n}\big\{T^{i_\ell,p}_n\,|\, T^{i_\ell,p}_n<t_k,\,\ell\leq k\big\}.\end{equation}
Then, given $(\ell,p,n)$ such that $t_{k+1}=T^{i_\ell,p}_n$, $i_{k+1}=p$ and $j_{k+1}=i_\ell$ so that $r_{k+1}=~(i_{k+1},j_{k+1})$. The procedure is stopped once the set on the right-hand side of \cref{eq:recursivejumpingtimes} is empty (it happens almost surely after a finite number of iterations). The resulting interaction graph $\mathcal{G}_i(\mathcal{T}_k,\mathcal{R}_k)$ is called the random interaction graph with rate $\Lambda$ rooted on $i$ at time $t$. The definition is extended similarly starting from a finite number of indexes $(i_0,i_1,\ldots,i_k)$ instead of just $i$.  
\end{definition}

As explained before, a realisation of a Boltzmann particle system immediately gives an interaction graph for each particle. More importantly, given an interaction graph, it is possible to construct a forward realisation of a Boltzmann particle. More precisely, when the interaction graph is sampled as a random interaction graph following, then the following straightforward lemma constructs a forward realisation of a stochastic process whose pathwise law is equal to $f^{1,N}_{[0,t]}$, the first marginal of the law $f^N_{[0,t]}$ of a Boltzmann particle system given by the generator~\cref{eq:boltzmanngenerator} on the time interval $[0,t]$. 


\begin{lemma}
Let us consider the Boltzmann setting given by Assumption \cref{assum:L2} together with the uniform bound \cref{eq:uniformboundlambda} on $\lambda$. Given a realisation of a random interaction graph sampled beforehand as in the previous definition, apply the following procedure:
\begin{enumerate}
    \item At time $t=0$, let the particles $Z^{i_\ell}_0$ be distributed according to the initial law.
    \item Between two collision times, the particles evolve according to $L^{(1)}$. 
    \item At a collision time $t_\ell$, with probability $\lambda(Z^{i_\ell}_{t_\ell^-},Z^{j_\ell}_{t_\ell^-})/\Lambda$, the new states of particles $i_\ell$ and $j_\ell$ are sampled according to \[{\left(Z^{i_\ell}_{t_\ell^+},Z^{j_\ell}_{t_\ell^+}\right)}\sim\Gamma^{(2)}{\left(Z^{i_\ell}_{t_\ell^-},Z^{j_\ell}_{t_\ell^-},\dd z_1, \dd z_2\right)}.\]
\end{enumerate}
Then the process $(Z^{i}_s)_{s\leq t}$ is distributed according to the one-particle marginal $f^{1,N}_{[0,t]}$ of the law $f^N_{[0,t]}$ of a Boltzmann particle system given by the generator \cref{eq:boltzmanngenerator} on the time interval $[0,t]$.
\end{lemma}

This trajectorial construction of a particle (sub-)system is used in \cite{graham_stochastic_1997,meleard_asymptotic_1996}. Using purely combinatorial arguments the authors prove a pathwise version of Kac's theorem with an explicit optimal convergence rate in total variation norm. The main idea is that it is possible to compute the probability of sampling a \emph{bad} interaction graph, that is a graph which would give a system of particles with too much correlation. When $N\to+\infty$, this probability goes to zero. As a consequence, when $N$ is large, with high probability, the particle system is close to a system of independent particles which are shown to be distributed according to the solution of the Boltzmann equation. This will be reviewed in Section \cref{sec:pathwisekactheorem}. As in the proof of Lanford's theorem (Section \cref{sec:lanford}), a fundamental idea is to reduce the problem to the estimation of the number of recollisions in a sampled graph or to the number of graphs in which two given particles are linked by a route of arbitrary size. Indeed, in the probabilistic setting, if a random binary \emph{tree} with branching rate $\Lambda$ is sampled first and then a particle system is constructed as above but starting from independent particles, then this gives a trajectorial representation of a process whose law is the solution of the Boltzmann equation \cref{eq:Boltzmannequationgeneral}. 

\section*{Acknowledgments} The authors wish to thank Pierre Degond for his precious advice and careful proofreading of this manuscript. The authors also thank Paul Thevenin for fruitful comments and discussions.

\appendix

\section{Probability reminders}\label[I]{appendix:probabilityreminders}

For the convenience reader, classical elements of stochastic analysis and probability theory which are used throughout this review are gathered in this section. The Sections \cref{appendix:convergenceproba}, \cref{appendix:skorokhod} and \cref{appendix:martingales} are devoted to the classical construction of stochastic processes on the Skorokhod space. Notions related to the theory of Markov processes and their links with linear and nonlinear PDEs can be found in Section \cref{appendix:timeinhomoegeneousmarkov}. Section \cref{appendix:largedeviations}, Section \cref{appendix:girsanov} and Section \cref{appendix:poisson} summarize probability results regarding respectively large deviations, the Girsanov theorem and Poisson random measures. 

\subsection{Convergence of probability measures}\label[I]{appendix:convergenceproba} 

The two following classical theorems complete the results of Section \cref{sec:distancesproba} and Section \cref{sec:convergenceproba} to study the limit of sequences of probability measures. 

The first theorem links weak convergence and almost sure convergence of random variables. 

\begin{theorem}[Skorokhod's representation theorem] Let $(f_n)_{n \in \N}$ be a sequence of probability measures on a Polish space $E$ which converges weakly towards $f \in \pb(E)$ as $n\to+\infty$. Then there exist a probability space $(\Omega,\mathscr{F},\mathbb{P})$ and some $E$-valued random variables $X,X_n$ defined on this space for all $n \in\N$, such that
\[ \mathrm{Law}(X_n) = f_n, \quad \mathrm{Law}(X) = f, \quad X_n(\omega) \underset{n\to+\infty}{\longrightarrow} X(\omega),\,\,\, \mathbb{P}\text{-a.s.} \]
\end{theorem}

The second theorem is the widely used Prokhorov's theorem which gives a helpful characterization of compactness for the weak convergence topology. The following results can be found for instance in \cite[Section 5]{billingsley_convergence_1999}. Compactness is linked to the notion of tightness defined below. 

\begin{definition}[Tightness] A family $(f_i)_{i \in I}$ of probability measures on a separable metric space $E$ (endowed with its Borel $\sigma$-field) is said to be \emph{tight} when for every $\varepsilon > 0$, there exists a compact set $K_{\varepsilon} \subset E$ such that
\[ \forall i \in I, \quad f_i(K_{\varepsilon}) > 1-\varepsilon. \]
A sequence of random variables is said to be tight when the sequence of their laws is tight. 
\end{definition}

\begin{theorem}[Prokhorov's Theorem] A tight sequence $(f_n)_{n \in \N}$ of probability measures on $E$ is weakly relatively compact. Conversely, if $E$ is also complete, any weakly relatively compact family $(f_n)_{n \in \N}$ is tight. 
\end{theorem}

\subsection{Skorokhod's topology and tightness on the Skorokhod space}\label[I]{appendix:skorokhod}

A stochastic process is a random function from a time interval to a state space $(E,\rho)$ assumed to be Polish. Throughout this article, the stochastic process are assumed to belong (at least) to the Skorokhod space of c\`adl\`ag functions. 

\begin{definition}[c\`adl\`ag] Let $T$ in $(0,+ \infty]$. A function $x: [0,T] \to E$ is said to belongs to the Skorokhod space $D([0,T],E)$ of c\`adl\`ag functions when $x$ is right-continuous and has a left-limit at any time $t\in[0,T]$:
\[ x(t^-) := \lim_{ \substack{s \to t \\ s < t} } x(s) \text{ exists,} \quad x(t) = x(t^+). \]
We recall that a c\`adl\`ag function admits an at most countable number of discontinuities.
\end{definition}

The law of a stochastic process is therefore an element of $\pb(D([0,T],E))$. In order to characterize the compact sets of this space, it is first necessary to precise the topology on $D([0,T],E)$. For a much more detailed study of the Skorokhod space, we refer to \cite[Section 12]{billingsley_convergence_1999}. 

\begin{definition}[Skorokhod $J1$ topology] \label[I]{def:SkorokhodTopoJ1}
Let $\Lambda$ denote the set of strictly increasing homeomorphisms from $[0,T]$ onto itself. The Skorokhod $J1$ metrics on $D([0,T],E)$ is defined by
\[ d(x,y) := \inf_{\lambda \in \Lambda} {\left\{ \sup_{0 \leq t \leq T} \rho\big( x(t) , y(\lambda(t)) \big) + \sup_{s<t} {\left| \log {\frac{\lambda(t) - \lambda(s)}{t-s}} \right|} \right\}}. \]
Endowed with this metric, the Skorokhod space $D([0,T],E)$ is complete and separable.
\end{definition}

This topology is weaker than the one of continuous functions, which does not make $D([0,T],E)$ a complete space. However, estimates in the strong $\| . \|_{\infty}$ metrics implies estimates a complete metric, up to a fixed multiplicative constant. In practise, working with the $\| . \|_{\infty}$ metrics is thus often sufficient. 

Let us recall that for a continuous function $x \in C([0,T],E)$, the continuity modulus is defined for $0 < \delta < T$ by
\[ w_x(\delta) := \sup_{0 \leq t \leq T - \delta} w_x [t,t+\delta], \quad w_x I = \sup_{s,t \in I} \rho(x(t) , x(s)). \]
A function $x$ belongs to $C([0,T],E)$ if and only if $\lim_{\delta \to 0} w_x(\delta) = 0$. For c\`adl\`ag functions, another notion of modulus is defined. 

\begin{definition}[c\`adl\`ag modulus]
The c\`adl\`ag-modulus on $D([0,T],E)$ is defined by
\[ w'_x(\delta) := \inf_{\{ t_i \}} \max_i w_x [t_{i-1},t_i), \]
where the infimum is taken over sub-divisions $\{ t_i \}$ of $[0,T]$ such that $\min_i t_{i+1} - t_i > \delta$. A function $x$ belongs to $D([0,T],E)$ if and only if $\lim_{\delta \to 0} w'_x(\delta) = 0$.
\end{definition}

An analog of the Ascoli-Arzela theorem in the space of c\`adl\`ag functions is given by \cite[Theorem 12.3]{billingsley_convergence_1999}. It states that a subset $A \subset D([0,T],E)$ is relatively compact if and only if
\begin{enumerate}[(1)]
    \item $\sup_{x \in A} \| x \|_{\infty} < \infty$;
    \item $\lim_{\delta \to 0} \sup_{x \in A} w'_x(\delta) = 0$.
\end{enumerate}
In some cases, it is easier to use the modulus
\[ w''_x(\delta) := \sup_{\substack{t_1 \leq t \leq t_2 \\ t_2 - t_1 \leq \delta}} |x(t_1) - x(t) | \land |x(t) - x(t_2) |,\]
see \cite[Theorem 12.4]{billingsley_convergence_1999}. Using these results the following tightness criterion for probability measures on $D([0,T],E)$ is proved in \cite[Chapter 3, Corollary 7.4]{ethier_markov_1986}).

\begin{theorem}[Basic tightness criterion in $D$] \label[I]{thm:kurtztight}
For each $n\in\N$, let $(X^n_t)^{}_t$ be an adapted $E$-valued c\`adl\`ag process on the filtered probability space $(\Omega,\mathscr{F},(\mathscr{F}_t)_t,\mathbb{P})$. The sequence $(X^n_t)^{}_t$ is tight if and only if the following two conditions hold.
\begin{enumerate}[(1)]
\item For every $\varepsilon>0$ and every rational number $t\geq 0$, there exists a compact set $K_{\varepsilon,t} \subset E$ such that 
\[ \liminf_{n \to + \infty} \mathbb{P}(X^n_t \in K_{\varepsilon,t}) \geq 1 - \varepsilon. \]
\item For every $\varepsilon>0$ and $T > 0$, there exists $\delta >0$ such that
\[\limsup_{n \to + \infty} \mathbb{P}\big( {w'_{X^n|_{\left[0,T\right]}}(\delta) \geq \varepsilon} \big) \leq \varepsilon. \]
\end{enumerate}
\end{theorem}

Still, this criterion remains abstract and requires to find a suitable partition of the time interval to evaluate the c\`adl\`ag modulus. The tightness criteria which are used to prove propagation of chaos will be detailed in Appendix \cref{appendix:tightness}. 

\subsection{Stochastic processes and martingales}\label[I]{appendix:martingales}
 
In the following, let $(\Omega,\mathscr{F},\mathbb{P})$ be a probability space. This section reminds the reader of the basic properties related to stochastic processes and martingales. An exhaustive rigorous study can be found in the classical books \cite[Chapter 2]{ethier_markov_1986}, \cite{jacod_limit_2003}, \cite{revuz_continuous_1999} or \cite{le_gall_brownian_2016}. The present section also summarizes elements of \cite{joffe_weak_1986}.

\subsubsection{Martingales} 

\begin{definition}[Stochastic process]
A \emph{stochastic process} with state space $E$ (a measurable space) and indexed by a set $I$ (often $I=[0,T]$) is a function $X: I\times \Omega  \to E$ such that $X_t = X(t,\cdot)$ is a $E$-valued random variable (that is is to say $\omega \mapsto X_t(\omega)$ is measurable) for every $t \in I$. 
\begin{enumerate}[(1)]
    \item The (pathwise) \emph{law} of $(X_t)_t$ is the push-forward measure $f_I = X_{\#}\mathbb{P}$ on a space $\mathcal{F}(I,E)$ of functions $I \to E$. 
    \item The \emph{time marginal laws} are the push-forward measures on $E$ defined for any $t\in I$ by $f_t = {\mathsf{X}^E_t}_{\#}f_I$, provided that the evaluation maps $\mathsf{X}^E_t : \mathcal{F}(I,E) \to E, \omega\mapsto \omega(t)$ are measurable.
\end{enumerate}
The stochastic process $X = (X_t)_{t \in I}$ can be seen as a $\mathcal{F}(I,E)$-valued random variable. For a given $\omega \in \Omega$, $(X_t(\omega))_{t \geq 0}$ is called a \emph{sample path} (or trajectory) of $X$.
\end{definition}

From now on let $I=[0,T]$ with $T\in(0,+\infty]$.

\begin{example}[Canonical stochastic process] \label[I]{example:canonicalprocess} Given a probability distribution on the space of functions $I\to E$, $f_I \in \pb( \mathcal{F}(I,E) )$, the \emph{canonical stochastic process with law $f_I$} is $\mathsf{X} = ( \mathsf{X}^E_t )_{t \in I}$ defined on the probability space $(\mathcal{F}(I,E),\mathscr{F},f_I)$: the law of $\mathsf{X}$ is indeed ${\mathsf{X}}_{\#}f_I = f_I$.
\end{example}

\begin{definition}[Filtration]
A filtration is an increasing family of $\sigma$-algebras $(\mathscr{F}_t)_{t \geq 0}$, i.e. such that $\mathscr{F}_s \subset \mathscr{F}_t$ for $s \leq t$. A filtration is said to be:
\begin{enumerate}[(1)]
    \item \emph{complete} when
    \[ \forall t \geq 0,\quad \{ \mathcal{A} \subset \Omega, \,\, \mathcal{A}\subset \mathcal{B}, \, \mathbb{P}(\mathcal{B}) = 0 \} \subset \mathscr{F}_t; \]
    \item \emph{right-continuous} when
    \[ \forall t \geq 0, \quad \mathscr{F}_t = \mathscr{F}_{t^+}, \quad \mathscr{F}_{t^+} := \bigcap_{\varepsilon > 0} \mathscr{F}_{t+\varepsilon}. \]
\end{enumerate}
In the following, let $(\Omega,\mathscr{F},(\mathscr{F}_t)_{t \geq 0},\mathbb{P})$ be a filtered probability space, whose filtration is assumed to be complete and right-continuous, together with $\mathscr{F} = \bigcup_{t \geq 0} \mathscr{F}_t$.
\end{definition}

\begin{definition}[Regularity for stochastic processes]
A stochastic process $X = (X_t)_{t \geq 0}$ is said to be 
\begin{enumerate}[(1)]
\item \emph{adapted} to the filtration $(\mathscr{F}_t)_{t \geq 0}$ when the random variable $X_t$ is $\mathscr{F}_t$-measurable;
\item \emph{predictable} when $X$ is measurable for the $\sigma$-algebra generated by the sets $(s,t] \times F$ for $0 \leq s \leq t$ and $F\in\mathscr{F}_s$ ; 
\item \emph{with finite variation paths} when a.s. $(X_t(\omega))_{0\leq t \leq T}$ has bounded variation on any finite time-interval $[0,T]$ ; 
\item \emph{continuous} (resp. c\`adl\`ag, right-continuous...) when for $\mathbb{P}$-almost every $\omega \in \Omega$, the sample path $(X_t(\omega))_{t \geq 0}$ is a continuous (resp. c\`adl\`ag, right-continuous...) function of $t$.
\end{enumerate}
\end{definition}

\begin{definition}[Martingale and sub-martingale]
An adapted real-valued process $X = (X_t)_{t \geq 0}$ such that $\E|X_t|<+\infty$ for every $t\geq0$ is a $(\mathscr{F}_t)_t$-martingale when
\[ \forall t,s\geq 0, \quad \E[X_{t+s} | \mathscr{F}_t ] = X_t. \]
It is a \emph{sub-martingale} when 
\[ \forall t,s\geq 0, \quad \E[X_{t+s} | \mathscr{F}_t ] \geq X_t. \]
These definitions are extended componentwise to $\R^N$-valued processes.
\end{definition}

Sub-martingales enjoy many useful properties, starting with the following proposition \cite[Chapter 2, Proposition 2.9]{ethier_markov_1986}.

\begin{proposition}[C\`adl\`ag modification]
If $X=(X_t)_{t \geq 0}$ is a sub-martingale, then there exists an adapted real-valued process $Y=(Y_t)_{t \geq 0}$ such that
\begin{equation} \label[I]{def:modification}
\forall t \geq 0, \quad P( X_t  = Y_t ) = 1, 
\end{equation}
and $Y$ is c\`adl\`ag outside a countable set of times $t$. In the following this set is assumed to be $\emptyset$, so (sub-)martingales will be considered as c\`adl\`ag processes. 
\end{proposition}

Equation \cref{def:modification} means that $Y$ is a \emph{modification} of $X$. The following inequality controls the growth of sub-martingales.

\begin{proposition}[Doob's maximal inequality]
Given a sub-martingale $(X_t)_{t \geq 0}$, $T>0$ and any $p>1$, it holds that
\[\E\Big|\sup_{0\leq t \leq T} X_t\Big|^p \leq {\left(\frac{p}{p-1}\right)}^p \E|X_T|^p.\]
\end{proposition}

\subsubsection{Local martingales and quadratic variation}

In order to define a stochastic integral, the notion of martingale needs to be weakened to the notion of local martingale. 

\begin{definition}[Local martingale]
A real-valued process $(M_t)_{t \geq 0}$ is a \emph{local martingale} when there exists a sequence $(\tau_n)_{n \in \N}$ of stopping times such that $\tau_n \to + \infty$ and $(M_{t \wedge \tau_n})_{t \geq 0}$ is a martingale for every $n \geq 0$. This definition is extended componentwise to $\R^N$-valued processes.
\end{definition}

\begin{definition}[Quadratic variation]
The quadratic variation $[ M ] = ( [ M ]_t )_{t \geq 0}$ of a (local) square integrable martingale $(M_t)_{t \geq 0}$ is defined as the limit in probability
\[ [M]_t := \lim_{\delta(\pi) \rightarrow 0} \sum_{t_n \in \pi} \big(M_{t_{n+1} \land t} - M_{t_{n} \wedge t}\big)^2, \]
where $\pi$ is a subdivision $0=t_0<t_1<\ldots$ of $\R_+$ with mesh $\delta(\pi)=\sup_n (t_{n+1}-t_n)$. It is the only (up to modification) adapted increasing process with jumps $\Delta[M]_t = \Delta M^2_t= M_t^2 - M^2_{t^-}$ such that $\left( {M^2_t - {\left[ {M} \right]}_t} \right)_t$ is a (local) martingale.
\end{definition}

\begin{example}[The case of finite variation paths] When $(M_t)_{t \geq 0}$ has finite variation paths, it is straightforward to check that
\[ [M]_t = \sum_{s \leq t} (\Delta M_s)^2.\]
\end{example}

\begin{definition}[Cross-variation]
The cross-variation $[ M,N ] = ( [ M,N ]_t )_{t \geq 0}$ of two (local) square integrable martingales $(M_t)_{t \geq 0}$ and $(N_t)_{t \geq 0}$ is defined by 
\[ [M,N]_t := \frac{1}{2}\big( [M+N]_t - [M]_t - [N]_t \big). \]
It is the only (up to modification) adapted increasing process with jump at time $t>0$, $\Delta[M,N]_t = \Delta (M N)_t= M_t N_t - M_{t^-} N_{t^-}$ such that $( {M_t N_t - [M,N]_t} )_t$ is a (local) martingale. 
\end{definition}

\begin{proposition}[BDG inequality]\label[I]{prop:bdg} For every $p \geq 1$, there exist two constants $c_p, C_p>0$ such that  for any local martingale $(M_t)_{t \geq 0}$, the supremum $M^{\star}_t := \sup_{0 \leq s \leq t} | M_s |$ can be controlled by the quadratic variation in $L^p$-norm:
\[ c_p \E [ {M} ]_t^{p/2} \leq \E ( {M^{\star}_t} )^p \leq C_p \E [ {M} ]_t^{p/2}. \]
\end{proposition}

\begin{definition}[Predictable quadratic and cross-variation]
The predictable quadratic variation $\langle M \rangle = (\langle M \rangle_t )_{t \geq 0}$ of a (local) square integrable martingale $(M_t)_{t \geq 0}$ is the unique \emph{predictable} c\`adl\`ag increasing process (with finite variation paths) such that
\[ \big(M^2_t - \langle M\rangle_t\big)_{t\geq 0} \]
is a (local) martingale. The existence of the predictable quadratic variation stems from the Doob-Meyer decomposition theorem for $DL$-supermartingales (see \cite[Chapter 2, Proposition 5]{ethier_markov_1986} for more details). The predictable cross-variation of two (local) square integrable martingales $M,N$ is defined the same way, setting
\[ \langle M,N\rangle_t := \frac{1}{2}\big(  \langle M+N\rangle_t -  \langle M\rangle_t -  \langle N\rangle_t \big). \]
\end{definition}

Note that by substracting the martingale characterizations of the quadratic variation and of the predictable quadratic variation, we get that $([M]_t - \langle M\rangle_t)_{t\geq 0}$ is a (local) square integrable martingale. The predictable quadratic variation is said to be the \emph{compensator} of the quadratic variation. Moreover, the equality $[M] = \langle M\rangle$ holds as soon as $(M_t)_{t \geq 0}$ is a continuous process.

\subsubsection{Semimartingales}

The notion of local martingales is then extended to the notion of semimartingale, which forms a large class of processes against which a stochastic integral can be defined (a notion that will not be detailed here). 

\begin{definition}[Semimartingale]
A real-valued process $(X_t)_{t \geq 0}$ is a semimartingale when it can be decomposed as 
\[ X_t = M_t + A_t, \]
where $(M_t)_{t \geq 0}$ is a local martingale and $(A_t)_{t \geq 0}$ is an adapted process with finite variation paths. This definition is extended component-wise to $\R^N$-valued processes.
\end{definition}

The cross-variation $[X,Y]$ and predictable cross-variation $\langle X,Y\rangle$ can be extended to semimartingales $(X_t)_{t \geq 0}$ and $(Y_t)_{t \geq 0}$ as the limit in probability 
\[ [X,Y]_t = \lim_{\delta(\pi) \rightarrow 0} \sum_{t_n \in \pi} (X_{t_{n+1} \wedge t} - X_{t_{n} \wedge t})(Y_{t_{n+1} \wedge t} - Y_{t_{n} \wedge t}). \]
Writing $X_t = M_t + A_t$, note that $[X] = [M]$ and $\langle X \rangle = \langle M \rangle$ because $A$ has finite variation paths. 

For vector-valued semimartingales, the cross-variation is defined componentwise as follows. 

\begin{definition}[Vectorial cross-variation] Let $(\mathbf{X}_t)_{t \geq 0}$ be a $\R^d$-valued semimartingale with the notation $\mathbf{X}_t = \left(X^{1}_t,\ldots,X^{d}_t\right)$, the matrix-valued cross-variations are defined by 
\[ \llbracket {\mathbf{X}} \rrbracket = \big([ {X^{i},X^{j}} ]\big)_{1 \leq i,j \leq d}, \quad \llangle {\mathbf{X}} \rrangle = \big(\langle {X^{i},X^{j}} \rangle\big)_{1 \leq i,j \leq d}, \]
and the related scalar quantities are the traces of theses matrices (defined as the sum of diagonal elements)
\[ [ {\mathbf{X}} ] = \Tr\big( \llbracket {\mathbf{X}} \rrbracket\big), \quad \langle {\mathbf{X}} \rangle = \Tr\big( \llangle {\mathbf{X}} \rrangle\big). \]
\end{definition}

The integration of locally bounded predictable processes $(H_t)_{t \geq 0}$ against finite variation processes is well-defined using Stieltjes formalism. The theory of stochastic integration extends this to integrate $(H_t)_{t \geq 0}$ against any square integrable semimartingale $(X_t)_{t \geq 0}$. Notable identities are the It\=o isometry
\[ \E{\left[{\left( \int_0^t H_{s^-} \dd X_s \right)}^2\right]} = \E{\left[ \int_0^t H^2_{s^-} \dd [ {X} ]_s\right]} =  \E{\left[ \int_0^t H^2_{s^-} \dd \langle {X} \rangle_s \right]}, \]
(the last equality holds because $[ {X} ] - \langle {X} \rangle$ is a local martingale) and the integration by parts formula
\[ X_t Y_t = X_0 Y_0 + \int_0^t X_{s^-} \dd Y_s + \int_0^t Y_{s^-} \dd X_s + [ X,Y ]_t, \]
for general semimartingales $(X_t)_{t \geq 0}$ and $(Y_t)_{t \geq 0}$. Given a (local) square integrable martingale $(M_t)_{t \geq 0}$, the process
\[ {\left( \int_0^t H_{s^-} \dd M_s \right)}_{t\geq 0} \]
is built in such a way that it is a (local) square integrable martingale.

\begin{example}[Dynkin's representation formula] Anticipating on the next section, a Markov process $(X_t)_{t \geq 0}$ with generator $L$ is such that for any $\varphi$ in $\mathcal{D}(L)$ (see Theorem \cref{thm:charactFullGene})
\[M^{\varphi}_t := \varphi(X_t) - \varphi(X_0) - \int_0^t L \varphi (X_s)\dd s,\]
is a martingale. In this case, the predictable quadratic variation $\langle {M^{\varphi}} \rangle$ can be computed as
\[ \langle {M^{\varphi}} \rangle_t = \int_0^t \Gamma_L(\varphi,\varphi)(X_s) \dd s, \]
involving the \emph{carr\'e du champ operator}
\[ \Gamma_L(\varphi,\psi) := L[\varphi\psi] - \varphi L \psi - \psi L \varphi, \]
for any $\varphi,\psi\in \mathcal{D}(L)$. 
These properties can be generalized to obtain a wider class of semimartingales, which are non-necessarily Markov processes.
\end{example}

\subsubsection{\texorpdfstring{$D$}{D}-semimartingales}\label[I]{appendix:Dsemimartingales}

Dynkin's formula can be turned into a definition to obtain a wider class of (non-necessarily Markov) $\R^d$-valued semimartingales, named class $D$ (in the honour of Dynkin). The following definition can be found in \cite{joffe_weak_1986}. The notation $\pi_i$ denotes the coordinate function $\mathbf{x}^d \mapsto x^i$.

\begin{definition}[Semimartingale of class $D$] \label[I]{def:semiMartD}
A $\R^d$-valued semimartingale $\mathbf{X}_t$ belongs to the class $D$ when there exist an increasing c\`adl\`ag function $t \mapsto A(t)$, a vector space $\mathcal{C}$ of $\R$-valued continuous functions on $\R^d$ and a mapping $L: \mathcal{C} \times \R^d \times \R_+ \times \Omega \to \R$ such that the following properties hold.
\begin{enumerate}[(1)]
    \item For every $1 \leq i,j \leq d$, $\pi_i$ and $\pi_i \pi_j$ belong to $\mathcal{C}$. 
    \item For every $(\mathbf{x}^d,t,\omega)\in \R^d \times \R_+ \times \Omega$, the map $\varphi \mapsto L( \varphi,\mathbf{x}^d,t,\omega )$ is linear and maps $\mathcal{C}$ to itself. Moreover for every $\varphi \in \mathcal{C}$, the map $(\mathbf{x}^d,t,\omega) \mapsto L( \varphi,\mathbf{x}^d,t,\omega )$ is measurable for the $\sigma$-algebra on $\R^d$ of Borel sets, and the one of predictable events on $\Omega$.
    \item For every $\varphi\in\mathcal{C}$,
    \[ M^{\varphi}_t := \varphi(\mathbf{X}_t) - \varphi(\mathbf{X}_0) - \int_0^t L( \varphi,\mathbf{X}_{s^-},s,\cdot ) \dd A(s), \]
    is a local square integrable martingale.
\end{enumerate}
To this process are associated the \emph{local coefficients} 
\[ b_i(\mathbf{x}^d , t, \omega) = L(\pi_i,\mathbf{x}^d , t,\omega), \quad a_{ij}( \mathbf{x}^d , t, \omega) = \Gamma_{L(\cdot,\mathbf{x}^d t,\omega)} (\pi_i,\pi_j)(\mathbf{x}^d), \]
and the \emph{drift vector} $\mathbf{b} = (b_i)_{1\leq i \leq d}$ and the \emph{diffusion matrix} $\mathbf{a} = (a_{ij})_{1\leq i,j \leq d}$. In the following we often omit to write the dependency in $\omega$.  
\end{definition}

The following lemma is proved in \cite[Lemma 3.1.3]{joffe_weak_1986}.

\begin{lemma}[Predictable variation] \label[I]{lemma:semiMartVariation}
Given a $D$-semimartingale $\mathbf{X}$, define
\[ \mathbf{M}_t := \mathbf{X}_t - \mathbf{X}_0 - \int_0^t \mathbf{b}(\mathbf{X}_{s^-}, s) \dd A(s). \]
Then $(\mathbf{M}_t)_{t \geq 0}$ is a $\R^d$-valued local square integrable martingale, whose (scalar) predictable quadratic variation reads
\[ \langle\mathbf{X}\rangle_t = \int_0^t \Tr\big(\mathbf{a}( \mathbf{X}_{s^-} , s)\big) \dd A(s) - \sum_{s \leq t} \| {\mathbf{b}(\mathbf{X}_{s^-}, s)} \|^2 | \Delta A (s) |^2. \]
\end{lemma}

\begin{remark}[Generalized SDE and intrinsic randomness] The previous notions naturally extend the notion of diffusion SDE, since at least formally
\[ \dd \mathbf{X}_t = \mathbf{b}(\mathbf{X}_{t^-}, t) \dd A(t) + \mathbf{a}( \mathbf{X}_{t^-} , t) \dd \mathbf{M}_t. \]
This parallel could be used to extend coupling and completeness methods to irregular processes or even non-Markov ones. However, remember that in this case, the generator $L$ depends on $\omega$ and has therefore \emph{its own source of randomness}, together with the local coefficients $\mathbf{b}$ and $\mathbf{a}$.
\end{remark}

Similarly to Markov processes, a convenient way to build a $D$-semimartingale on $\Omega = D([0,T],\R^d)$ (for $T\in(0,+\infty])$ is to solve a martingale problem.

\begin{definition}[Martingale problem]
A probability distribution on the path space $f_I\in\pb(D([0,T],\R^d))$ is a solution to the martingale problem issued from $f_0\in\pb(\R^d)$ whenever for all all $\varphi\in\mathcal{C}$,
\[M^\varphi_t := \varphi(\mathsf{X}_t)-\varphi(\mathsf{X}_0) - \int_0^t L(\varphi,\mathsf{X}_{s^-},s,\cdot) \dd A(s),\]
is a $f_I$-martingale, where $\mathsf{X}_t$ is the canonical process $D([0,T],\R^d) \to \R^d$. 
\end{definition}

\subsection{Markov processes and Markov representation for PDEs}\label[I]{appendix:timeinhomoegeneousmarkov}

The purpose of this section is to briefly review the probabilistic framework for linear and nonlinear Markov processes. Classical references and review articles on the subject include \cite{ethier_markov_1986, revuz_continuous_1999, bottcher_levy_2013, cont_markov_2010}. The prototypical nonlinear Markov process is the solution of the McKean-Vlasov SDE:
\[ \dd X_t = b(X_t,f_t)\dd t + \sigma(X_t,f_t)\dd B_t, \quad X_t \sim f_t, \]
where $(B_t)_{t\geq 0}$ is a standard Brownian motion and $f_t$ the law of the random variable $X_t$ at time $t$. Such a process is said to be nonlinear (in the sense of McKean) since its definition depends on its own law. This type of nonlinearity has been introduced in the seminal \cite{mckean_propagation_1969}. Unlike, classical (linear) Markov process, the law $f_t$ of the nonlinear Markov process satisfies the nonlinear the Fokker-Planck equation
\[ \partial_t f_t( x ) = -\nabla_x\cdot\{b( x ,f_t)f_t\}+\frac{1}{2}\sum_{i,j=1}^d \partial_{x_i}\partial_{x_j}\{a_{ij}(x,f_t)f_t\}, \]
where the matrix $a = (a_{ij})$ is defined by $a = \sigma\sigma^\mathrm{T}$. This kind of equations is derived in the limit $N \to +\infty$ in the Kolmogorov equations associated to large interacting particle systems defined by $N$ linear Markov processes. While the theory of linear Markov processes is well-established, nonlinear Markov processes are not so classical in the literature and require specific tools to be built (for instance, the well-posedness result proved in Proposition \cref{prop:wellposednessmckeanlipschitz} for the McKean-Vlasov SDE). 

The classical theory of time homogeneous linear Markov processes is presented in Section \cref{sec:timehomogeneouslinearmarkov}. Elements of the theory of time-inhomogeneous Markov processes can be found in Section \cref{sec:timeinhomogeneouslinearmarkov}. Using these concepts and following \cite{johnson_class_1968} and \cite{mckean_class_1966}, a theoretical framework is presented in Section \cref{sec:nonlinearmarkov} for nonlinear Markov processes (in the sense of McKean). Before that, we first recall the basic notions about Markov process (linear or not).  Following \cite[Chapter 4]{ethier_markov_1986}, the very general definition of Markov processes is the following.

\begin{definition}[Markov process]
A stochastic process $(X_t)_{t \geq 0}$ is a Markov process on $(\Omega,\mathscr{F},\mathbb{P})$ when for every $s, t \geq 0$
\begin{equation} \label[Ieq]{eq:Markovprop}
\forall \mathscr{A} \in \mathcal{B}(E), \quad \mathbb{P}{\left(X_{s+t} \in \mathscr{A} | (X_r)_{0 \leq r \leq t}\right)} = \mathbb{P}(X_{s+t} \in \mathscr{A} | X_t). 
\end{equation}
On the left-hand side the probability is conditioned by the filtration generated by the process, but stronger definitions could involve wider filtrations. The relation~\cref{eq:Markovprop} will be referred as the \emph{Markov property}. When the time $t$ is replaced by a random \emph{stopping time}, this relation is called the \emph{strong Markov property}. 
\end{definition}

This means the law of $X_{s+t}$ at time $s+t$ conditionally on the past history up to time~$t$, is the same as the law at time $s+t$ conditionally on the state at time $t$ only. The definition can be equivalently written in terms of bounded measurable test functions $\varphi \in B_b(E)$ as
\begin{equation}
\E {\left[\varphi(X_{s+t})|(X_r)_{0 \leq r \leq t}\right]} = \E [\varphi(X_{t+s})|X_t].
\end{equation}

This definition remains abstract and does not tell how to \emph{build} a Markov process. If the law $f_I$ is known (built from given processes or e.g. by solving a martingale process), a Markov process with law $f_I$ is given by the canonical process $\mathsf{X}^E =(\mathsf{X}^E_t)^{}_{t \geq 0}$ on the probability space $\Omega =D(I,E)$ endowed with $f_I$ and an adequate filtration (see Example \cref{example:canonicalprocess}). Further constructions which are closer to the PDE point of view are presented in the next subsections.

\subsubsection{Time-homogeneous Markov processes and linear PDEs}\label[I]{sec:timehomogeneouslinearmarkov}

The content of this section is quite classical and can be found in the classical reference \cite{ethier_markov_1986}. 

\subsubsection*{Transition functions.} A usual way to present Markov processes is to first think about time discrete Markov chains with jumps which are given by a transition kernel (or transition matrix when $E$ is discrete). This definition can be generalized to the time continuous framework using the notion of transition function. 

\begin{definition}[Transition functions, homogeneous case]
A family of maps $(P_t)_{t\geq0}$ where $P_t : E\to\pb(E)$ is a family of transition functions when the following properties hold.  
\begin{itemize}
    \item The map $(t,x) \mapsto P_t\left(x,\cdot\right)$ is a measurable map $[0,\infty) \times E \rightarrow \pb\left(E\right)$.
    \item For all $x\in E$, $P_0(x, \cdot) = \delta_x$.
    \item For all $s,t\geq0$ and all $\mathscr{A}\subset\mathcal{B}(E)$, $P_{s+t}(x,\mathscr{A}) = \int_E P_{s}(y,\mathscr{A}) P_{t}(x,\dd y)$.
\end{itemize}
This last relation is called the \emph{Chapman-Kolmogorov property}.
\end{definition}

A family of transition functions $(P_t)_{t\geq0}$ is said to be \emph{adapted to the Markov process} $(X_t)_{t \geq 0}$ when for all $s,t\geq0$, and all $\mathscr{A}\subset\mathcal{B}(E)$,
\[ \mathbb{P}(X_{s+t} \in \mathscr{A} | (X_r)_{0 \leq r \leq t}) = P_s(X_{t},\mathscr{A}),\]
which is equivalent to  \[\E[\varphi(X_{t+s})|(X_r)_{0 \leq r \leq t}] = \int_E \varphi(y) P_{s}(X_t,\dd y), \]
for any bounded Borel measurable test function $\varphi\in B_p(E)$. Note that this relation and the Markov property \cref{eq:Markovprop} imply the Chapman-Kolmogorov property with $x = X_r$ for $r\geq0$. Moreover given the initial distribution $X_0\sim f_0$, the finite dimensional distributions of $(X_t)_t$ can be computed by
\begin{align} 
&\mathbb{P}(X_0 \in \mathscr{A}_0, X_{t_1} \in \mathscr{A}_1, \ldots, X_{t_n} \in \mathscr{A}_n) \nonumber\\  
&=\int_{\mathscr{A}_0} \ldots \int_{\mathscr{A}_{n-1}} P_{t_n - t_{n-1}}(y_{n-1},\mathscr{A}_n) P_{t_{n-1} - t_{n-2}}(y_{n-2},\dd y_{n-1}) \nonumber \\
&\phantom{=\int_{\mathscr{A}_0} \ldots \int_{\mathscr{A}_{n-1}} P_{t_n - t_{n-1}}(y_{n-1},\mathscr{A}_n) P_{t_{n-1} - t_{n-2}}(y_{n-2}}\ldots P_{t_1}(y_0,\dd y_1)f_0(\dd y_0), \label[Ieq]{eq:Markovmarginals}
\end{align}
for any $\mathscr{A}_0,\ldots, \mathscr{A}_n\in \mathcal{B}(E)$ and $0\leq t_1\leq\ldots\leq t_n$.  In fact transition functions are sufficient to build a Markov process; this is a consequence of the Kolmogorov extension theorem (or more generally of the Ionescu-Tulcea theorem), see for instance \cite[Chapter 4, Theorem 1.1]{ethier_markov_1986} or \cite[Chapter 3, Theorem 1.5]{revuz_continuous_1999}. 

\begin{theorem}[Markov process built from transition functions]
Let $E$ be a Polish space and $f_0\in \pb(E)$. Given the transition functions $(P_t)_{t\geq0}$, there exists a Markov process whose finite dimensional distributions are uniquely determined by \cref{eq:Markovmarginals}. Its law on the path space is a probability measure $f_I \in \pb(D( I,E))$.
\end{theorem}

If $f_0 = \delta_x$, then $f_I$ is denoted by $f_I^x$ and \cite[Chapter 4,Propositon 1.2]{ethier_markov_1986} proves that the map $x \mapsto f_I^x(\mathscr{B})$ is measurable for any Borel set $\mathscr{B} \subset D( I,E)$. 

\begin{example}[Brownian motion]
Transition functions are not explicit in general. A notable exception is the $d$-dimensional Brownian motion for which
\[ P_t(x,\mathscr{A}) = (2 \pi t)^{-d/2}\int_\mathscr{A} \exp {\left(-\frac{|x-y|^2}{2t}\right)} \dd y. \]
Note that the map $t \mapsto P_t(x,\dd y)$ is the measure solution of the 1D heat equation $\partial_t u = \partial_{xx}^2 u$ with initial condition $u(t=0,\cdot)=\delta_x$. 
\end{example}

\subsubsection*{Semigroup representation.} The connection between Markov processes and linear PDEs is given by the semigroup representation as explained below. For $t\geq0$, let the linear operator $T_t$ acting on  be defined by:
\begin{equation} \label[Ieq]{transitionsemi}
T_t \varphi ( x ) := \int_E \varphi (y) P_t(x,\dd y),
\end{equation}
for any test functions $\varphi \in B_b(E)$. Thanks to the Chapman-Kolmogorov property, this defines a positive measurable semigroup of contractions on $B_b(E)$, where we recall that a family of bounded operators $(T_t)_{t\geq0}$ on a closed subspace $\mathcal{D} \subset B_b(E)$ is a semigroup when $T_{0} = \mathrm{Id}$ and $T_{t+s} = T_{t} T_{s}$ for all $s,t\geq0$. It is said to be a contraction semigroup when the operators are bounded with norm smaller or equal to 1. The semi-group $(T_t)_t$ is said to correspond to a Markov process $(X_t)_{t \geq 0}$ when
\[ \forall s,t \geq0,\,\, \forall \varphi \in \mathcal{D}, \quad T_t \varphi(X_{t+s}) = \E {\left[\varphi(X_{s+t})|(X_r)_{0 \leq r \leq t}\right]} \]
The semi-group representation characterises a Markov process, as stated in \cite[Chapter 4, Proposition 1.6]{ethier_markov_1986}.

\begin{theorem}
Let $E$ be a Polish space and let $\mathcal{D} \subset B_b(E)$ be a closed subspace assumed to be separating. Let $f_0$ in $\pb(E)$ and let $(T_t)_{t\geq0}$ be a semigroup on $\mathcal{D}$ corresponding to a Markov process $(X_t)_t$. Then the finite dimensional distributions of $(X_t)_t$ are determined by $(T_t)_{t\geq0}$ and $f_0$.
\end{theorem}

In the following, starting from a semigroup $(T_t)_t$, the goal is to construct a corresponding Markov process. With the notable exception of jump processes, stronger assumptions on $E$ and $(T_t)_t$ are often needed, as the ones given in the following definition. 

\begin{definition}[Feller semigroup]\label[I]{def:feller} Let $E$ be a locally compact Polish space. A semi-group $(T_t)_{t\geq0}$ is a Feller semigroup when its elements satisfy the following properties. 
\begin{enumerate}[(1)]
    \item \textbf{(Feller).} For all $t\geq0$, $T_t$ maps $C_0(E)$ to $C_0(E)$, where $C_0(E)$ is the space of continuous functions vanishing at infinity. It means that
    \[\forall \varphi\in C_0(E),\quad T_t\varphi \in C_0(E).\]
    \item \textbf{(Contraction).} For all $t\geq0$ and all  $\varphi \in C_0(E)$, $\| T_t \varphi \|_{\infty} \leq \| \varphi \|_{\infty}$.
    \item \textbf{(Mass preserving).} For all $t\geq0$, $T_{t} 1 = 1 $. 
    \item \textbf{(Positivity).} For all $t\geq0$ and all  $\varphi \in C_0(E)$ such that  $\varphi \geq 0$, then $T_{t} \varphi \geq 0 $.
    \item \textbf{(Strongly continuous).} For all $\varphi\in C_0(E)$, $\| T_t \varphi - \varphi \|_{\infty} \to 0$ when $t \to 0$.
\end{enumerate}
\end{definition}
The properties (2) to (4) hold when $T_t$ is defined from \cref{transitionsemi}. Also note that every semigroup defined on $C_0(E)$ can be uniquely extended to the whole space $B_b(E)$ (see \cite[Theoem 1.5]{bottcher_levy_2013}). A Markov process corresponding to a Feller semigroup is called a \emph{Feller process}. There is a one-to-one correspondence between Feller semigroups and Feller processes; this is a consequence of the Riesz representation theorem, see for instance \cite[Section 1.2]{bottcher_levy_2013}, \cite[Chapter 3, Proposition 2.2]{revuz_continuous_1999} or \cite[Chapter 4, Theorem 2.7]{ethier_markov_1986}. 

\begin{theorem}[C\`adl\`ag Markov process from a Feller semigroup]
Let $E$ be a locally compact Polish space and let $f_0\in\pb(E)$. Let $(T_t)_t$ be a Feller semigroup on $E$. Then there exists a unique transition function $(P_t)_t$ on $E$ such that \cref{transitionsemi} holds. As a consequence, there exists a Markov process corresponding to $(T_t)_t$ with initial distribution $f_0$, whose finite dimensional distributions are uniquely determined by $(T_t)_{t}$. Moreover this process has a c\`adl\`ag modification and satisfies the strong Markov property with respect to the right-continuous filtration $\mathscr{G}_t = \cap_{\varepsilon>0} \sigma((X_s)_{s\leq t+\varepsilon})$. 
\end{theorem}

\begin{example}[Brownian motion] The $d$-dimensional Brownian motion is a Feller process. More generally diffusion processes are Feller processes under mild assumptions on the diffusion coefficients.    
\end{example}

\begin{example}[Markov jump processes and $C_b$-Feller processes] The definition of Feller process is not universal in the literature. An important variant is the notion of $C_b$-Feller process for which the space $C_0(E)$ is replaced by $C_b(E)$ in Definition \cref{def:feller}. In this case, the local compactness assumption on $E$ can be dropped. Diffusion processes are \emph{not} $C_b$-Feller processes because the diffusion semigroup is not strongly continuous on $C_b(E)$ (see \cite[Example 1.7d]{bottcher_levy_2013}). The main class of $C_b$-Feller process that are considered in this review is the class of Markov jump processes, defined by the transition function: 
\[\forall x\in E,\quad P_t(x,\dd y) = \sum_{k=0}^{+\infty} \frac{\e^{-t}t^k}{k!} P^k(x,\dd y),\]
where $P:E\times \mathcal{B}(E)\to\pb(E)$ is a transition probability and for $k\in \N$, $x\in E$ and $\mathscr{A}\in \mathcal{B}(E)$,  
\[P^k(x,\mathscr{A}) = \int_E P(x,\dd y) P^{k-1}(y, \mathscr{A}).\]
An explicit construction of a Markov jump process can be found in \cite[Chapter 4, Section 2]{ethier_markov_1986}, see also \cite[Chapter 3, Exercise 1.8]{revuz_continuous_1999}. Markov jump processes are more easily understood through their generator as defined below. More generally, the strong continuity property holds on $C_b(E)$ for semigroups which have a bounded generator (for the $\|\cdot\|_\infty$ topology). Further links between Feller semigroups, $C_b$-Feller semigroups and other notions of Feller semigroups can be found in \cite[Section~1.1]{bottcher_levy_2013} and the references therein. 
\end{example}

\subsubsection*{Infinitesimal generator and PDE.} Strongly continuous contraction semigroups are determined by their infinitesimal generator.

\begin{definition}[Infinitesimal generator] \label{def:infgene}
The infinitesimal generator $L$ of a strongly continuous contraction semigroup $(T_t)_{t\geq0}$ on a closed subspace $\mathcal{D}\subset B_p(E)$ is the linear operator defined by
\[ L \varphi = \lim_{\substack{t \to 0 \\ t > 0}} \frac{T_t \varphi - \varphi}{t} \]
on the domain $\mathcal{D}(L) \subset \mathcal{D}$ of functions $\varphi$ such that the limit exists (for the topology on $\mathcal{D}$).
\end{definition}
The first connection between Markov processes and PDE comes from the observation that for all $\varphi \in \mathcal{D}(L)$, 
\begin{equation}\label[Ieq]{eq:backkolmogorovappendix}\frac{\dd}{\dd t} T_t \varphi = T_t L\varphi,\quad  \frac{\dd}{\dd t} T_t \varphi = LT_t\varphi,\end{equation}
which are called the \emph{Kolmogorov equations}. The first equation is called the forward Kolmogorov equation and the second, the backward Kolmogorov equation. The terminology will appear more clearly in the time-inhomogeneous setting below. Since $T_t$ and $L$ commute, the two equations are of course equivalent but each one has its own physical interpretation. The backward equation gives a Markov representation of the solution of the linear PDE $\partial_t u = Lu$ with initial condition $\varphi$ as the conditional expectation $ u(t,x) = T_t \varphi (x) = \E [\varphi(X_t)| X_0=x]$. The backward equation thus describes the evolution of an observable of the Markov process. A more general version when source terms or boundary conditions are added is the Feynman-Kac formula. More on the forward equation is given in the next paragraph. Note that the generator $L$ can include differential and jump terms in which case the (backward) Kolmogorov equation is an integro-differential PDE. 

\begin{example} The generator of the $d$-dimensional Brownian motion is $L\varphi = \frac{1}{2}\Delta\varphi$ and $\mathcal{D}(L)\subset C^2_0(\R^d)$. The generator of a Markov jump process on a space $E$ is 
\[L\varphi(x) = \int_E \{\varphi(y)-\varphi(x)\}P(x,\dd y),\]
and $L$ is bounded on $(C_b(E), \|\cdot\|_\infty)$. 
\end{example}

In stochastic analysis, we often use the notion of \emph{full generator} which turns the (forward) Kolmogorov equation into a definition. 

\begin{definition}[Full generator] The full generator of a strongly continuous contraction semigroup $(T_t)_t$ on $\mathcal{D}$ is the subset 
\[\widehat{L} := {\left\{(\varphi,\psi)\in \mathcal{D}\times \mathcal{D},\,\, T_t\varphi - \varphi = \int_0^t T_s \psi\, \dd s\right\}}.\]
\end{definition}
The (forward) Kolmogorov equation says that $\{(\varphi, L\varphi),\,\,\varphi\in \mathcal{D}(L)\}\subset \widehat{L}$. The full generator is often used in connection with the martingale characterisation of a Markov process due to Stroock and Varhadan, see \cite[Corollary 1.37]{bottcher_levy_2013} and \cite[Chapter 4, Section 3]{ethier_markov_1986}. 

\begin{theorem} \label[I]{thm:charactFullGene}
Let $(X_t)_t$ be a strong Markov process with full generator $\widehat{L}$. For $\varphi,\psi\in B_p(E)$, let us define for $t\geq 0$, 
\[M^{\varphi,\psi}_t := \varphi(X_t) - \varphi(X_0) - \int_0^t \psi(X_s)\dd s.\]
Then the full generator is characterised by 
\[\widehat{L} = {\left\{ (\varphi,\psi)\in B_p(E)\times B_p(E),\,\, (M^{\varphi,\psi}_t)_t\,\,\,\text{is a} \,\,\sigma(X)\text{-martingale}\right\}}.\]
\end{theorem}

Taking $\psi = L\varphi$ shows that 
\[ M^\varphi_t = \varphi(X_t) - \varphi(X_0) - \int_0^t L \varphi\left(X_s\right) \dd s \]
is a martingale. The (forward) Kolmogorov equation is retrieved by simply taking the expectation. 

Conversely, a linear operator (bounded or unbounded) $L$ with domain $\mathcal{D}\subset B_p(E)$ is the generator of a strongly continuous contraction semigroup (and thus is the generator of a Markov process) if and only if it satisfies the Hille-Yosida theorem. In the context of Feller processes, the result is stated in \cite[Theorem 3.1]{bottcher_levy_2013} or \cite[Chapter~4, Theorem 2.2]{ethier_markov_1986}. The following example is a classical and important application. Further examples using various points of view (SDE, martingale problem\ldots) can be found in \cite[Chapter 3]{bottcher_levy_2013}. 

\begin{example}\label[I]{example:diffusionprocess} The second order differential operator on $\R^d$,  
\[L\varphi(x) = c(x)\varphi(x)+\sum_{i=1}^ d b_i(x)\frac{\partial\varphi}{\partial x_i} + \frac{1}{2}\sum_{i,j=1}^d a_{ij}(x)\frac{\partial^2\varphi}{\partial x_i\partial x_j},\]
is the generator of a Feller semigroup when $a_{ij}$, $b_i$ and $c$ are respectively $C^3_b(\R^d)$, $C^2_b(\R^d)$ and $C^1_b(\R^d)$ with $c(x)\leq 0$ and the matrix $a=(a_{ij})$ is uniformly elliptic in the sense that there exists $\lambda>0$ such that 
\[\forall x,\xi\in \R^d,\quad \langle a(x)\xi,\xi\rangle \geq \lambda |x|^2.\]
The corresponding Feller process is said to be a diffusion process. 
\end{example}

\subsubsection*{The dual semi-group.} Let $E$ be locally compact and let $(T_t)_t$ be a Feller semigroup. To obtain a forward (or strong) version of the PDE representation \cref{eq:backkolmogorovappendix}, let the dual version of \cref{transitionsemi} be given by: 
\[\forall \nu\in \pb(E),\quad \varphi \in C_0(E)\mapsto \langle S_t\nu,\varphi\rangle  := \langle \nu, T_t\varphi\rangle\in\R.\]
By the Riesz representation theorem, this defines a family of operator $(S_t)_t$ on $\pb(E)$. For all $\nu\in\pb(E)$, $S_t\nu\in \pb(E)$ and by the Chapman-Kolmogorov property, the family of operators $(S_t)_t$ also forms a semigroup. When transition functions $(P_t)_{t\geq0}$ are available, this duality relation implies
\[ S_t \nu(\dd y) = \int_E P_t(x,\dd y)\nu(\dd x). \]
Given the generator $L$ of the semigroup $(T_t)_t$, the dual semigroup $(S_t)_t$ satisfies the dual Kolmogorov equations: 
\[ \partial_t S_t = L^{\star} S_t,\quad \partial_t S_t = S_t L^\star, \]
where $L^\star$ is the dual operator of $L$. This time, the forward Kolmogorov equation (the first one) can be interpreted as the initial value problem $\partial_t f_t = L^\star f_t$.  According to \cref{eq:Markovmarginals}, the solution of the forward equation $f_t = S_t f_0$ is the law at time $t\geq 0$ of the Feller process with initial distribution $f_0$.   

\begin{example} The law $f_t$ of a diffusion process (see Example \cref{example:diffusionprocess}) satisfies the forward Kolmogorov equation, also called the Fokker-Planck equation in this context:
\[\partial f_t(x) = c(x)f_t(x) - \sum_{i=1}^d \partial_{x_i}\{b_i(x)f_t(x)\} + \frac{1}{2}\sum_{i,j=1}^d \partial_{x_i}\partial_{x_j} \{a_{ij}(x) f_t(x)\}.\]
\end{example}

\subsubsection{Time-inhomogeneous Markov processes}\label[I]{sec:timeinhomogeneouslinearmarkov} 

One of the goal of this subsection is to extend this formalism to cases where the generator $L_t$ has a time-dependence. The analog of transition functions is the following time-inhomogeneous version. 

\begin{definition}[Transition functions, inhomogeneous case]
A family of maps $(P_{s,t})_{0 \leq s \leq t}$ where $P_{s,t}:E\to\pb(E)$ is a family of time-inhomogeneous transition functions when the following properties hold.
\begin{itemize}
    \item The map $(s,t,x) \mapsto P_{s,t}(x,\cdot)$ is a measurable map $[0,\infty)\times[0,\infty)\times E\to\pb(E)$. 
    \item For all $t \in \R_+$ and for all $x \in E$, $P_{t,t}(x, \cdot) = \delta_x$.
    \item For all $0 \leq r \leq s \leq t$, for all $x \in E$ and for all $\mathscr{A} \in \mathcal{B}(E)$, \[P_{r,t}(x,A) = \int_E P_{s,t}(y,\mathscr{A}) P_{r,s}(x,\dd y).\]
\end{itemize}
This last relation is the equivalent of the Chapman-Kolmogorov property.
\end{definition}

Similarly to the time-homogeneous case, the transition functions $(P_{s,t})_{0 \leq s \leq t}$ are said to be adapted to the time-inhomogeneous Markov process $(X_t)_{t \geq 0}$ when
\[ \mathbb{P}{\left(X_{t} \in \mathscr{A} |(X_r)_{0 \leq r \leq s} \right)} = P_{s,t}(X_{s}\in\mathscr{A}),\]
or equivalently, 
\[\E {\left[\varphi(X_{t})|(X_r)_{0 \leq r \leq s}\right]} = \int_E \varphi(y) P_{s,t}(X_s,\dd y), \]
for every $0 \leq s \leq t$ and any $\varphi\in B_p(E)$. The time-homogeneous setting is recovered when $P_{s,t}\equiv P_{t-s}$ depends on $t-s$ only. Moreover given the initial distribution $f_0\in \pb(E)$, the finite dimensional distributions can be computed as
\begin{align*} 
&\mathbb{P}(X_0 \in  \mathscr{A}_0, X_{t_1} \in \mathscr{A}_1,X_{t_1}, \ldots, X_{t_n} \in \mathscr{A}_n) \\
&= \int_{\mathscr{A}_0} \ldots \int_{\mathscr{A}_{n-1}} P_{t_{n-1},t_n}(y_{n-1},\mathscr{A}_n) P_{t_{n-2} , t_{n-1}}(y_{n-2},\dd y_{n-1})\\
&\phantom{= \int_{\mathscr{A}_0} \ldots \int_{\mathscr{A}_{n-1}} P_{t_{n-1},t_n}(y_{n-1},\mathscr{A}_n) P_{t_{n-2} , t_{n-1}}(y_{n-2}}\ldots P_{0,t_1}(y_0,\dd y_1)f_0(\dd y_0),
\end{align*}
so that given $f_0$, a time inhomogeneous Markov process is fully characterised by the time-inhomogeneous transition functions. The link with PDEs is retrieved similarly by introducing the operators
\[\forall \varphi\in B_p(E),\quad T_{s,t} \varphi (x) := \int_E \varphi(y) P_{s,t}(x,\dd y). \]
The family $(T_{s,t})_{s\leq t}$ is an evolution system in the following sense, see \cite{bottcher_feller_2014}, \cite{rueschendorf_comparison_2016} and the references therein. 

\begin{definition}[Evolution system]
An evolution system $(T_{s,t})_{0 \leq s \leq t}$ is a family of bounded linear operators on a closed subspace $\mathcal{D} \subset B_p(E)$ such that
\begin{itemize}
    \item for all $ t \geq 0, T_{t,t} = \mathrm{Id}$;
    \item for all $ 0 \leq s \leq t$, $T_{s,t} 1 = 1$;
    \item for all $0 \leq r \leq s \leq t$, $T_{r,t} = T_{r,s} T_{s,t}$.
\end{itemize}
\end{definition}
These properties can be directly checked in the case of a system defined by time inhomogeneous transition functions.

An evolution system is said to correspond to a Markov process $(X_t)_{t \geq 0}$ when
\[ \forall 0 \leq s \leq t,\,\, \forall \varphi \in \mathcal{D}, \quad T_{s,t} \varphi(X_{s}) = \E {\left[\varphi(X_{t})|(X_r)_{0 \leq r \leq s}\right]}. \]
The previous notion of Feller semigroup readily extends to evolution systems in a time-inhomogeneous setting with the strong continuity property being replaced by 
\[\forall \varphi\in C_0(E),\quad \|T_{s,t}\varphi-\varphi\|_\infty\underset{(s,t)\to(0,0)}{\longrightarrow}0.\]

In the time inhomogeneous case, there are two notions of infinitesimal generators, depending on if the derivative is taken from the left of from the right. The left and right generators are defined respectively by 
\[L_t^-\varphi := \lim_{\varepsilon\to0^+}\frac{T_{t-\varepsilon,t}\varphi-\varphi}{\varepsilon},\quad L^+_t \varphi = \lim_{\varepsilon\to0^+} \frac{T_{t,t+\varepsilon} \varphi - \varphi}{\varepsilon}.\]
They are respectively defined on the domains denoted by $\mathcal{D}(L_t^-)$ and $\mathcal{D}(L_t^+)$. Note that in both cases, the generators depend on a time variable. In general, the left and right generators do not coincide but they do under stronger uniform continuity assumptions with respect to the time variable, see \cite[Lemma 2.2]{bottcher_feller_2014}. They also coincide for time homogeneous systems. 

If $(T_{s,t})_{s\leq t}$ is strongly continuous then the forward and backward Kolmogorov equations reads, respectively:
\[\frac{\dd^\pm}{\dd t} T_{s,t} = T_{s,t}L_t^\pm, \quad \frac{\dd^\pm}{\dd s} T_{s,t} = -L_s^\pm T_{s,t},\]
where $\frac{\dd^+}{\dd t}$ (resp. $\frac{\dd^-}{\dd t}$) denotes the right (resp.) left derivative. Given $\varphi\in \mathcal{D}(L_s^-)$ and a fixed $t>0$, the solution of the backward Kolmogorov equation is 
\[u(s,x) = T_{s,t}\varphi(x) = \E[\varphi(X_t)|X_s = x].\]
The terminology \emph{backward} refers to the fact that $u$ satisfies a final value problem with terminal condition $\varphi$ at $s=t$. In the time homogeneous setting,
\[T_{s,t} = T_{0,t-s} \equiv T_{t-s},\]
and consequently for $\tau\leq t$, the quantity
\[ U(\tau,x) := u(t-\tau,x) = T_{\tau}\varphi(x),\]
satisfies
\[\partial_{\tau} U = L U,\quad U(\tau=0) = \varphi,\]
which is the backward Kolmogorov equation previously obtained (since in this case $L^+_s = L^-_s = L$ does not depend on $s$).

As before, the forward equation is better understood with the dual formulation. Given $s\leq t$, the operator $S_{s,t}$ acting on $\pb(E)$ is defined by 
\[\forall \nu\in \pb(E),\,\,\forall \varphi\in B_p(E),\quad \langle S_{s,t}\nu, \varphi\rangle = \langle \nu, T_{s,t}\varphi\rangle.\]
This duality relation implies 
\[S_{s,t}\nu(\dd y) = \int_{E} P_{s,t}(x,\dd y) \nu(\dd x) \]
Let $f_0(x \dd x)$ be the initial distribution. Then
\[f_t(\dd x) := S_{0,s}f_0(\dd x) = \int_E P_{0,t}(y,\dd x) f_0(\dd y),\]
is the law of the associated Markov process at time $t$. In this last equality, note the change of variables which exchanges the roles of the $y$ and $x$ variables. Since the family of operators $(S_{s,t})_{s\leq t}$ satisfies the dual Chapman-Kolmogorov property $S_{s,t}S_{r,s}=S_{r,t}$ for $r\leq s\leq t$, then $S_{s,t}f_s = f_t$ for all $0\leq s\leq t$. Let us assume that the left and right generators coincide and let $L^\star_t$ be the formal adjoint of $L_t=L_t^+$. Then, the forward Kolmogorov equation becomes: 
\[\partial_t S_{s,t} = L_t^\star S_{s,t},\]
which is an initial value problem for the density $f_t$ with initial condition $f_s$ at $s$. 

\begin{example}\label[I]{example:timeinhomgeneousfp} The time inhomogeneous setting allows to consider diffusion processes with time variable coefficients (see Example \cref{example:diffusionprocess}). The backward Kolmogorov equation reads: 
\[-\partial_s u(s,x) = \sum_{i=1}^d b_i(s,x) \partial_{x_i}u(s,x) + \frac{1}{2}\sum_{i,j=1}^d a_{ij}(s,x)\partial_{x_i}\partial_{x_j} u(s,x),\quad u(t,x) = \varphi(x),\]
and the forward Kolmogorov equation (or Fokker-Planck equation) is: 
\[\partial_t f_t(x) = - \sum_{i=1}^d \partial_{x_i}\{b_i(t,x)f_t(x)\} + \frac{1}{2}\sum_{i,j=1}^d \partial_{x_i}\partial_{x_j}\{a_{ij}(t,x)f_t(x)\},\quad f_{t=0}(x) = f_0(x).\]
\end{example}

\subsubsection{Non-linear Markov processes}\label[I]{sec:nonlinearmarkov}

The previous steps have given Markov representations for linear PDEs (with and without time varying coefficients). The class of \emph{nonlinear} PDEs studied in this review can also be seen as the forward Kolmogorov equation associated to a particular class of time-inhomogeneous Markov processes. This point of view was introduced by the seminal work of McKean \cite{mckean_class_1966, mckean_propagation_1969} and Johnson \cite{johnson_class_1968}. 

Let us first introduce an extension of the notion of transition functions.

\begin{definition}[Transition functions, nonlinear case]
A family of maps $(P^{\nu}_t)_{t\geq0}$ from $E\to\pb(E)$ defined for any $\nu\in\pb(E)$ is a family of non-linear transition functions when it satisfies the following properties.  
\begin{itemize}
    \item The map $(t,x,\nu) \mapsto P^{\nu}_t(x,\cdot)$ is a measurable map $[0,\infty) \times E \times \pb(E) \rightarrow \pb(E)$.
    \item For all $x \in E$ and for all  $\nu \in \pb(E)$, $P^{\nu}_0(x, \cdot) = \delta_x$.
    \item For all $s,t \geq 0$, for all $x \in E$, for all $\nu \in \pb\left(E\right)$ and for all $\mathscr{A} \in \mathcal{B}(E)$, it holds that \[P^{\nu}_{s+t}(x,\mathscr{A}) = \int_E P^{\int_E P^{\nu}_{t}(x,\dd y) \nu(\dd x)}_{s}(y,\mathscr{A}) P^{\nu}_{t}(x,\dd y).\]
\end{itemize}
This last relation is a nonlinear version of the Chapman-Kolmogorov property; the linear case is recovered when $P^{\nu}_t$ does not depend on $\nu$.
\end{definition}

At this point, it is natural to introduce the nonlinear operator $\overline{S}_t:\pb(E)\to\pb(E)$ defined for $\nu\in\pb(E)$ by 
\[ \overline{S}_t (\nu)(\dd y) = \int_E P^{\nu}_{t}(x,\dd y)\nu(\dd x), \]
so that the non-linear Chapman-Kolmogorov relation reads
\begin{equation} 
P^{\nu}_{s+t}(x,\cdot) = \int_E P^{\overline{S}_t (\nu)}_{s}(y,\cdot) P^{\nu}_{t}(x,\dd y).
\end{equation}
Integrating against $\nu$, this gives
\begin{equation} 
\overline{S}_{t+s} (\nu) = \int_E P^{\overline{S}_t (\nu)}_{s} (y,\cdot) {\left(\int_E P^{\nu}_{t}(x,\dd y) \nu(\dd x)\right)} =  \overline{S}_{s} (\overline{S}_{t}(\nu)). 
\end{equation}
In particular $\overline{S}_0 = \mathrm{Id}$, and $(\overline{S}_{t})_{t \geq 0}$ appears to be a \emph{non-linear semi-group}: the bar notation reminds of the non-linearity when writing $\overline{S}_{t}$ alone. The semigroup $(\overline{S}_{t})_{t \geq 0}$ is the analog of the previous dual semi-group.

Finally a nonlinear Markov process in the sense of McKean with initial distribution $f_0\in\pb(E)$ is a time-inhomogeneous Markov process with transition functions of the form
\[\overline{P}_{s,t} = P^{\overline{S}_s(f_0)}_{t-s},\]
for $s\leq t$. This corresponds to the construction of Johnson \cite{johnson_class_1968} established when the state space $E$ is the 2-state space $E=\{-1,+1\}$. 

The nonlinear evolution system is defined for any $s\leq t$ by 
\[\forall \varphi\in B_p(E),\quad \overline{T}_{s,t}\varphi(x) = \int_{E} \varphi(y)\overline{P}_{s,t}(x,\dd y) = \int_E \varphi(y) P_{t-s}^{f_s}(x,\dd y).\] 
In the above expression, $f_s:=\overline{S}_s(f_0)$ is the law at time $s$ of the associated nonlinear Markov process in the sense of Mckean.  Its right generator is given by: 
\[ L_{f_t}\varphi(x) := \lim_{\varepsilon\to0} \frac{\overline{T}_{t,t+\varepsilon}\varphi(x) - \varphi(x)}{\varepsilon} = \lim_{\varepsilon\to0} \frac{\int_E \varphi(y) P^{f_t}_\varepsilon(x,\dd y) - \varphi(x)}{\varepsilon},\]
from which it can be seen that it depends on $f_t$ only. With this particular form for the dependence in time, the forward Kolmogorov equation given in Example \cref{example:timeinhomgeneousfp} thus appears to be the nonlinear Fokker-Planck equation satisfied by the law of the solution of the McKean-Vlasov diffusion SDE 
\[ \dd \overline{X}_t = b(\overline{X}_t,f_t)\dd t + \sigma(\overline{X}_t,f_t)\dd B_t, \quad \overline{X}_t \sim f_t,\]
the wellposedness of which has been studied in Proposition \cref{prop:wellposednessmckeanlipschitz}. 

\subsection{Large Deviation Principles and Sanov theorem}\label[I]{appendix:largedeviations}

\begin{definition}[Large Deviation Principle]
Given a sequence sequence $(a_N)_N$ of positive numbers $a_N \rightarrow 0$ and a non-negative lower-semicontinuous function $I$ on $E$, a sequence $( \mu_N )_N$ in $\mathcal{P} ( E )$ satisfies a Large Deviation Principle (LDP) with speed $a_N$ and rate function $I$ when for any Borel set $A\subset E$, it holds that
\[- \inf_{\mathring{A}} I \leq \liminf_{N \rightarrow \infty} a_N \log \mu_N ( A ) \leq \limsup_{N \rightarrow \infty} a_N \log \mu_N ( A ) \leq - \inf_{\overline{A}} I,\]
where $\mathring{A}$ and $\overline{A}$ denote respectively the interior and closure of $A$. 
\end{definition}

Given a sequence $(X^i)_{i\in\N}$ of i.i.d. real-valued random variables, Cram\'er's theorem states that the sequence of the laws $\mathrm{Law}(\frac{1}{N}\sum_{i=1}^N X^i)$, $N\in\N$, satisfies a LDP with rate function $\Lambda^*(x) = \sup_{t\in\R}(tx-\log\E[\exp(tX^1)])$ and speed $a_N = 1/N$. By taking the image random variables $\varphi(X^1),\varphi(X^2),\ldots$ for a fixed test function $\varphi$, a LDP can be obtained for the sequence of laws $\mathrm{Law}(\langle \mu_{\mathcal{X}^N},\varphi\rangle)$, where $\mathcal{X}^N=(X^1,\ldots,X^N)$. However, in this case, the rate function depends on the choice of the test function $\varphi$. A more precise theorem which gives a LDP for the laws of the sequence of empirical measures is Sanov theorem. 

\begin{theorem}[Sanov]
Let $\mu$ be a probability measure on a Polish space $E$, and let $( X^i )_{i\in\N}$ be a sequence of independent $\mu$-distributed random variables. For $N\in\N$, we recall the notation $\mathcal{X}^N = (X^1,\ldots,X^N)\in E^N$. Then the laws in $\mathcal{P} ( \mathcal{P} ( E ) )$ of the measure-valued random variables $\mu_{\mathcal{X}^N}$ satisfies a large deviation principle with speed $N^{-1}$ and rate function the relative entropy $\nu \mapsto H ( \nu | \mu )$.
\end{theorem}

\subsection{Girsanov transform}\label[I]{appendix:girsanov}

There are many versions of Girsanov theorem. We do not give the most general result (which can be found for instance in \cite[Theorem~5.22]{le_gall_brownian_2016}) but only the one which will is used in this review and which can be found in \cite[Chapter 3, Theorem 5.1]{karatzas_brownian_1998}. 

\begin{theorem}[Girsanov]
Let $(\Omega,\mathscr{F}, (\mathscr{F}_t)_t,\mathbb{P})$ be a filtered probability space and let $(B_t)_t$ be a $d$-dimensional Brownian motion on this space with $\mathbb{P}(B_0=0)=1$. Let $(X_t)_t$ be a $\R^d$-valued adapted measurable process and let the process be defined (whenever it exists) for $t<+\infty$ by: 
\[H_t := \int_0^t X_s\cdot\dd B_s - \frac{1}{2}\int_0^t |X_s|^2\dd s.\]
Let $\mathbb{Q}$ the probability measure on $(\Omega,\mathscr{F})$ defined by its Radon-Nikodym derivative on each $\mathscr{F}_T$, $T<+\infty$ :
\[\frac{\dd\mathbb{Q}}{\dd\mathbb{P}}\Big|_{\mathscr{F}_T} = \exp(H_T).\]
Assume that $\exp(H)$ is a martingale and let us define the process 
\[\widetilde{B}_t = B_t - \int_0^t X_s\dd s.\]
Then for each fixed $T\in[0,+\infty)$, $(\widetilde{B}_t)_{t\leq T}$ is a Brownian motion on $(\Omega,\mathscr{F}_T,\mathbb{Q}|_{\mathscr{F}_T})$.
\end{theorem}

\subsection{Poisson random measures}\label[I]{appendix:poisson}

This section briefly explains how to model jump processes using Poisson random measures. The theory of random measures is explained in great details in \cite{jacod_limit_2003}. Another classical reference on stochastic integration with respect to random measures is \cite{ikeda_stochastic_1989}. The following presentation is also inspired by \cite[Appendix A]{bansaye_stochastic_2015} and \cite[Section 3]{murata_propagation_1977}.  

Let us fix a filtered probability space $(\Omega,\mathscr{F},(\mathscr{F}_t)_t,\mathbb{P})$. 

\begin{definition}[Poisson random measure] Let $(\mathscr{E},\mu)$ be a measurable Polish space endowed with a $\sigma$-finite measure $\mu$. Let $\widehat{\mathcal{M}}(\mathscr{E})$ be the set of all measures $\lambda$ on $\mathscr{E}$ which are expressed as a countable sum of Dirac measures on $\mathscr{E}$ and such that $\lambda(\mathscr{A})<+\infty$ for any $\mu$-finite set $\mathscr{A}$. A Poisson random measure with intensity $\mu$ is a mapping  $\mathcal{N}:\Omega\to \widehat{\mathcal{M}}(\mathscr{E})$ with the following properties. 
\begin{enumerate}[(i)]
    \item The mapping $\omega\in\Omega\mapsto\mathcal{N}(\omega,\mathscr{A})$ is measurable for any $\mu$-finite set $\mathscr{A}$. 
    \item For every disjoints $\mu$-finite sets $\mathscr{A}_1,\ldots,\mathscr{A}_k$, the random variables $\mathcal{N}(\mathscr{A}_j)$, $j\in\{1,\ldots,k\}$, are independent and $\mathcal{N}(\mathscr{A}_j)$ follows a Poisson distribution on $\N$ with parameter $\mu(\mathscr{A}_j)$. 
\end{enumerate}
\end{definition}

We will only consider the case $\mathscr{E}=\R_+\times \Theta$ where $\Theta$ is a Polish space and $\mu$ is of the form $\mu(\dd t, \dd x)=\dd t\otimes \nu(\dd \theta)$. The results below also extend to the case where $\nu$ is replaced by a family of $\sigma$-finite measures $(\nu_t)_t$ on $\Theta$ which depend on the time parameter. The Poisson random measure $\mathcal{N}$ is assumed to be \emph{adapted} which means that it satisfies the following properties. 
\begin{enumerate}[(i)]
\item $\mathcal{N}(\mathscr{A})$ is $\mathscr{F}_t$-measurable for each Borel measurable set $\mathscr{A}\in\mathscr{B}([0,t]\times \Theta)$ with $t>0$ 
\item The $\sigma$-field generated by $\{\mathcal{N}(\mathscr{A}),\,\,\mathscr{A}\in\mathscr{B}((t,+\infty)\times \Theta)\}$ is independent of $\mathscr{F}_t$. 
\end{enumerate}

Some ideas on the construction of the stochastic integral against $\mathcal{N}$ are gathered below. 

To better understand what a Poisson random measure does, it is useful to consider first the case $\nu(\Theta)<+\infty$. In this case, on every finite time interval $[0,T]$, $\mathcal{N}_t(\Theta):=\mathcal{N}((0,t]\times \Theta)$ defines a classical Poisson process. The Poisson random measure $\mathcal{N}$ can be shown to admit the representation: 
\[\mathcal{N}(\dd t, \dd \theta) = \sum_{n=1}^{\gamma} \delta_{(T_n,\theta_n)}(\dd t,\dd\theta),\]
where $T_1,\ldots,T_\gamma$ are the jump times of $\mathcal{N}_t(\Theta)$ and $\theta_n$ are i.i.d. random variables with distribution $\nu(\dd \theta)/\nu(\Theta)$. For any measurable function $a\equiv a(\omega,t,\theta)$ on $\Omega\times\R_+\times \Theta$ with values in $\R$, the integral with respect to $\mathcal{N}$ is defined by: 
\[\int_0^T\int_\Theta a(\omega,s,\theta)\mathcal{N}(\dd s, \dd \theta) := \sum_{n=1}^\gamma a(\omega,T_n,\theta_n).\]
That is to say, it is the sum of the random amplitudes $a(\omega,T_n,\theta_n)$ added at each jumping time $T_n$. 

To extend the previous construction to the case $\nu(\Theta)=+\infty$, let us consider a predictable real-valued function $a\equiv a(\omega,t,\theta)$ on $\Omega\times\R_+\times \Theta$. We do not write the dependence in $\omega$ in the following. We recall that $(\Theta,\nu(\dd\theta))$ is $\sigma$-finite, so there exists an increasing sequence of subsets $(\Theta_p)_p$ such that $\nu(\Theta_p)<+\infty$ and $\Theta = \cup_p \Theta_p$. By the previous construction, the integral of $a$ against $\mathcal{N}$ is well-defined on each subset $[0,T]\times\Theta_p$. There are two cases to distinguish. 

\begin{enumerate}
    \item When $a$ satisfies the $L^1$ condition
    \[\E{\left[\int_0^t\int_\Theta|a(s,\theta)|\nu(\dd \theta)\dd s\right]}<+\infty,\]
    it is possible to show that the sequence ${\left(\int_0^T\int_{\Theta_p}a(s,\theta)\mathcal{N}(\dd s, \dd\theta)\right)_p}$ is Cauchy in $L^1$. Its limit is denoted by $\int_0^T\int_\Theta a(s,\theta)\mathcal{N}(\dd s, \dd\theta)$. In this case, the process
    \begin{equation}\label[Ieq]{eq:martingalePoissonmeasure} M_t = \int_0^t\int_\Theta a(s,\theta)\mathcal{N}(\dd s,\dd \theta)-\int_0^t\int_\Theta a(s,\theta)\nu(\dd \theta)\dd s,\end{equation}
    is a $\mathscr{F}_t$-martingale (in fact, $M_t$  characterises $\mathcal{N}$). 
    \item When $a$ satisfies the $L^2$ condition 
    \[\E{\left[\int_0^T\int_\Theta|a(s,\theta)|^2\nu(\dd \theta)\dd s\right]}<+\infty,\]
    then it is possible to prove that \cref{eq:martingalePoissonmeasure} still defines a square integrable martingale with quadratic variation 
    \[\langle M\rangle_t = \int_0^t\int_\Theta |a(s,\theta)|^2 \widetilde{\mathcal{N}}(\dd s,\dd\theta),\]
    where 
    \[\widetilde{\mathcal{N}}(\dd t,\dd\theta) := \mathcal{N}(\dd t,\dd\theta)-\nu(\dd\theta)\dd t,\]
    is called the compensated measure of $\mathcal{N}$. Note however that the quantity  $\int_0^T\int_\Theta a(s,\theta)\mathcal{N}(\dd s, \dd\theta)$ may not be defined. 
\end{enumerate}

The next step is to make sense of the jump-diffusion SDE: 
\begin{multline}\label[Ieq]{eq:SDEpoisson}
    X_t = X_0 + \int_0^t b(X_s)\dd s + \int_0^t \sigma(X_s)\dd B_s + \int_0^t\int_\Theta \alpha(X_{s^-},\theta)\mathcal{N}(\dd s,\dd\theta)\\ 
    + \int_0^t\int_\Theta \widetilde{\alpha}(X_{s^-},\theta)\widetilde{\mathcal{N}}(\dd s,\dd\theta), 
\end{multline}  
where this time $\alpha,\widetilde{\alpha}:E\times \Theta\to E$ for $E=\R^d$. We always assume the following Lipschitz integrability conditions: 
\begin{enumerate}[(i)]
    \item For all $x\in E$ and $T>0$, it holds that $\int_0^T\int_\Theta |\alpha(x,\theta)|\nu(\dd\theta)<+\infty $ and $\int_0^T\int_\Theta |\widetilde{\alpha}(x,\theta)|^2\nu(\dd\theta)<+\infty $. 
    \item There exists $C>0$ such that for any $x,y\in E$, 
    \[|\sigma(x)-\sigma(y)|^2 + |b(x)-b(y)|^2 + \int_\Theta |\widetilde{\alpha}(x,\theta)-\widetilde{\alpha}(y,\theta)|^2\nu(\dd\theta)\leq C|x-y|^2.\]
\end{enumerate}
In the classical theory of SDE (see \cite[Chapter IV, Section 9]{ikeda_stochastic_1989}) the Lipschitz integrability condition 
\[\int_\Theta |\alpha(x,\theta)-\alpha(y,\theta)|^p\nu(\dd\theta)\leq C|x-y|^p\]
with $p=2$ is also assumed. However, from a modelling point of view, it makes more sense in the context of this review to assume this condition with $p=1$ (see \cite[Remark 2.1]{andreis_mckeanvlasov_2018} or the introduction of \cite{graham_mckean-vlasov_1992}). In this $L^1$ setting, strong existence and uniqueness for the SDE \cref{eq:SDEpoisson} is proved in \cite[Theorem 1.2]{graham_mckean-vlasov_1992}. Moreover, the generator of the process is the sum of the three generators
\begin{align*}
\mathcal{L}\varphi(x) &= \sum_{i=1}^d b(x)\cdot\nabla\varphi(x)+\frac{1}{2}\sum_{i,j=1}^d (\sigma\sigma^\mathrm{T})_{ij}(x)\partial_{x_i}\partial_{x_j}\varphi(x),\\
\mathcal{J}\varphi(x) &= \int_\Theta \big\{\varphi\big(x+\alpha(x,\theta)\big)-\varphi(x)\big\}\nu(\dd\theta),\\
\widetilde{\mathcal{J}}\varphi(x) &= \int_\Theta \big\{\varphi\big(x+\widetilde{\alpha}(x,\theta)\big)-\varphi(x)-\widetilde{\alpha}(x,\theta)\cdot\nabla\varphi(x)\big\}\nu(\dd\theta). 
\end{align*}

\section{A strengthened version of Hewitt-Savage theorem and its partial system version}\label[I]{appendix:hewittsavage}

The following computation gives the Cauchy-estimate of Theorem \cref{thm:quantitativehewitt} using stronger metrics.

\begin{corollary}[Extending Cauchy-estimates to other metrics] \label[I]{CauchyStrongHewit} 
The Cauchy-estimate of Theorem \cref{thm:quantitativehewitt} can be obtained in $W_{\delta,p}$-distance for $M \leq N$ with rate $\varepsilon ( M )$, where $\delta$ is any metric on $\pb ( E )$ such that a uniform $\delta^p$-law of large numbers holds, i.e. for every $N \geq 1$, there exists $\varepsilon ( N ) > 0$ with $\varepsilon ( N ) \to 0$ as $N \to \infty$, such that for every $f$ in $\pb_p ( E )$ and every $f^{\otimes N}$-distributed vector $\overline{\mathcal{X}}^N$, 
\[ \E [ \delta^p ( \mu_{\overline{\mathcal{X}}^N} , f ) ]^{1/p} \leq \varepsilon(N) . \] 
Thanks to \cite{fournier_rate_2015}, this includes all the Wassertein-$p$ metrics on a compact space $E$ or on any $E$ endowed with a bounded distance. 
\end{corollary}

\begin{proof}
In order to couple $\mathcal{X}^M$ to the sub-vector $\mathcal{X}^{M,N}$ of $\mathcal{X}^N$, we choose the transference plane $( \boldsymbol{\mu}_{M,N} , \boldsymbol{\mu}_N )_{\#} \pi^{N} \otimes \pi^N$: this is well defined thanks to the compatibility property for $\pi^{M,N}$ and $\pi^N$. Thus
\[ W_{\delta,p}^p \big( \mathrm{Law} ( \mu_{\mathcal{X}^M} ) , \mathrm{Law} ( \mu_{\mathcal{X}^N} ) \big) \leq \E_{\mathcal{X}^N \sim \pi^N} \big[ \delta^p ( \mu_{\mathcal{X}^{N,M}} , \mu_{\mathcal{X}^N} ) \big] = \big\langle \pi^N , \delta^p \big( \boldsymbol{\mu}_{M,N} , \boldsymbol{\mu}_N \big) \big\rangle.\]
Theorem \cref{thm:quantitativehewitt} showed that $\pi^N$ is the $N$-th moment measure of the limit $\pi$ in $\mathcal{W}_{H^{-s}}$-distance, allowing to take advantage of superposing i.i.d. states:
\begin{align*}
\big\langle \pi^N , \delta^p \big( \boldsymbol{\mu}_{M,N} , \boldsymbol{\mu}_N \big) \big\rangle &= \int_{\pb ( E )} \langle \nu^{\otimes N} , \delta^p ( \boldsymbol{\mu}_{M,N} , \boldsymbol{\mu}_N ) \rangle \mathrm{Law}(\mu_{\mathcal{X}^N}) ( \dd \nu ) \\
&\leq C(p) \int_{\pb ( E )} \langle \nu^{\otimes N} , \delta^p ( \boldsymbol{\mu}_{M,N} , \nu ) \rangle \mathrm{Law}(\mu_{\mathcal{X}^N}) ( \dd \nu ) \\
&\phantom{abcd} + C(p) \int_{\pb ( E )} \langle \nu^{\otimes N} , \delta^p ( \nu , \boldsymbol{\mu}_N ) \rangle \mathrm{Law}(\mu_{\mathcal{X}^N}) ( \dd \nu ), 
\end{align*}
for a constant $C(p)$ which only depends on $p$. Since $\nu^{\otimes N}$ is the law of a vector of i.i.d. particles
\[ \langle \nu^{\otimes N} , \delta^p ( \boldsymbol{\mu}_{M,N} , \nu ) \rangle = \langle \nu^{\otimes M} , \delta^p ( \boldsymbol{\mu}_M , \nu ) \rangle = \E_{\overline{\mathcal{X}}^{M} \sim \nu^{\otimes M}} \big[ \delta^p \big( \mu_{\overline{\mathcal{X}}^{M}} , \nu \big) \big], \]
and the quantitative assumption on the $\delta^p$-law of large numbers concludes. 
\end{proof}

The following proposition gives a useful Cauchy-estimate for the empirical measure of a sub-system of a finite exchangeable particle system. More precisely, let $M < N$ and let $\mathcal{X}^N\sim f^N$ be a finite exchangeable particle system. One wish to compare $\mu_{\mathcal{X}^N}$ and $\mu_{\mathcal{X}^{M,N}}$. This is not possible by a direct coupling argument since the two empirical measures do not have the same size. We thus use an alternative argument based on the special polynomial structure of the $H^{-s}$ norm (see Lemma~\cref{lemma:sobolevpolynomial}). This argument is also used in the proof of the Hewitt-Savage Theorem \cref{thm:quantitativehewitt}.

\begin{proposition}[Block empirical measures approximation] \label[I]{mckeanG} Let $s>d/2$. 
Let $M < N$ and let $\mathcal{X}^N\sim f^N$ be a finite exchangeable particle system. It holds that
\[ \mathcal{W}_{H^{-s}} \big( \mathrm{Law}(\mu_{\mathcal{X}^{M,N} }) , \mathrm{Law}(\mu_{\mathcal{X}^M }) \big) \leq 2 \| \Phi_s \|_{\infty} {\left( \frac{1}{M} - \frac{1}{N} \right)}, \]
where $\Phi_s(z):=\int_{\R^d} \e^{-iz\cdot\xi}(1+|\xi|^2)^{-s}\dd\xi$.
\end{proposition}

\begin{proof}
Thanks to the identity \cref{eq:polynomialsobolev}, we get: 
\begin{align*}
&\E {\left[ \| \mu_{\mathcal{X}^{M,N}} - \mu_{\mathcal{X}^N} \|^2_{H^{-s}} \right]} \\
&= \int_{(\R^{d})^N} {\left( \int_{\R^d\times\R^d} \Phi_s(x-y) {\left[ \mu_{\mathbf{x}^{M,N}}^{\otimes 2} - \mu_{\mathbf{x}^{M,N}} \otimes \mu_{\mathbf{x}^N} \right]} (\dd x,\dd y) \right)} f^N \big( \dd \mathbf{x}^N \big) \\
&\quad + \int_{(\R^{d})^N} {\left( \int_{\R^d\times\R^d} \Phi_s(x-y) {\left[ \mu_{\mathbf{x}^N}^{\otimes 2} - \mu_{\mathbf{x}^N} \otimes \mu_{\mathbf{x}^{M,N}} \right]} (\dd x,\dd y) \right)} f^N \big( \dd \mathbf{x}^N \big) \\
&= \int_{(\R^{d})^N} {\left( \frac{1}{M^2} \sum_{k,\ell=1}^M  \Phi_s ( x^k - x^\ell ) - \frac{1}{NM} \sum_{k=1}^M\sum_{j=1}^N  \Phi_s ( x^k - x^j ) \right)} f^N \big( \dd \mathbf{x}^N \big) \\
&\quad + \int_{(\R^{d})^N} {\left( \frac{1}{N^2} \sum_{i,j=1}^N  \Phi_s ( x^i -x^j ) - \frac{1}{NM} \sum_{i=1}^N\sum_{\ell=1}^M  \Phi_s ( x^i -x^\ell ) \right)} f^N \big( \dd \mathbf{x}^N \big) \\
&= {\left(\frac{M}{M^2} - \frac{M}{NM} \right)} \Phi_s ( 0) \\
&\quad+ {\left( \frac{M^2 - M}{M^2} - \frac{NM - M}{NM} \right)} \int_{\R^d\times\R^d} \Phi_s ( x - y ) f^{2,N} ( \dd x , \dd y ) \\
&\quad + {\left(\frac{N}{N^2} - \frac{M}{NM} \right)} \Phi_s ( 0) \\
&\quad+ {\left( \frac{N^2 - N}{N^2} - \frac{NM - M}{NM} \right)} \int_{\R^d\times\R^d} \Phi_s ( x - y ) f^{2,N} ( \dd x , \dd y ) \\
&= {\left(\frac{1}{M} - \frac{1}{N} \right)}{\left( \Phi_s ( 0) - \int_{\R^d\times\R^d} \Phi_s ( x - y ) f^{2,N} ( \dd x , \dd y ) \right)}.
\end{align*}
Since $\Phi_s$ is bounded, this gives the desired estimate.
\end{proof}

In order to compare the empirical measures of two sub-systems with different sizes but which come from the same exchangeable finite particle system, one can use the same method as in Corollary \cref{CauchyStrongHewit}. However, we need to replace the marginals $f^{k,N}$ by their moment measures approximation. Thanks to Lemma \cref{lemma:grunbaumlemma}, this error term is quantitative which gives an explicit maximal size for the subsystems.

\begin{proposition}[Strong Cauchy-estimates for block empirical measures] Under the assumptions of Corollary \cref{CauchyStrongHewit}, let $M<N$ and let $k<\ell$ such that $\varepsilon(\ell) \geq \ell^2 N^{-1}$. Then it holds that
\[ W^p_{\delta,p} \big( \mathrm{Law}(\mu_{\mathcal{X}^{k,N} }) , \mathrm{Law}(\mu_{\mathcal{X}^{\ell,N} }) \big) \leq \varepsilon(k). \]
\end{proposition}

\begin{proof}
The transference plane is the same as in Corollary \cref{CauchyStrongHewit}. Let us consider for simplicity the Wassertein-2 distance (other adaptations are straightforward). Assuming the distance on $E=\R^d$ to be bounded, the condition on $\ell$ allows to write
\begin{align*}
&\E {\left[ W_2^2 ( \mu_{\mathcal{X}^{k,N}} , \mu_{\mathcal{X}^{\ell,N}} ) \right]} = \big\langle f^{\ell,N} , W_2^2 \big( \boldsymbol{\mu}_{k,\ell} , \boldsymbol{\mu}_{\ell} \big) \big\rangle \\
&= \big\langle \E_{\mathcal{X}^N} \mu_{\mathcal{X}^N}^{\otimes \ell} , W_2^2 \big( \boldsymbol{\mu}_{k,\ell} , \boldsymbol{\mu}_{\ell} \big) \big\rangle + \mathcal{O} \big( \ell^2 N^{-1} \big), 
\end{align*}
where $\mu_{\mathcal{X}^N}^{\otimes \ell}$ is the $\ell$-th moment measure of $\mathrm{Law}(\mu_{\mathcal{X}^N})$. This leads to
\begin{align*}
&\big\langle \E_{\mathcal{X}^N} \mu_{\mathcal{X}^N}^{\otimes \ell} , W_2^2 \big( \boldsymbol{\mu}_{k,\ell} , \boldsymbol{\mu}_{\ell} \big) \big\rangle \\
&= \int_{\pb ( E )} \big\langle \nu^{\otimes \ell} , W_2^2 \big( \boldsymbol{\mu}_{k,\ell} , \boldsymbol{\mu}_{\ell} \big) \big\rangle \mathrm{Law}(\mu_{\mathcal{X}^N_t}) ( \dd \nu ) \\
&\leq \int_{\pb ( E )} \big\langle \nu^{\otimes \ell} , W_2^2 \big( \boldsymbol{\mu}_{k,\ell} , \nu \big) \big\rangle \mathrm{Law}(\mu_{\mathcal{X}^N}) ( \dd \nu ) \\
&\quad+ \int_{\pb ( E )} \big\langle \nu^{\otimes N} , W_2^2 ( \nu , \boldsymbol{\mu}_{\ell} ) \big\rangle \mathrm{Law}(\mu_{\mathcal{X}^N}) ( \dd \nu ). 
\end{align*}
Since $\nu^{\otimes \ell}$ is the law of a vector of i.i.d. particles
\begin{align*} \big\langle \nu^{\otimes \ell} , W_2^2 \big( \boldsymbol{\mu}_{k,\ell} , \nu \big) \big\rangle &= \big\langle \nu^{\otimes k} , W_2^2 \big( \boldsymbol{\mu}_{k} , \nu \big) \big\rangle \\&= \E_{\overline{\mathcal{X}}^{k} \sim \nu^{\otimes k}} W_2^2 \big( \mu_{\overline{\mathcal{X}}^{k}} , \nu \big), \end{align*}
and the quantitative law of large numbers for the Wassertein-2 distance in \cite{fournier_rate_2015} concludes.
\end{proof}

\begin{remark}[Recovering Cauchy-estimates on finite marginals]
If the Wassertein-1 distance is considered instead, this leads to Cauchy estimates on $f^{\ell,N}$ thanks to Proposition \cref{prop:W1measure}
\begin{equation*} W_1 ( f^{k,M} , f^{\ell,N} ) = \mathcal{W}_1 \big( \mathrm{Law}(\mu_{\mathcal{X}^{k,M}}) , \mathrm{Law}(\mu_{\mathcal{X}^{\ell,N}}) \big)  \leq \E \big[ W_1 ( \mu_{\mathcal{X}^{k,M}} , \mu_{\mathcal{X}^{\ell,N}}) \big]. \end{equation*}
\end{remark}

\providecommand{\href}[2]{#2}
\providecommand{\arxiv}[1]{\href{http://arxiv.org/abs/#1}{arXiv:#1}}
\providecommand{\url}[1]{\texttt{#1}}
\providecommand{\urlprefix}{}


\begin{thebibliography}{100}

\bibitem{acebron_kuramoto_2005}
\newblock J.~A. Acebr{\'o}n, L.~L. Bonilla, C.~J. P{\'e}rez~Vicente, F.~Ritort
  and R.~Spigler,
\newblock The {Kuramoto} model: {A} simple paradigm for synchronization
  phenomena,
\newblock \emph{Rev. Modern Phys.}, \textbf{77} (2005), 137--185,
\newblock \urlprefix\url{https://link.aps.org/doi/10.1103/RevModPhys.77.137}.

\bibitem{albi_vehicular_2019}
\newblock G.~Albi, N.~Bellomo, L.~Fermo, S.-Y. Ha, J.~Kim, L.~Pareschi,
  D.~Poyato and J.~Soler,
\newblock Vehicular traffic, crowds, and swarms: {From} kinetic theory and
  multiscale methods to applications and research perspectives,
\newblock \emph{Math. Models Methods Appl. Sci.}, \textbf{29} (2019),
  1901--2005,
\newblock
  \urlprefix\url{https://www.worldscientific.com/doi/abs/10.1142/S0218202519500374}.

\bibitem{ambrosio_gradient_2008}
\newblock L.~Ambrosio, N.~Gigli and G.~Savar{\'e},
\newblock \emph{Gradient {Flows} in {Metric} {Spaces} and in the {Space} of
  {Probability} {Measures}},
\newblock 2nd edition,
\newblock Lectures in {Mathematics} {ETH} {Z{\"u}rich}, Birkh{\"a}user, Basel,
  2008.

\bibitem{andreis_mckeanvlasov_2018}
\newblock L.~Andreis, P.~Dai~Pra and M.~Fischer,
\newblock {McKean}{\textendash}{Vlasov} limit for interacting systems with
  simultaneous jumps,
\newblock \emph{Stoch. Anal. Appl.}, \textbf{36} (2018), 960--995.

\bibitem{bakry_hypercontractivite_1994}
\newblock D.~Bakry,
\newblock L'hypercontractivit{\'e} et son utilisation en th{\'e}orie des
  semigroupes,
\newblock in \emph{Lectures on {Probability} {Theory}: {Ecole} d'{Et{\'e}} de
  {Probabilit{\'e}s} de {Saint}-{Flour} {XXII}-1992} (eds. D.~Bakry, R.~D.
  Gill, S.~A. Molchanov and P.~Bernard),
\newblock Lecture {Notes} in {Mathematics}, Springer Berlin Heidelberg, 1994,
\newblock 1--114,
\newblock \urlprefix\url{https://doi.org/10.1007/BFb0073872}.

\bibitem{azema_diffusions_1985}
\newblock D.~Bakry and M.~{\'E}mery,
\newblock Diffusions hypercontractives,
\newblock in \emph{S{\'e}minaire de {Probabilit{\'e}s} {XIX} 1983/84} (eds.
  J.~Az{\'e}ma and M.~Yor),
\newblock no. 1123 in Lecture {Notes} in {Mathematics}, Springer Berlin
  Heidelberg, 1985,
\newblock 177--206,
\newblock \urlprefix\url{http://link.springer.com/10.1007/BFb0075847},
\newblock Series Title: Lecture Notes in Mathematics.

\bibitem{bansaye_stochastic_2015}
\newblock V.~Bansaye and S.~M{\'e}l{\'e}ard,
\newblock \emph{Stochastic {Models} for {Structured} {Populations}},
\newblock no. 1.4 in Mathematical {Biosciences} {Institute} {Lecture} {Series},
  Springer International Publishing Switzerland, 2015.

\bibitem{bellomo_active_2017}
\newblock Nicola Bellomo, Pierre Degond and Eitan Tadmor (eds.),
\newblock \emph{Active {Particles}, {Volume} 1: {Advances} in {Theory},
  {Models}, and {Applications}},
\newblock Modeling and {Simulation} in {Science}, {Engineering} and
  {Technology}, Springer International Publishing, 2017,
\newblock \urlprefix\url{http://link.springer.com/10.1007/978-3-319-49996-3}.

\bibitem{bellomo_active_2019}
\newblock Nicola Bellomo, Pierre Degond and Eitan Tadmor (eds.),
\newblock \emph{Active {Particles}, {Volume} 2: {Advances} in {Theory},
  {Models}, and {Applications}},
\newblock Modeling and {Simulation} in {Science}, {Engineering} and
  {Technology}, Springer International Publishing, 2019,
\newblock \urlprefix\url{http://link.springer.com/10.1007/978-3-030-20297-2}.

\bibitem{ben_arous_methode_1990}
\newblock G.~Ben~Arous and M.~Brunaud,
\newblock Methode de {Laplace}: {\'e}tude variationnelle des fluctuations de
  diffusions de type ``champ moyen",
\newblock \emph{Stochastics and Stochastics Reports}, \textbf{31} (1990),
  79--144.

\bibitem{ben_arous_increasing_1999}
\newblock G.~Ben~Arous and O.~Zeitouni,
\newblock Increasing propagation of chaos for mean field models,
\newblock \emph{Ann. Inst. Henri Poincar{\'e} Probab. Stat.}, \textbf{35}
  (1999), 85--102.

\bibitem{berlyand_continuum_2019}
\newblock L.~Berlyand, R.~Creese, P.-E. Jabin and M.~Potomkin,
\newblock Continuum {Approximations} to {Systems} of {Correlated} {Interacting}
  {Particles},
\newblock \emph{J. Stat. Phys.}, \textbf{174} (2019), 808--829,
\newblock \urlprefix\url{http://link.springer.com/10.1007/s10955-018-2205-8}.

\bibitem{bertini_dynamical_2009}
\newblock L.~Bertini, G.~Giacomin and K.~Pakdaman,
\newblock Dynamical {Aspects} of {Mean} {Field} {Plane} {Rotators} and the
  {Kuramoto} {Model},
\newblock \emph{J. Stat. Phys.}, \textbf{138} (2009), 270--290.

\bibitem{bertini_synchronization_2014}
\newblock L.~Bertini, G.~Giacomin and C.~Poquet,
\newblock Synchronization and random long time dynamics for mean-field plane
  rotators,
\newblock \emph{Probab. Theory Related Fields}, \textbf{160} (2014), 593--653,
\newblock \urlprefix\url{http://link.springer.com/10.1007/s00440-013-0536-6}.

\bibitem{bhatnagar_model_1954}
\newblock P.~L. Bhatnagar, E.~P. Gross and M.~Krook,
\newblock A {Model} for {Collision} {Processes} in {Gases}. {I}. {Small}
  {Amplitude} {Processes} in {Charged} and {Neutral} {One}-{Component}
  {Systems},
\newblock \emph{Phys. Rev.}, \textbf{94} (1954), 511--525,
\newblock \urlprefix\url{https://link.aps.org/doi/10.1103/PhysRev.94.511}.

\bibitem{billingsley_convergence_1999}
\newblock P.~Billingsley,
\newblock \emph{Convergence of {Probability} {Measures}},
\newblock 2nd edition,
\newblock Wiley {Series} in {Probability} and {Statistics}., Wiley, New York,
  1999.

\bibitem{bodineau_brownian_2016}
\newblock T.~Bodineau, I.~Gallagher and L.~Saint-Raymond,
\newblock The {Brownian} motion as the limit of a deterministic system of
  hard-spheres,
\newblock \emph{Invent. math.}, \textbf{203} (2016), 493--553,
\newblock \urlprefix\url{http://link.springer.com/10.1007/s00222-015-0593-9}.

\bibitem{bolley_quantitative_2010}
\newblock F.~Bolley,
\newblock Quantitative concentration inequalities on sample path space for mean
  field interaction,
\newblock \emph{ESAIM Probab. Stat.}, \textbf{14} (2010), 192--209.

\bibitem{bolley_convergence_2012}
\newblock F.~Bolley, I.~Gentil and A.~Guillin,
\newblock Convergence to equilibrium in {Wasserstein} distance for
  {Fokker}{\textendash}{Planck} equations,
\newblock \emph{J. Funct. Anal.}, \textbf{263} (2012), 2430--2457,
\newblock
  \urlprefix\url{https://linkinghub.elsevier.com/retrieve/pii/S0022123612002777}.

\bibitem{bolley_uniform_2013}
\newblock F.~Bolley, I.~Gentil and A.~Guillin,
\newblock Uniform {Convergence} to {Equilibrium} for {Granular} {Media},
\newblock \emph{Arch. Ration. Mech. Anal.}, \textbf{208} (2013), 429--445,
\newblock \urlprefix\url{http://link.springer.com/10.1007/s00205-012-0599-z}.

\bibitem{bolley_quantitative_2006}
\newblock F.~Bolley, A.~Guillin and C.~Villani,
\newblock Quantitative concentration inequalities for empirical measures on
  non-compact spaces,
\newblock \emph{Probab. Theory Related Fields}, \textbf{137} (2006), 541--593.

\bibitem{bolley_weighted_2005}
\newblock F.~Bolley and C.~Villani,
\newblock Weighted {Csisz{\'a}r}-{Kullback}-{Pinsker} inequalities and
  applications to transportation inequalities,
\newblock \emph{Ann. Fac. Sci. Toulouse Math. (6)}, \textbf{14} (2005),
  331--352,
\newblock \urlprefix\url{http://www.numdam.org/item/AFST_2005_6_14_3_331_0/}.

\bibitem{bolthausen_laplace_1986}
\newblock E.~Bolthausen,
\newblock Laplace approximations for sums of independent random vectors,
\newblock \emph{Probab. Theory Related Fields}, \textbf{72} (1986), 305--318,
\newblock \urlprefix\url{https://doi.org/10.1007/BF00699109}.

\bibitem{boltzmann_weitere_1872}
\newblock L.~Boltzmann,
\newblock Weitere {Studien} {\"u}ber das {W{\"a}rmegleichgewicht} unter
  {Gasmolek{\"u}len},
\newblock \emph{Sitzungsberichte der Akademie der Wissenschaften}, \textbf{66}
  (1872), 275--370,
\newblock Translation: Further studies on the thermal equilibrium of gas
  molecules, in \textit{Kinetic Theory 2}, 88{\textendash}174, Ed. S.G. Brush,
  Pergamon, Oxford (1966).

\bibitem{cont_markov_2010}
\newblock M.~Bossy and N.~Champagnat,
\newblock Markov {Processes},
\newblock in \emph{Encyclopedia of {Quantitative} {Finance}} (ed. R.~Cont),
\newblock John Wiley \& Sons, Ltd, 2010,
\newblock \urlprefix\url{http://doi.wiley.com/10.1002/9780470061602.eqf02016}.

\bibitem{bottcher_feller_2014}
\newblock B.~B{\"o}ttcher,
\newblock Feller evolution systems: {Generators} and approximation,
\newblock \emph{Stoch. Dyn.}, \textbf{14} (2014), 1350025,
\newblock
  \urlprefix\url{https://www.worldscientific.com/doi/abs/10.1142/S0219493713500251}.

\bibitem{bottcher_levy_2013}
\newblock B.~B{\"o}ttcher, R.~Schilling and J.~Wang,
\newblock \emph{L{\'e}vy {Matters} {III}. {L{\'e}vy}-type processes:
  construction, approximation and sample path properties},
\newblock no.~3 in L{\'e}vy {Matters}, Springer International Publishing, 2013.

\bibitem{bresch_mean-field_2019}
\newblock D.~Bresch, P.-E. Jabin and Z.~Wang,
\newblock On mean-field limits and quantitative estimates with a large class of
  singular kernels: {Application} to the
  {Patlak}{\textendash}{Keller}{\textendash}{Segel} model,
\newblock \emph{C. R. Math. Acad. Sci. Paris}, \textbf{357} (2019), 708--720.

\bibitem{cardaliaguet_notes_2010}
\newblock P.~Cardaliaguet,
\newblock Notes on mean field games (from {P}.-{L}. {Lions}' lectures at
  {Coll{\`e}ge} de {France}),
\newblock in \emph{Lecture given at {Tor} {Vergata}}, 2010,
\newblock 1--59.

\bibitem{cardaliaguet_master_2019}
\newblock P.~Cardaliaguet, F.~Delarue, J.-M. Lasry and P.-L. Lions,
\newblock \emph{The {Master} {Equation} and the {Convergence} {Problem} in
  {Mean} {Field} {Games}},
\newblock no. 201 in Annals of {Mathematics} {Studies}, Princeton University
  Press, 2019,
\newblock \urlprefix\url{http://www.jstor.org/stable/10.2307/j.ctvckq7qf}.

\bibitem{carlen_central_2000}
\newblock E.~Carlen, M.~C. Carvalho and E.~Gabetta,
\newblock Central limit theorem for {Maxwellian} molecules and truncation of
  the {Wild} expansion,
\newblock \emph{Commun. Pure Appl. Math.}, \textbf{53} (2000), 370--397.

\bibitem{carlen_entropy_2008}
\newblock E.~Carlen, M.~C. Carvalho, J.~Le~Roux, M.~Loss and C.~Villani,
\newblock Entropy and chaos in the {Kac} model,
\newblock \emph{Kinet. Relat. Models}, \textbf{3} (2008), 85--122.

\bibitem{carlen_kinetic_2013-1}
\newblock E.~Carlen, R.~Chatelin, P.~Degond and B.~Wennberg,
\newblock Kinetic hierarchy and propagation of chaos in biological swarm
  models,
\newblock \emph{Phys. D}, \textbf{260} (2013), 90--111,
\newblock
  \urlprefix\url{https://linkinghub.elsevier.com/retrieve/pii/S0167278912001492}.

\bibitem{carlen_kinetic_2013}
\newblock E.~Carlen, P.~Degond and B.~Wennberg,
\newblock Kinetic limits for pair-interaction driven master equations and
  biological swarm models,
\newblock \emph{Math. Models Methods Appl. Sci.}, \textbf{23} (2013),
  1339--1376.

\bibitem{carmona_lectures_2016}
\newblock R.~Carmona,
\newblock \emph{Lectures on {BSDEs}, {Stochastic} {Control}, and {Stochastic}
  {Differential} {Games} with {Financial} {Applications}},
\newblock SIAM, 2016.

\bibitem{carmona_probabilistic_2018}
\newblock R.~Carmona and F.~Delarue,
\newblock \emph{Probabilistic {Theory} of {Mean} {Field} {Games} with
  {Applications} {I}, {Mean} {Field} {FBSDEs}, {Control}, and {Games}},
\newblock no.~83 in Probability {Theory} and {Stochastic} {Modelling}, Springer
  International Publishing, 2018.

\bibitem{carmona_probabilistic_2018-1}
\newblock R.~Carmona and F.~Delarue,
\newblock \emph{Probabilistic {Theory} of {Mean} {Field} {Games} with
  {Applications} {II}, {Mean} {Field} {Games} with {Common} {Noise} and
  {Master} {Equations}},
\newblock no.~84 in Probability {Theory} and {Stochastic} {Modelling}, Springer
  International Publishing, 2018.

\bibitem{naldi_particle_2010}
\newblock J.~A. Carrillo, M.~Fornasier, G.~Toscani and F.~Vecil,
\newblock Particle, kinetic, and hydrodynamic models of swarming,
\newblock in \emph{Mathematical {Modeling} of {Collective} {Behavior} in
  {Socio}-{Economic} and {Life} {Sciences}} (eds. G.~Naldi, L.~Pareschi and
  G.~Toscani),
\newblock Birkh{\"a}user Boston, 2010,
\newblock 297--336,
\newblock
  \urlprefix\url{http://link.springer.com/10.1007/978-0-8176-4946-3_12}.

\bibitem{carrillo_consensus-based_2021}
\newblock J.~A. Carrillo, S.~Jin, L.~Li and Y.~Zhu,
\newblock A consensus-based global optimization method for high dimensional
  machine learning problems,
\newblock \emph{ESAIM Control Optim. Calc. Var.}, \textbf{27} (2021), 1--22,
\newblock \urlprefix\url{https://www.esaim-cocv.org/10.1051/cocv/2020046}.

\bibitem{carrillo_-convexity_2020}
\newblock J.~A. Carrillo, M.~Delgadino and G.~Pavliotis,
\newblock A $\lambda$-convexity based proof for the propagation of chaos for
  weakly interacting stochastic particles,
\newblock \emph{J. Funct. Anal.}, \textbf{279} (2020), 108734.

\bibitem{carrillo_double_2009}
\newblock J.~A. Carrillo, M.~R. D'Orsogna and V.~Panferov,
\newblock Double milling in self-propelled swarms from kinetic theory,
\newblock \emph{Kinet. Relat. Models}, \textbf{2} (2009), 363--378,
\newblock
  \urlprefix\url{http://aimsciences.org//article/doi/10.3934/krm.2009.2.363}.

\bibitem{cattiaux_stochastic_2018}
\newblock P.~Cattiaux, F.~Delebecque and L.~P{\'e}d{\`e}ches,
\newblock Stochastic {Cucker}{\textendash}{Smale} models: {Old} and new,
\newblock \emph{Ann. Appl. Probab.}, \textbf{28} (2018), 3239--3286,
\newblock
  \urlprefix\url{https://projecteuclid.org/journals/annals-of-applied-probability/volume-28/issue-5/Stochastic-CuckerSmale-models-Old-and-new/10.1214/18-AAP1400.full}.

\bibitem{cercignani_boltzmann_1988}
\newblock C.~Cercignani,
\newblock \emph{The {Boltzmann} {Equation} and {Its} {Applications}},
\newblock no.~67 in Applied {Mathematical} {Sciences}, Springer New York, 1988,
\newblock \urlprefix\url{http://link.springer.com/10.1007/978-1-4612-1039-9}.

\bibitem{cercignani_ludwig_2006}
\newblock C.~Cercignani,
\newblock \emph{Ludwig {Boltzmann}, the {Man} {Who} {Trusted} {Atoms}},
\newblock Oxford University Press, 2006,
\newblock
  \urlprefix\url{https://oxford.universitypressscholarship.com/view/10.1093/acprof:oso/9780198570646.001.0001/acprof-9780198570646}.

\bibitem{cercignani_mathematical_1994}
\newblock C.~Cercignani, R.~Illner and M.~Pulvirenti,
\newblock \emph{The {Mathematical} {Theory} of {Dilute} {Gases}},
\newblock no. 106 in Applied {Mathematical} {Sciences}, Springer-Verlag New
  York, 1994.

\bibitem{champion_monge_2011}
\newblock T.~Champion and L.~De~Pascale,
\newblock The {Monge} problem in ?$^{\textrm{\textit{d}}}$,
\newblock \emph{Duke Math. J.}, \textbf{157} (2011), 551--572.

\bibitem{chizat_global_2018}
\newblock L.~Chizat and F.~Bach,
\newblock On the {Global} {Convergence} of {Gradient} {Descent} for
  {Over}-parameterized {Models} using {Optimal} {Transport},
\newblock in \emph{Advances in {Neural} {Information} {Processing} {Systems} 31
  ({NeurIPS} 2018)} (eds. S.~Bengio, H.~Wallach, H.~Larochelle, K.~Grauman,
  N.~Cesa-Bianchi and R.~Garnett),
\newblock Curran Associates, Inc., Montreal, Canada, 2018,
\newblock 3040--3050.

\bibitem{cortez_quantitative_2016}
\newblock R.~Cortez and J.~Fontbona,
\newblock Quantitative propagation of chaos for generalized {Kac} particle
  systems,
\newblock \emph{Ann. Appl. Probab.}, \textbf{26} (2016), 892--916,
\newblock \urlprefix\url{http://projecteuclid.org/euclid.aoap/1458651823}.

\bibitem{cortez_quantitative_2018}
\newblock R.~Cortez and J.~Fontbona,
\newblock Quantitative {Uniform} {Propagation} of {Chaos} for {Maxwell}
  {Molecules},
\newblock \emph{Commun. Math. Phys.}, \textbf{357} (2018), 913--941,
\newblock \urlprefix\url{http://link.springer.com/10.1007/s00220-018-3101-4}.

\bibitem{csiszar_sanov_1984}
\newblock I.~Csisz{\'a}r,
\newblock Sanov {Property}, {Generalized} {I}-{Projection} and a {Conditional}
  {Limit} {Theorem},
\newblock \emph{Ann. Probab.}, \textbf{12} (1984), 768--793,
\newblock \urlprefix\url{https://projecteuclid.org/euclid.aop/1176993227}.

\bibitem{cucker_mathematics_2007}
\newblock F.~Cucker and S.~Smale,
\newblock On the mathematics of emergence,
\newblock \emph{Jpn. J. Math.}, \textbf{2} (2007), 197--227,
\newblock \urlprefix\url{http://link.springer.com/10.1007/s11537-007-0647-x}.

\bibitem{dai_pra_mckean-vlasov_1996}
\newblock P.~Dai~Pra and F.~den Hollander,
\newblock {McKean}-{Vlasov} limit for interacting random processes in random
  media,
\newblock \emph{J. Stat. Phys.}, \textbf{84} (1996), 735--772.

\bibitem{pajot_lecture_2014}
\newblock S.~Danieri and G.~Savar{\'e},
\newblock Lecture notes on gradient flows and optimal transport,
\newblock in \emph{Optimal {Transportation}} (eds. H.~Pajot, Y.~Ollivier and
  C.~Villani),
\newblock Cambridge University Press, Cambridge, 2014,
\newblock 100--144,
\newblock
  \urlprefix\url{https://www.cambridge.org/core/product/identifier/CBO9781107297296A015/type/book_part}.

\bibitem{dawson_measure-valued_1993}
\newblock D.~Dawson,
\newblock Measure-valued {Markov} processes,
\newblock in \emph{{\'E}cole d'{\'E}t{\'e} de {Probabilit{\'e}s} de
  {Saint}-{Flour} {XXI}-1991} (ed. P.~Hennequin),
\newblock no. 1541 in Lecture {Notes} in {Mathematics}, Springer Berlin
  Heidelberg, 1993.

\bibitem{dawson_large_1987}
\newblock D.~Dawson and J.~G{\"a}rtner,
\newblock Large deviations from the {McKean}-{Vlasov} limit for weakly
  interacting diffusions,
\newblock \emph{Stochastics}, \textbf{20} (1987), 247--308.

\bibitem{de_bortoli_quantitative_2020}
\newblock V.~De~Bortoli, A.~Durmus and X.~Fontaine,
\newblock Quantitative {Propagation} of {Chaos} for {SGD} in {Wide} {Neural}
  {Networks},
\newblock in \emph{Advances in {Neural} {Information} {Processing} {Systems} 33
  ({NeurIPS} 2020)}, 2020,
\newblock 278--288,
\newblock
  \urlprefix\url{https://proceedings.neurips.cc/paper/2020/file/02e74f10e0327ad868d138f2b4fdd6f0-Paper.pdf}.

\bibitem{de_masi_hydrodynamic_2015}
\newblock A.~De~Masi, A.~Galves, E.~L{\"o}cherbach and E.~Presutti,
\newblock Hydrodynamic {Limit} for {Interacting} {Neurons},
\newblock \emph{J. Stat. Phys.}, \textbf{158} (2015), 866--902,
\newblock \urlprefix\url{http://link.springer.com/10.1007/s10955-014-1145-1}.

\bibitem{bellomo_macroscopic_2004}
\newblock P.~Degond,
\newblock Macroscopic limits of the {Boltzmann} equation: a review,
\newblock in \emph{Modeling and {Computational} {Methods} for {Kinetic}
  {Equations}} (eds. N.~Bellomo, P.~Degond, L.~Pareschi and G.~Russo),
\newblock Birkh{\"a}user Boston, Boston, MA, 2004,
\newblock 3--57,
\newblock \urlprefix\url{http://link.springer.com/10.1007/978-0-8176-8200-2_1},
\newblock Series Title: Modeling and Simulation in Science, Engineering and
  Technology.

\bibitem{degond_mathematical_2018}
\newblock P.~Degond,
\newblock Mathematical models of collective dynamics and self-organization,
\newblock in \emph{Proceedings of the {International} {Congress} of
  {Mathematicians} {ICM} 2018}, vol.~4,
\newblock Rio de Janeiro, Brazil, 2018,
\newblock 3943--3964.

\bibitem{degond_phase_2015}
\newblock P.~Degond, A.~Frouvelle and J.-G. Liu,
\newblock Phase {Transitions}, {Hysteresis}, and {Hyperbolicity} for
  {Self}-{Organized} {Alignment} {Dynamics},
\newblock \emph{Arch. Ration. Mech. Anal.}, \textbf{216} (2015), 63--115,
\newblock \urlprefix\url{http://link.springer.com/10.1007/s00205-014-0800-7}.

\bibitem{giacomin_alignment_2019}
\newblock P.~Degond, A.~Frouvelle, S.~Merino-Aceituno and A.~Trescases,
\newblock Alignment of {Self}-propelled {Rigid} {Bodies}: {From} {Particle}
  {Systems} to {Macroscopic} {Equations},
\newblock in \emph{Stochastic {Dynamics} {Out} of {Equilibrium}, {Institut}
  {Henri} {Poincar{\'e}}, {Paris}, {France}, 2017} (eds. G.~Giacomin, S.~Olla,
  E.~Saada, H.~Spohn and G.~Stoltz),
\newblock no. 282 in Springer {Proceedings} in {Mathematics} \& {Statistics},
  Springer, Cham, 2019,
\newblock 28--66,
\newblock \urlprefix\url{http://link.springer.com/10.1007/978-3-030-15096-9_2}.

\bibitem{degond_continuum_2008}
\newblock P.~Degond and S.~Motsch,
\newblock Continuum limit of self-driven particles with orientation
  interaction,
\newblock \emph{Math. Models Methods Appl. Sci.}, \textbf{18} (2008),
  1193--1215,
\newblock
  \urlprefix\url{https://www.worldscientific.com/doi/abs/10.1142/S0218202508003005}.

\bibitem{del_moral_measure-valued_1998}
\newblock P.~Del~Moral,
\newblock Measure-valued processes and interacting particle systems.
  {Application} to nonlinear filtering problems,
\newblock \emph{Ann. Appl. Probab.}, \textbf{8} (1998), 438--495.

\bibitem{del_moral_feynman-kac_2004}
\newblock P.~Del~Moral,
\newblock \emph{Feynman-{Kac} {Formulae}, {Genealogical} and {Interacting}
  {Particle} {Systems} with {Applications}},
\newblock Probability and {Its} {Applications}, Springer-Verlag New York, 2004.

\bibitem{del_moral_mean_2013}
\newblock P.~Del~Moral,
\newblock \emph{Mean field simulation for {Monte} {Carlo} integration},
\newblock no. 126 in Monographs on {Statistics} and {Applied} {Probability},
  CRC Press, Taylor \& Francis Group, 2013.

\bibitem{del_moral_stochastic_2017}
\newblock P.~Del~Moral and S.~Penev,
\newblock \emph{Stochastic {Processes}: {From} {Applications} to {Theory}},
\newblock Texts in {Statistical} {Science}, CRC Press : Taylor \& Francis
  Group, 2017.

\bibitem{del_moral_uniform_2019}
\newblock P.~Del~Moral and J.~Tugaut,
\newblock Uniform propagation of chaos and creation of chaos for a class of
  nonlinear diffusions,
\newblock \emph{Stoch. Anal. Appl.}, \textbf{37} (2019), 909--935.

\bibitem{diaconis_finite_1980}
\newblock P.~Diaconis and D.~Freedman,
\newblock Finite {Exchangeable} {Sequences},
\newblock \emph{Ann. Probab.}, \textbf{8} (1980), 745--764,
\newblock \urlprefix\url{https://projecteuclid.org/euclid.aop/1176994663}.

\bibitem{diez_propagation_2020}
\newblock A.~Diez,
\newblock Propagation of chaos and moderate interaction for a piecewise
  deterministic system of geometrically enriched particles,
\newblock \emph{Electron. J. Probab.}, \textbf{25} (2020), 1--38.

\bibitem{dorsogna_self-propelled_2006}
\newblock M.~R. D{\textquoteright}Orsogna, Y.~L. Chuang, A.~L. Bertozzi and
  L.~S. Chayes,
\newblock Self-{Propelled} {Particles} with {Soft}-{Core} {Interactions}:
  {Patterns}, {Stability}, and {Collapse},
\newblock \emph{Phys. Rev. Lett.}, \textbf{96} (2006), 104302,
\newblock
  \urlprefix\url{https://link.aps.org/doi/10.1103/PhysRevLett.96.104302}.

\bibitem{durmus_elementary_2020}
\newblock A.~Durmus, A.~Eberle, A.~Guillin and R.~Zimmer,
\newblock An elementary approach to uniform in time propagation of chaos,
\newblock \emph{Proc. Amer. Math. Soc.}, \textbf{148} (2020), 5387--5398.

\bibitem{eberle_reflection_2016}
\newblock A.~Eberle,
\newblock Reflection couplings and contraction rates for diffusions,
\newblock \emph{Probab. Theory Related Fields}, \textbf{166} (2016), 851--886.

\bibitem{eberle_quantitative_2019}
\newblock A.~Eberle, A.~Guillin and R.~Zimmer,
\newblock Quantitative {Harris}-type theorems for diffusions and
  {McKean}-{Vlasov} processes,
\newblock \emph{Trans. Amer. Math. Soc.}, \textbf{371} (2019), 7135--7173.

\bibitem{el_karoui_martingale_1990}
\newblock N.~El~Karoui and S.~M{\'e}l{\'e}ard,
\newblock Martingale measures and stochastic calculus,
\newblock \emph{Probab. Theory Related Fields}, \textbf{84} (1990), 83--101,
\newblock \urlprefix\url{http://link.springer.com/10.1007/BF01288560}.

\bibitem{ethier_markov_1986}
\newblock S.~N. Ethier and T.~G. Kurtz,
\newblock \emph{Markov processes: characterization and convergence},
\newblock Wiley series in probability and mathematical statistics, Wiley, New
  York, 1986.

\bibitem{feldman_monges_2001}
\newblock M.~Feldman and R.~J. McCann,
\newblock Monge{\textquoteright}s transport problem on a {Riemannian} manifold,
\newblock \emph{Trans. Amer. Math. Soc.}, \textbf{354} (2001), 1667--1697.

\bibitem{fontbona_measurability_2009}
\newblock J.~Fontbona, H.~Gu{\'e}rin and S.~M{\'e}l{\'e}ard,
\newblock Measurability of optimal transportation and convergence rate for
  {Landau} type interacting particle systems,
\newblock \emph{Probab. Theory Related Fields}, \textbf{143} (2009), 329--351,
\newblock \urlprefix\url{http://link.springer.com/10.1007/s00440-007-0128-4}.

\bibitem{fournier_rate_2015}
\newblock N.~Fournier and A.~Guillin,
\newblock On the rate of convergence in {Wasserstein} distance of the empirical
  measure,
\newblock \emph{Probab. Theory Related Fields}, \textbf{162} (2015), 707--738,
\newblock Publisher: Springer.

\bibitem{fournier_toy_2016}
\newblock N.~Fournier and E.~L{\"o}cherbach,
\newblock On a toy model of interacting neurons,
\newblock \emph{Ann. Inst. Henri Poincar{\'e} Probab. Stat.}, \textbf{52}
  (2016), 1844--1876,
\newblock
  \urlprefix\url{https://projecteuclid.org/journals/annales-de-linstitut-henri-poincare-probabilites-et-statistiques/volume-52/issue-4/On-a-toy-model-of-interacting-neurons/10.1214/15-AIHP701.full}.

\bibitem{fournier_markov_2001}
\newblock N.~Fournier and S.~M{\'e}l{\'e}ard,
\newblock A {Markov} {Process} {Associated} with a {Boltzmann} {Equation}
  {Without} {Cutoff} and for {Non}-{Maxwell} {Molecules},
\newblock \emph{J. Stat. Phys.}, \textbf{104} (2001), 359--385,
\newblock \urlprefix\url{http://link.springer.com/10.1023/A:1010322130480}.

\bibitem{fournier_rate_2016}
\newblock N.~Fournier and S.~Mischler,
\newblock Rate of convergence of the {Nanbu} particle system for hard
  potentials and {Maxwell} molecules,
\newblock \emph{Ann. Probab.}, \textbf{44} (2016), 589--627,
\newblock \urlprefix\url{http://projecteuclid.org/euclid.aop/1454423051}.

\bibitem{funaki_certain_1984}
\newblock T.~Funaki,
\newblock A {Certain} {Class} of {Diffusion} {Processes} {Associated} with
  {Nonlinear} {Parabolic} {Equations},
\newblock \emph{Zeitschrift f{\"u}r Wahrscheinlichkeitstheorie und Verwandte
  Gebiete}, \textbf{67} (1984), 331--348.

\bibitem{gallagher_newton_2014}
\newblock I.~Gallagher, L.~Saint-Raymond and B.~Texier,
\newblock \emph{From {Newton} to {Boltzmann}: {Hard} {Spheres} and
  {Short}-{Range} {Potentials}},
\newblock no.~18 in Zurich {Lectures} in {Advanced} {Mathematics}, European
  Mathematical Society, 2014.

\bibitem{gartner_mckean-vlasov_1988}
\newblock J.~G{\"a}rtner,
\newblock On the {McKean}-{Vlasov} limit for {Interacting} {Diffusions},
\newblock \emph{Math. Nachr.}, \textbf{137} (1988), 197--248.

\bibitem{golse_dynamics_2016}
\newblock F.~Golse,
\newblock On the {Dynamics} of {Large} {Particle} {Systems} in the {Mean}
  {Field} {Limit},
\newblock Lecture notes, \arxiv{1301.5494}.

\bibitem{grad_asymptotic_1963}
\newblock H.~Grad,
\newblock Asymptotic theory of {Boltzmann} equation,
\newblock \emph{Phys. Fluids}, \textbf{6} (1963), 147--181.

\bibitem{graham_mckean-vlasov_1992}
\newblock C.~Graham,
\newblock {McKean}-{Vlasov} {It{\=o}}-{Skorohod} equations, and nonlinear
  diffusions with discrete jump sets,
\newblock \emph{Stochastic Process. Appl.}, \textbf{40} (1992), 69--82.

\bibitem{graham_nonlinear_1992}
\newblock C.~Graham,
\newblock Nonlinear diffusion with jumps,
\newblock \emph{Ann. Inst. Henri Poincar{\'e} Probab. Stat.}, \textbf{28}
  (1992), 393--402.

\bibitem{graham_stochastic_1997}
\newblock C.~Graham and S.~M{\'e}l{\'e}ard,
\newblock Stochastic particle approximations for generalized {Boltzmann} models
  and convergence estimates,
\newblock \emph{Ann. Probab.}, \textbf{25} (1997), 115--132.

\bibitem{grassi_particle_2021}
\newblock S.~Grassi and L.~Pareschi,
\newblock From particle swarm optimization to consensus based optimization:
  stochastic modeling and mean-field limit,
\newblock \emph{Math. Models Methods Appl. Sci.}, \textbf{31} (2021),
  1625--1657.

\bibitem{grunbaum_propagation_1971}
\newblock F.~A. Gr{\"u}nbaum,
\newblock Propagation of chaos for the {Boltzmann} equation,
\newblock \emph{Arch. Ration. Mech. Anal.}, \textbf{42} (1971), 323--345.

\bibitem{guillin_uniform_2020}
\newblock A.~Guillin and P.~Monmarch{\'e},
\newblock Uniform long-time and propagation of chaos estimates for mean field
  kinetic particles in non-convex landscapes,
\newblock \emph{J. Stat. Phys.}, \textbf{185} (2020), 1--20.

\bibitem{ha_emergence_2009}
\newblock S.-Y. Ha, K.~Lee and D.~Levy,
\newblock Emergence of time-asymptotic flocking in a stochastic
  {Cucker}-{Smale} system,
\newblock \emph{Commun. Math. Sci.}, \textbf{7} (2009), 453--469.

\bibitem{hauray_kacs_2014}
\newblock M.~Hauray and S.~Mischler,
\newblock On {Kac}'s chaos and related problems,
\newblock \emph{J. Funct. Anal.}, \textbf{266} (2014), 6055--6157.

\bibitem{holding_propagation_2016}
\newblock T.~Holding,
\newblock Propagation of chaos for {H{\"o}lder} continuous interaction kernels
  via {Glivenko}-{Cantelli},
\newblock preprint, \arxiv{1608.02877}.

\bibitem{ikeda_stochastic_1989}
\newblock N.~Ikeda and S.~Watanabe,
\newblock \emph{Stochastic differential equations and diffusion processes},
\newblock no.~24 in North-{Holland} {Mathematical} {Library}, Elsevier, 1989,
\newblock
  \urlprefix\url{http://qut.eblib.com.au/patron/FullRecord.aspx?p=1888657}.

\bibitem{jabin_review_2014}
\newblock P.-E. Jabin,
\newblock A review of the mean field limits for {Vlasov} equations,
\newblock \emph{Kinet. Relat. Models}, \textbf{7} (2014), 661--711.

\bibitem{jabin_mean_2017}
\newblock P.-E. Jabin and Z.~Wang,
\newblock Mean {Field} {Limit} for {Stochastic} {Particle} {Systems},
\newblock in \emph{Active {Particles}, {Volume} 1 : {Advances} in {Theory},
  {Models}, and {Applications}} (eds. N.~Bellomo, P.~Degond and E.~Tadmor),
\newblock Modeling and {Simulation} in {Science}, {Engineering} and
  {Technology}, Birkh{\"a}user Basel, 2017,
\newblock 379--402,
\newblock \urlprefix\url{https://doi.org/10.1007/978-3-319-49996-3_10}.

\bibitem{jabin_quantitative_2018}
\newblock P.-E. Jabin and Z.~Wang,
\newblock Quantitative estimates of propagation of chaos for stochastic systems
  with \textit{{W}}$^{\textrm{-1,$\infty$}}$ kernels,
\newblock \emph{Invent. Math.}, \textbf{214} (2018), 523--591.

\bibitem{jacod_limit_2003}
\newblock J.~Jacod and A.~N. Shiryaev,
\newblock \emph{Limit {Theorems} for {Stochastic} {Processes}},
\newblock Second edition edition,
\newblock no. 288 in Grundlehren der mathematischen {Wissenschaften}, Springer
  Berlin Heidelberg, 2003,
\newblock \urlprefix\url{http://link.springer.com/10.1007/978-3-662-05265-5}.

\bibitem{joffe_weak_1986}
\newblock A.~Joffe and M.~M{\'e}tivier,
\newblock Weak convergence of sequences of semimartingales with applications to
  multitype branching processes,
\newblock \emph{Adv. in Appl. Probab.}, \textbf{18} (1986), 20--65,
\newblock
  \urlprefix\url{https://www.cambridge.org/core/product/identifier/S0001867800015585/type/journal_article}.

\bibitem{johnson_class_1968}
\newblock D.~P. Johnson,
\newblock On a class of stochastic processes and its relationship to infinite
  particle gases,
\newblock \emph{Trans. Amer. Math. Soc.}, \textbf{132} (1968), 275--275,
\newblock
  \urlprefix\url{http://www.ams.org/jourcgi/jour-getitem?pii=S0002-9947-1968-0256452-X}.

\bibitem{jourdain_propagation_1998}
\newblock B.~Jourdain and S.~M{\'e}l{\'e}ard,
\newblock Propagation of chaos and fluctuations for a moderate model with
  smooth initial data,
\newblock \emph{Ann. Inst. Henri Poincar{\'e} Probab. Stat.}, \textbf{34}
  (1998), 727--766,
\newblock \urlprefix\url{http://www.numdam.org/item/AIHPB_1998__34_6_727_0},
\newblock Publisher: Gauthier-Villars.

\bibitem{kac_foundations_1956}
\newblock M.~Kac,
\newblock Foundations of kinetic theory,
\newblock in \emph{Proceedings of the {Third} {Berkeley} {Symposium} on
  {Mathematical} {Statistics} and {Probability}}, vol.~3,
\newblock University of California Press Berkeley and Los Angeles, California,
  1956,
\newblock 171--197.

\bibitem{cohen_probabilistic_1973}
\newblock M.~Kac,
\newblock Some {Probabilistic} {Aspects} of the {Boltzmann} {Equation},
\newblock in \emph{The {Boltzmann} {Equation}. {Acta} {Physica} {Austriaca}
  ({Supplementum} {X} {Proceedings} of the {International} {Symposium}
  {\textquotedblleft}100 {Years} {Boltzmann} {Equation}{\textquotedblright} in
  {Vienna} 4th{\textendash}8th {September} 1972)} (eds. E.~G.~D. Cohen and
  W.~Thirring),
\newblock Springer Vienna, 1973,
\newblock 379--400,
\newblock
  \urlprefix\url{http://link.springer.com/10.1007/978-3-7091-8336-6_17}.

\bibitem{karatzas_brownian_1998}
\newblock I.~Karatzas and S.~Shreve,
\newblock \emph{Brownian {Motion} and {Stochastic} {Calculus}},
\newblock 2nd edition,
\newblock no. 113 in Graduate {Texts} in {Mathematics}, Springer-Verlag New
  York, 1998.

\bibitem{kusuoka_gibbs_1984}
\newblock S.~Kusuoka and Y.~Tamura,
\newblock Gibbs measures for mean field potentials,
\newblock \emph{J. Fac. Sci. Univ. Tokyo}, \textbf{31} (1984), 223--245.

\bibitem{lacker_strong_2018}
\newblock D.~Lacker,
\newblock On a strong form of propagation of chaos for {McKean}-{Vlasov}
  equations,
\newblock \emph{Electron. Commun. Probab.}, \textbf{23} (2018), 1--11.

\bibitem{lacker_hierarchies_2021}
\newblock D.~Lacker,
\newblock Hierarchies, entropy, and quantitative propagation of chaos for mean
  field diffusions,
\newblock preprint, \arxiv{2105.02983}.

\bibitem{lanford_time_1975}
\newblock O.~E. Lanford,
\newblock Time evolution of large classical systems,
\newblock in \emph{Dynamical {Systems} {Theory} and {Application}, {Battelle}
  {Seattle} 1974 {Rencontres}} (ed. J.~Moser),
\newblock Springer-Verlag Berlin Heidelberg, 1975.

\bibitem{le_gall_brownian_2016}
\newblock J.-F. Le~Gall,
\newblock \emph{Brownian {Motion}, {Martingales}, and {Stochastic} {Calculus}},
\newblock no. 274 in Graduate {Texts} in {Mathematics}, Springer International
  Publishing, 2016.

\bibitem{ledoux_concentration_1999}
\newblock M.~Ledoux,
\newblock Concentration of measure and logarithmic {Sobolev} inequalities,
\newblock \emph{S{\'e}minaire de probabilit{\'e}s de Strasbourg}, \textbf{33}
  (1999), 120--216,
\newblock \urlprefix\url{http://www.numdam.org/item/SPS_1999__33__120_0/}.

\bibitem{leonard_loi_1986}
\newblock C.~L{\'e}onard,
\newblock Une loi des grands nombres pour des syst{\`e}mes de diffusions avec
  interaction et {\`a} coefficients non born{\'e}s,
\newblock \emph{Ann. Inst. Henri Poincar{\'e} Probab. Stat.}, \textbf{22}
  (1986), 237--262.

\bibitem{leonard_large_1995}
\newblock C.~L{\'e}onard,
\newblock Large deviations for long range interacting particle systems with
  jumps,
\newblock \emph{Ann. Inst. Henri Poincar{\'e} Probab. Stat.}, \textbf{31}
  (1995), 289--323,
\newblock \urlprefix\url{http://www.numdam.org/item/AIHPB_1995__31_2_289_0/}.

\bibitem{leonard_large_1995-1}
\newblock C.~L{\'e}onard,
\newblock On large deviations for particle systems associated with spatially
  homogeneous {Boltzmann} type equations,
\newblock \emph{Probab. Theory Related Fields}, \textbf{101} (1995), 1--44,
\newblock \urlprefix\url{https://doi.org/10.1007/BF01192194}.

\bibitem{leonard_girsanov_2012}
\newblock C.~L{\'e}onard,
\newblock Girsanov theory under a finite entropy condition,
\newblock in \emph{S{\'e}minaire de {Probabilit{\'e}s} {XLIV}},
\newblock no. 2046 in Lecture {Notes} in {Mathematics}, Springer, Berlin,
  Heidelberg, 2012,
\newblock 429--465,
\newblock \urlprefix\url{https://doi.org/10.1007/978-3-642-27461-9_20}.

\bibitem{letta_sur_1989}
\newblock G.~Letta,
\newblock Sur les th{\'e}or{\`e}mes de {Hewitt}-{Savage} et de de {Finetti},
\newblock \emph{S{\'e}minaire de probabilit{\'e}s de Strasbourg}, \textbf{23}
  (1989), 531--535,
\newblock \urlprefix\url{http://www.numdam.org/item/SPS_1989__23__531_0/}.

\bibitem{lindvall_coupling_1986}
\newblock T.~Lindvall and L.~C.~G. Rogers,
\newblock Coupling of {Multidimensional} {Diffusions} by {Reflection},
\newblock \emph{Ann. Probab.}, \textbf{14},
\newblock
  \urlprefix\url{https://projecteuclid.org/journals/annals-of-probability/volume-14/issue-3/Coupling-of-Multidimensional-Diffusions-by-Reflection/10.1214/aop/1176992442.full}.

\bibitem{liu_long-time_2021}
\newblock W.~Liu, L.~Wu and C.~Zhang,
\newblock Long-time behaviors of mean-field interacting particle systems
  related to {McKean}-{Vlasov} equations,
\newblock \emph{Commun. Math. Phys.}, \textbf{387} (2021), 179--214.

\bibitem{goncalves_large_2015}
\newblock E.~Lu{\c c}on,
\newblock Large {Population} {Asymptotics} for {Interacting} {Diffusions} in a
  {Quenched} {Random} {Environment},
\newblock in \emph{From {Particle} {Systems} to {Partial} {Differential}
  {Equations} {II}} (eds. P.~Gon{\c c}alves and A.~J. Soares),
\newblock no. 129 in Springer {Proceedings} in {Mathematics} \& {Statistics},
  Springer, Cham, 2015,
\newblock 231--251,
\newblock \urlprefix\url{http://link.springer.com/10.1007/978-3-319-16637-7_8}.

\bibitem{malrieu_logarithmic_2001}
\newblock F.~Malrieu,
\newblock Logarithmic {Sobolev} inequalities for some nonlinear {PDE}'s,
\newblock \emph{Stochastic Process. Appl.}, \textbf{95} (2001), 109--132.

\bibitem{matthes_steady_2008}
\newblock D.~Matthes and G.~Toscani,
\newblock On {Steady} {Distributions} of {Kinetic} {Models} of {Conservative}
  {Economies},
\newblock \emph{J. Stat. Phys.}, \textbf{130} (2008), 1087--1117,
\newblock \urlprefix\url{http://link.springer.com/10.1007/s10955-007-9462-2}.

\bibitem{mckean_class_1966}
\newblock H.~P. McKean,
\newblock A class of {Markov} processes associated with nonlinear parabolic
  equations,
\newblock \emph{Proc. Nat. Acad. Sci.}, \textbf{56} (1966), 1907,
\newblock Publisher: National Academy of Sciences.

\bibitem{mckean_exponential_1967}
\newblock H.~P. McKean,
\newblock An exponential formula for solving {Boltzmann}'s equation for a
  {Maxwellian} gas,
\newblock \emph{Journal of Combinatorial Theory}, \textbf{2} (1967), 358--382,
\newblock
  \urlprefix\url{http://www.sciencedirect.com/science/article/pii/S0021980067800358}.

\bibitem{mckean_propagation_1969}
\newblock H.~P. McKean,
\newblock Propagation of chaos for a class of non-linear parabolic equations,
\newblock in \emph{Lecture {Series} in {Differential} {Equations}, {Volume} 2}
  (ed. A.~K. Aziz),
\newblock no.~19 in Van {Nostrand} {Mathematical} {Studies}, Van Nostrand
  Reinhold Company, 1969,
\newblock 177--194.

\bibitem{mei_mean_2018}
\newblock S.~Mei, A.~Montanari and P.-M. Nguyen,
\newblock A mean field view of the landscape of two-layer neural networks,
\newblock \emph{Proc. Natl. Acad. Sci. USA}, \textbf{115} (2018), E7665--E7671,
\newblock
  \urlprefix\url{http://www.pnas.org/lookup/doi/10.1073/pnas.1806579115}.

\bibitem{meleard_asymptotic_1996}
\newblock S.~M{\'e}l{\'e}ard,
\newblock Asymptotic {Behaviour} of some interacting particle systems;
  {McKean}-{Vlasov} and {Boltzmann} models,
\newblock in \emph{Probabilistic {Models} for {Nonlinear} {Partial}
  {Differential} {Equations}} (eds. D.~Talay and L.~Tubaro),
\newblock no. 1627 in Lecture {Notes} in {Mathematics}, Springer-Verlag Berlin
  Heidelberg, 1996.

\bibitem{meleard_systemes_1988}
\newblock S.~M{\'e}l{\'e}ard and S.~Roelly-Coppoletta,
\newblock Syst{\`e}mes de particules et mesures-martingales : un
  th{\'e}or{\`e}me de propagation du chaos,
\newblock \emph{S{\'e}minaire de probabilit{\'e}s (Strasbourg)}, \textbf{22}
  (1988), 438--448.

\bibitem{mischler_kacs_2012}
\newblock S.~Mischler,
\newblock Kac's chaos and {Kac}'s program,
\newblock in \emph{S{\'e}minaire {Laurent} {Schwartz} - {EDP} et applications},
  vol. 2012-2013,
\newblock Institut des hautes {\'e}tudes scientifiques \& Centre de
  math{\'e}matiques Laurent Schwartz, 2012,
\newblock Expos{\'e} no XXII, 1--17.

\bibitem{mischler_kacs_2013}
\newblock S.~Mischler and C.~Mouhot,
\newblock Kac's program in kinetic theory,
\newblock \emph{Invent. Math.}, \textbf{193} (2013), 1--147,
\newblock Publisher: Springer.

\bibitem{mischler_new_2015}
\newblock S.~Mischler, C.~Mouhot and B.~Wennberg,
\newblock A new approach to quantitative propagation of chaos for drift,
  diffusion and jump processes,
\newblock \emph{Probab. Theory Related Fields}, \textbf{161} (2015), 1--59.

\bibitem{muntean_collective_2014}
\newblock Adrian Muntean and Federico Toschi (eds.),
\newblock \emph{Collective {Dynamics} from {Bacteria} to {Crowds}: {An}
  {Excursion} {Through} {Modeling}, {Analysis} and {Simulation}},
\newblock no. 553 in {CISM} {International} {Centre} for {Mechanical}
  {Sciences}, Springer, Vienna, 2014,
\newblock \urlprefix\url{http://link.springer.com/10.1007/978-3-7091-1785-9}.

\bibitem{murata_propagation_1977}
\newblock H.~Murata,
\newblock Propagation of chaos for {Boltzmann}-like equation of non-cutoff type
  in the plane,
\newblock \emph{Hiroshima Math. J.}, \textbf{7} (1977), 479--515,
\newblock \urlprefix\url{https://projecteuclid.org/euclid.hmj/1206135751}.

\bibitem{naldi_mathematical_2010}
\newblock Giovanni Naldi, Lorenzo Pareschi and Giuseppe Toscani (eds.),
\newblock \emph{Mathematical {Modeling} of {Collective} {Behavior} in
  {Socio}-{Economic} and {Life} {Sciences}},
\newblock Modeling and {Simulation} in {Science}, {Engineering} and
  {Technology}, Birkh{\"a}user Boston, 2010,
\newblock \urlprefix\url{http://link.springer.com/10.1007/978-0-8176-4946-3}.

\bibitem{nanbu_direct_1980}
\newblock K.~Nanbu,
\newblock Direct {Simulation} {Scheme} {Derived} from the {Boltzmann}
  {Equation}. {I}. {Monocomponent} {Gases},
\newblock \emph{Journal of the Physical Society of Japan}, \textbf{49} (1980),
  2042--2049,
\newblock \urlprefix\url{http://journals.jps.jp/doi/10.1143/JPSJ.49.2042}.

\bibitem{oelschlager_martingale_1984}
\newblock K.~Oelschl{\"a}ger,
\newblock A {Martingale} {Approach} to the {Law} of {Large} {Numbers} for
  {Weakly} {Interacting} {Stochastic} {Processes},
\newblock \emph{Ann. Probab.}, \textbf{12} (1984), 458--479,
\newblock \urlprefix\url{https://projecteuclid.org/euclid.aop/1176993301}.

\bibitem{otto_generalization_2000}
\newblock F.~Otto and C.~Villani,
\newblock Generalization of an {Inequality} by {Talagrand} and {Links} with the
  {Logarithmic} {Sobolev} {Inequality},
\newblock \emph{J. Funct. Anal.}, \textbf{173} (2000), 361--400.

\bibitem{parthasarathy_probability_1967}
\newblock K.~R. Parthasarathy,
\newblock \emph{Probability {Measures} on {Metric} {Spaces}},
\newblock Academic Press, 1967,
\newblock
  \urlprefix\url{https://linkinghub.elsevier.com/retrieve/pii/C20130081078}.

\bibitem{pinnau_consensus-based_2017}
\newblock R.~Pinnau, C.~Totzeck, O.~Tse and S.~Martin,
\newblock A consensus-based model for global optimization and its mean-field
  limit,
\newblock \emph{Math. Models Methods Appl. Sci.}, \textbf{27} (2017), 183--204,
\newblock
  \urlprefix\url{https://www.worldscientific.com/doi/abs/10.1142/S0218202517400061}.

\bibitem{pulvirenti_kinetic_1996}
\newblock M.~Pulvirenti,
\newblock Kinetic limits for stochastic particle systems,
\newblock in \emph{Probabilistic {Models} for {Nonlinear} {Partial}
  {Differential} {Equations}} (eds. D.~Talay and L.~Tubaro),
\newblock no. 1627 in Lecture {Notes} in {Mathematics}, Springer-Verlag Berlin
  Heidelberg, 1996.

\bibitem{revuz_continuous_1999}
\newblock D.~Revuz and M.~Yor,
\newblock \emph{Continuous {Martingales} and {Brownian} {Motion}},
\newblock 3rd edition,
\newblock no. 293 in Grundlehren der mathematischen {Wissenschaften},
  Springer-Verlag Berlin Heidelberg, 1999,
\newblock \urlprefix\url{http://link.springer.com/10.1007/978-3-662-06400-9}.

\bibitem{rotskoff_trainability_2019}
\newblock G.~M. Rotskoff and E.~Vanden-Eijnden,
\newblock Trainability and {Accuracy} of {Neural} {Networks}: {An}
  {Interacting} {Particle} {System} {Approach},
\newblock preprint, \arxiv{1805.00915}.

\bibitem{rougerie_finetti_2015}
\newblock N.~Rougerie,
\newblock De {Finetti} theorems, mean-field limits and {Bose}-{Einstein}
  condensation,
\newblock \emph{Lectures notes from a course at the LMU, Munich, 2015},
\newblock \arxiv{1506.05263}.

\bibitem{rueschendorf_comparison_2016}
\newblock L.~Rueschendorf, A.~Schnurr and V.~Wolf,
\newblock Comparison of time-inhomogeneous {Markov} processes,
\newblock \emph{Adv. in Appl. Probab.}, \textbf{48} (2016), 1015--1044.

\bibitem{salem_gradient_2020}
\newblock S.~Salem,
\newblock A gradient flow approach to propagation of chaos,
\newblock \emph{Discrete Contin. Dyn. Syst.}, \textbf{40} (2020), 5729--5754.

\bibitem{serfaty_systems_2019}
\newblock S.~Serfaty,
\newblock Systems of points with {Coulomb} interactions,
\newblock in \emph{Proceedings of the {International} {Congress} of
  {Mathematicians} ({ICM} 2018)},
\newblock World Scientific, Rio de Janeiro, Brazil, 2019,
\newblock 935--977,
\newblock
  \urlprefix\url{https://www.worldscientific.com/doi/abs/10.1142/9789813272880_0033}.

\bibitem{sirignano_mean_2020}
\newblock J.~Sirignano and K.~Spiliopoulos,
\newblock Mean {Field} {Analysis} of {Neural} {Networks}: {A} {Law} of {Large}
  {Numbers},
\newblock \emph{SIAM J. Appl. Math.}, \textbf{80} (2020), 725--752,
\newblock \urlprefix\url{https://epubs.siam.org/doi/10.1137/18M1192184}.

\bibitem{stroock_multidimensional_1997}
\newblock D.~W. Stroock and S.~R.~S. Varadhan,
\newblock \emph{Multidimensional {Diffusion} {Processes}},
\newblock Classics in {Mathematics}, Springer Berlin Heidelberg, 1997,
\newblock \urlprefix\url{http://link.springer.com/10.1007/3-540-28999-2}.

\bibitem{sznitman_equations_1984}
\newblock A.-S. Sznitman,
\newblock {\'E}quations de type de {Boltzmann}, spatialement homog{\`e}nes,
\newblock \emph{Zeitschrift f{\"u}r Wahrscheinlichkeitstheorie und Verwandte
  Gebiete}, \textbf{66} (1984), 559--592.

\bibitem{sznitman_nonlinear_1984}
\newblock A.-S. Sznitman,
\newblock Nonlinear {Reflecting} {Diffusion} {Process}, and the {Propagation}
  of {Chaos} and {Fluctuations} {Associated},
\newblock \emph{J. Funct. Anal.}, \textbf{56} (1984), 311--336.

\bibitem{sznitman_topics_1991}
\newblock A.-S. Sznitman,
\newblock Topics in propagation of chaos,
\newblock in \emph{{\'E}c. {\'E}t{\'e} {Probab}. {St}.-{Flour}
  {XIX}{\textemdash}1989},
\newblock Springer, 1991,
\newblock 165--251.

\bibitem{graham_probabilistic_1996}
\newblock Denis Talay and Luciano Tubaro (eds.),
\newblock \emph{Probabilistic {Models} for {Nonlinear} {Partial} {Differential}
  {Equations}},
\newblock no. 1627 in Lecture {Notes} in {Mathematics}, Springer-Verlag Berlin
  Heidelberg, 1996,
\newblock \urlprefix\url{http://link.springer.com/10.1007/BFb0093175}.

\bibitem{tanaka_probabilistic_1978}
\newblock H.~Tanaka,
\newblock Probabilistic {Treatment} of the {Boltzmann} {Equation} of
  {Maxwellian} {Molecules},
\newblock \emph{Zeitschrift f{\"u}r Wahrscheinlichkeitstheorie und Verwandte
  Gebiete}, \textbf{46} (1978), 67--105.

\bibitem{tanaka_probabilistic_1983}
\newblock H.~Tanaka,
\newblock Some probabilistic problems in the spatially homogeneous {Boltzmann}
  equation,
\newblock in \emph{Theory and {Application} of {Random} {Fields}, {Proceedings}
  of the {IFIP}-{WG} 7/1 {Working} {Conference}, {Bangalore} 1982} (ed.
  G.~Kallianpur),
\newblock Lecture {Notes} in {Control} and {Information} {Sciences},
  Springer-Verlag Berlin Heidelberg, 1983,
\newblock 258--267.

\bibitem{totzeck_trends_2021}
\newblock C.~Totzeck,
\newblock Trends in {Consensus}-based optimization,
\newblock preprint, \arxiv{2104.01383}.

\bibitem{trotter_approximation_1958}
\newblock H.~F. Trotter,
\newblock Approximation of semi-groups of operators,
\newblock \emph{Pacific J. Math.}, \textbf{8} (1958), 887--919.

\bibitem{vicsek_collective_2012}
\newblock T.~Vicsek and A.~Zafeiris,
\newblock Collective motion,
\newblock \emph{Phys. Rep.}, \textbf{517} (2012), 71--140,
\newblock
  \urlprefix\url{https://linkinghub.elsevier.com/retrieve/pii/S0370157312000968}.

\bibitem{villani_limite_2001}
\newblock C.~Villani,
\newblock Limite de champ moyen,
\newblock \emph{Cours de DEA}.

\bibitem{villani_review_2002}
\newblock C.~Villani,
\newblock A {Review} of {Mathematical} {Topics} in {Collisional} {Kinetic}
  {Theory},
\newblock in \emph{Handbook of {Mathematical} {Fluid} {Dynamics}} (eds.
  S.~Friedlander and D.~Serre), vol.~1,
\newblock Elsevier Science, 2002,
\newblock 71--74,
\newblock
  \urlprefix\url{https://linkinghub.elsevier.com/retrieve/pii/S1874579202800040}.

\bibitem{villani_topics_2003}
\newblock C.~Villani,
\newblock \emph{Topics in {Optimal} {Transportation}},
\newblock no.~58 in Graduate {Studies} in {Mathematics}, American Mathematical
  Society, 2003.

\bibitem{villani_optimal_2009}
\newblock C.~Villani,
\newblock \emph{Optimal {Transport}, {Old} and {New}},
\newblock no. 338 in Grundlehren der mathematischen {Wissenschaften},
  Springer-Verlag Berlin Heidelberg, 2009,
\newblock \urlprefix\url{http://link.springer.com/10.1007/978-3-540-71050-9}.

\bibitem{wagner_functional_1996}
\newblock W.~Wagner,
\newblock A functional law of large numbers for {Boltzmann} type stochastic
  particle systems,
\newblock \emph{Stoch. Anal. Appl.}, \textbf{14} (1996), 591--636.

\bibitem{wild_boltzmanns_1951}
\newblock E.~Wild,
\newblock On {Boltzmann}'s equation in the kinetic theory of gases,
\newblock \emph{Math. Proc. Camb. Phil. Soc.}, \textbf{47} (1951), 602--609,
\newblock
  \urlprefix\url{https://www.cambridge.org/core/product/identifier/S0305004100026992/type/journal_article}.

\end{thebibliography}

\medskip
Received xxxx 20xx; revised xxxx 20xx.
\medskip

\end{document}